\documentclass{article}

\usepackage{arxiv}

\usepackage[cp1252]{inputenc}
\usepackage[T1]{fontenc}    
\usepackage{hyperref}       
\usepackage{url}            
\usepackage{booktabs}       
\usepackage{amsfonts}       
\usepackage{nicefrac}       
\usepackage{microtype}      
\usepackage{lipsum}
\usepackage{graphicx}

\usepackage{textcomp}
\usepackage{amsmath}
\usepackage{amssymb}
\usepackage{adjustbox}
\usepackage[abs]{overpic}
\usepackage{hyphenat}
\usepackage[ruled,vlined]{algorithm2e}
\usepackage{amsthm}
\usepackage{bbold}
\usepackage{enumitem}
\setenumerate[1]{label = \arabic*.}
\setenumerate[2]{label*= \arabic*}
\setenumerate[3]{label*=.\arabic*}
\setenumerate[4]{label*=.\arabic*}
\setenumerate[4]{label*=.\arabic*}
\usepackage{xcolor}
\usepackage{multirow}
\usepackage[titletoc]{appendix}
\usepackage{anyfontsize}
\usepackage{framed}
\usepackage{footnote}
\usepackage{mathbbol}
\usepackage{subcaption}
\usepackage{placeins}

\usepackage{adjustbox}
\usepackage{algorithmic}

\usepackage{comment}
\usepackage{xcolor}

\newcommand{\cblack}{\color{black}}

\newcommand{\cgray}{\color{gray}}
\usepackage{marginnote}

\usepackage{tabulary}
\usepackage{bbm} 
\newtheoremstyle{dotless}{}{}{\itshape}{}{\bfseries}{}{ }{}
\theoremstyle{dotless}

\makeatletter

\renewcommand{\mathbf}{\boldsymbol}

\graphicspath{{figures/}}

\usepackage{color}
\usepackage{soul}
\usepackage{amsthm}

\usepackage{units}

\newcommand{\bb}{\mbox{\boldmath{$b$}}}

\newcommand{\bd}{\mbox{\boldmath{$d$}}}

\newcommand{\bg}{\mbox{\boldmath{$g$}}}
\newcommand{\bI}{\mbox{\boldmath{$I$}}}
\newcommand{\bK}{\mbox{\boldmath{$K$}}}

\newcommand{\bn}{\mbox{\boldmath{$n$}}}
\newcommand{\bN}{\mbox{\boldmath{$N$}}}
\newcommand{\bp}{\mbox{\boldmath{$p$}}}

\newcommand{\bs}{\mbox{\boldmath{$s$}}}
\newcommand{\bt}{\mbox{\boldmath{$t$}}}

\newcommand{\bu}{\mbox{\boldmath{$u$}}}
\newcommand{\bU}{\mbox{\boldmath{$U$}}}

\newcommand{\bv}{\mbox{\boldmath{$v$}}}

\newcommand{\bx}{\mbox{\boldmath{$x$}}}

\newcommand{\bepsilon}{\mbox{\boldmath{$\varepsilon$}}}

\newcommand{\bsigma}{\mbox{\boldmath{$\sigma$}}}

\newcommand{\DeltaVarepsilon}{\mbox{$\Delta\text{\large{\mbox{$\varepsilon$}}}$}}

\newcommand{\bzero}{\mbox{$\bf 0$}}

\newfont{\twelvemsb}{msbm10 at 11.6pt}

\renewcommand{\div}{\mathop{\rm div}\nolimits}

\newcommand{\tr}{\mathop{\rm tr}}

\title{On mesh refinement procedures for polygonal virtual elements}

\author{
	Daniel van Huyssteen \\
	Institute of Applied Mechanics, 
	Friedrich-Alexander-Universität Erlangen-Nürnberg
	Erlangen, 91058, Germany \\
	\texttt{daniel.van.huyssteen@fau.de} \\
	\And
	Felipe Lopez Rivarola \\
	Facultad de Ingeniería,
	Universidad de Buenos Aires,
	Buenos Aires, C1127AAR, Argentina,
	CONICET\\
	\And
	Guillermo Etse \\
	Facultad de Ingeniería,
	Universidad de Buenos Aires,
	Buenos Aires, C1127AAR, Argentina,
	CONICET\\
	\AND
	Paul Steinmann \\
	Institute of Applied Mechanics, 
	Friedrich-Alexander-Universität Erlangen-Nürnberg
	Erlangen, 91058, Germany \\
}

\begin{document}
	\maketitle
	
	\begin{abstract}
		This work concerns adaptive refinement procedures for the virtual element method. Specifically, adaptive refinement procedures previously proposed by the authors, and investigated for meshes of structured quadrilateral elements with compressible material behaviour, are studied in more general applications.	
		The performance of the proposed refinement procedures is studied through a broad numerical campaign considering meshes of both structured quadrilateral elements and irregular unstructured/Voronoi elements. Furthermore, the robustness of the proposed procedures with respect to compressibility is investigated through the consideration of problems involving both compressible and nearly incompressible material behaviour.
		The results demonstrate that the high level of efficacy and efficiency exhibited by the proposed procedures on meshes of structured quadrilateral elements with compressible material behaviour, as seen in previous work by the authors, is also achieved in the case of meshes of irregular unstructured/Voronoi meshes. Furthermore, the proposed procedures exhibit robust behaviour with respect to compressibility and exhibit good performance in the cases of both compressibility and near-incompressibility.
		The versatility and efficacy of the proposed refinement procedures demonstrated over a variety of mesh types, for varying levels of compressibility, and over a large range of problems indicates that the procedures are well-suited to virtual element applications.
	\end{abstract}

	
	\section{Introduction}
	\label{S:Introduction}
	The virtual element method (VEM) \cite{VEIGA2012,Veiga2014} is a popular extension of the finite element method (FEM) and has been applied to a growing range of problems in solid mechanics \cite{VEMContactWriggers2016,Aldakheel2018,Aldakheel2019a,Chi2017,WriggersPlastic2017,HudobivnikPlastic2018,Aldakheel2019}.
	An attractive feature of the VEM is its robust performance under challenging numerical conditions. For example, the low-order displacement-based VEM exhibits locking-free behaviour in cases of limiting material properties, such as near-incompressibility and near-inextensibility \cite{Reddy2019,Tang2020}, even under large deformations \cite{Huyssteen2020,Huyssteen2021}.
	However, the most noted feature of the VEM is its permission of arbitrary polygonal/polyhedral element geometries, both convex and non-convex. 
	Additionally, the method has been extended to accommodate elements with arbitrarily curved edges \cite{ArtioliBeiraoDassi2020,WriggersHudobivnikAldakheel2020,Bellis2019}. 
	
	The mesh flexibility offered by the VEM means that it is well suited to problems involving complex geometries \cite{Antonietti2017} and material features such as inclusions or grains \cite{Rivarola2019}. The discretization of complex geometries can result in highly distorted and/or stretched element geometries which, in the case of both finite and virtual elements, could significantly reduce the accuracy of the method \cite{knupp2007remarks,Lowrie2011,Sorgente2021,Sorgente2022}. However, the VEM stabilization term can be easily tuned to improve the method's performance in cases of challenging element geometries \cite{ReddyvanHuyssteen2021,DvH_BDR_MeshQuality}. Thus, further motivating application of the VEM to problems involving complex geometries.
	The geometric freedom offered by the VEM, together with it suitability for complex geometries, naturally lends the method to problems involving adaptive remeshing, a topic of rapidly growing interest.
	
	%
	
	In the context of adaptive remeshing, the VEM provides significant advantages over the FEM as new elements may have arbitrary geometries. Furthermore, additional nodes may be inserted arbitrarily along element edges with no consideration or treatment of hanging nodes required, and anisotropic mesh refinement poses no challenge to the VEM \cite{vanHuyssteen2022}. 
	However, since the virtual element basis functions are not known explicitly, computation of well-known mesh refinement indicators for finite elements, such as the $Z^{2}$ and Kelly estimators \cite{Zienkiewicz1987,Kelly1983}, cannot be performed trivially in a virtual element setting.
	These advantages and considerations have motivated the study and development of various approaches related to the calculation of error estimators for the VEM \cite{Veiga2015a,Cangiani2017,Berrone2017,Mora2017,Veiga2019a,DAltri2020,Artioli2020a}.
	These approaches have largely focused on scalar problems such as Poisson's problem or Steklov's eigenvalue problem \cite{Veiga2015a,Cangiani2017,Berrone2017,Mora2017,Veiga2019a}, however, there is growing interest in the computation of error estimators and adaptive meshing procedures for elastic problems \cite{DAltri2020,Artioli2020a}.

	In \cite{vanHuyssteen2022} a variety of adaptive mesh refinement procedures, constructed specifically for application to problems involving the VEM, were proposed and numerically evaluated in the context of linear elasticity. The procedures involved mesh refinement indicators computed from the displacement and strain field approximations and were studied both separately and in combination. The convergence behaviour of the various procedures was studied in the $\mathcal{H}^{1}$ error norm and compared to a reference refinement procedure. It was found that the best performance was achieved when using a combination of refinement procedures based on both the displacement and strain fields. Furthermore, all of the proposed procedures provided significant improvements in computational efficiency and were able to generate solutions of equivalent accuracy to a reference procedure while using significantly fewer degrees of freedom and much less run time. The scope of the study was limited to the consideration of only structured quadrilateral meshes and problems involving compressible material behaviour. Since the VEM permits arbitrary element geometries, and is locking-free in the case of nearly incompressible materials, questions naturally arise about the performance of the procedures proposed in \cite{vanHuyssteen2022} when applied to meshes of irregular elements and in the presence of near-incompressibility.
	
	This work represents the extension of the study presented in \cite{vanHuyssteen2022} to more general applications. In this work the mesh refinement procedures presented in \cite{vanHuyssteen2022} are applied to problems involving regular structured meshes and irregular unstructured/Voronoi meshes to investigate the geometric versatility of the proposed refinement procedures. Additionally, problems involving both compressible and nearly incompressible material behaviour are considered to investigate the robustness of the proposed refinement procedures with respect to compressibility.
	The proposed refinement procedures are comparatively evaluated over a wide range of numerical problems to study their efficacy and efficiency. The efficacy is measured in the $\mathcal{H}^{1}$ error norm, while the efficiency is evaluated in terms of both the number of degrees of freedom involved in a problem and the run time compared to a reference procedure. Additionally, the convergence behaviour in the $\mathcal{H}^{1}$ error norm with respect to mesh size, and the convergence of the energy contribution from the stabilization term are studied as in \cite{vanHuyssteen2022}. Finally, the investigation of the convergence behaviour of the displacement and strain field approximations in the $\mathcal{L}^{2}$ error norm represents a further extension to \cite{vanHuyssteen2022}.
	Through this investigation the most effective of the proposed mesh refinement procedures is identified.
	
	
	The structure of the rest of this work is as follows. The governing equations of linear elasticity are set out in Section~\ref{sec:GovEq}. This is followed in Section~\ref{sec:VEM} by a description of the first-order virtual element method. 
	The procedures used to generate and refine meshes are presented in Section~\ref{sec:MeshGenerationAndRefinement}. This is followed, in Section~\ref{sec:MeshRefinementIndicators}, by a description of the various mesh refinement indicators along with the procedure used to identify elements qualifying for refinement.
	Section~\ref{sec:Results} comprises a set of numerical results through which the performance of the various refinement procedures is evaluated.
	Finally, the work concludes in Section~\ref{sec:Conclusion} with a discussion of the results.
	
	\section{Governing equations of linear elasticity} 
	\label{sec:GovEq}
\begin{figure}[ht!]
		\centering
		\includegraphics[width=0.25\textwidth]{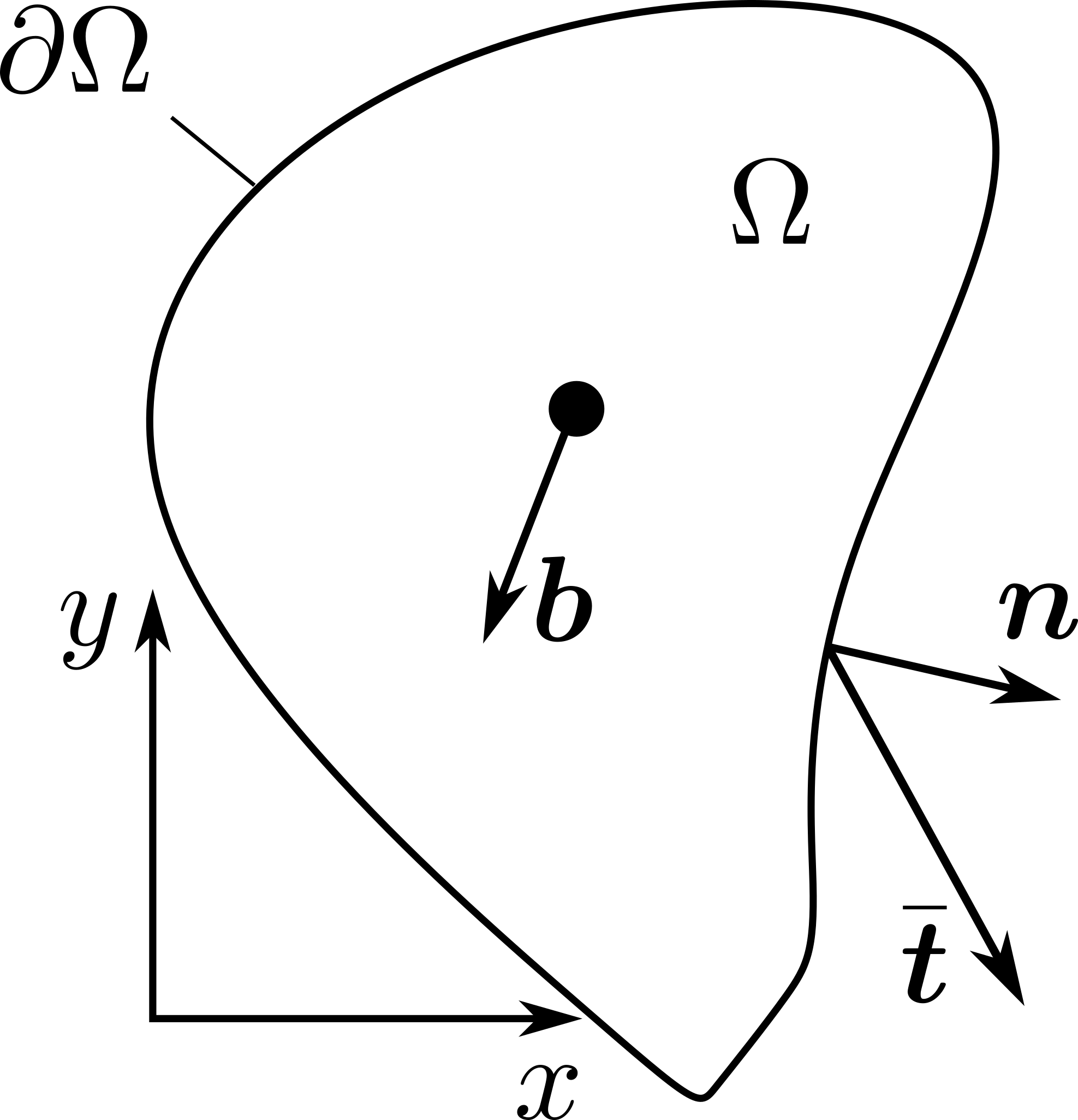}
		\caption{Arbitrary elastic body subject to body force and traction.
			\label{fig:ElasticBody}}
	\end{figure}

	Consider an arbitrary elastic body occupying a plane, bounded, domain ${\Omega \subset \mathbb{R}^{2} }$ subject to a traction ${\bar{\bt}}$ and body force ${\bb}$ (see Figure~\ref{fig:ElasticBody}).
	The boundary ${\partial \Omega}$ has an outward facing normal denoted by $\bn$ and comprises a non-trivial Dirichlet part $\Gamma_{D}$ and a Neumann part $\Gamma_{N}$ such that ${\Gamma_{D} \cap \Gamma_{N} = \emptyset}$ and ${\overline{\Gamma_{D} \cup \Gamma_{N}}=\partial \Omega}$.
	
	\FloatBarrier
	\FloatBarrier
	
	In this work small displacements are assumed and the strain-displacement is relation given by 
	\begin{equation}
		\bepsilon \left( \bu \right) = \frac{1}{2} \left[\nabla\,\bu + \left[ \nabla\,\bu \right]^{T} \right] \, . \label{eqn:DisplacementStrain}
	\end{equation}
	Here the displacement is denoted by $\bu$, ${\bepsilon}$ is the symmetric infinitesimal strain tensor and ${\nabla \left( \bullet \right) = \frac{\partial \left( \bullet \right)_{i} }{\partial \, x_{j}} \, \boldsymbol{e}_{i} \otimes \boldsymbol{e}_{j} }$ is the gradient of a vector quantity. Additionally, the body is assumed to be linear elastic with the stress-strain relation given by
	\begin{equation}
		\bsigma = \mathbb{C} : \bepsilon \,. \label{eqn:StressStrain1}
	\end{equation}
	Here, ${\bsigma}$ is the Cauchy stress tensor and ${\mathbb{C}}$ is a fourth-order constitutive tensor. For a linear elastic and isotropic material (\ref{eqn:StressStrain1}) is given by
	\begin{eqnarray}
		\bsigma = \lambda \tr \left( \bepsilon \right)  \bI + 2\mu \, \bepsilon \, , \label{eqn:StressStrain2}
	\end{eqnarray}
	where ${\tr \left( \bullet \right)}$ denotes the trace, $\bI$ is the second-order identity tensor, and $\lambda$ and $\mu$ are the Lam\'{e} parameters.
	
	For equilibrium it is required that 
	\begin{equation}
		\div \, \bsigma + \bb = \bzero \, , \label{eqn:Equilibrium}
	\end{equation}
	where ${\div \left( \bullet \right) = \frac{\partial \left( \bullet \right)_{ij} }{\partial \, x_{j}} \boldsymbol{e}_{i} }$ is the divergence of a tensor quantity.
	The Dirichlet and Neumann boundary conditions are given by 
	\begin{align}
		&\bu = \bg \quad \text{on } \Gamma_{D} \, , \text{ and} \label{eqn:DirichletBC} \\
		&\bsigma \cdot \bn = \bar{\bt} \quad \text{on } \Gamma_{N} \, , \label{eqn:NeumannBC}
	\end{align}
	respectively, with $\bg$ and $\bar{\bt}$ denoting prescribed displacements and tractions respectively.
	Equations (\ref{eqn:StressStrain2})-(\ref{eqn:NeumannBC}), together with the displacement-strain relationship (\ref{eqn:DisplacementStrain}), constitute the boundary-value problem for a linear elastic isotropic body.
	
	\subsection{Weak form}
	\label{subsec:WeakForm}
	The space of square-integrable functions on $\Omega$ is hereinafter denoted by ${\mathcal{L}^{2}\left(\Omega\right)}$. The Sobolev space of functions that, together with their first derivatives, are square-integrable on $\Omega$ is hereinafter denoted by ${\mathcal{H}^{1}\left(\Omega\right)}$. Additionally, the function space $\mathcal{V}$ is introduced and defined as
	\begin{align}
		\mathcal{V} = \left[ \mathcal{H}^{1}_{D} \left(\Omega\right) \right]^{d} 
		=
		\left\{ \bv \, | \, v_{i} \in \mathcal{H}^{1}\left(\Omega\right), \, \bv = \boldsymbol{0} \,\, \text{on} \,\, \Gamma_{D}  \right\} \, 
	\end{align}
	where ${d=2}$ is the dimension.
	Furthermore, the function ${\bu_{g}\in \left[ \mathcal{H}^{1} \left(\Omega\right) \right]^{d} }$ is introduced satisfying (\ref{eqn:DirichletBC}) such that ${\bu_{g}|_{\Gamma_{D}}=\bg}$.
	
	The bilinear form ${a\left(\cdot,\cdot\right)}$, where ${a : \left[ \mathcal{H}^{1} \left(\Omega\right) \right]^{d}  \times \left[ \mathcal{H}^{1} \left(\Omega\right) \right]^{d} \rightarrow \mathbb{R}}$, and the linear functional ${\ell \left(\cdot\right)}$, where ${\ell : \left[ \mathcal{H}^{1} \left(\Omega\right) \right]^{d} \rightarrow \mathbb{R}}$, are defined respectively by
	\begin{equation}
		a\left(\bu, \, \bv \right) = \int_{\Omega} \bsigma \left( \bu \right) : \bepsilon \left( \bv \right) \, dx \, , \label{eqn:BilinearForm}
	\end{equation}
	and 
	\begin{equation}
		\ell \left( \bv \right) = \int_{\Omega} \bb \cdot \bv \, dx + \int_{\Gamma_{N}} \bar{\bt} \cdot \bv \, ds - a\left(\bu_{g},\, \bv \right) \, . \label{eqn:LinearFuntional} 
	\end{equation}
	
	The weak form of the problem is then: given ${\bb \in \left[ \mathcal{L}^{2}\left(\Omega\right) \right]^{d} }$ and ${ \bar{\bt} \in \left[ \mathcal{L}^{2}\left(\Gamma_{N}\right) \right]^{d} }$, find ${\bU \in \left[ \mathcal{H}^{1}\left(\Omega\right) \right]^{d} }$ such that 
	\begin{equation}
		\bU = \bu + \bu_{g} \, , \quad \bu \in \mathcal{V} \, ,
	\end{equation}
	and
	\begin{equation}
		a\left(\bu , \, \bv \right) = \ell \left( \bv \right) \, , \quad \forall \bv \in \mathcal{V} \, .
	\end{equation}
	
	\section{The virtual element method}
	\label{sec:VEM}
	
\begin{figure}[ht!]
		\centering
		\includegraphics[width=0.45\textwidth]{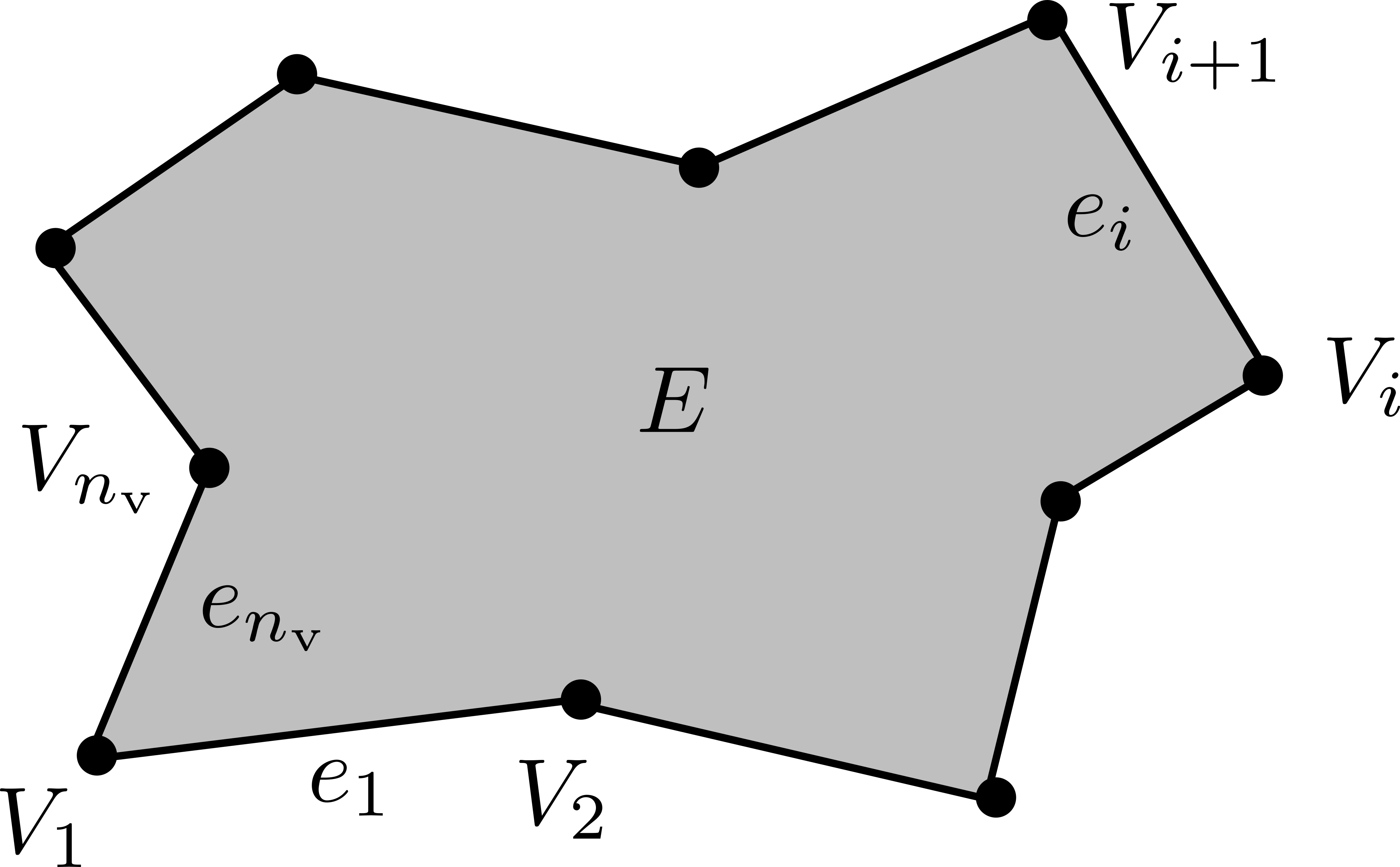}
		\caption{Sample virtual element.
			\label{fig:SampleElement}}
	\end{figure}

	The domain $\Omega$ is partitioned into a mesh of non-overlapping arbitrary polygonal elements\footnote{If $\Omega$ is not polygonal the mesh will be an approximation of $\Omega$} $E$ with $\overline{\cup E}=\overline{\Omega}$. Here $E$ denotes the element domain and $\partial E$ its boundary, with ${\overline{(\, \bullet \,)}}$ denoting the closure of a set.
	An example of a typical element is depicted in Figure~\ref{fig:SampleElement} with edge $e_{i}$ connecting vertices $V_{i}$ and $V_{i+1}$. Here ${i=1,\dots,n_{\rm v}}$ with $n_{\rm v}$ denoting the total number of vertices of element $E$. 
	
	\FloatBarrier
	\FloatBarrier
	A conforming approximation is constructed in a space ${\mathcal{V}^{h} \subset \mathcal{V}}$ where $\mathcal{V}^{h}$ comprises functions that are continuous on the domain $\Omega$, piecewise linear on the boundary ${\partial E}$, and with ${\nabla^{2}\, \bv_{h} }$ vanishing on $E$ (see \cite{VEIGA2012}). 
	The space $\mathcal{V}^{h}$ is built-up element-wise with 
	\begin{equation}
		\mathcal{V}^{h}|_{E} = \left\{ \bv_{h} \in \mathcal{V} \, | \, \bv_{h} \in \left[ \mathcal{C}(E) \right]^{2} \, , \, \nabla^{2} \, \bv_{h}  = \boldsymbol{0} \text{ on } E \, , \, \bv_{h}|_{e} \in \mathcal{P}_{1}(e)  \right\} \,.
	\end{equation}
	Here ${\mathcal{P}_{k}(X)}$ is the space of polynomials of degree ${\leq \, k}$ on the set ${X \, \subset \, \mathbb{R}^{d} }$ with ${d=1,\,2}$ and ${\nabla^{2}=\nabla\cdot\nabla}$ is the Laplacian operator. 
	
	Since all computations are performed on element edges it is convenient to write, for element $E$,
	\begin{equation}
		\bv_{h}|_{\partial E} = \bN \cdot \bd^{E} \,. \label{eqn:DisplacementTrace}
	\end{equation}
	Here, $\bN$ is a matrix of standard linear Lagrangian basis functions and $\bd^{E}$ is a ${{2}n_{\rm v} \times 1}$ vector of the degrees of freedom associated with $E$. The virtual basis functions are not known, nor required to be known on $E$; their traces, however, are known and are simple Lagrangian functions.
	
	
	The virtual element projection ${\Pi \, : \, \mathcal{V}^{h}|_{E} \rightarrow \mathcal{P}_{0}(E)  }$ is required to satisfy
	\begin{equation}
		\int_{E} \Pi \, \bv_{h} \cdot \bepsilon\left( \bp \right) \, dx = \int_{E} \bepsilon\left(\bv_{h}\right) \cdot \bepsilon\left( \bp \right) \, dx \quad \forall \bp \in \mathcal{P}_{1} \,, \label{eqn:Projection}
	\end{equation}
	where ${\Pi \, \bv_{h}}$ represents the ${\mathcal{L}^{2}}$ projection of the symmetric gradient of ${\bv_{h}}$ onto constants (see \cite{Artioli2017}). Since the projection is constant at element-level, after applying integration by parts to (\ref{eqn:Projection}), and considering (\ref{eqn:DisplacementTrace}), the components of the projection can be computed as
	\begin{align}
		\left(\Pi\,\bv_{h}\right)_{ij} &= \frac{1}{2}\frac{1}{|E|}  \sum_{e\in\partial E} \int_{e} \left[ N_{iA} \, d_{A}^{E} \, n_{j} + N_{jA} \, d_{A}^{E} \, n_{i}\right] ds \,, \label{eqn:ComputeProjection}
	\end{align}
	where summation is implied over repeated indices. Additionally, the condition ${\nabla^{2}\, \bv_{h}  = \boldsymbol{0} \text{ on } E}$ ensures computability of the projection \cite{Artioli2017}.
	
	The virtual element approximation of the bilinear form (\ref{eqn:BilinearForm}) is constructed by writing 
	\begin{align}
		a^{E}\left(\bu,\,\bv\right) :&= a\left(\bu,\,\bv\right)|_{E} 
		= \int_{E} \bepsilon\left(\bv_{h}\right) : \left[ \mathbb{C} : \bepsilon\left(\bu_{h}\right) \right] dx \, , \label{eqn:ElementBilinearForm}
	\end{align}
	where ${a^{E}\left(\cdot,\cdot\right)}$ is the contribution of element $E$ to the bilinear form ${a\left(\cdot,\cdot\right)}$. Consideration of (\ref{eqn:ComputeProjection}) allows (\ref{eqn:ElementBilinearForm}) to be written as (see \cite{Reddy2019})
	\begin{align}
		a^{E}\left(\bu_{h},\,\bv_{h}\right) 
		&= \underbrace{\int_{E} \Pi\,\bv_{h} : \left[ \mathbb{C} : \Pi\,\bu_{h} \right] dx }_{\text{Consistency term}} 
		+ \underbrace{\int_{E} \left[ \bepsilon\left( \bv_{h} \right) : \left[ \mathbb{C} : \bepsilon \left( \bu_{h} \right) \right] - \Pi\,\bv_{h} : \left[ \mathbb{C} : \Pi \, \bu_{h} \right] \right] dx }_{\text{Stabilization term}} \,, \label{eqn:ExpandedBilinearForm}
	\end{align}
	where the remainder term is discretized by means of a stabilization.
	
	\subsection{The consistency term} 
	\label{subsec:ConsistencyTerm}
	
	The projection (\ref{eqn:ComputeProjection}), and thus the consistency term, can be computed exactly and yields
	\begin{equation}
		a_{\rm c}^{E}\left(\bu_{h},\,\bv_{h}\right) \, = \, \int_{E} \Pi\,\bv_{h} : \left[ \mathbb{C} : \Pi\,\bu_{h} \right] dx \, = \, \widehat{\bd}^{E} \cdot \left[ \bK_{\rm c}^{E} \cdot \bd^{E} \right] \,.
	\end{equation}
	%
	Here $\bK_{\rm c}^{E}$ is the consistency part of the stiffness matrix of element $E$ with ${\widehat{\bd}^{E}}$ and ${\bd^{E}}$ the degrees of freedom of $\bv_{h}$ and $\bu_{h}$ respectively that are associated with element $E$.
	
	\subsection{The stabilization term} 
	\label{subsec:Stab}
	
	The remainder term cannot be computed exactly and is approximated by means of a discrete stabilization term \cite{Gain2014,Veiga2015}.
	The approximation employed in this work is motivated by seeking to approximate the difference between the element degrees of freedom $\bd^{E}$ and the nodal values of a linear function that is closest to $\bd^{E}$ in some way (see \cite{Artioli2017,Reddy2019,vanHuyssteen2022}). 
	The nodal values of the linear function are given by
	\begin{equation}
		\widetilde{\bd} = \boldsymbol{\mathcal{D}} \cdot \bs \,. \label{eqn:LinearApprox}
	\end{equation}
	Here $\bs$ is a vector of the degrees of freedom of the linear function and $\boldsymbol{\mathcal{D}}$ is a matrix relating $\widetilde{\bd}$ to $\bs$ with respect to a scaled monomial basis. For the full expression of $\boldsymbol{\mathcal{D}}$ see \cite{Artioli2017,Reddy2019}.
	After some manipulation (see, again, \cite{Reddy2019}) the stabilization term of the bilinear form can be approximated as
	\begin{equation}
		a_{\text{stab}}^{E}\left(\bu_{h},\,\bv_{h}\right) \, = \,
		\int_{E} \left[ \bepsilon\left( \bv_{h} \right) : \left[ \mathbb{C} : \bepsilon \left( \bu_{h} \right) \right] - \Pi\,\bv_{h} : \left[ \mathbb{C} : \Pi \, \bu_{h} \right]  \right] dx \, \approx \, \widehat{\bd}^{E} \cdot \bK_{\rm s}^{E} \cdot \bd^{E} \, ,
	\end{equation}
	where $\bK_{\rm s}^{E}$ is the stabilization part of the stiffness matrix of element $E$ and is defined as
	\begin{equation}
		\bK_{\rm s}^{E} = \mu \left[ \bI - \boldsymbol{\mathcal{D}} \cdot \left[ \boldsymbol{\mathcal{D}}^{T} \cdot \boldsymbol{\mathcal{D}}\right]^{-1} \cdot  \boldsymbol{\mathcal{D}}^{T} \right] \, .
	\end{equation}
	The total element stiffness matrix ${\bK^{E}}$ is then computed as the sum of the consistency and stabilization matrices, i.e. ${\bK^{E} = \bK_{\rm c}^{E} + \bK_{\rm s}^{E} \,}$.

	\section{Mesh generation and refinement} 
	\label{sec:MeshGenerationAndRefinement}
	In this section the procedures used to generate meshes and refine elements are described.
	
	\subsection{Mesh generation} 
	\label{subsec:MeshGeneration}
\begin{figure}[ht!]
		\centering
		\includegraphics[width=0.95\textwidth]{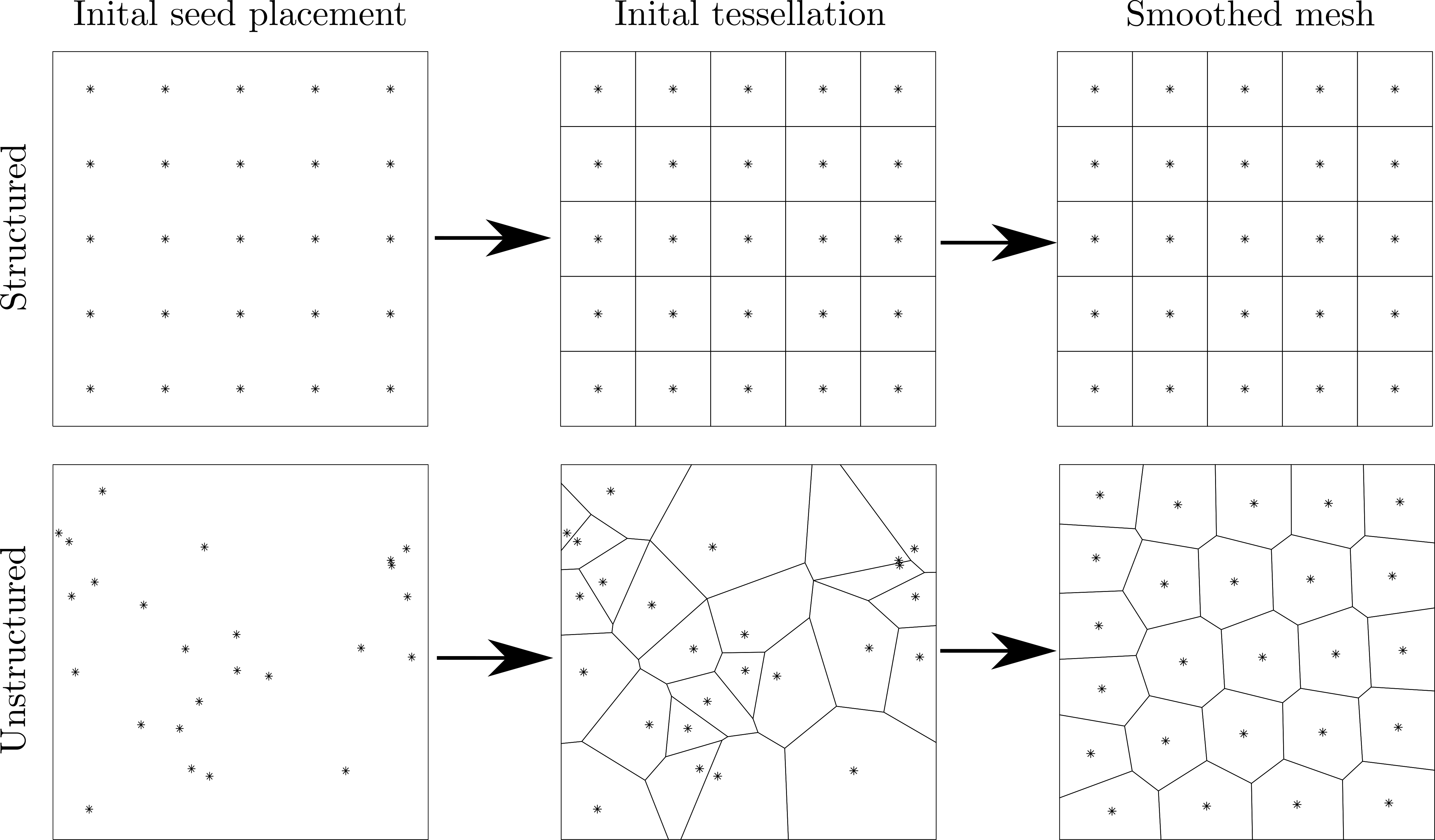}
		\caption{Mesh generation procedure for structured and unstructured/Voronoi meshes.
			\label{fig:MeshGeneration}}
	\end{figure}

	In this work all meshes are created by Voronoi tessellation of a set of seed points. These seeds will be generated in both structured and unstructured sets to create structured and unstructured meshes respectively. 
	In the case of structured meshes seeds points are placed to form a structured grid, while in the case of unstructured/Voronoi meshes seeds are placed arbitrarily within the problem domain. 
	After placement of the seeds an initial Voronoi tessellation is created using PolyMesher \cite{PolyMesher}. Finally, a smoothing algorithm in PolyMesher is used to modify the locations of the seed points to create a mesh in which all elements have (approximately) equal areas.
	The mesh generation procedure is illustrated in Figure~\ref{fig:MeshGeneration} where the top and bottom rows depict the generation of structured and unstructured/Voronoi meshes respectively. The left-hand figures depict the initial sets of seed points, and the central figures depict the initial Voronoi tessellations. Finally, the right-hand figures depict the mesh after the smoothing algorithm has been applied. It is noted that in the case of structured meshes the smoothing procedure is trivial since the seeds are placed such that structured elements of equal areas are created during the initial tessellation.  
	
	\FloatBarrier
	\FloatBarrier
	
	\subsection{Mesh refinement} 
	\label{subsec:MeshRefinement}
\begin{figure}[ht!]
		\centering
		\includegraphics[width=0.95\textwidth]{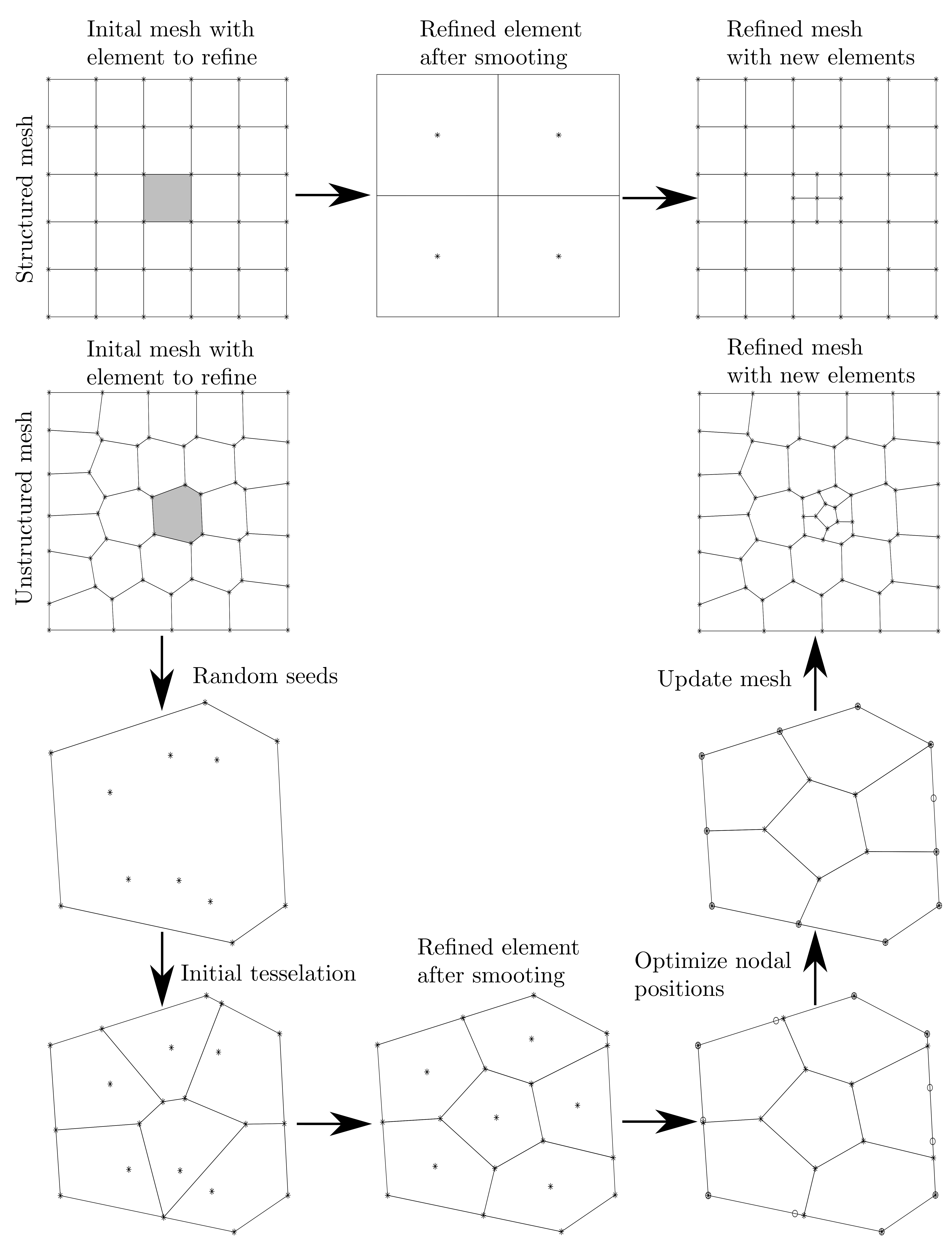}
		\caption{Refinement procedure for structured and unstructured/Voronoi elements.
			\label{fig:MeshRefinement}}
	\end{figure}

	The elements qualifying for refinement are identified using the procedure described in the next section. 
	Once an element has been marked for refinement the refinement procedure is performed using a modified version of PolyMesher \cite{PolyMesher}. The element refinement procedure is illustrated in Figure~\ref{fig:MeshRefinement} for structured and unstructured/Voronoi meshes.
	The element refinement process is similar for both structured and unstructured/Voronoi elements. However, as was the case during mesh generation, a number of steps in the refinement process are trivial in the case of structured elements. 
	The element marked for refinement is shown in grey within the initial mesh (see Figure~\ref{fig:MeshRefinement}). The first step of the refinement process involves generating a set of seed points within the marked element. For simplicity the number of seeds to be generated is chosen to be equal to the number of nodes of the element. In the case of structured meshes the seeds are placed in a structured grid, while in the case of unstructured/Voronoi elements the seeds are placed randomly. Then, as was the case during mesh generation, an initial Voronoi tessellation of the seeds is created and then smoothed using the PolyMesher algorithm.
	After smoothing, a procedure is used to optimize the positions of the newly created nodes that lie on the edges of the original element. This procedure involves looping over the original edges and identifying how many new nodes lie on a specific edge. For example, in the case of the unstructured element depicted in Figure~\ref{fig:MeshRefinement}, and focusing specifically on the middle figure on the bottom row, one new node exists on the left-hand edge, while two new nodes lie on the right-hand edge. For each edge a set of optimal node positions is created. This set contains ${n_{\text{new}}+2}$ node positions with ${n_{\text{new}}}$ the number of new nodes lying on the specific edge, and where the optimal node positions are spaced linearly along the edge. The first and last optimal node positions coincide with the original nodes at the ends of the edge and the remaining positions are evenly distributed along the edge. For clarity, the optimal node positions are indicated as circles in the bottom right-hand figure. The positions of the newly created nodes are then shifted to coincide with the nearest optimal node position. Finally, the new elements are incorporated into the mesh in the place of the original element.
	The node optimization procedure is used to; improve the compatibility of groups of refined elements when neighboring elements are marked for refinement, prevent the creation of very short edges, and to reduce the total number of new nodes introduced during refinement.
	
	\FloatBarrier
	\FloatBarrier
	
	\section{Mesh refinement indicators} 
	\label{sec:MeshRefinementIndicators}
	In this section the various mesh refinement indicators are presented along with the procedure used to determine the elements qualifying for refinement. The mesh refinement indicators are the same as those presented in \cite{vanHuyssteen2022}. 
	
	
	\subsection{Displacement-based refinement indicator} 
	\label{subsec:DispBasedIndicator}
	The displacement-based refinement indicator is motivated by seeking to quantify the deviation from planar of the nodal values of the displacement $\bu_{h}$ on an element $E$.
	A reference planar/linear displacement $\bu_{\rm{lin}}$ is computed on element $E$ following a similar procedure to that presented in \cite{WriggersIsotropic2017} for the projection of $\bu$ onto a linear ansatz space. The $x$- and $y$-components of $\bu_{\rm{lin}}$ are given by
	\begin{equation}
		u_{\rm{lin}}^{x}= a_{1} + a_{3}\,x + a_{5}\,y \, , \quad \text{and} \quad
		u_{\rm{lin}}^{y}= a_{2} + a_{4}\,x + a_{6}\,y \, ,
	\end{equation}
	respectively, with $a_{i}$ the degrees of freedom of $\bu_{\rm{lin}}$.
	The parameters $a_{3,\,\dots,\,6}$ represent components of the average non-symmetric gradient of $\bu_{h}$ on element $E$. The parameters $a_{3,\,6}$ are automatically determined during the computation of the projection (see (\ref{eqn:ComputeProjection})) and are respectively the components $\left(\Pi\,\bu_{h}\right)_{11}$ and $\left(\Pi\,\bu_{h}\right)_{22}$. The parameters $a_{4,\,5}$ relate specifically to the non-symmetric gradient and, depending on the implementation of (\ref{eqn:ComputeProjection}), may require additional computations with 
	\begin{align}
		a_{4}=\frac{1}{|E|} \int_{E}  \frac{\partial \, u_{y}^{h}}{\partial x} \, d\Omega\, , \text{ and} \quad 
		a_{5}&=\frac{1}{|E|} \int_{E}  \frac{\partial \, u_{x}^{h}}{\partial y} \, d\Omega \, .
	\end{align}
	The remaining parameters $a_{1,\,2}$ are computed by summation over element vertices with 
	\begin{align}
		a_{1}=\frac{1}{n_{\rm{v}}}\sum_{i=1}^{n_{\rm{v}}}\left[u_{x}^{h}(V_{i})-a_3\,x_{i} - a_5\,y_{i}\right]
		\, , \quad  \text{and} \quad  
		a_{2}=\frac{1}{n_{\rm{v}}}\sum_{i=1}^{n_{\rm{v}}}\left[u_{y}^{h}(V_{i})-a_4\,x_{i} - a_6\,y_{i}\right] \, .
	\end{align}
	Here $u_{x}^{h}(V_{i})$ and $u_{y}^{h}(V_{i})$ are respectively the $x$- and $y$-displacement components of the $i$-th vertex of element $E$. Similarly, $x_{i}$ and $y_{i}$ are respectively the $x$- and $y$-coordinates of the $i$-th vertex of $E$.
	
	The total displacement-based refinement indicator $\Upsilon^{t}_{\rm{DB}}$ comprises $x$ and $y$ components from each node defined as
	\begin{align}
		\Upsilon_{\text{DB},\,i}^{x} = u_{\rm{lin}}^{x}(V_{i}) - u_{x}^{h}(V_{i}) \, , \quad
		\text{and} \quad 
		\Upsilon_{\text{DB},\,i}^{y} = u_{\rm{lin}}^{y}(V_{i}) - u_{y}^{h}(V_{i}) \, , 
	\end{align}
	respectively. The $x$ and $y$ components are combined into a nodal contribution $\Upsilon_{\text{DB},\,i}^{t}$ by means of an $\mathcal{L}^{2}$/Pythagorean norm with ${\Upsilon_{\text{DB},\,i}= \| \Upsilon_{\text{DB},\,i}^{x},\, \Upsilon_{\text{DB},\,i}^{y}  \|_{\mathcal{L}^{2}}}$.
	Finally, the total indicator is computed as the $\mathcal{L}^{2}$ norm of the nodal contributions with
	\begin{equation}
		\Upsilon_{\text{DB},\,E}^{t} = \| \Upsilon_{\text{DB},\,i} \|_{\mathcal{L}^{2}}
		= \sqrt{\sum_{i=1}^{n_{\rm v}} \left[  \Upsilon_{\text{DB},\,i}\right]^{2} } \, .
	\end{equation}
	It is noted that the motivation behind the displacement-based indicator is similar to that of the stabilization term. It is possible to formulate a qualitatively similar, and computationally cheaper, indicator in terms of the stabilization term (see Section~\ref{subsec:PlateWithHoleTraction}). However, since a variety of approaches to the formulation of the stabilization exist, the displacement-based indicator is presented in this way so that it can be easily implemented into existing VEM formulations/codes.
	Additionally, the presented approach can be easily extended to higher-order and hyperelastic VEM formulations such as those presented in \cite{Bellis2019,WriggersIsotropic2017}.
	Finally, it is noted that the displacement-based indicator degenerates in the case of three-noded triangular elements on which the displacement approximation is inherently planar. This could be resolved by exploiting the flexibility of the VEM and inserting additional nodes along element edges. Alternatively, a combination of mesh refinement indicators (see Section~\ref{subsec:SelectElementsForRefinement}), or an indicator designed for P1 finite elements could be used.
	
	\subsection{Strain jump-based indicator} 
	\label{subsec:StrainJumpBasedIndicator}
\begin{figure}[ht!]
		\centering
		\includegraphics[width=0.4\textwidth]{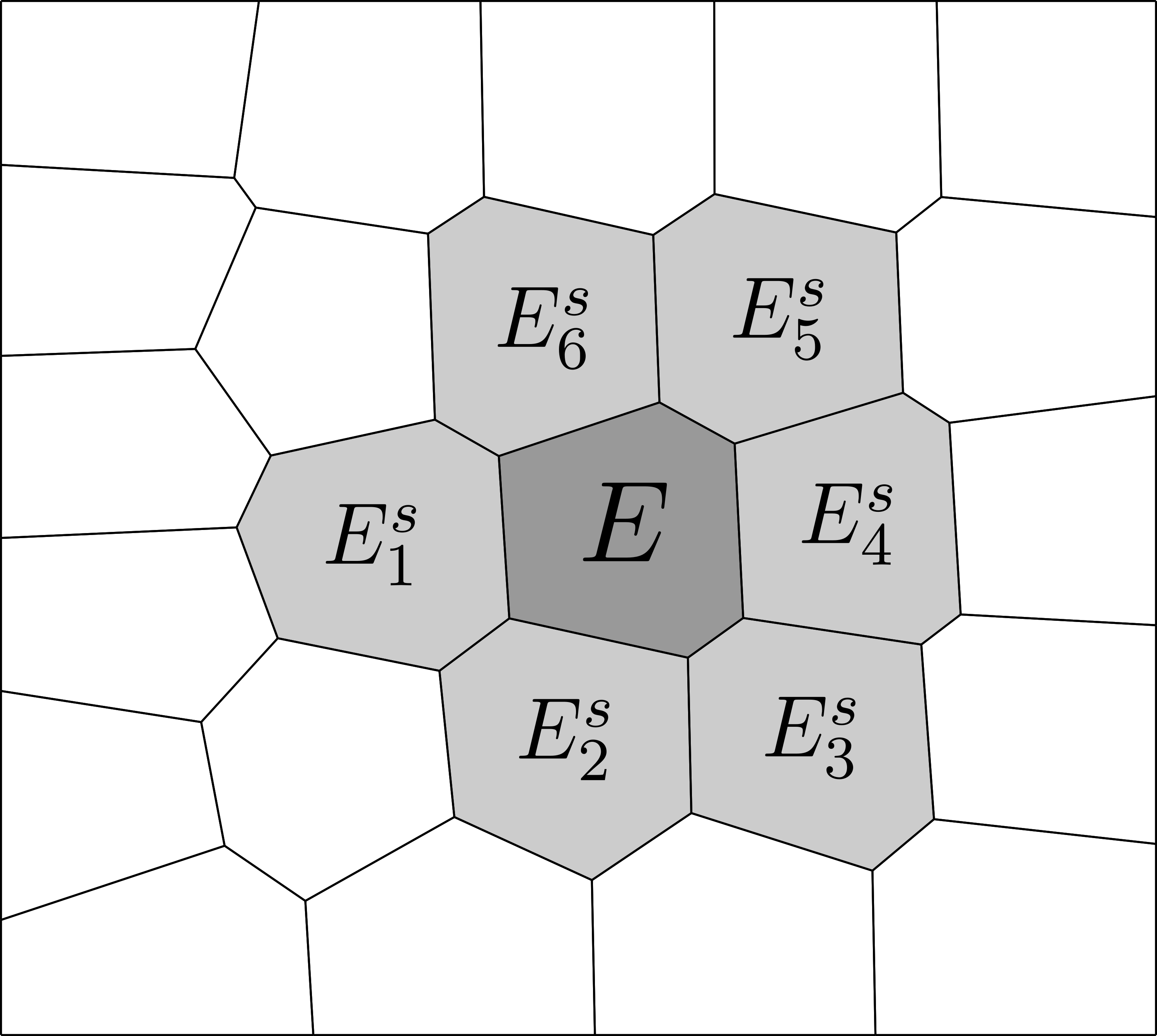}
		\caption{Example of elements $E_{i}^{\rm s}$ surrounding/connected to element $E$. 
			\label{fig:StrainJump}}
	\end{figure}

	The strain jump-based refinement indicator, inspired by minimum residual-based FEM approaches (see for example \cite{chan2013minimum,axelsson1998minimum}), is motivated by seeking to reduce the differences in strains between elements and, thus, to create a smoother approximation of the strain field over the problem domain\footnote{It is noted that a similar indicator in terms stress jumps could be formulated. Such an indicator has been implemented and tested and yields very similar results to the strain jump version. For consistency with \cite{vanHuyssteen2022} it was chosen to present here the strain jump version.}.
	To calculate the components of the strain jump-based indicator $\Upsilon^{jk}_{\rm{SJ},\,E}$ for element $E$ it is first necessary to identify the elements $E_{i}^{\rm s}$ surrounding/connected to $E$.
	An element $E_{i}$ is considered connected to $E$ if the elements share any vertices.
	Figure~\ref{fig:StrainJump} shows a sample mesh with element $E$, depicted in dark grey, connected to 6 surrounding elements $E_{i}^{\rm s}$, depicted in light grey. Elements not connected to $E$ are shown in white.
	
	\FloatBarrier
	\FloatBarrier
	
	The components of the refinement indicator are computed as
	\begin{align}
		\Upsilon_{\text{SJ},\,E}^{jk} &= |E|\, \text{Max} \left[  |\DeltaVarepsilon_{jk}(E_{i}^{\text{s}})| \right] \, , 
	\end{align} 
	where ${\DeltaVarepsilon_{jk}(E_{i}^{\text{s}}) = \varepsilon_{jk}(E)-\varepsilon_{jk}(E_{i}^{\rm s})}$. Furthermore, the strains ${\varepsilon_{jk}(E_{i})}$ are those computed via the projection operator (see (\ref{eqn:ComputeProjection})).
	Finally, the total strain jump indicator is computed as
	\begin{align}
		\Upsilon_{\text{SJ},\,E}^{t} = \| \Upsilon_{\text{SJ},\,E}^{xx} , \, \Upsilon_{\text{SJ},\,E}^{yy} , \, \Upsilon_{\text{SJ},\,E}^{xy} \|_{\mathcal{L}^{2}} \, .
	\end{align}

	\subsection{$Z^{2}$-like indicator} 
	\label{subsec:Z2Indicator}
	The $Z^{2}$-like refinement indicator is inspired by the well-known $Z^{2}$ indicator originally presented in \cite{Zienkiewicz1987}. The indicator involves computing the difference between the element strains (computed via (\ref{eqn:ComputeProjection})) and some higher-order/smoothed approximation of the strains\footnote{It is noted that a similar indicator in terms stresses could be formulated. Such an indicator has been implemented and tested and yields very similar results to the strain version. For consistency with \cite{vanHuyssteen2022} it was chosen to present the strain version here.}.
	
\begin{figure}[ht!]
		\centering
		\includegraphics[width=\textwidth]{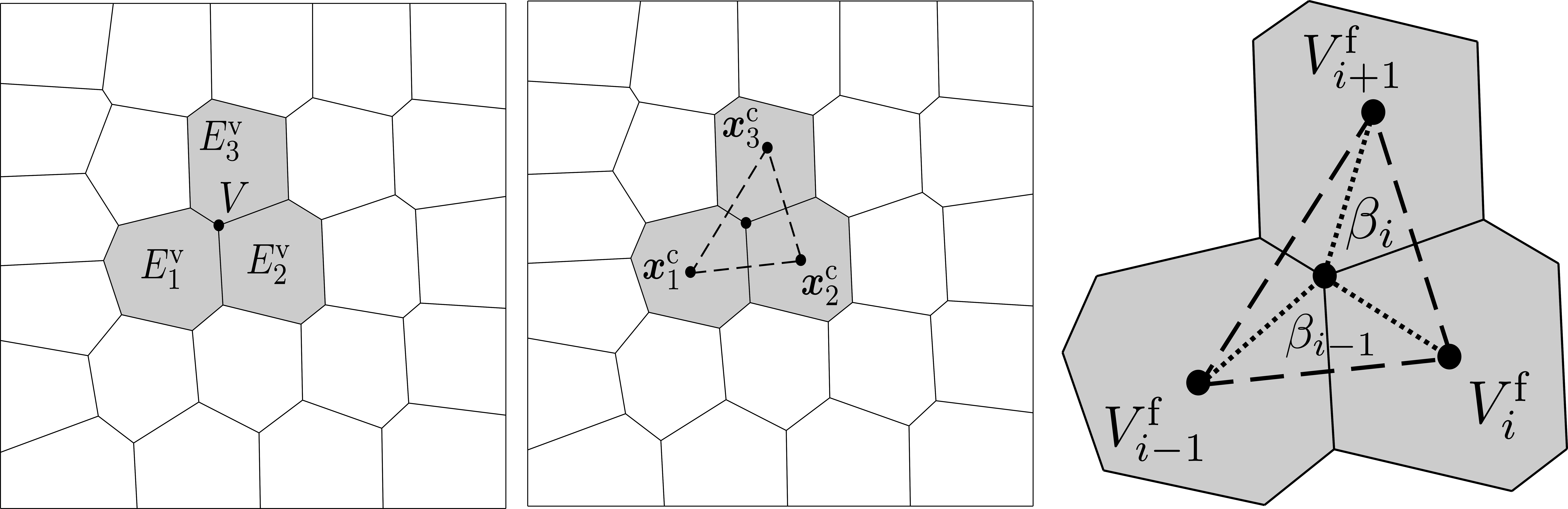}
		\caption{Example of elements $E_{i}^{\rm{v}}$ connected to vertex $V$, and depiction of mean value coordinates.
			\label{fig:Z2}}
	\end{figure}

	Smoothed strains $\bepsilon^{*}$ are computed at each node/vertex in the domain using mean value coordinates \cite{Floater2003}.
	The smoothed strain at vertex $V$ is computed by considering the elements $E^{\rm v}_{i}$ connected to $V$ (see Figure~\ref{fig:Z2}). The centroids $\bx^{\rm c}_{i}$ of the elements $E^{\rm v}_{i}$ are then treated as `fictitious vertices' $V^{\rm f}$ and the strains of element $E^{\rm v}_{i}$, computed via (\ref{eqn:ComputeProjection}), are considered as the degrees of freedom of $V^{\rm f}_{i}$. The smoothed strain at vertex $V$ is then a weighted sum of the strains associated with the fictitious vertices.
	The weight assigned to $V^{\rm f}_{i}$ is given by
	\begin{equation}
		w_{i} = \frac{\theta_{i}}{\sum_{j=1}^{n_{\rm{vf}}} \theta_{j}} \, , 
		\quad \text{where} \quad
		\theta_{i} = \frac{\tan \left(\beta_{i-1} / 2 \right) + \tan \left(\beta_{i} / 2 \right)}{|V - V_{i}^{\rm f}|} \, .
	\end{equation}
	Here $n_{\rm{vf}}$ denotes the number of fictitious vertices, $\beta_{i}$ denotes the angle between the line segments connecting vertices $V$ and $V_{i}^{\rm{f}}$ and vertices $V$ and $V_{i+1}^{\rm{f}}$. Similarly, $\beta_{i-1}$ denotes the angle between the line segments connecting vertices $V$ and $V_{i-1}^{\rm{f}}$ and vertices $V$ and $V_{i}^{\rm{f}}$.
	Additionally, ${|V - V_{i}^{\rm{f}}|}$ is the distance between vertices $V$ and $V_{i}^{\rm{f}}$. The smoothed strain components are then computed as
	\begin{eqnarray}
		\varepsilon_{ij}^{*} = \sum_{k=1}^{n_{\rm{vf}}} w_{k} \, \varepsilon_{ij}(V_{k}^{\rm f}) \, .
	\end{eqnarray}
	
	\FloatBarrier
	\FloatBarrier
	
	The components $\Upsilon_{\rm{Z2}}^{ij}$ of the indicator are computed as 
	\begin{align}
		\Upsilon_{\rm{Z2}}^{ij} &= |E| \sqrt{ \sum_{i=1}^{n_{\rm v}} \left[  \varepsilon_{ij}^{*}(V_{i}^{E}) - (\Pi\,\bu_{h})_{ij}^{E} \right]^{2} } \, .
	\end{align}
	Finally, the total indicator $\Upsilon_{\rm{Z2}}^{t} $ is computed as
	\begin{equation}
		\Upsilon_{\rm{Z2}}^{t} =  \| \Upsilon_{\rm{Z2}}^{xx} , \, \Upsilon_{\rm{Z2}}^{yy}, \, \Upsilon_{\rm{Z2}}^{xy} \|_{\mathcal{L}^{2}} \,.
	\end{equation}
	
	\subsection{Selecting elements to refine} 
	\label{subsec:SelectElementsForRefinement}
	
	The procedure for identifying elements to refine is similar to that presented in \cite{vanHuyssteen2022}.
	Specifically, a refinement threshold percentage ${T=X\%}$ is introduced from which an allowable threshold value ${T_{\text{val}}}$ is determined. 
	To determine ${T_{\text{val}}}$ a domain-level list of the total refinement indicators ${\boldsymbol{\Upsilon}_{\Omega}^{t}}$ is assembled from the element-level indicators such that
	\begin{equation}
		\boldsymbol{\Upsilon}_{\Omega}^{t} = \left\{ \Upsilon_{E_{1}}^{t} , \, \Upsilon_{E_{2}}^{t} , \, \Upsilon_{E_{3}}^{t} ,\, \dots ,\, \Upsilon_{E_{n_{\text{el}}}}^{t}  \right\} \, ,
	\end{equation}
	where $n_{\text{el}}$ is the number of elements in the domain. A reduced list is created by removing duplicate values in $\boldsymbol{\Upsilon}_{\Omega}^{t}$ and this reduced list is sorted in descending order. Finally, the value of the entry ${X\%}$ of the way down the reduced list is found and set as ${T_{\text{val}}}$. Then any element whose indicator $\Upsilon_{E_{i}}^{t}$ is greater than or equal to ${T_{\text{val}}}$ is marked for refinement.
	This approach ensures that ${T_{\text{val}}}$ is always scaled to suit the specific problem.
	Additionally, it is chosen to consider only unique values of $\Upsilon_{E_{i}}^{t}$ to promote meshes that are symmetric when the problem and initial mesh are both also symmetric.
	
	Hereinafter, the term refinement procedure refers to the use of a specific refinement indicator in conjunction with the presented approach to identifying elements qualifying for refinement.
	
	In addition to the displacement-based, strain jump-based, and $Z^{2}$-like refinement indicators, cases of combined refinement indicators will be considered. Specifically, the combination of the displacement-based and strain jump-based indicators, and the combination of the displacement-based and $Z^{2}$-like indicators will be considered. In the case of combined refinement procedures lists of elements marked for refinement are generated for each of the indicators considered as described previously. For example the lists are named List~DB and List~SJ. The total number of elements to refine is chosen to be equal to the length of the shorter of the two lists. The process of creating the combined list of elements to refine (for example, named List~C) comprises two parts. First, the elements that appear on both List~DB and List~SJ are identified and added to List~C. Second, the remaining elements on List~DB and List~SJ are sorted in descending order based on the values of their refinement indicators. Elements are then selected from List~DB and List~SJ in an alternating fashion, and in descending order, and added to List~C until List~C contains the required number of elements.
	
	\section{Numerical Results}
	\label{sec:Results}
	
	In this section numerical results are presented for a range of example problems to demonstrate the efficacy of the various proposed refinement procedures.
	
	To facilitate a fair investigation, a wide range of example problems is considered as some procedures may be more or less suited to a particular problem. The supplementary material contains the full range of six example problems that have been considered (see \cite{DanielHuyssteen2023}). However, for compactness, only three of the problems are presented here. When discussing results over the range of example problems this refers to the full range of six example problems.
	For consistency with the numbering of the problems in the supplementary material, and in \cite{vanHuyssteen2022}, the problems in this work are numbered alphanumerically. For example, the second problem presented is labelled as B(4). The `B' indicates it is the second problem presented in this work, while the (4) indicates that it corresponds to the fourth problem in the supplementary material.
	
	The efficacy of the refinement procedures is evaluated in terms of an $\mathcal{H}^{1}$ error norm. The results generated using the refinement procedures are compared against those generated using a reference procedure in which every element is refined at every refinement step\footnote{To facilitate a fair comparison in terms of run time the same mesh refinement function is used for the adaptive and reference refinement procedures. Additionally, when using the reference procedure no refinement indicators are computed, all elements are automatically marked for refinement, and no element selection process is performed.}. 
	The $\mathcal{H}^{1}$ error norm is defined by 
	\begin{align}
		||\widetilde{\bu}-\bu_{h}||_{1}=&\left[ \int_{\Omega}\left[ |\widetilde{\bu}-\bu_{h}|^{2} + |\nabla \widetilde{\bu} - \nabla \bu_{h}|^{2} \right] \, d\Omega \right]^{0.5} \, , \label{eqn:H1Error}
	\end{align}
	in which integration of $\bu_{h}$ is required over the domain. Since in the case of VEM formulations $\bu_{h}$ is only known on element boundaries, it suffices to approximate (\ref{eqn:H1Error}) as (see \cite{Huyssteen2020,Huyssteen2021})
	\begin{align}
		\begin{split}
			||\widetilde{\bu}-\bu_{h}||_{1} \approx & \left[ \sum_{i=1}^{n_{el}}\frac{|E_{i}|}{n_{\rm v}^{i}} \sum_{j=1}^{n_{\rm v}} \Bigl[ \left[ \widetilde{\bu}(\bx_{j}) - \bu_{h}^{i}(\bx_{j}) \right] \cdot \left[ \widetilde{\bu}(\bx_{j}) -  \bu_{h}^{i}(\bx_{j}) \right] \Bigr. \right.  \\
			&\Biggl. \Bigl. + \left[ \nabla \widetilde{\bu}(\bx_{j}) - \Pi \bu_{h}^{i}(\bx_{j}) \right] : \left[ \nabla \widetilde{\bu}(\bx_{j}) - \Pi \bu_{h}^{i}(\bx_{j}) \right] \Bigr]  \Biggr]^{0.5} \, .
		\end{split} 
	\end{align}
	Here $\widetilde{\bu}$ is a reference solution generated using an overkill mesh of biquadratic finite elements, the location of the $j$-th vertex is denoted by $\bx_{j}$, and ${\Pi \bu_{h}}$ is the gradient of $\bu_{h}$ computed via the projection operator (see (\ref{eqn:ComputeProjection})).
	
	
	In the examples that follow the material is isotropic with a Young's modulus of ${E=1~\rm{Pa}}$. For consistency with \cite{vanHuyssteen2022} a Poisson's ratio corresponding to the case of a compressible material is considered with ${\nu=0.3}$. Additionally, to investigate the performance of the refinement procedures in the case of a nearly incompressible material, results are also presented for a Poisson's ratio of ${\nu=0.49995}$. Finally, the shear modulus is computed as ${\mu = E /2 \left[1+\nu\right]}$. 
	
	For each problem results have been generated for the cases of both structured and unstructured/Voronoi meshes for compressible and nearly incompressible Poisson's ratios. For compactness, a subset of these results will be presented here. The full set of results can be found in the supplementary material \cite{DanielHuyssteen2023}.
	

	\subsection{Problem A(1): Plate with square hole - Prescribed traction \\ \cgray Displacement-based refinement procedure \cblack}
	\label{subsec:PlateWithHoleTraction}
	
\begin{figure}[ht!]
		\centering
		\begin{subfigure}[t]{0.45\textwidth}
			\centering
			\includegraphics[width=0.95\textwidth]{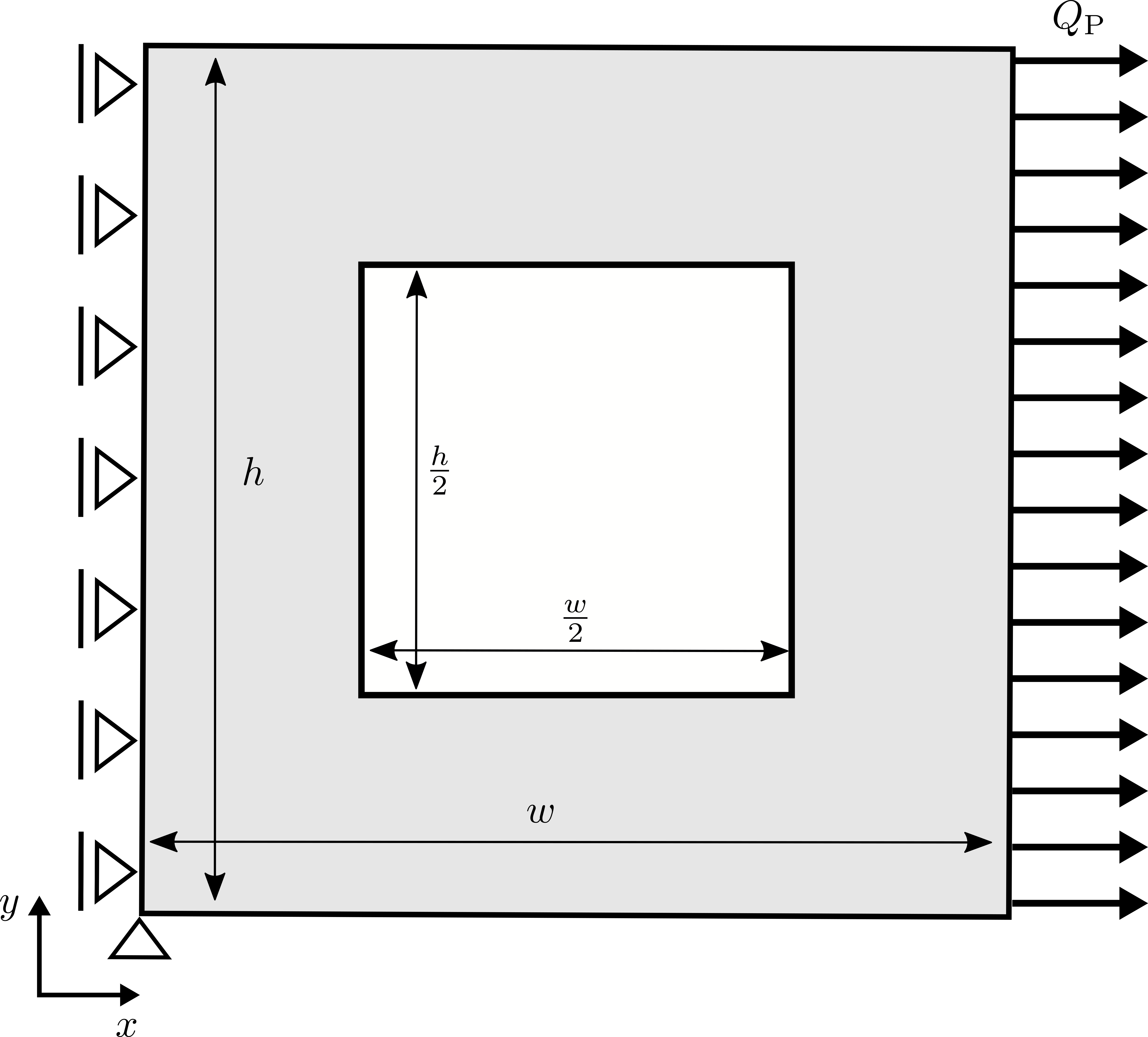}
			\caption{Problem geometry}
		\end{subfigure}%
		\begin{subfigure}[t]{0.55\textwidth}
			\centering
			\includegraphics[width=0.95\textwidth]{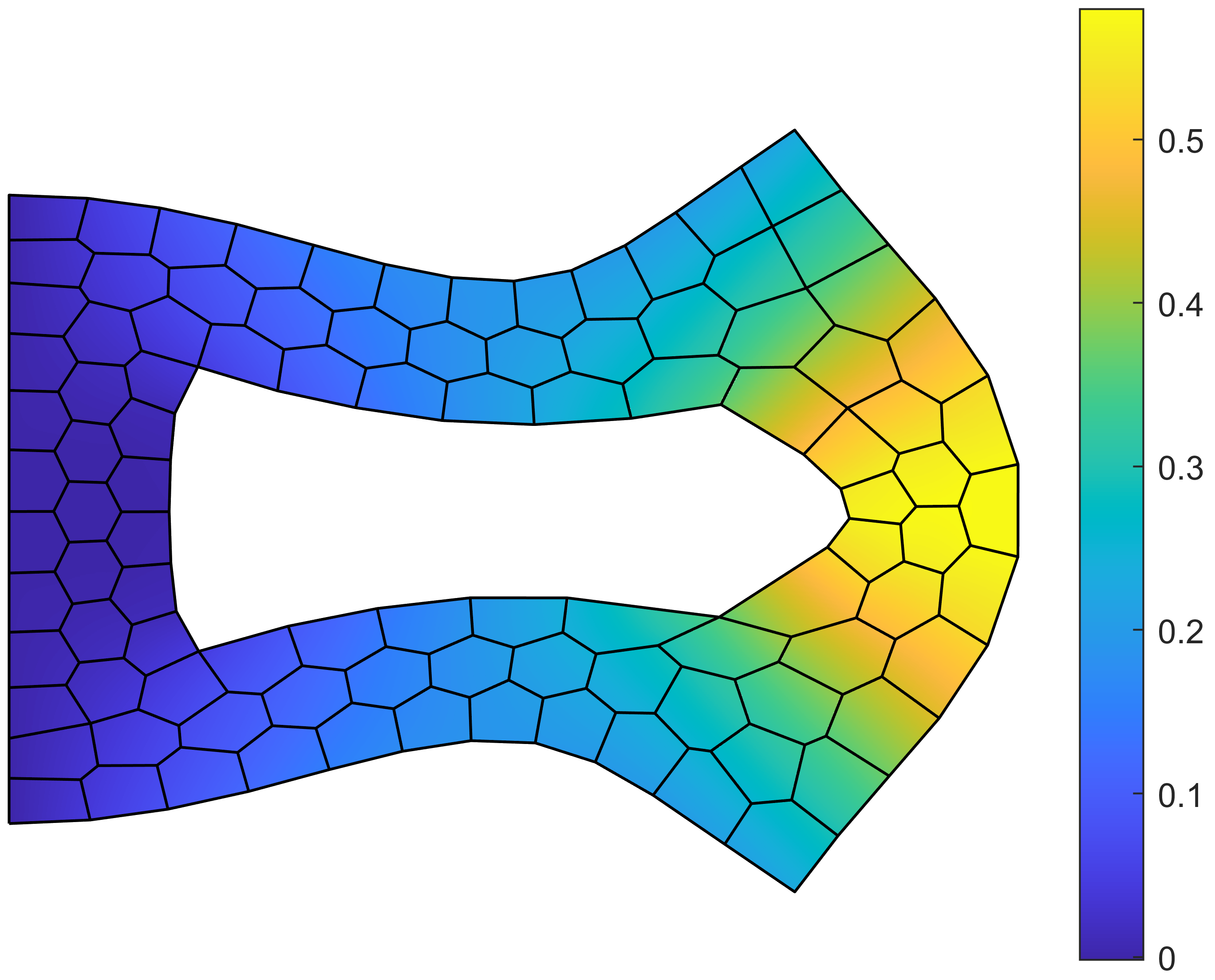}
			\caption{Deformed configuration}
		\end{subfigure}
		\caption{Problem 1 (a) geometry, and (b) deformed configuration of a Voronoi mesh with ${\nu=0.3}$. 
			\label{fig:PlateWithHoleTractionGeometry}}
	\end{figure}

	Problem~A(1) comprises a plate of width ${w=1~\rm{m}}$ and height ${h=1~\rm{m}}$ with a centrally located hole (see Figure~\ref{fig:PlateWithHoleTractionGeometry}(a)). The left-hand edge of the plate is constrained horizontally and the bottom left-hand corner is fully constrained. The right-hand edge is subject to a prescribed traction of ${Q_{\rm P}=0.2~\frac{\rm N}{\rm m}}$. The results presented for this problem were generated using the displacement-based refinement procedure. Figure~\ref{fig:PlateWithHoleTractionGeometry}(b) depicts a sample deformed configuration of the plate with a Voronoi mesh and a compressible Poisson's ratio with ${\nu=0.3}$. The horizontal displacement ${u_{x}}$ is plotted on the colour axis.
	
	\FloatBarrier
	\FloatBarrier
	
\begin{figure}[ht!]
		\centering
		\begin{subfigure}[t]{0.33\textwidth}
			\centering
			\includegraphics[width=0.95\textwidth]{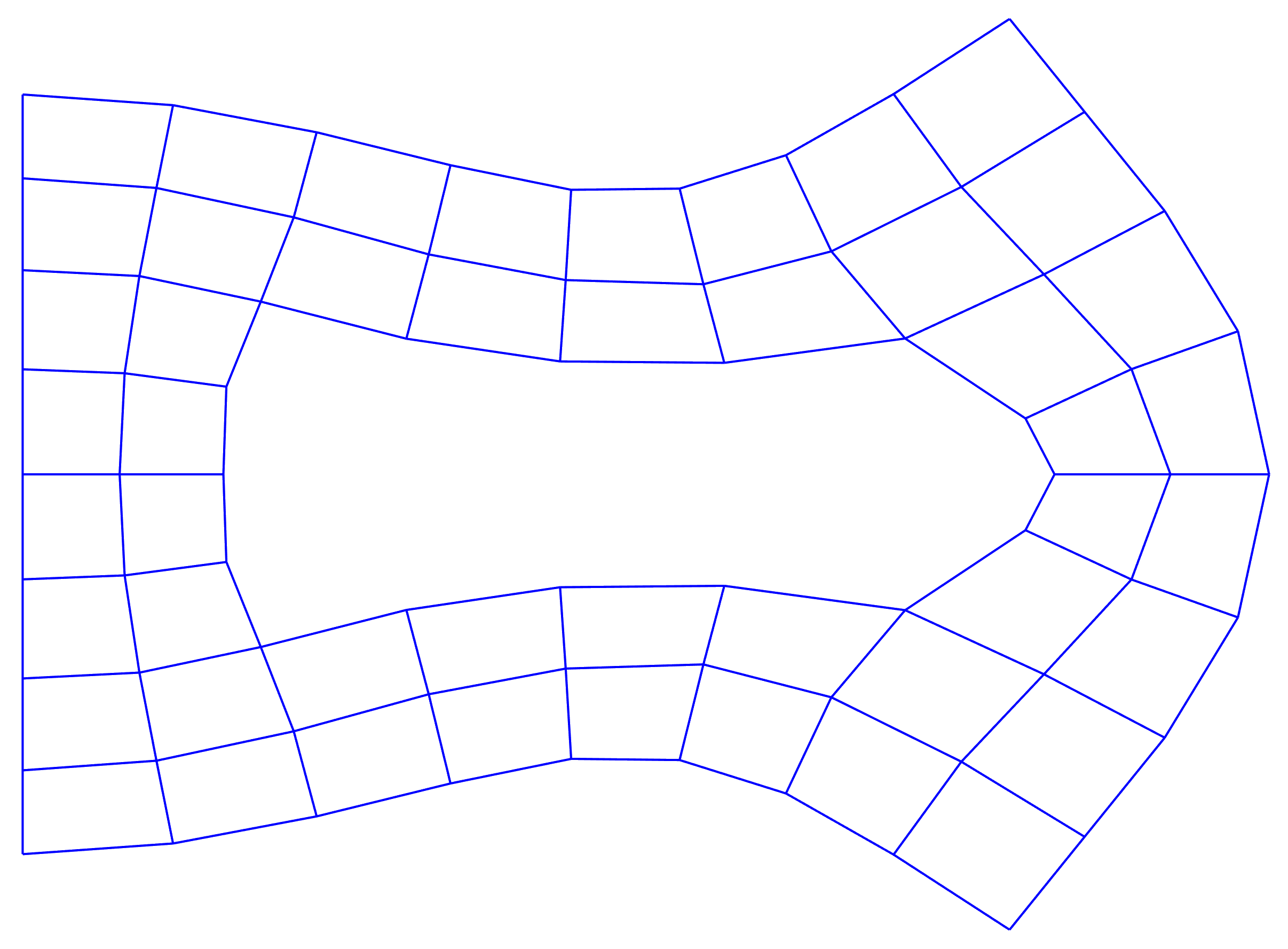}
			\caption{Structured - Step 1}
		\end{subfigure}%
		\begin{subfigure}[t]{0.33\textwidth}
			\centering
			\includegraphics[width=0.95\textwidth]{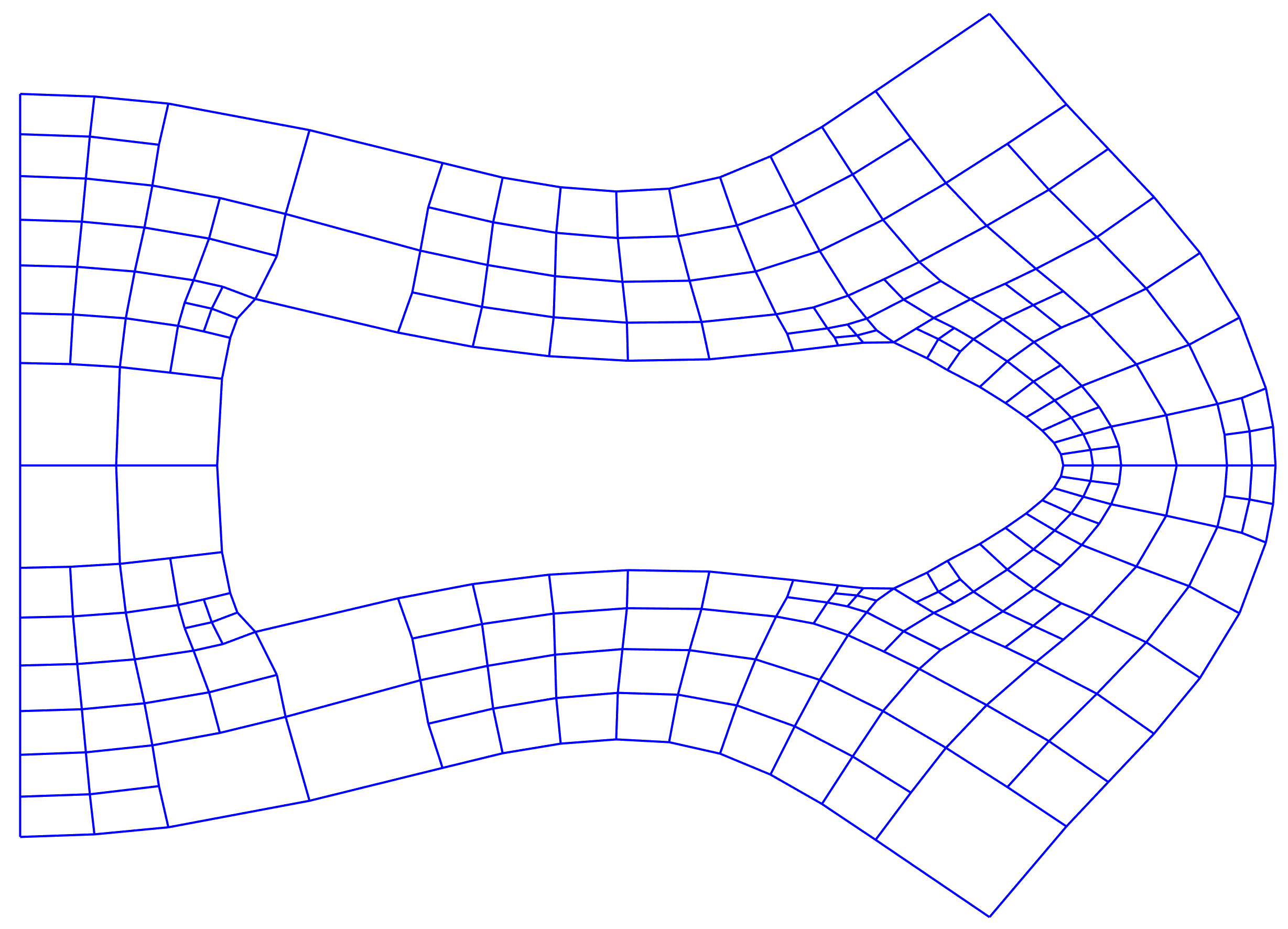}
			\caption{Structured - Step 4}
		\end{subfigure}%
		\begin{subfigure}[t]{0.33\textwidth}
			\centering
			\includegraphics[width=0.95\textwidth]{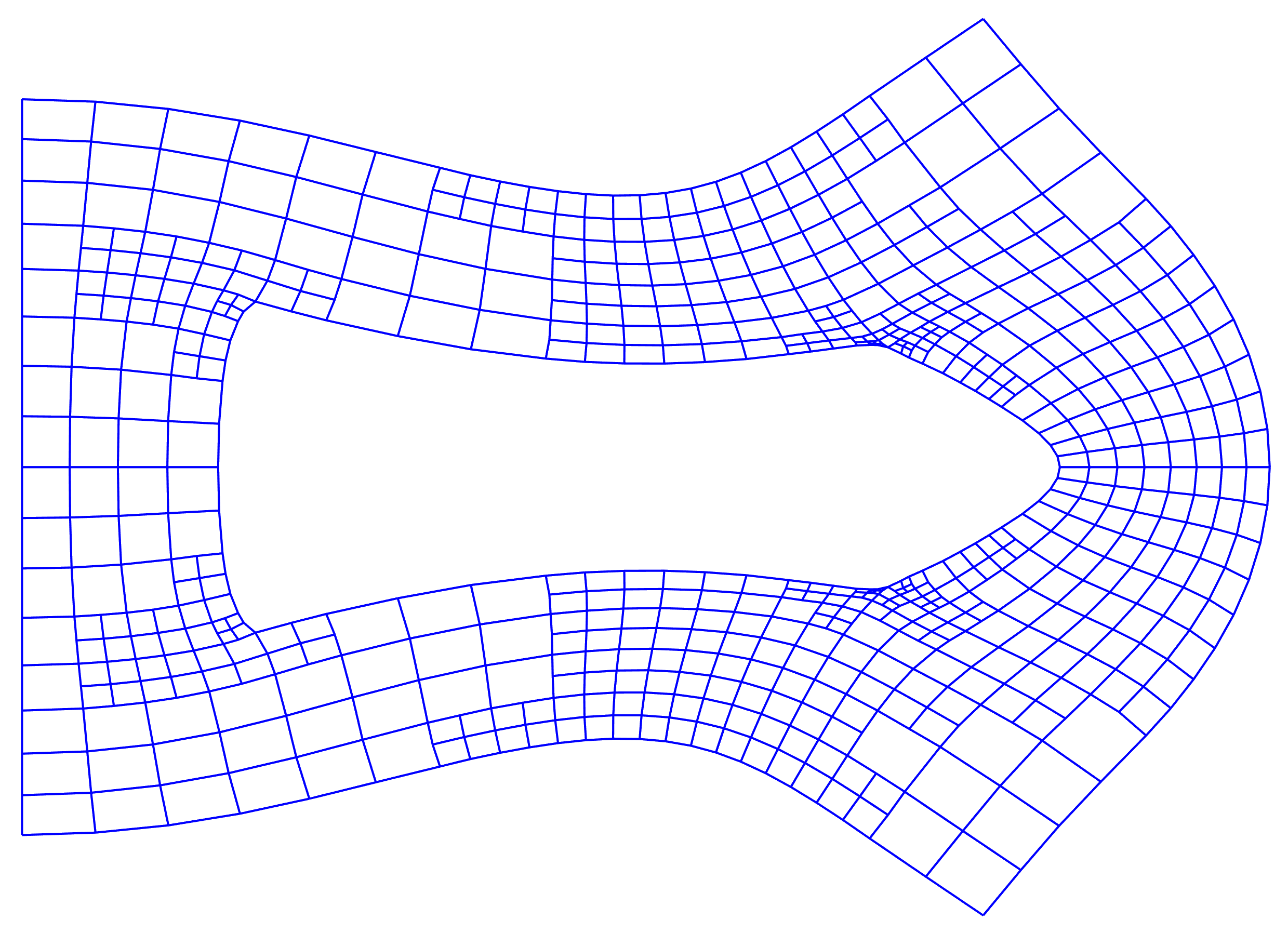}
			\caption{Structured - Step 6}
		\end{subfigure}
		\vskip \baselineskip 
		\begin{subfigure}[t]{0.33\textwidth}
			\centering
			\includegraphics[width=0.95\textwidth]{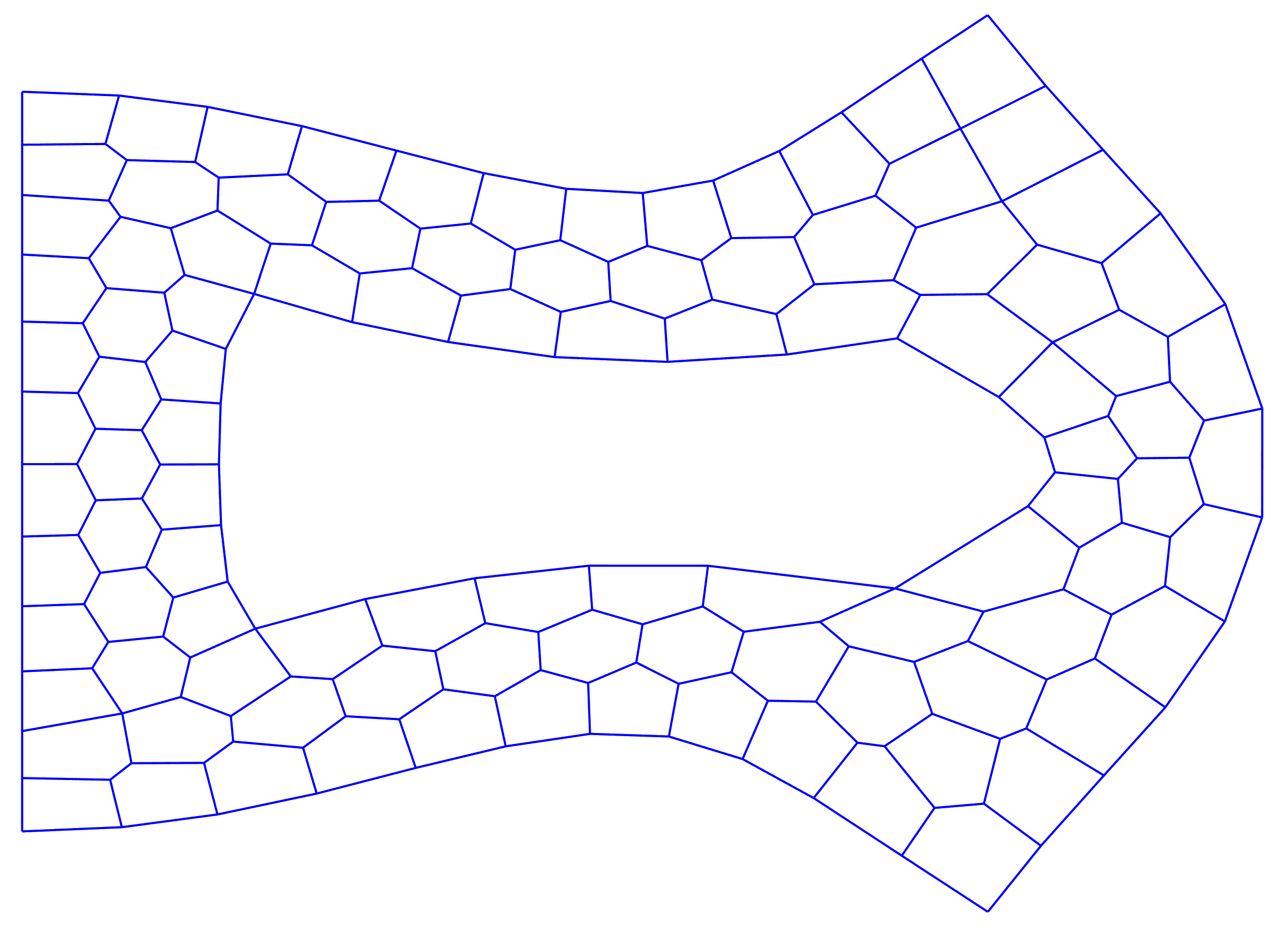}
			\caption{Voronoi - Step 1}
		\end{subfigure}%
		\begin{subfigure}[t]{0.33\textwidth}
			\centering
			\includegraphics[width=0.95\textwidth]{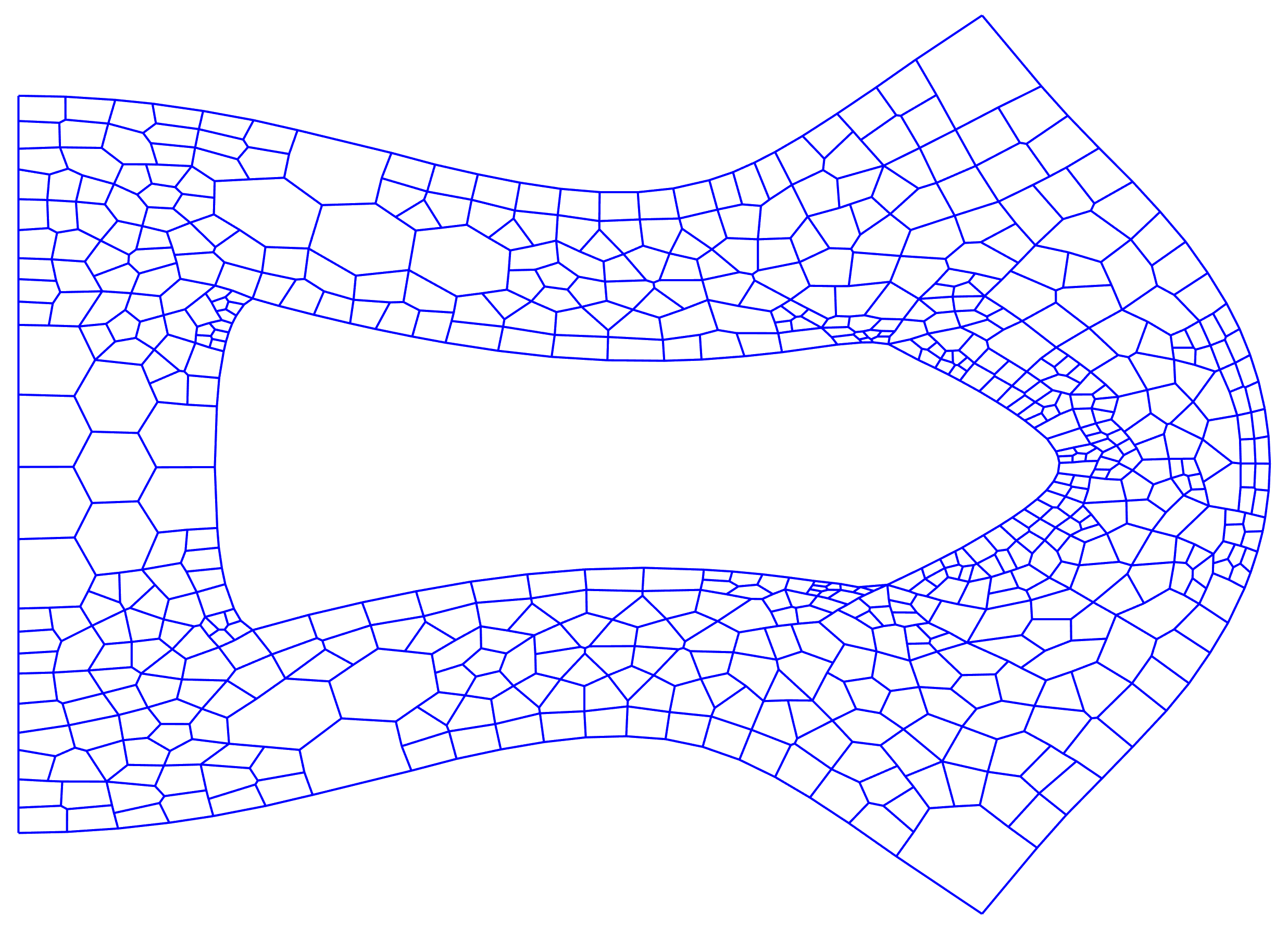}
			\caption{Voronoi - Step 4}
		\end{subfigure}%
		\begin{subfigure}[t]{0.33\textwidth}
			\centering
			\includegraphics[width=0.95\textwidth]{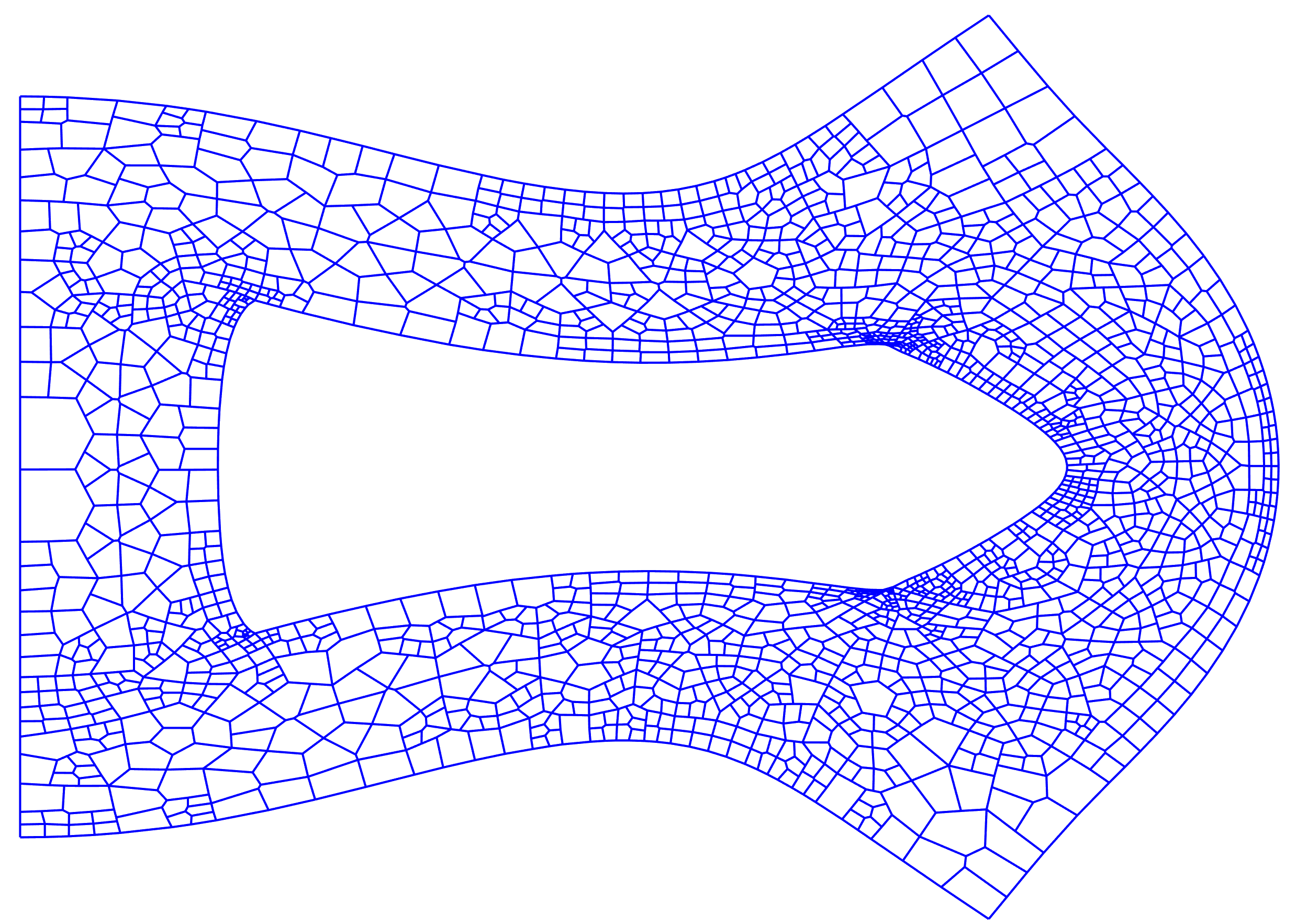}
			\caption{Voronoi - Step 6}
		\end{subfigure}
		\caption{Mesh refinement process for problem~A(1) using the displacement-based refinement procedure with ${T=20\%}$ for structured and Voronoi meshes with a compressible Poisson's ratio.
			\label{fig:PlateWithHoleTractionMeshes}}
	\end{figure}

	Figure~\ref{fig:PlateWithHoleTractionMeshes} depicts the mesh refinement process for problem~A(1) using the displacement-based refinement procedure with ${T=20\%}$ for structured and Voronoi meshes with a compressible Poisson's ratio. Meshes are shown at various refinement steps with step~1 corresponding to the initial mesh.
	The evolution of the structured and Voronoi meshes are qualitatively similar. For both mesh types regions of increased refinement form around the corners of the hole and in regions undergoing the most severe deformation, i.e. the right-hand portion of the domain which experiences large compressions and tensions. Furthermore, the meshes remain coarser in the regions experiencing simpler deformations, i.e. the bottom and top right-hand corners of the domain, in which the deformation is largely due to translation, and the middle of the left-hand edge which undergoes very little deformation. Thus, the evolution of both mesh types is sensible for this example problem.
	
	\FloatBarrier
	\FloatBarrier
	
\begin{figure}[ht!]
		\centering
		\begin{subfigure}[t]{0.5\textwidth}
			\centering
			\includegraphics[width=0.95\textwidth]{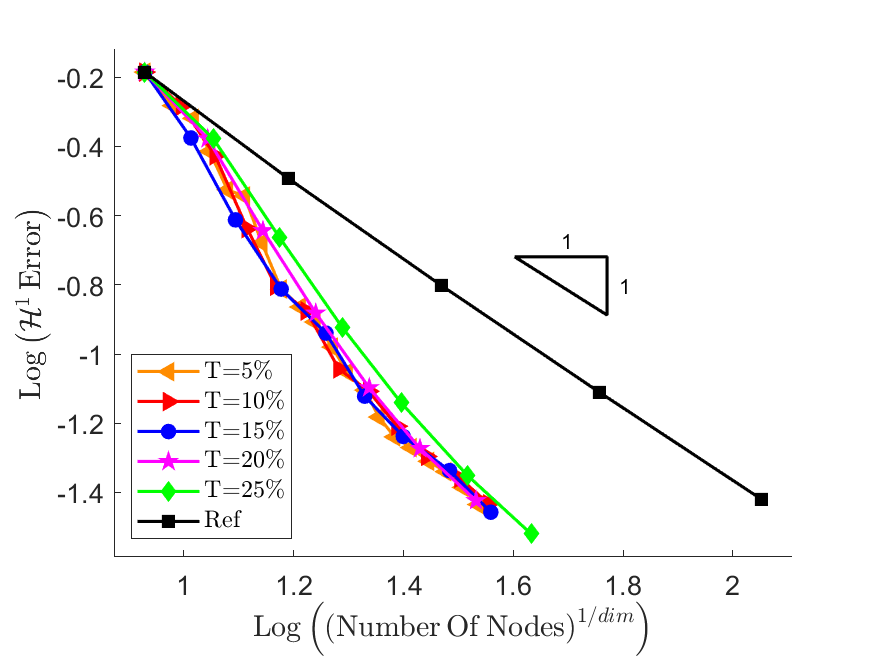}
			\caption{Structured mesh}
		\end{subfigure}%
		\begin{subfigure}[t]{0.5\textwidth}
			\centering
			\includegraphics[width=0.95\textwidth]{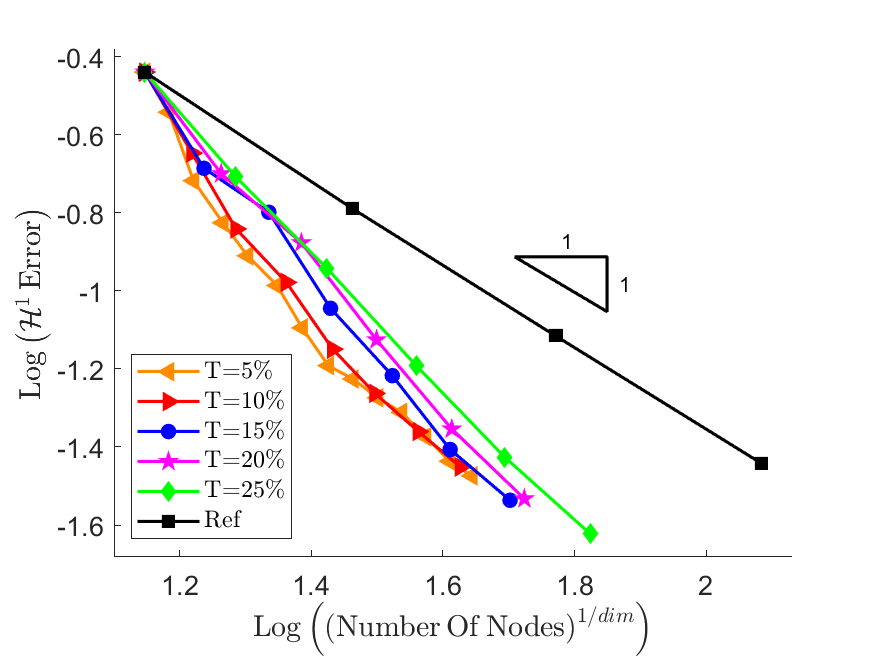}
			\caption{Voronoi mesh}
		\end{subfigure}
		\caption{$\mathcal{H}^{1}$ error vs $n_{\rm v}$ for problem~A(1) using the displacement-based refinement procedure with a variety of choices of $T$ for structured and Voronoi meshes with a compressible Poisson's ratio.
			\label{fig:PlateWithHoleTractionConvergenceNumberOfNodes}}
	\end{figure}

	The convergence behaviour in the ${\mathcal{H}^{1}}$ error norm of the VEM for problem~A(1) using the displacement-based refinement procedure is depicted in Figure~\ref{fig:PlateWithHoleTractionConvergenceNumberOfNodes} on a logarithmic scale for a variety of choices of $T$. Here, the ${\mathcal{H}^{1}}$ error is plotted against the number of vertices/nodes in the discretization for cases of structured and Voronoi meshes with a compressible Poisson's ratio. Additionally, each marker on the curves depicted in Figure~\ref{fig:PlateWithHoleTractionConvergenceNumberOfNodes} corresponds to a refinement step.
	The convergence behaviour is qualitatively similar for all choices of $T$. However, particularly in the case of Voronoi meshes, for lower choices of $T$ the convergence rate is initially slightly faster and decreases marginally as the number of nodes increases and matches the convergence rate exhibited by the other choices of $T$. For larger choices of $T$ the convergence rate is more consistent throughout the figure domain. For all choices of $T$ the adaptive refinement procedure exhibits a superior convergence rate to, and significantly outperforms, the reference refinement procedure. That is, the same level of accuracy as the reference procedure is achieved by the adaptive procedure while using significantly fewer nodes than the reference procedure.
	
	\FloatBarrier
	\FloatBarrier
	
	The performance of the VEM in terms of its convergence behaviour in the $\mathcal{H}^{1}$ error norm with respect to the number of vertices/nodes in the discretization when using the displacement-based refinement procedure for problem~A(1) is summarized in Table~\ref{tab:PerformancePlateWithHoleTractionNumberOfNodes}. Here, the performance of the adaptive procedure, as measured by the percentage relative effort, is presented for the cases of structured and Voronoi meshes with compressible and nearly incompressible Poisson's ratios for a variety of choices of $T$. 
	
	The performance of the adaptive procedure is compared against that of the reference procedure by determining, via interpolation, how many nodes are required to reach the same accuracy (i.e. $\mathcal{H}^{1}$ error) as the reference procedure with its finest mesh.
	The percentage relative effort (PRE) quantifies what proportion of the computational resources (in this case, number of nodes) used by the reference procedure is required by the adaptive procedure to reach the same level of accuracy. 
	For clarity, an example calculation is presented for the case of a structured mesh with a compressible Poisson's ratio for ${T=20\%}$.
	\begin{align}
		\text{PRE}\left(nNodes \right)|_{T=20\%} &= \frac{nNodes|_{T=20\%}}{nNodes|_{Ref}} =  \frac{1142.64}{12672} = 9.02\% \, .
	\end{align}
	Thus, the same level of accuracy as the reference procedure is achieved by the displacement-based refinement procedure with ${T=20\%}$ while using ${9.02\%}$ of the number of vertices/nodes used by the reference procedure.
	
	The general trend observed in Table~\ref{tab:PerformancePlateWithHoleTractionNumberOfNodes} for the case of a compressible Poisson's ratio is that a lower choice of $T$ yields improved performance. This behaviour is expected as refining a smaller number of elements at each step allows the refinement to be more targeted/localized in critical areas of a domain, for example at a corner/notch. However, a lower choice of $T$ has the consequence of requiring more refinement steps to achieve a particular level of accuracy (see the number of markers on the curves depicted in Figure~\ref{fig:PlateWithHoleTractionConvergenceNumberOfNodes}). 
	
	In the case of near-incompressibility the choice of $T$ has significantly less influence on the PRE than in the compressible case with all choices of $T$ exhibiting a similar PRE. This is because in the case of near-incompressibility the $\mathcal{H}^{1}$ error is distributed over a larger area than in the compressible case. Thus, to reach the desired level of accuracy, refinement is required over a larger area and the more targeted/localized refinement offered by low choices of $T$ is less effective. Nevertheless, for both mesh types, both choices of Poisson's ratio, and for all choices of $T$, the adaptive procedure significantly outperforms the reference procedure.
	
	\FloatBarrier
	\begin{table}[ht!]
		\centering 
		\begin{adjustbox}{max width=\textwidth}
			\begin{tabular}{|c|cccc|cccc|}
				\hline
				\multirow{3}{*}{Threshold} & \multicolumn{4}{c|}{Compressible}                                                                  & \multicolumn{4}{c|}{Nearly-incompressible}                                                         \\ \cline{2-9} 
				& \multicolumn{2}{c|}{Structured}                            & \multicolumn{2}{c|}{Voronoi}          & \multicolumn{2}{c|}{Structured}                            & \multicolumn{2}{c|}{Voronoi}          \\ \cline{2-9} 
				& \multicolumn{1}{c|}{nNodes}   & \multicolumn{1}{c|}{PRE}   & \multicolumn{1}{c|}{nNodes}   & PRE   & \multicolumn{1}{c|}{nNodes}   & \multicolumn{1}{c|}{PRE}   & \multicolumn{1}{c|}{nNodes}   & PRE   \\ \hline
				T=5\%                      & \multicolumn{1}{c|}{1130.77}  & \multicolumn{1}{c|}{8.92}  & \multicolumn{1}{c|}{1696.59}  & 11.56 & \multicolumn{1}{c|}{4109.83}  & \multicolumn{1}{c|}{32.43} & \multicolumn{1}{c|}{3878.80}  & 26.67 \\ \hline
				T=10\%                     & \multicolumn{1}{c|}{1233.28}  & \multicolumn{1}{c|}{9.73}  & \multicolumn{1}{c|}{1736.73}  & 11.83 & \multicolumn{1}{c|}{4161.85}  & \multicolumn{1}{c|}{32.84} & \multicolumn{1}{c|}{3842.04}  & 26.42 \\ \hline
				T=15\%                     & \multicolumn{1}{c|}{1176.53}  & \multicolumn{1}{c|}{9.28}  & \multicolumn{1}{c|}{1866.84}  & 12.72 & \multicolumn{1}{c|}{4135.27}  & \multicolumn{1}{c|}{32.63} & \multicolumn{1}{c|}{3847.68}  & 26.46 \\ \hline
				T=20\%                     & \multicolumn{1}{c|}{1142.64}  & \multicolumn{1}{c|}{9.02}  & \multicolumn{1}{c|}{2162.95}  & 14.74 & \multicolumn{1}{c|}{4081.88}  & \multicolumn{1}{c|}{32.21} & \multicolumn{1}{c|}{3862.79}  & 26.56 \\ \hline
				T=25\%                     & \multicolumn{1}{c|}{1338.61}  & \multicolumn{1}{c|}{10.56} & \multicolumn{1}{c|}{2559.68}  & 17.44 & \multicolumn{1}{c|}{4026.12}  & \multicolumn{1}{c|}{31.77} & \multicolumn{1}{c|}{3877.02}  & 26.66 \\ \hline
				Ref                        & \multicolumn{1}{c|}{12672.00} & \multicolumn{1}{c|}{}      & \multicolumn{1}{c|}{14675.00} &       & \multicolumn{1}{c|}{12672.00} & \multicolumn{1}{c|}{}      & \multicolumn{1}{c|}{14542.00} &       \\ \hline
			\end{tabular}
		\end{adjustbox}
		\caption{Performance summary of the VEM in terms of its convergence behaviour in the $\mathcal{H}^{1}$ error norm with respect to the number of vertices/nodes in the discretization when using the displacement-based refinement procedure for problem~A(1).
			\label{tab:PerformancePlateWithHoleTractionNumberOfNodes}}
	\end{table}
	\FloatBarrier
	
\begin{figure}[ht!]
		\centering
		\begin{subfigure}[t]{0.5\textwidth}
			\includegraphics[width=0.95\textwidth]{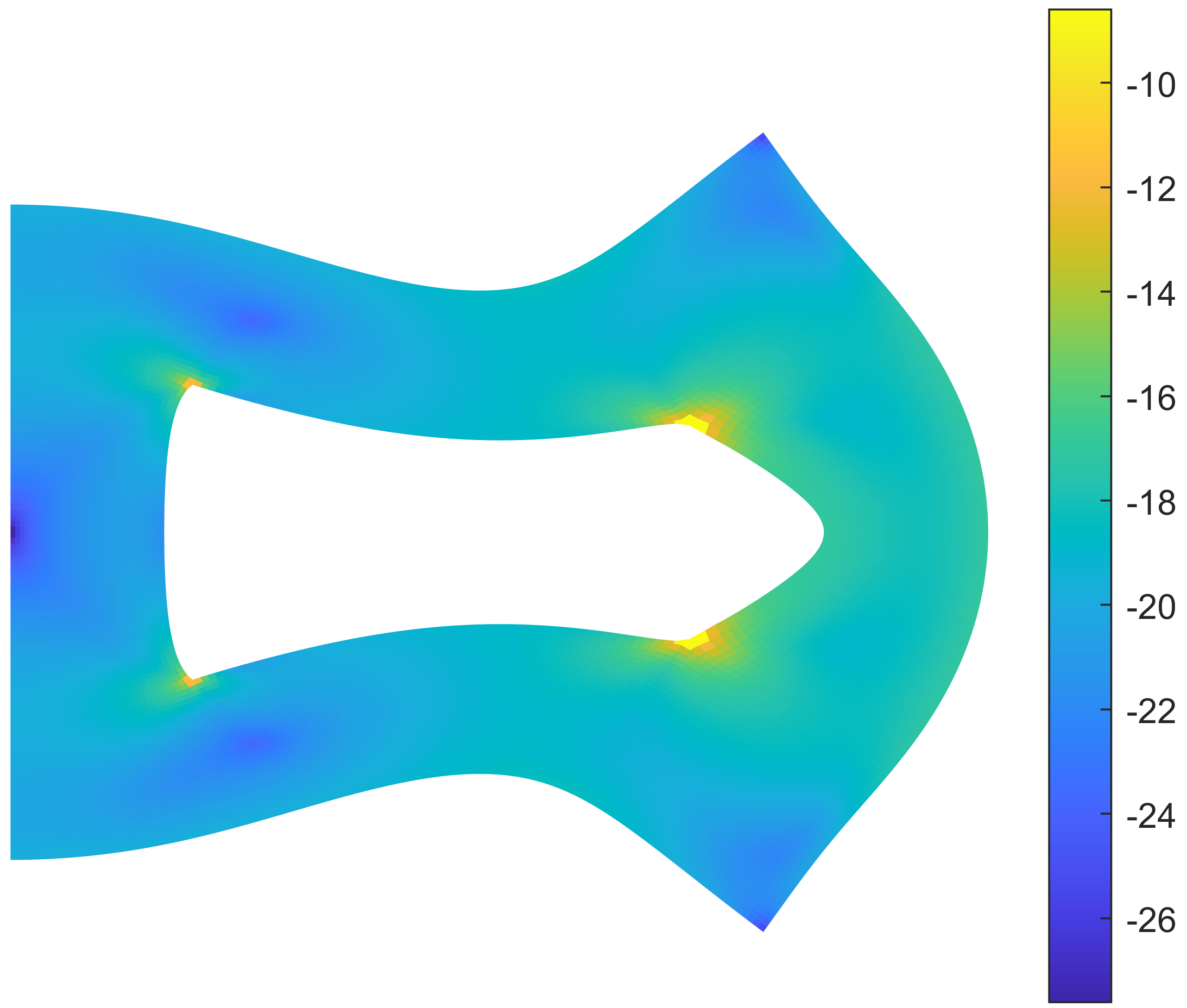}
			\caption{$\mathcal{H}^{1}$ error - Compressible - Reference procedure: Step 5}
		\end{subfigure}%
		\begin{subfigure}[t]{0.5\textwidth}
			\includegraphics[width=0.95\textwidth]{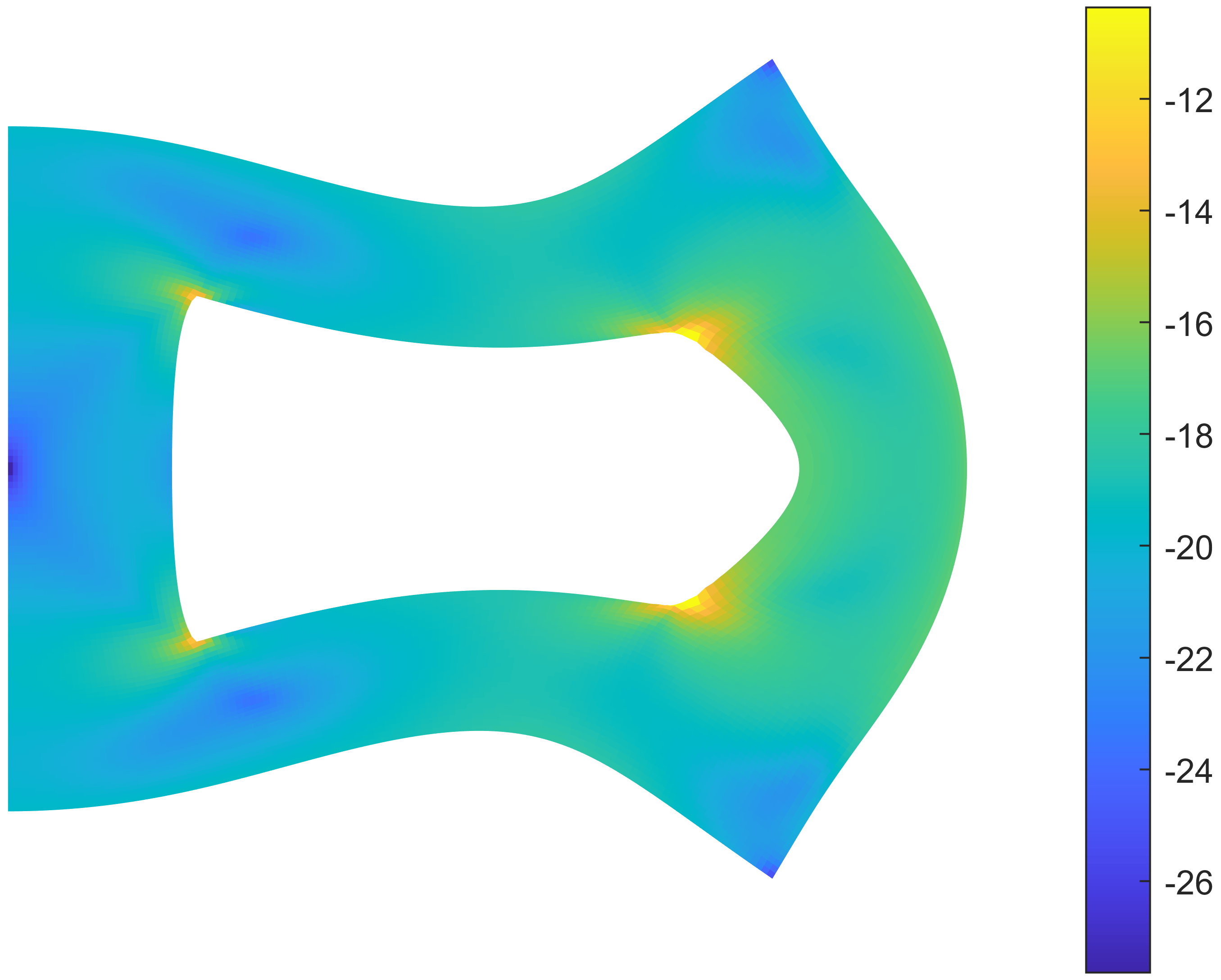}
			\caption{$\mathcal{H}^{1}$ error - Nearly incompressible - Reference procedure: Step 5}
		\end{subfigure}
		\vskip \baselineskip 
		\begin{subfigure}[t]{0.5\textwidth}
			\includegraphics[width=0.79\textwidth]{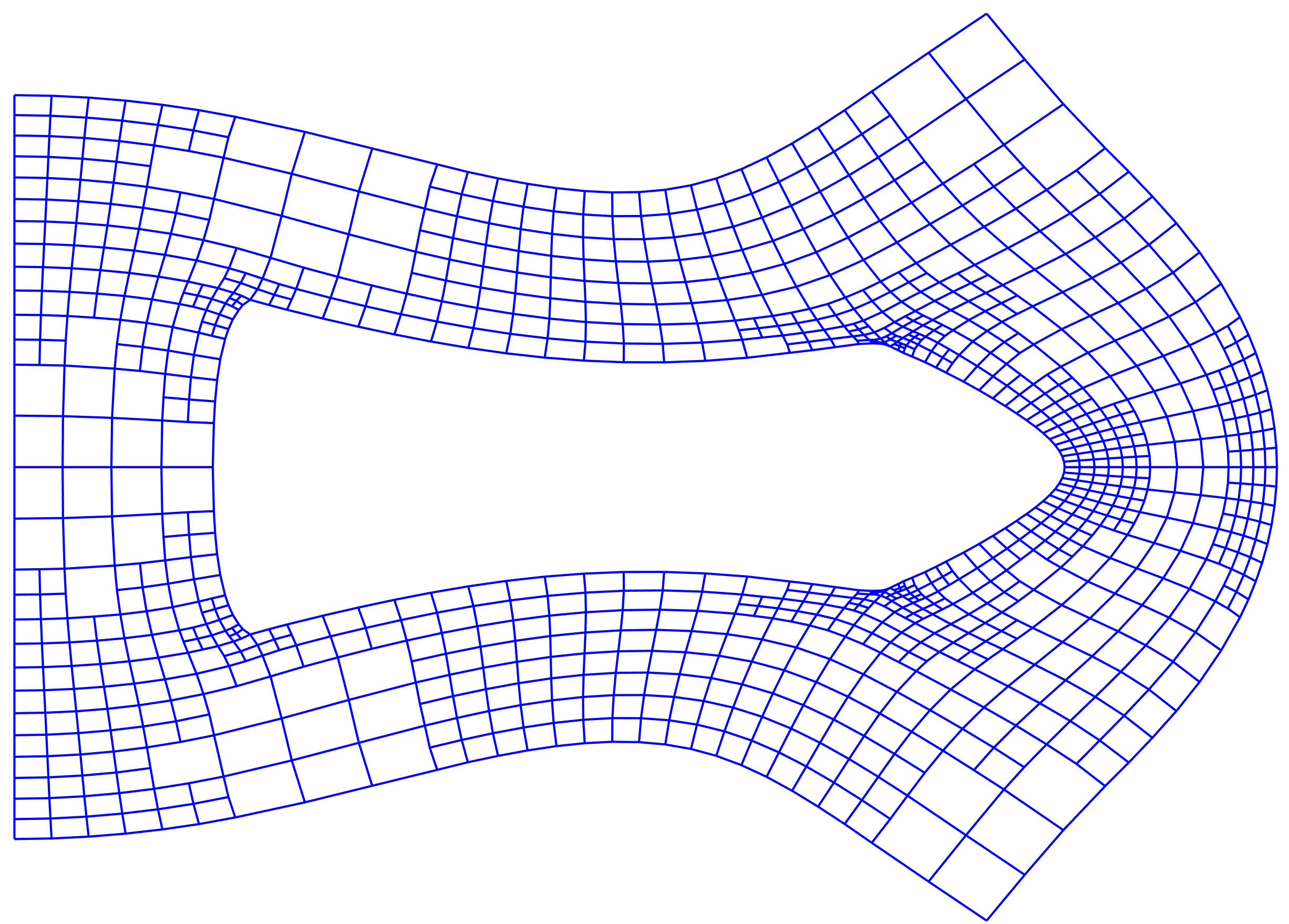}
			\caption{Adaptively refined mesh - Compressible - Step 7}
		\end{subfigure}%
		\begin{subfigure}[t]{0.5\textwidth}
			\includegraphics[width=0.79\textwidth]{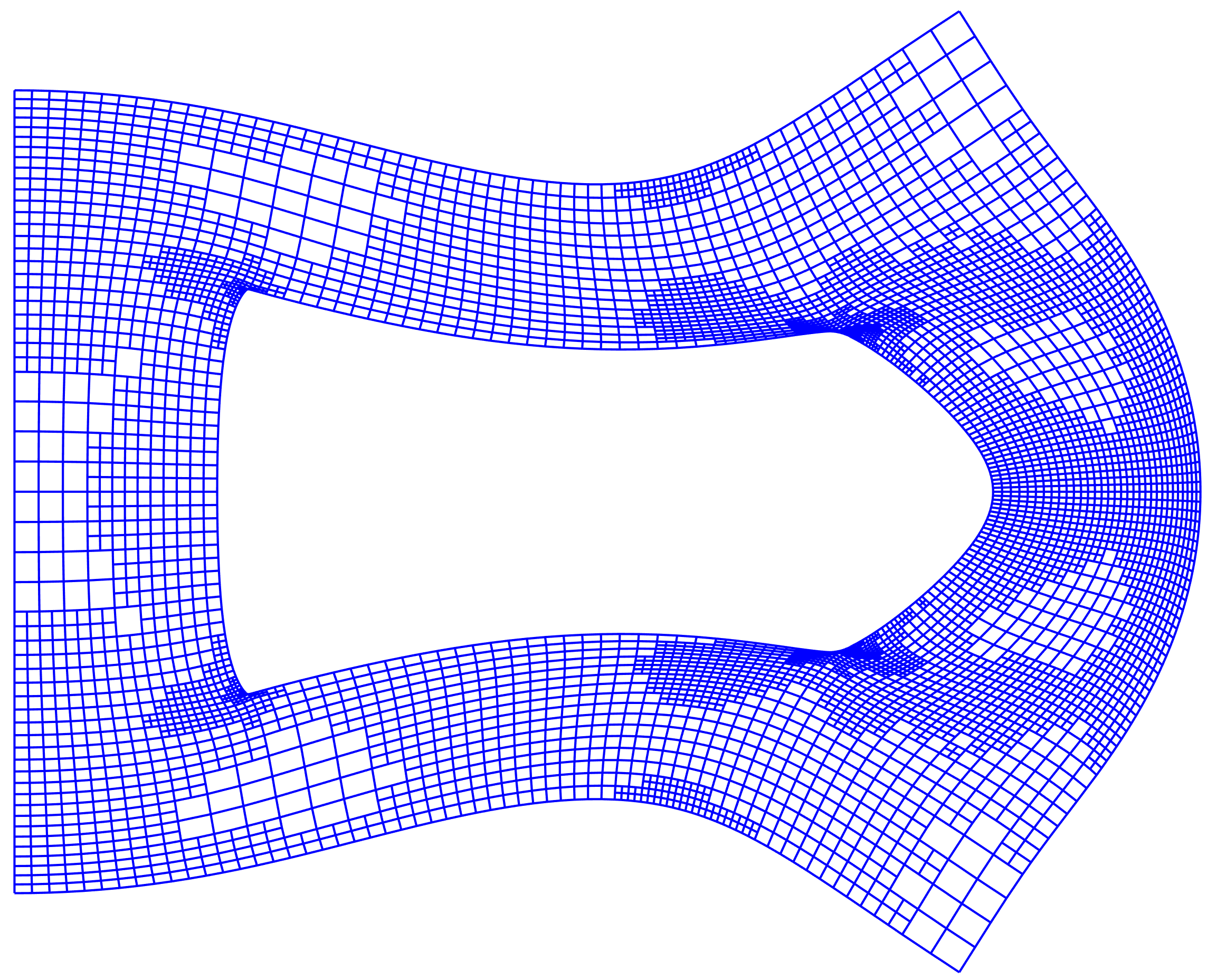}
			\caption{Adaptively refined mesh - Nearly incompressible - Step 11}
		\end{subfigure}
		\caption{$\mathcal{H}^{1}$ error distribution for the reference procedure and final deformed meshes (i.e. the adaptively refined mesh that has the same accuracy as the reference solution) for problem~A(1) using the displacement-based refinement procedure with ${T=20\%}$ for structured meshes with compressible and nearly incompressible Poisson's ratios.
			\label{fig:PlateWithHoleTractionEffectOfPoisson}}
	\end{figure}

	The influence of compressibility on the error distribution and on the adaptively generated meshes is demonstrated in Figure~\ref{fig:PlateWithHoleTractionEffectOfPoisson} for structured meshes. 
	The top row depicts the distribution of the $\mathcal{H}^{1}$ error (in log scale) for the reference procedure with its finest mesh. 
	The bottom row depicts the final adaptively generated meshes (i.e. the adaptively refined mesh that has the same accuracy as the reference solution) using the displacement-based refinement procedure with ${T=20\%}$.
	The left-hand and right-hand columns correspond respectively to compressible and nearly incompressible Poisson's ratios.
	The difference in the error distribution in the cases of compressibility and near-incompressibility can be seen in the top row of figures. Here, the regions of high error around the corners of the hole extend slightly further in the nearly incompressible case than in the compressible case. Furthermore, the (relative) error in the nearly incompressible case is higher than that of the compressible case in the regions experiencing larger deformations, such as the right-hand portion of the domain. Overall, a larger portion of the domain experiences relatively higher error in the nearly incompressible case than in the compressible case. That is, the heat map is on average more green (indicating higher error) in the nearly incompressible case and more blue (indicating lower error) in the compressible case.
	The differences in the adaptively generated meshes in the cases of compressibility and near-incompressibility closely reflect the differences in the error distribution. In both cases there is increased refinement around the corners of the hole and in the areas of the domain experiencing relatively large deformations. However, although the adaptively generated meshes are qualitatively similar, in the nearly incompressible case the overall level of refinement is greater than that of the compressible case as a result of the larger distribution, and on average relatively higher magnitude, of the error.
	In both cases the regions of increased and decreased refinement coincide respectively with regions of higher and lower error. Thus, and as discussed in Figure~\ref{fig:PlateWithHoleTractionMeshes}, the adaptive refinement procedure generates sensible meshes for the problem.
	
	\FloatBarrier
	\FloatBarrier
	
\begin{figure}[ht!]
		\centering
		\begin{subfigure}[t]{0.5\textwidth}
			\centering
			\includegraphics[width=0.95\textwidth]{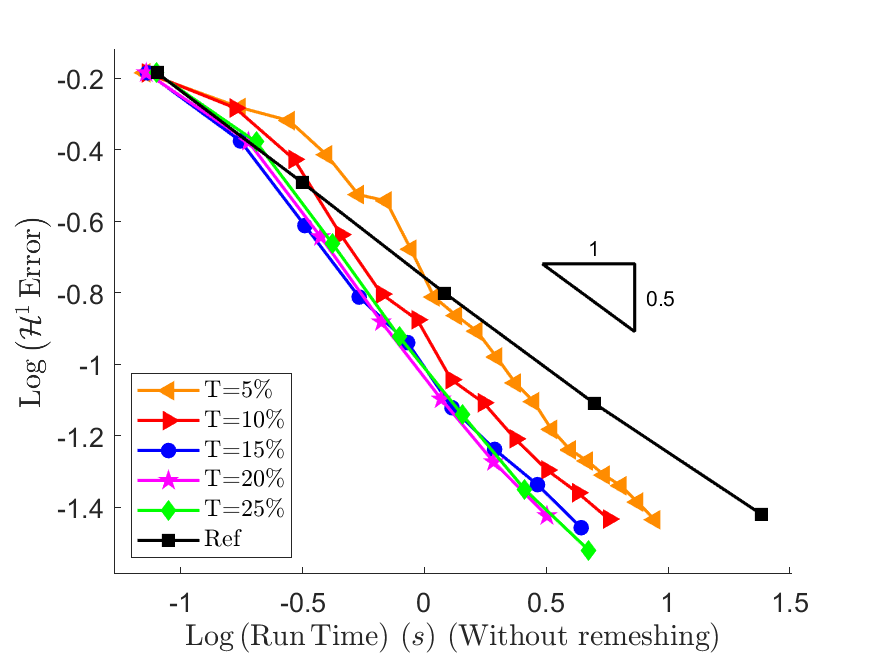}
			\caption{Structured mesh}
		\end{subfigure}%
		\begin{subfigure}[t]{0.5\textwidth}
			\centering
			\includegraphics[width=0.95\textwidth]{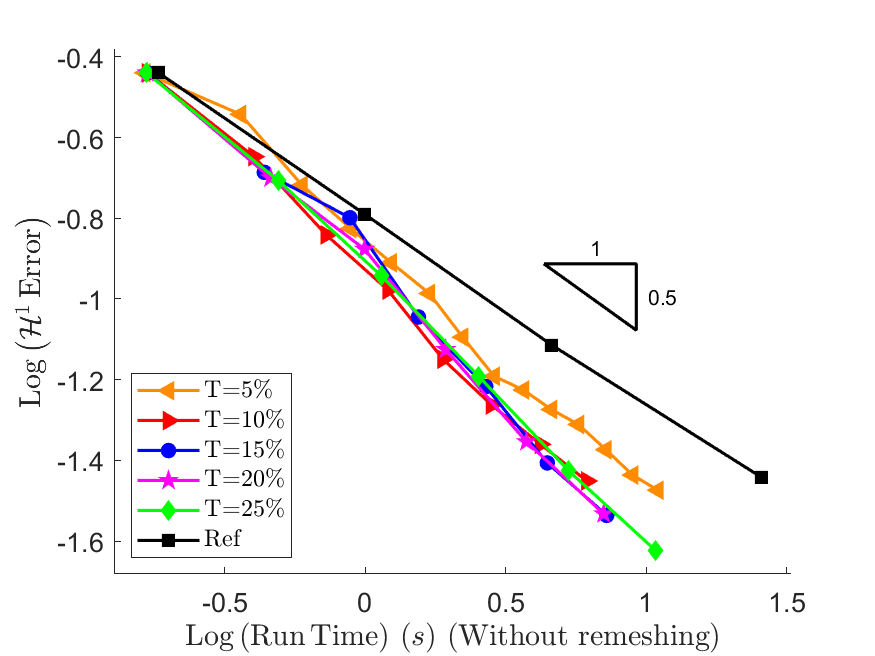}
			\caption{Voronoi mesh}
		\end{subfigure}
		\caption{$\mathcal{H}^{1}$ error vs run time (excluding remeshing time) for problem~A(1) using the displacement-based refinement procedure with a variety of choices of $T$ for structures and Voronoi meshes with a compressible Poisson's ratio.
			\label{fig:PlateWithHoleTractionConvergenceRunTime}}
	\end{figure}

	The convergence behaviour in the ${\mathcal{H}^{1}}$ error norm of the VEM for problem~A(1) using the displacement-based refinement procedure is depicted in Figure~\ref{fig:PlateWithHoleTractionConvergenceRunTime} on a logarithmic scale for a variety of choices of $T$. Here, the ${\mathcal{H}^{1}}$ error is plotted against run time for structured and Voronoi meshes with a compressible Poisson's ratio. In this figure the remeshing time (that is, time to perform element refinement, see Figure~\ref{fig:MeshRefinement}) is not included in the run time. However, the run time does include the time taken to compute the refinement indicators and mark elements for refinement. Furthermore, the run time is cumulative and, thus, includes that of all preceding remeshing steps. For clarity, each marker in Figure~\ref{fig:PlateWithHoleTractionConvergenceRunTime} represents a remeshing step.
	The rationale behind excluding remeshing time is discussed in Section~\ref{subsec:ThresholdComp}.
	It is clear from Figure~\ref{fig:PlateWithHoleTractionConvergenceRunTime} that lower choices of $T$ require significantly more remeshing steps, and consequently more run time, than larger values of $T$. Thus, lower choices of $T$ are not particularly efficient in terms of the run time. It follows then, that to determine an optimum choice of $T$ some balance of performance in terms of the number of nodes/vertices in the discretization and run time must be considered.
	However, for both mesh types and for all choices of $T$ the adaptive procedure exhibits a superior convergence rate to, and significantly outperforms, the reference refinement procedure.
	
	\FloatBarrier
	\FloatBarrier
	
	The performance of the VEM in terms of its convergence behaviour in the $\mathcal{H}^{1}$ error norm with respect to run time (excluding remeshing time) when using the displacement-based refinement procedure for problem~A(1) is summarized in Table~\ref{tab:PerformancePlateWithHoleTractionRunTime}. Here, the performance, as measured by the PRE, is presented for the cases of structured and Voronoi meshes with compressible and nearly incompressible Poisson's ratios for a variety of choices of $T$. The relative inefficiencies exhibited by lower choices of $T$ in Figure~\ref{fig:PlateWithHoleTractionConvergenceRunTime} are again evident in the table. Furthermore, these inefficiencies are exacerbated in the case of near-incompressibility where, as discussed previously, the more local/targeted refinement offered by lower choices of $T$ is less beneficial as a result in the greater distribution of error over the domain.
	This further motivates determination of a choice of $T$ that provides a good balance of performance in terms of the number of nodes/vertices and run time. Furthermore, it is clear that attention must also be paid to the balance of performance in the cases of compressibility and near-incompressibility. 
	While the lower choices of $T$ exhibit somewhat poor performance, particularly in the nearly incompressible case, the higher choices of $T$ exhibit good performance and represent significant improvements in efficiency compared to the reference procedure.
	
	\FloatBarrier
	\begin{table}[ht!]
		\centering 
		\begin{adjustbox}{max width=\textwidth}
			\begin{tabular}{|c|cccc|cccc|}
				\hline
				\multirow{3}{*}{Threshold} & \multicolumn{4}{c|}{Compressible}                                                                  & \multicolumn{4}{c|}{Nearly-incompressible}                                                           \\ \cline{2-9} 
				& \multicolumn{2}{c|}{Structured}                            & \multicolumn{2}{c|}{Voronoi}          & \multicolumn{2}{c|}{Structured}                             & \multicolumn{2}{c|}{Voronoi}           \\ \cline{2-9} 
				& \multicolumn{1}{c|}{Run time} & \multicolumn{1}{c|}{PRE}   & \multicolumn{1}{c|}{Run time} & PRE   & \multicolumn{1}{c|}{Run time} & \multicolumn{1}{c|}{PRE}    & \multicolumn{1}{c|}{Run time} & PRE    \\ \hline
				T=5\%                      & \multicolumn{1}{c|}{8.37}     & \multicolumn{1}{c|}{34.68} & \multicolumn{1}{c|}{9.20}     & 35.83 & \multicolumn{1}{c|}{37.24}    & \multicolumn{1}{c|}{154.21} & \multicolumn{1}{c|}{26.97}    & 105.42 \\ \hline
				T=10\%                     & \multicolumn{1}{c|}{5.41}     & \multicolumn{1}{c|}{22.41} & \multicolumn{1}{c|}{5.91}     & 23.01 & \multicolumn{1}{c|}{21.54}    & \multicolumn{1}{c|}{89.21}  & \multicolumn{1}{c|}{16.01}    & 62.58  \\ \hline
				T=15\%                     & \multicolumn{1}{c|}{3.86}     & \multicolumn{1}{c|}{16.01} & \multicolumn{1}{c|}{5.07}     & 19.75 & \multicolumn{1}{c|}{16.34}    & \multicolumn{1}{c|}{67.65}  & \multicolumn{1}{c|}{11.94}    & 46.69  \\ \hline
				T=20\%                     & \multicolumn{1}{c|}{3.14}     & \multicolumn{1}{c|}{12.99} & \multicolumn{1}{c|}{5.11}     & 19.90 & \multicolumn{1}{c|}{13.30}    & \multicolumn{1}{c|}{55.08}  & \multicolumn{1}{c|}{10.24}    & 40.04  \\ \hline
				T=25\%                     & \multicolumn{1}{c|}{3.30}     & \multicolumn{1}{c|}{13.68} & \multicolumn{1}{c|}{5.60}     & 21.82 & \multicolumn{1}{c|}{11.67}    & \multicolumn{1}{c|}{48.32}  & \multicolumn{1}{c|}{9.41}     & 36.80  \\ \hline
				Ref                        & \multicolumn{1}{c|}{24.13}    & \multicolumn{1}{c|}{}      & \multicolumn{1}{c|}{25.68}    &       & \multicolumn{1}{c|}{24.15}    & \multicolumn{1}{c|}{}       & \multicolumn{1}{c|}{25.58}    &        \\ \hline
			\end{tabular}
		\end{adjustbox}
		\caption{Performance summary of the VEM in terms of its convergence behaviour in the $\mathcal{H}^{1}$ error norm with respect to run time (excluding remeshing time) when using the displacement-based refinement procedure for problem~A(1).
			\label{tab:PerformancePlateWithHoleTractionRunTime}}
	\end{table}
	\FloatBarrier

\begin{figure}[ht!]
		\centering
		\begin{subfigure}[t]{0.5\textwidth}
			\centering
			\includegraphics[width=0.95\textwidth]{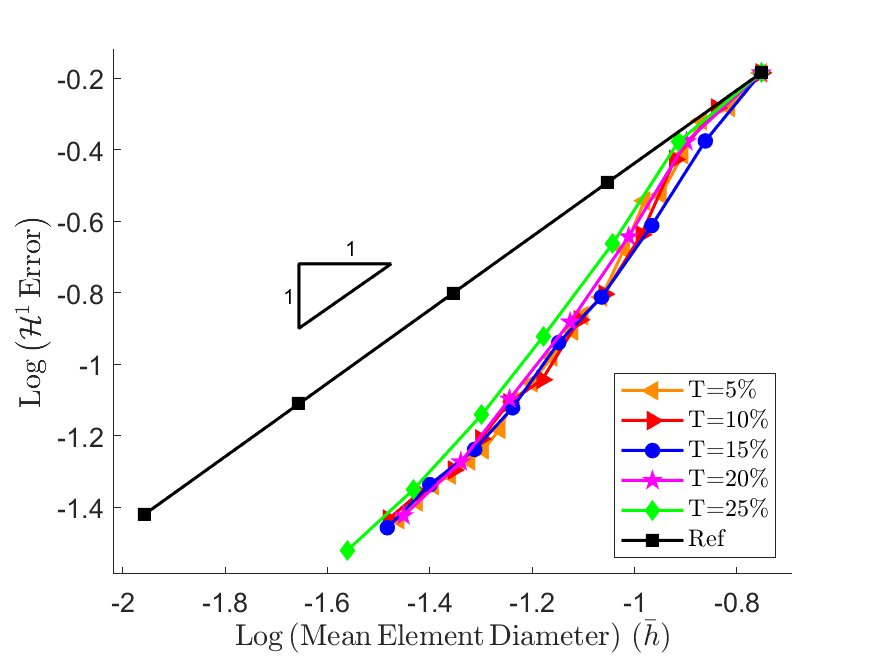}
			\caption{Structured mesh}
		\end{subfigure}%
		\begin{subfigure}[t]{0.5\textwidth}
			\centering
			\includegraphics[width=0.95\textwidth]{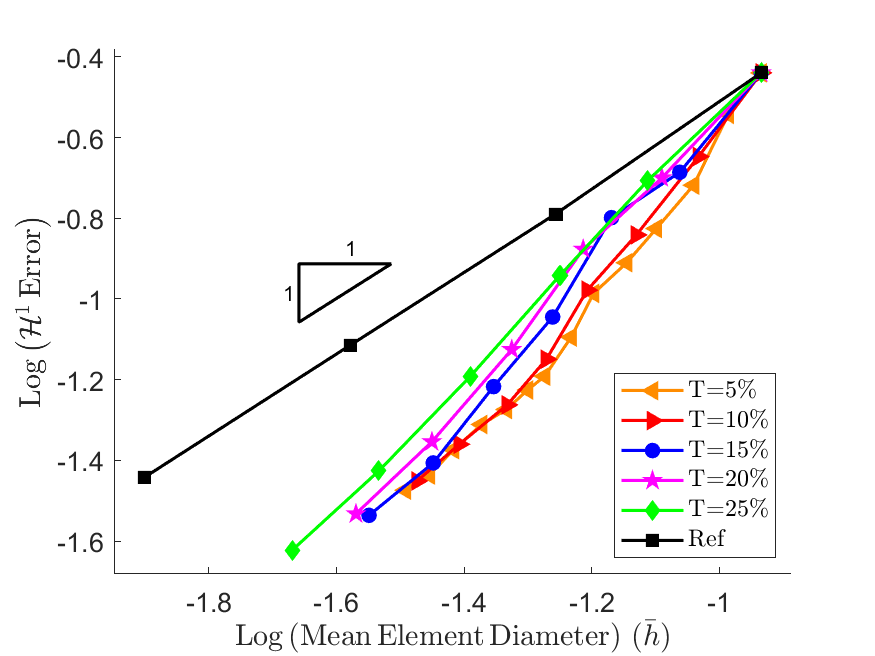}
			\caption{Voronoi mesh}
		\end{subfigure}
		\caption{$\mathcal{H}^{1}$ error vs mesh size for problem~A(1) using the displacement-based refinement procedure with a variety of choices of $T$ for structured and Voronoi meshes with a compressible Poisson's ratio.
			\label{fig:PlateWithHoleTractionConvergenceMeshSize}}
	\end{figure}

	The convergence behaviour in the ${\mathcal{H}^{1}}$ error norm of the VEM for problem~A(1) using the displacement-based refinement procedure is depicted in Figure~\ref{fig:PlateWithHoleTractionConvergenceMeshSize} on a logarithmic scale for a variety of choices of $T$. Here, the ${\mathcal{H}^{1}}$ error is plotted against mesh size, as measured by the mean element diameter, for structured and Voronoi meshes with a compressible Poisson's ratio. 
	The convergence in the ${\mathcal{H}^{1}}$ error norm with respect to mesh size is qualitatively similar to that observed in Figure~\ref{fig:PlateWithHoleTractionConvergenceNumberOfNodes} with respect to the number of nodes/vertices. In the coarse mesh range, particularly for the case of unstructured/Voronoi meshes, lower choices of $T$ initially exhibit a slightly superior convergence rate to larger values of $T$, while larger values of $T$ exhibit more consistent convergence behaviour throughout the domain. 
	However, as the level of mesh refinement increases all choices of $T$ exhibit similar convergence behaviour. For both mesh types, and for all choices of $T$, the adaptive refinement procedure exhibits a superior convergence rate to, and significantly outperforms, the reference refinement procedure.
	
	\FloatBarrier
	\FloatBarrier
	
	The performance of the VEM in terms of its convergence behaviour in the $\mathcal{H}^{1}$ error norm with respect to mesh size when using the displacement-based refinement procedure for problem~A(1) is summarized in Table~\ref{tab:PerformancePlateWithHoleTractionMeshSize}. Here, the performance, as measured by the PRE, is presented for the cases of structured and Voronoi meshes with compressible and nearly incompressible Poisson's ratios for a variety of choices of $T$. 
	The PRE is calculated as described previously, however, in the case of mesh size a larger PRE is desirable as it indicates that the adaptive procedure has on average used larger elements than the reference procedure to reach the same accuracy.
	For compatibility, a modified PRE is introduced and denoted by PRE*. The PRE* is computed as
	\begin{equation}
		PRE^{*} = \frac{10000}{PRE} \,.
	\end{equation}
	Thus, a lower PRE* indicates better performance.
	The PRE* results presented in Table~\ref{tab:PerformancePlateWithHoleTractionMeshSize} are qualitatively similar to the results in terms of the number of nodes/vertices. Specifically, the degree of compressibility has a significant influence on the performance of the adaptive refinement procedure as a greater level of refinement is required throughout the domain in the case of near-incompressibility.
	However, for both mesh types and for all choices of $T$ the adaptive procedure significantly outperforms the reference procedure.
	
	\FloatBarrier
	\begin{table}[ht!]
		\centering 
		\begin{adjustbox}{max width=\textwidth}
			\begin{tabular}{|c|cccc|cccc|}
				\hline
				\multirow{3}{*}{Threshold} & \multicolumn{4}{c|}{Compressible}                                                                    & \multicolumn{4}{c|}{Nearly-incompressible}                                                           \\ \cline{2-9} 
				& \multicolumn{2}{c|}{Structured}                             & \multicolumn{2}{c|}{Voronoi}           & \multicolumn{2}{c|}{Structured}                             & \multicolumn{2}{c|}{Voronoi}           \\ \cline{2-9} 
				& \multicolumn{1}{c|}{Mesh size} & \multicolumn{1}{c|}{PRE*}  & \multicolumn{1}{c|}{Mesh size} & PRE*  & \multicolumn{1}{c|}{Mesh size} & \multicolumn{1}{c|}{PRE*}  & \multicolumn{1}{c|}{Mesh size} & PRE*  \\ \hline
				T=5\%                      & \multicolumn{1}{c|}{0.035546}  & \multicolumn{1}{c|}{31.08} & \multicolumn{1}{c|}{0.034720}  & 36.25 & \multicolumn{1}{c|}{0.017815}  & \multicolumn{1}{c|}{62.02} & \multicolumn{1}{c|}{0.022563}  & 55.81 \\ \hline
				T=10\%                     & \multicolumn{1}{c|}{0.033932}  & \multicolumn{1}{c|}{32.56} & \multicolumn{1}{c|}{0.034147}  & 36.85 & \multicolumn{1}{c|}{0.017721}  & \multicolumn{1}{c|}{62.35} & \multicolumn{1}{c|}{0.022643}  & 55.61 \\ \hline
				T=15\%                     & \multicolumn{1}{c|}{0.034927}  & \multicolumn{1}{c|}{31.63} & \multicolumn{1}{c|}{0.033454}  & 37.62 & \multicolumn{1}{c|}{0.017817}  & \multicolumn{1}{c|}{62.01} & \multicolumn{1}{c|}{0.022706}  & 55.46 \\ \hline
				T=20\%                     & \multicolumn{1}{c|}{0.035655}  & \multicolumn{1}{c|}{30.99} & \multicolumn{1}{c|}{0.030985}  & 40.62 & \multicolumn{1}{c|}{0.018041}  & \multicolumn{1}{c|}{61.24} & \multicolumn{1}{c|}{0.022757}  & 55.34 \\ \hline
				T=25\%                     & \multicolumn{1}{c|}{0.032793}  & \multicolumn{1}{c|}{33.69} & \multicolumn{1}{c|}{0.028565}  & 44.06 & \multicolumn{1}{c|}{0.018207}  & \multicolumn{1}{c|}{60.68} & \multicolumn{1}{c|}{0.022781}  & 55.28 \\ \hline
				Ref                        & \multicolumn{1}{c|}{0.011049}  & \multicolumn{1}{c|}{}      & \multicolumn{1}{c|}{0.012585}  &       & \multicolumn{1}{c|}{0.011049}  & \multicolumn{1}{c|}{}      & \multicolumn{1}{c|}{0.012593}  &       \\ \hline
			\end{tabular}
		\end{adjustbox}
		\caption{Performance summary of the VEM in terms of its convergence behaviour in the $\mathcal{H}^{1}$ error norm with respect to mesh size when using the displacement-based refinement procedure for problem~A(1).
			\label{tab:PerformancePlateWithHoleTractionMeshSize}}
	\end{table}
	\FloatBarrier
	
\begin{figure}[ht!]
		\centering
		\begin{subfigure}[t]{0.5\textwidth}
			\centering
			\includegraphics[width=0.95\textwidth]{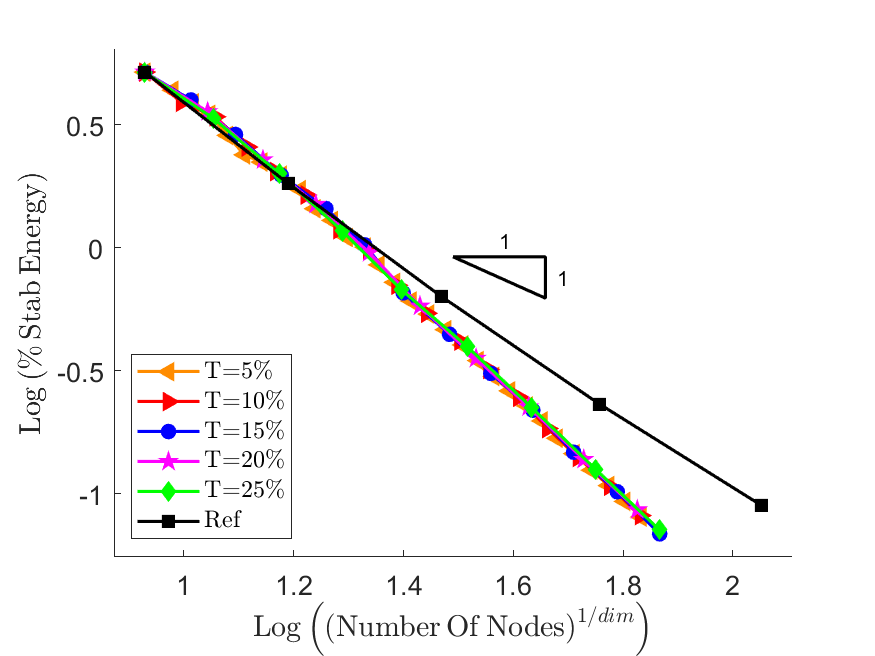}
			\caption{Structured mesh}
		\end{subfigure}%
		\begin{subfigure}[t]{0.5\textwidth}
			\centering
			\includegraphics[width=0.95\textwidth]{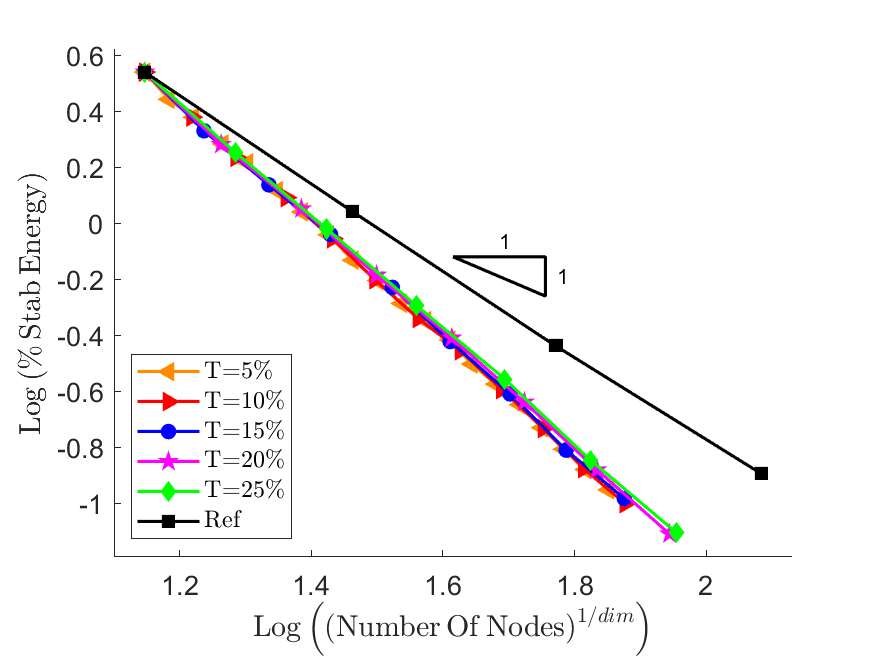}
			\caption{Voronoi mesh}
		\end{subfigure}
		\caption{PSE vs number of nodes/vertices in the discretization for problem~A(1) using the displacement-based refinement procedure with a variety of choices of $T$ for structured and Voronoi meshes with a compressible Poisson's ratio.
			\label{fig:PlateWithHoleTractionConvergencePSE}}
	\end{figure}

	The convergence behaviour in the percentage stabilization energy (PSE) of the VEM for problem~A(1) using the displacement-based refinement procedure is depicted in Figure~\ref{fig:PlateWithHoleTractionConvergencePSE} on a logarithmic scale for a variety of choices of $T$. Here, the PSE is plotted against the number of nodes/vertices in the discretization for the cases of structured and unstructured/Voronoi meshes with a compressible Poisson's ratio. 
	The PSE indicates the contribution of the total stored elastic energy in a deformed body arising from the stabilization term. The PSE is used as a measure of a mesh's suitability for modelling a specific problem, with a lower PSE corresponding to a more suitable mesh (see [\cite{vanHuyssteen2022}]). 
	The PSE is computed as
	\begin{equation}
		\text{PSE} = \frac{\text{E}|_{\text{s}}}{\text{E}|_{\text{t}}} \, ,
	\end{equation}
	where $\text{E}|_{\text{t}}$ is the total stored energy in a body and $\text{E}|_{\rm{s}}$ is the energy contribution from the stabilization term. Additionally, $\text{E}|_{\text{t}}$ and $\text{E}|_{\rm{s}}$ are computed as the sum of their element-level contributions defined respectively as 
	\begin{equation}
		\text{E}|_{\text{t}}^{E} = \frac{1}{2} \, \bd^{E} \cdot \left[ \bK^{E} \cdot \bd^{E} \right] \, 
		\quad \text{and} \quad 
		\text{E}|_{\text{s}}^{E} =  \frac{1}{2} \, \bd^{E} \cdot \left[ \bK_{\text{s}}^{E} \cdot \bd^{E} \right] \, .
	\end{equation}
	The convergence of the PSE is very similar for all choices of $T$. For structured meshes in the coarse mesh range the adaptive procedure exhibits similar convergence behaviour to the reference procedure. However, as the level of mesh refinement increases the adaptive procedure exhibits superior convergence behaviour and outperforms the reference procedure. For unstructured/Voronoi meshes the adaptive procedure exhibits superior convergence behaviour and outperforms the reference procedure throughout the domain.
	These results indicate that meshes generated by the adaptive procedure are well-suited to the specific example problem.
	
	\FloatBarrier
	\FloatBarrier
	
	The performance of the VEM in terms of its convergence behaviour in the PSE with respect to the number of nodes/vertices in the discretization when using the displacement-based refinement procedure for problem~A(1) is summarized in Table~\ref{tab:PerformancePlateWithHoleTractionPSE}. Here, the performance, as measured by the PRE, is presented for the cases of structured and Voronoi meshes with compressible and nearly incompressible Poisson's ratios for a variety of choices of $T$. 
	In general, the results presented in the table expectedly show that lower choices of $T$ perform slightly better than larger choices of $T$. However, the influence of the choice of $T$ is small. 
	Most notably, the results do not show any significant dependence on the degree of compressibility. This indicates that the meshes generated by the adaptive procedure are equally suited to the problem in the cases of compressibility and near incompressibility.
	
	\FloatBarrier
	\begin{table}[ht!]
		\centering 
		\begin{adjustbox}{max width=\textwidth}
			\begin{tabular}{|c|cccc|cccc|}
				\hline
				\multirow{3}{*}{Threshold} & \multicolumn{4}{c|}{Compressible}                                                                  & \multicolumn{4}{c|}{Nearly-incompressible}                                                         \\ \cline{2-9} 
				& \multicolumn{2}{c|}{Structured}                            & \multicolumn{2}{c|}{Voronoi}          & \multicolumn{2}{c|}{Structured}                            & \multicolumn{2}{c|}{Voronoi}          \\ \cline{2-9} 
				& \multicolumn{1}{c|}{nNodes}   & \multicolumn{1}{c|}{PRE}   & \multicolumn{1}{c|}{nNodes}   & PRE   & \multicolumn{1}{c|}{nNodes}   & \multicolumn{1}{c|}{PRE}   & \multicolumn{1}{c|}{nNodes}   & PRE   \\ \hline
				T=5\%                      & \multicolumn{1}{c|}{4170.08}  & \multicolumn{1}{c|}{32.91} & \multicolumn{1}{c|}{4435.69}  & 30.23 & \multicolumn{1}{c|}{3755.49}  & \multicolumn{1}{c|}{29.64} & \multicolumn{1}{c|}{4151.42}  & 28.55 \\ \hline
				T=10\%                     & \multicolumn{1}{c|}{4206.29}  & \multicolumn{1}{c|}{33.19} & \multicolumn{1}{c|}{4419.26}  & 30.11 & \multicolumn{1}{c|}{3785.78}  & \multicolumn{1}{c|}{29.88} & \multicolumn{1}{c|}{4168.32}  & 28.66 \\ \hline
				T=15\%                     & \multicolumn{1}{c|}{4247.38}  & \multicolumn{1}{c|}{33.52} & \multicolumn{1}{c|}{4570.95}  & 31.15 & \multicolumn{1}{c|}{3780.25}  & \multicolumn{1}{c|}{29.83} & \multicolumn{1}{c|}{4443.69}  & 30.56 \\ \hline
				T=20\%                     & \multicolumn{1}{c|}{4303.95}  & \multicolumn{1}{c|}{33.96} & \multicolumn{1}{c|}{4795.70}  & 32.68 & \multicolumn{1}{c|}{3799.42}  & \multicolumn{1}{c|}{29.98} & \multicolumn{1}{c|}{4597.87}  & 31.62 \\ \hline
				T=25\%                     & \multicolumn{1}{c|}{4341.72}  & \multicolumn{1}{c|}{34.26} & \multicolumn{1}{c|}{4972.25}  & 33.88 & \multicolumn{1}{c|}{3915.12}  & \multicolumn{1}{c|}{30.90} & \multicolumn{1}{c|}{4625.47}  & 31.81 \\ \hline
				Ref                        & \multicolumn{1}{c|}{12672.00} & \multicolumn{1}{c|}{}      & \multicolumn{1}{c|}{14675.00} &       & \multicolumn{1}{c|}{12672.00} & \multicolumn{1}{c|}{}      & \multicolumn{1}{c|}{14542.00} &       \\ \hline
			\end{tabular}
		\end{adjustbox}
		\caption{Performance summary of the VEM in terms of its convergence behaviour in the PSE with respect to the number of nodes/vertices in the discretization when using the displacement-based refinement procedure for problem~A(1).
			\label{tab:PerformancePlateWithHoleTractionPSE}}
	\end{table}
	\FloatBarrier
	
\begin{figure}[ht!]
		\centering
		\begin{subfigure}[t]{0.5\textwidth}
			\centering
			\includegraphics[width=0.8\textwidth]{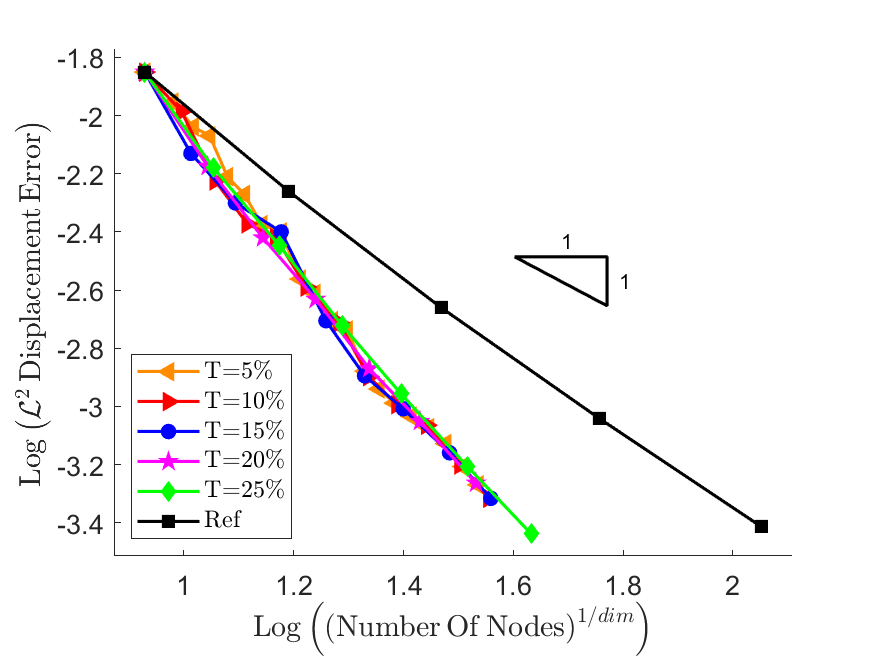}
			\caption{Displacement error - Structured meshes}
		\end{subfigure}%
		\begin{subfigure}[t]{0.5\textwidth}
			\centering
			\includegraphics[width=0.8\textwidth]{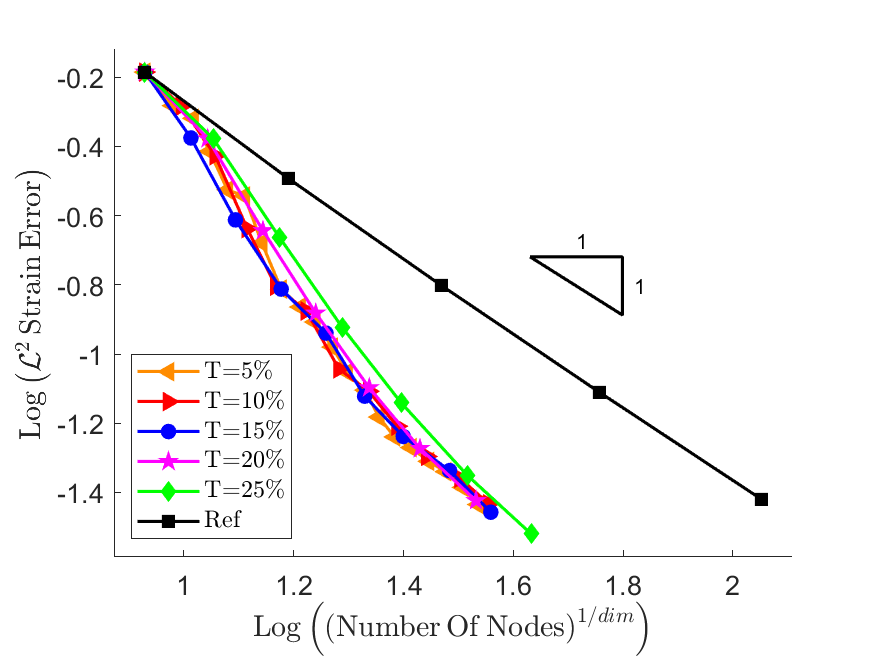}
			\caption{Strain error - Structured meshes}
		\end{subfigure}
		\vskip \baselineskip 
		\vspace*{-3mm}
		\begin{subfigure}[t]{0.5\textwidth}
			\centering
			\includegraphics[width=0.8\textwidth]{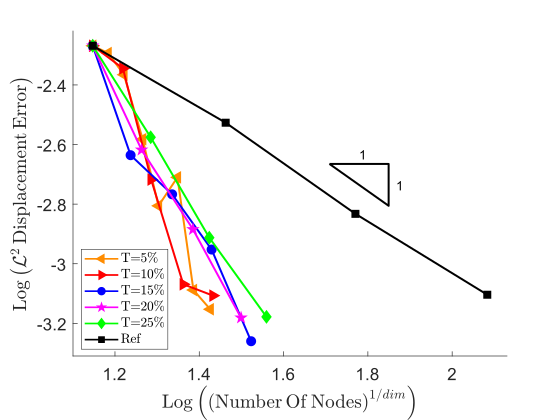}
			\caption{Displacement error - Voronoi meshes}
		\end{subfigure}%
		\begin{subfigure}[t]{0.5\textwidth}
			\centering
			\includegraphics[width=0.8\textwidth]{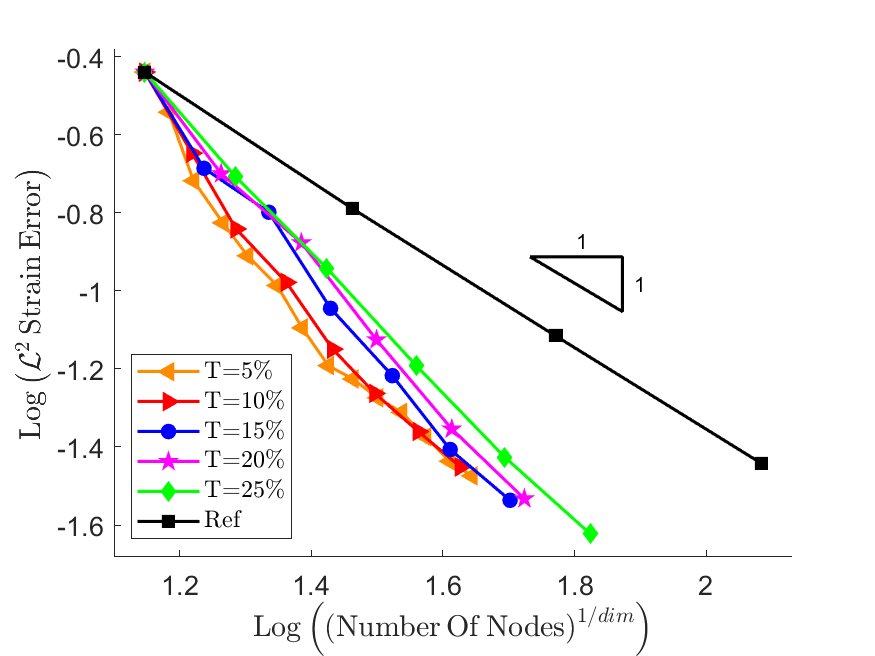}
			\caption{Strain error - Voronoi meshes}
		\end{subfigure}
		\caption{Displacement and strain $\mathcal{L}^{2}$ error components vs $n_{\rm v}$ for problem~A(1) using the displacement-based refinement procedure with a variety of choices of $T$ for structured and Voronoi meshes with a compressible Poisson's ratio.
			\label{fig:PlateWithHoleTractionConvergenceComponents}}
	\end{figure}

	In addition to the convergence behaviour in the ${\mathcal{H}^{1}}$ error norm, attention is also paid to the convergence of the different components/contributions to the ${\mathcal{H}^{1}}$ error. Specifically, the convergence behaviour in the ${\mathcal{L}^{2}}$ error norm of the displacement and strain field approximations when using the displacement-based refinement procedure for problem~A(1) is plotted in Figure~\ref{fig:PlateWithHoleTractionConvergenceComponents} against the number of vertices/nodes in the discretization. 
	Here, the first and second rows correspond to structured and Voronoi meshes respectively, while the left-hand and right-hand columns respectively correspond to the displacement and strain error components for the case of a compressible Poisson's ratio. 
	The convergence behaviour, in both error components, for both mesh types is qualitatively similar to the convergence behaviour in the ${\mathcal{H}^{1}}$ error norm presented in Figure~\ref{fig:PlateWithHoleTractionConvergenceNumberOfNodes}. For structured meshes the choice of $T$ has a comparatively small influence with similar convergence behaviour exhibited by all choices of $T$. In the case of unstructured/Voronoi meshes the convergence behaviour exhibited by lower choices of $T$ is slightly erratic, particularly in the case of the displacement error, while larger values of $T$ exhibit much more consistent convergence behaviour. Overall, for all choices of $T$ the adaptive procedure exhibits superior convergence behaviour to the reference procedure in both the displacement and strain error components for both mesh types. Furthermore, even though the formulation of displacement-based refinement indicator does not involve explicit consideration of the strain approximation, the indicator very effectively reduces the error in the approximation of the strain field. 
	
	\FloatBarrier
	\FloatBarrier
	
\begin{figure}[ht!]
		\centering
		\begin{subfigure}[t]{0.33\textwidth}
			\centering
			\includegraphics[width=0.95\textwidth]{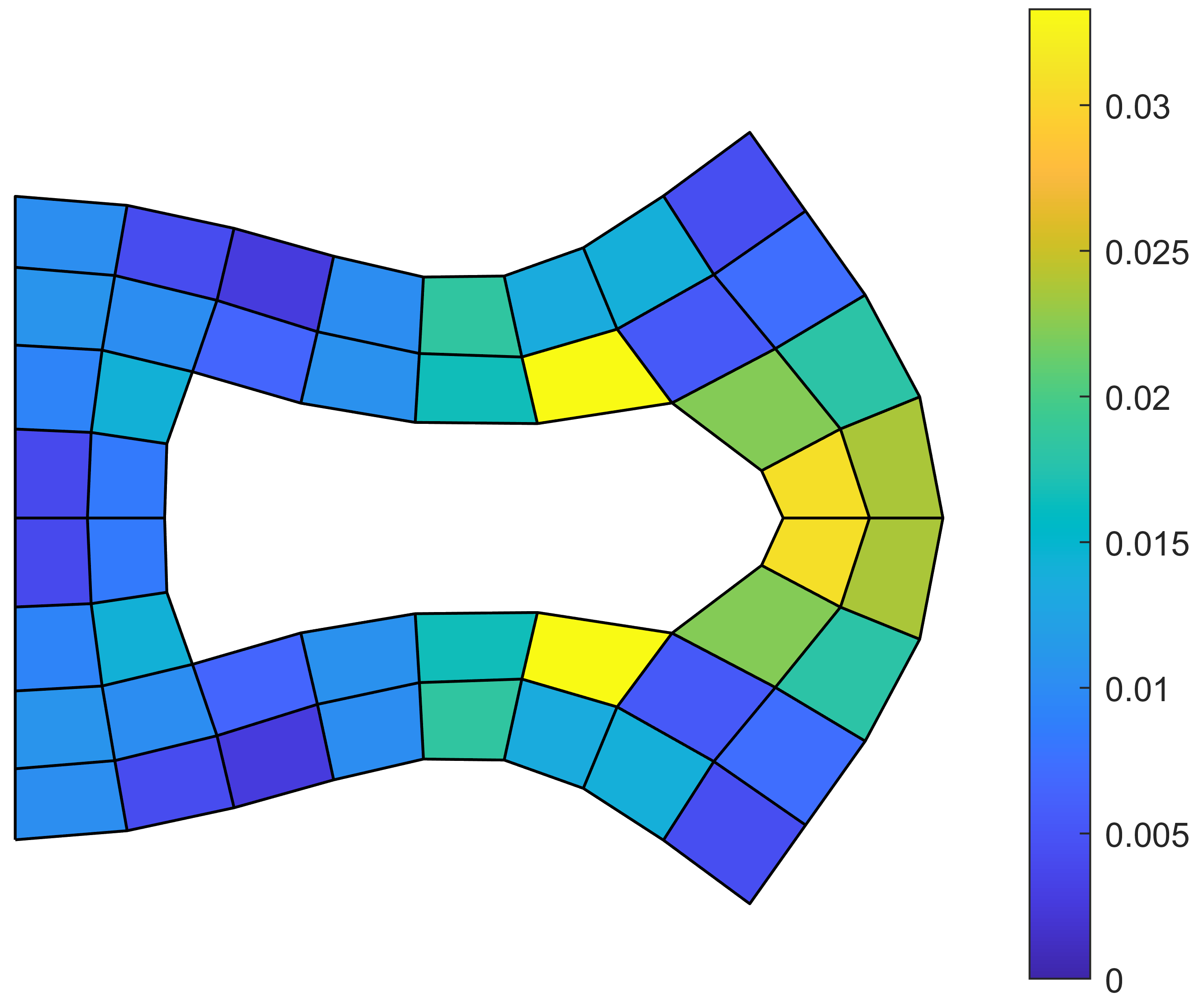}
			\caption{Structured - Step 1}
		\end{subfigure}%
		\begin{subfigure}[t]{0.33\textwidth}
			\centering
			\includegraphics[width=0.95\textwidth]{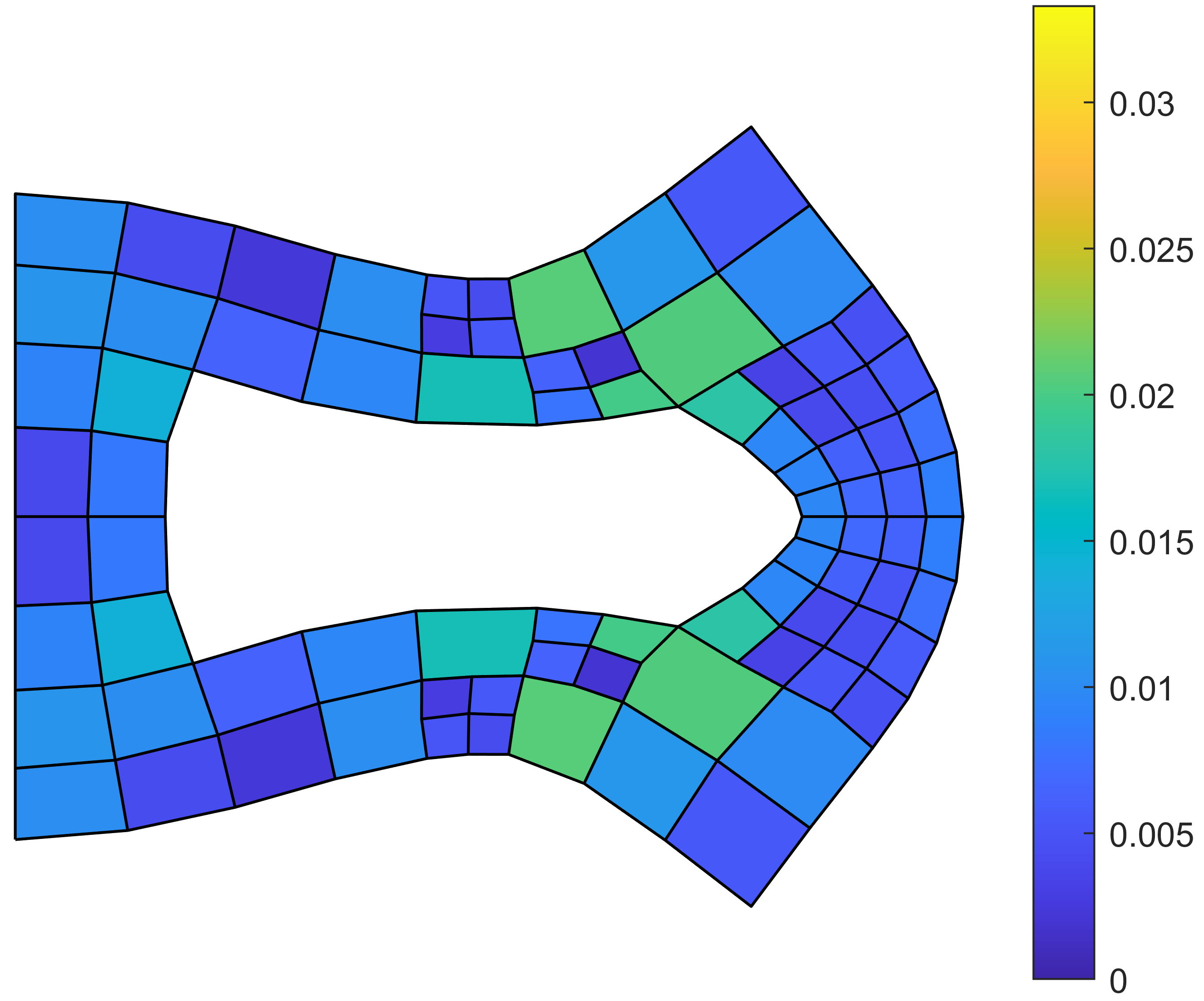}
			\caption{Structured - Step 2}
		\end{subfigure}%
		\begin{subfigure}[t]{0.33\textwidth}
			\centering
			\includegraphics[width=0.95\textwidth]{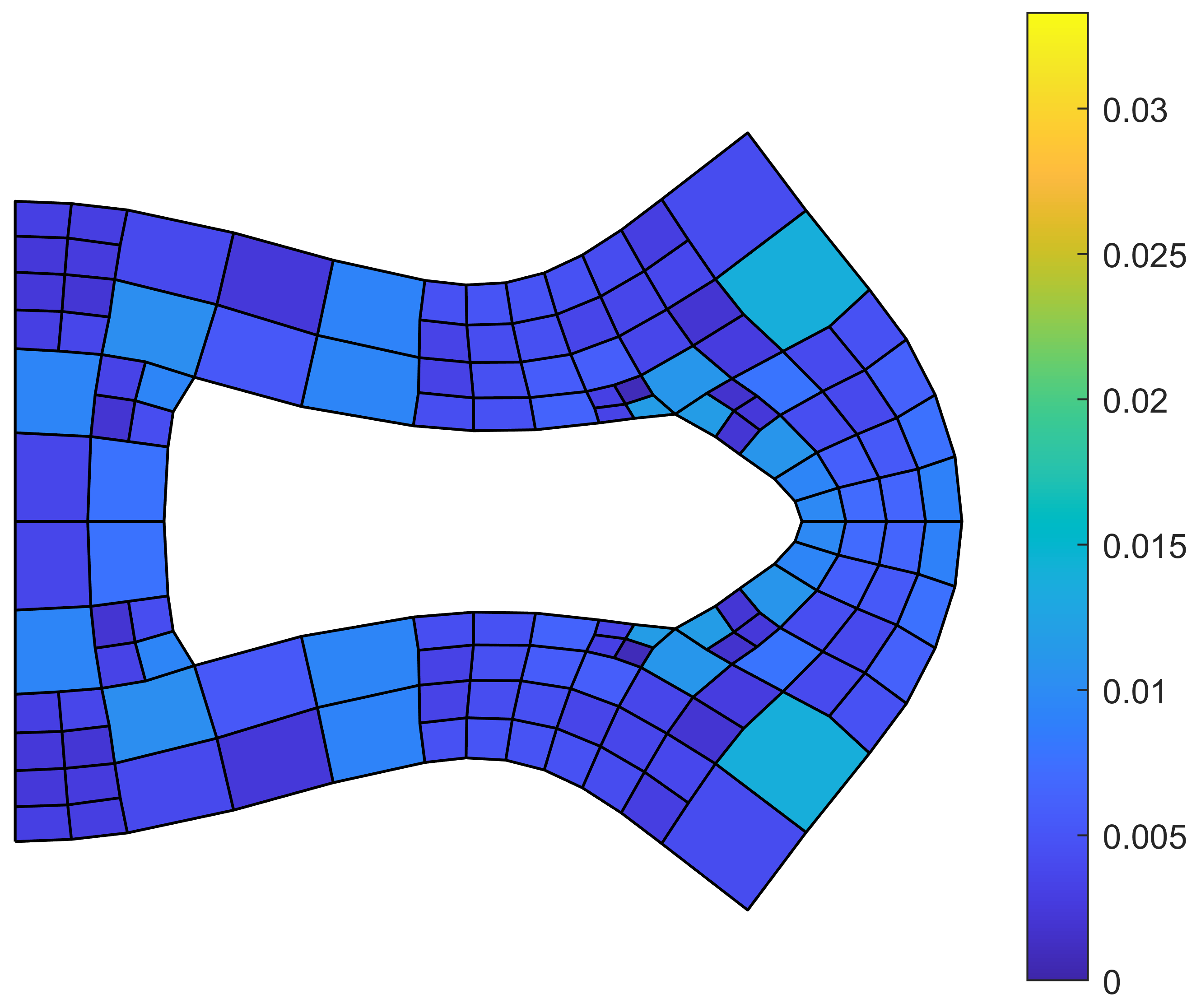}
			\caption{Structured - Step 3}
		\end{subfigure}
		\vskip \baselineskip 
		\begin{subfigure}[t]{0.33\textwidth}
			\centering
			\includegraphics[width=0.95\textwidth]{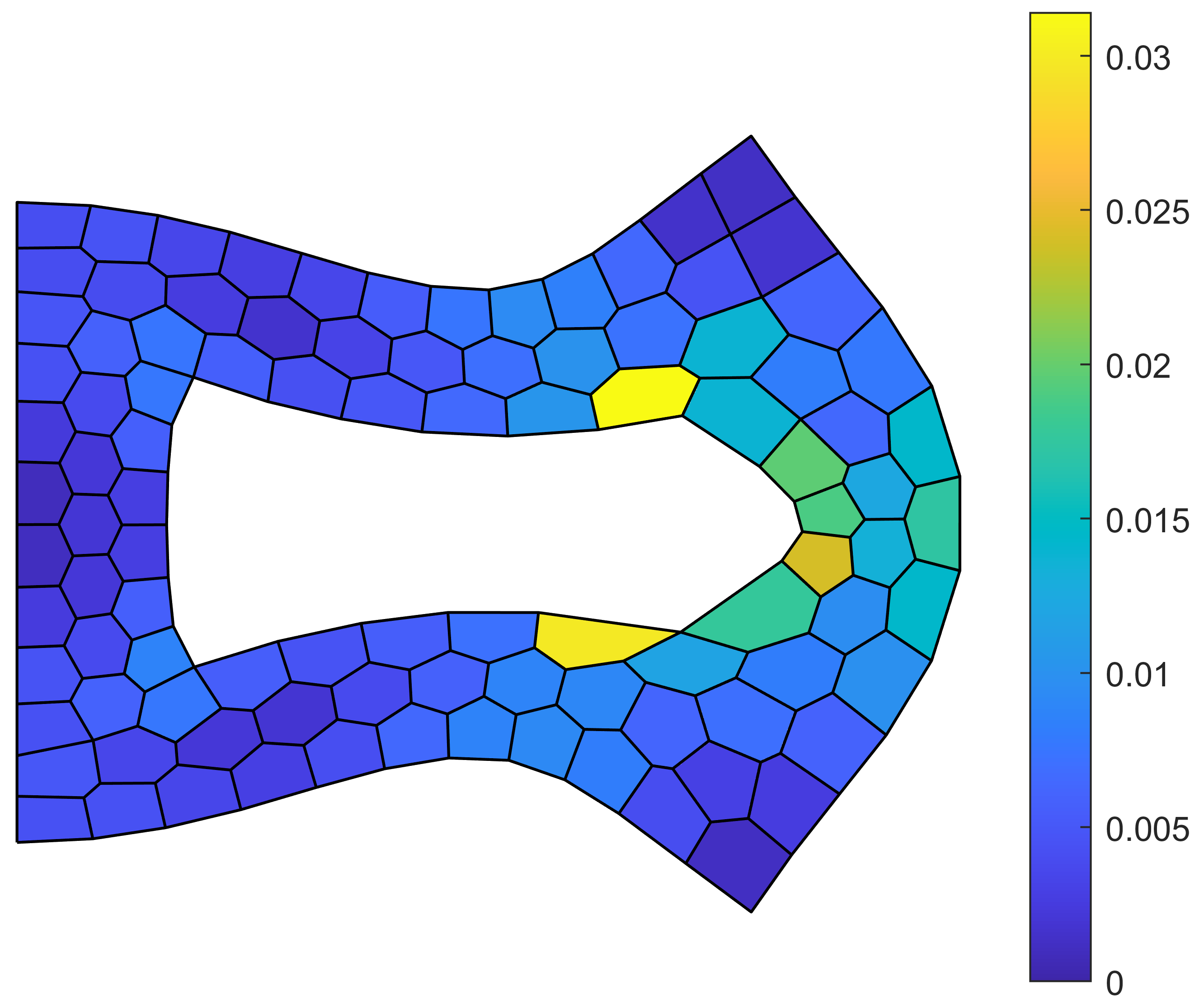}
			\caption{Voronoi - Step 1}
		\end{subfigure}%
		\begin{subfigure}[t]{0.33\textwidth}
			\centering
			\includegraphics[width=0.95\textwidth]{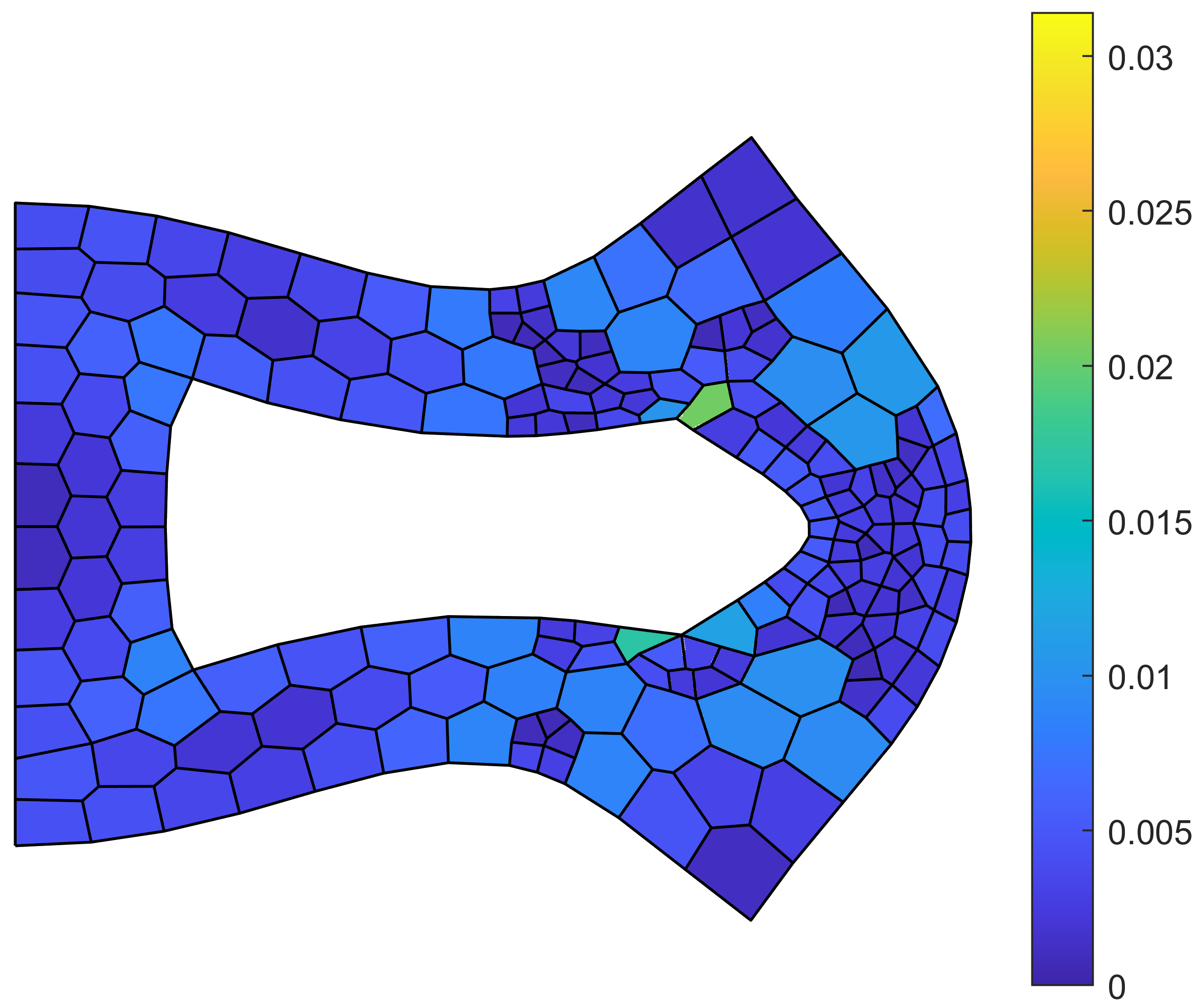}
			\caption{Voronoi - Step 2}
		\end{subfigure}%
		\begin{subfigure}[t]{0.33\textwidth}
			\centering
			\includegraphics[width=0.95\textwidth]{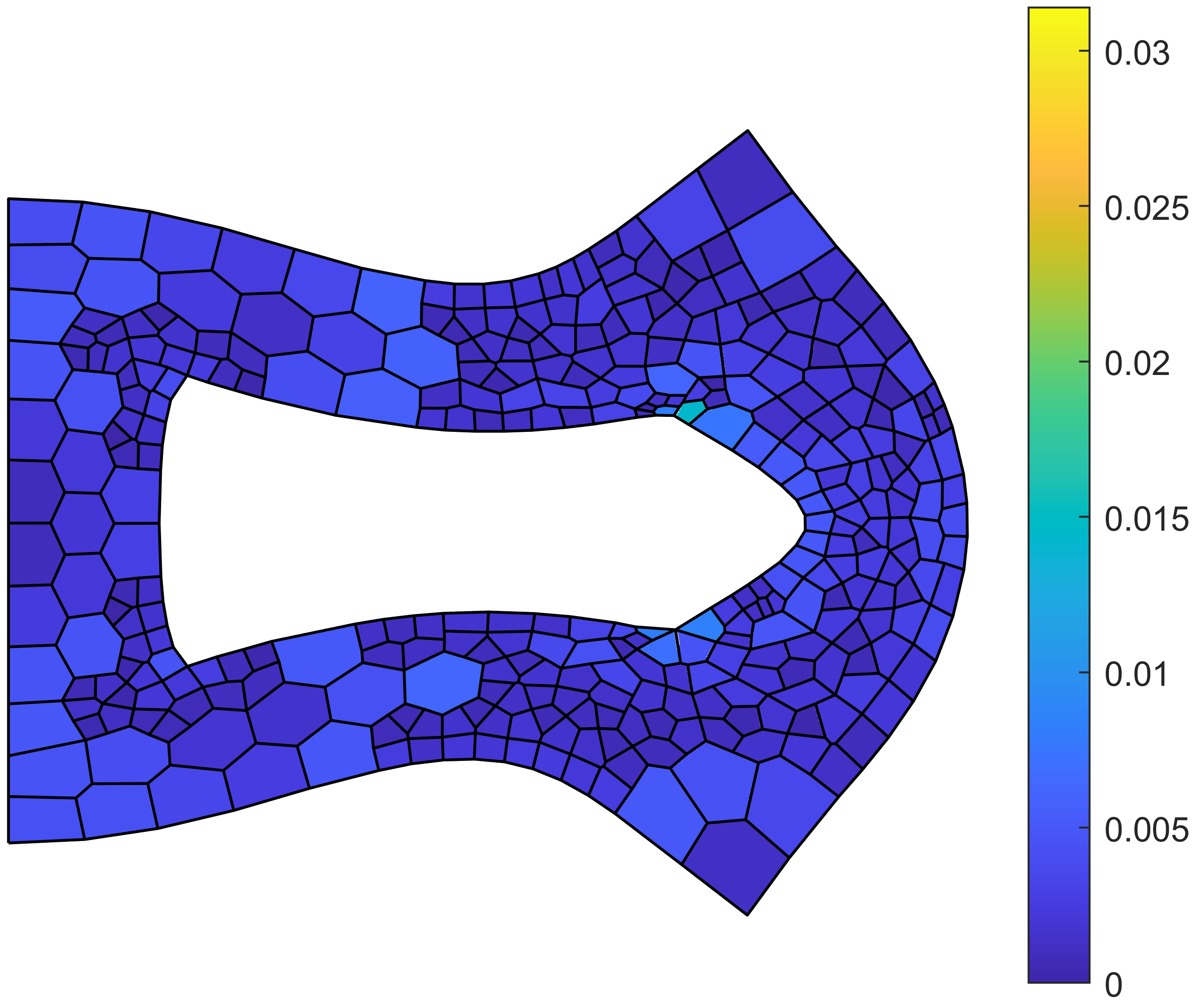}
			\caption{Voronoi - Step 3}
		\end{subfigure}
		\caption{Displacement-based refinement indicator ${\Upsilon_{\text{DB}}^{t}}$ over the domain for problem~A(1) using the displacement-based refinement procedure with ${T=20\%}$ for structured and Voronoi meshes with a compressible Poisson's ratio.
			\label{fig:PlateWithHoleTractionIndicatorMaps}}
	\end{figure}

	Figure~\ref{fig:PlateWithHoleTractionIndicatorMaps} shows plots of the displacement-based refinement indicator ${\Upsilon_{\text{DB}}^{t}}$ over the domain for problem~A(1) using the displacement-based refinement procedure with ${T=20\%}$. Plots are presented at successive refinement steps for structured and unstructured/Voronoi meshes with a compressible Poisson's ratio\footnote{For compactness, in the following sections refinement indicator maps are not presented in this work for the other refinement indicators considered. These maps, as well as maps corresponding the the case of a nearly incompressible Poisson's ratio, can be found in the supplementary material.}. 
	The refinement indicator is, expectedly, initially largest in the right-hand portion of the domain, particularly near the corners of the hole. It is clear that even after a single refinement step the magnitude of the indicator decreases significantly as the elements with the largest values of the indicator are refined. As the refinement procedure progresses the magnitude of ${\Upsilon_{\text{DB}}^{t}}$ continues to decrease and localizations/regions of larger values of ${\Upsilon_{\text{DB}}^{t}}$ are significantly reduced resulting in a smoother distribution of ${\Upsilon_{\text{DB}}^{t}}$ over the domain.
	
	\FloatBarrier
	\FloatBarrier
	
\begin{figure}[ht!]
		\centering
		\begin{subfigure}[t]{0.5\textwidth}
			\centering
			\includegraphics[width=0.95\textwidth]{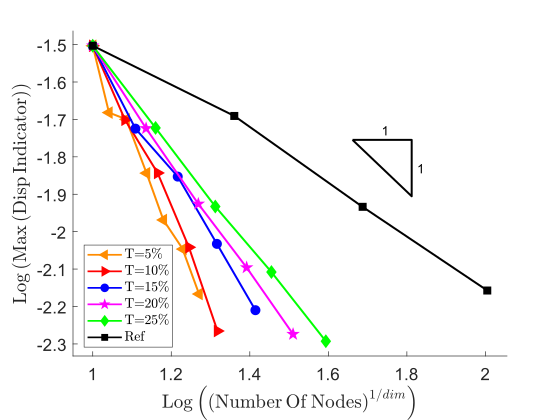}
			\caption{Max value of ${\Upsilon_{\text{DB}}^{t}}$}
		\end{subfigure}%
		\begin{subfigure}[t]{0.5\textwidth}
			\centering
			\includegraphics[width=0.95\textwidth]{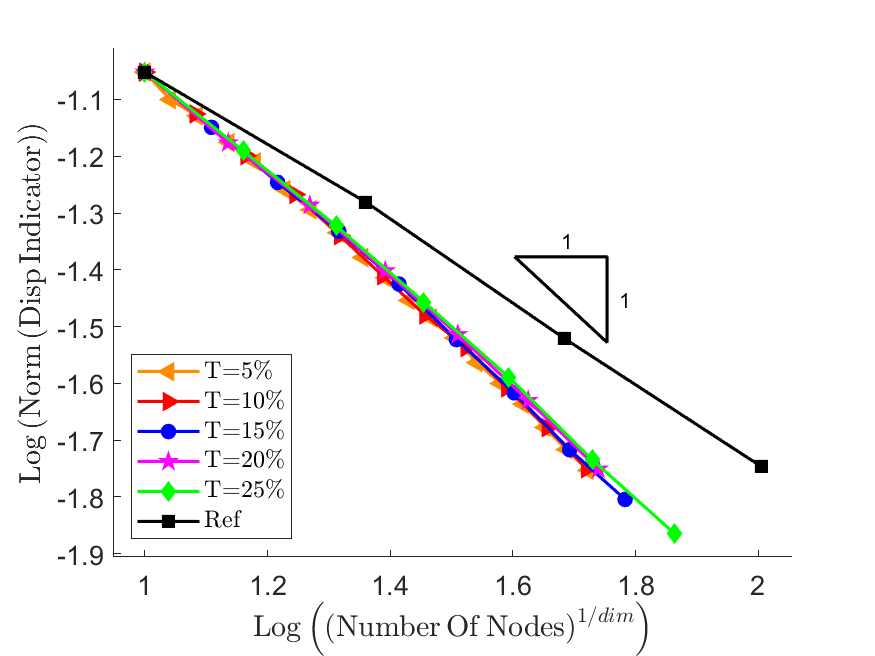}
			\caption{Norm of ${\Upsilon_{\text{DB}}^{t}}$}
		\end{subfigure}
		\caption{Maximum value and norm of the displacement-based indicator vs number of nodes/vertices in the discretization for problem~A(1) using the displacement-based refinement procedure with a variety of choices of $T$ for a compressible Poisson's ratio with unstructured/Voronoi meshes.
			\label{fig:PlateWithHoleTractionIndicatorConvergence}}
	\end{figure}

	The convergence behaviour of the maximum value and norm of the displacement-based indicator for problem~A(1) using the displacement-based refinement procedure is depicted in Figure~\ref{fig:PlateWithHoleTractionIndicatorConvergence} on a logarithmic scale for a variety of choices of $T$. Here, the values of the displacement-based indicator are plotted against the number of nodes/vertices in the discretization for unstructured/Voronoi meshes with a compressible Poisson's ratio. 
	In the case of the maximum value of ${\Upsilon_{\text{DB}}^{t}}$, the choice of $T$ has a notable influence on the convergence behaviour with lower choices of $T$ exhibiting superior convergence rates. This behaviour is expected as the lower choices of $T$ allow for more targeted/local refinement which is, of course, well suited to reducing the maximum values of the refinement indicator. 
	In terms of reducing the magnitude of the refinement indicator over the entire problem domain, as measured by the norm, all choices of $T$ are equally effective. The difference lies in the number of refinement steps required to reduce the norm of ${\Upsilon_{\text{DB}}^{t}}$. Lower values of $T$ result in more targeted/local refinement and, thus, require more steps than larger values of $T$ to reduce ${\Upsilon_{\text{DB}}^{t}}$ over the entire domain. 
	
	\FloatBarrier
	\FloatBarrier
	
\begin{figure}[ht!]
		\centering
		\begin{subfigure}[t]{0.33\textwidth}
			\centering
			\includegraphics[width=0.95\textwidth]{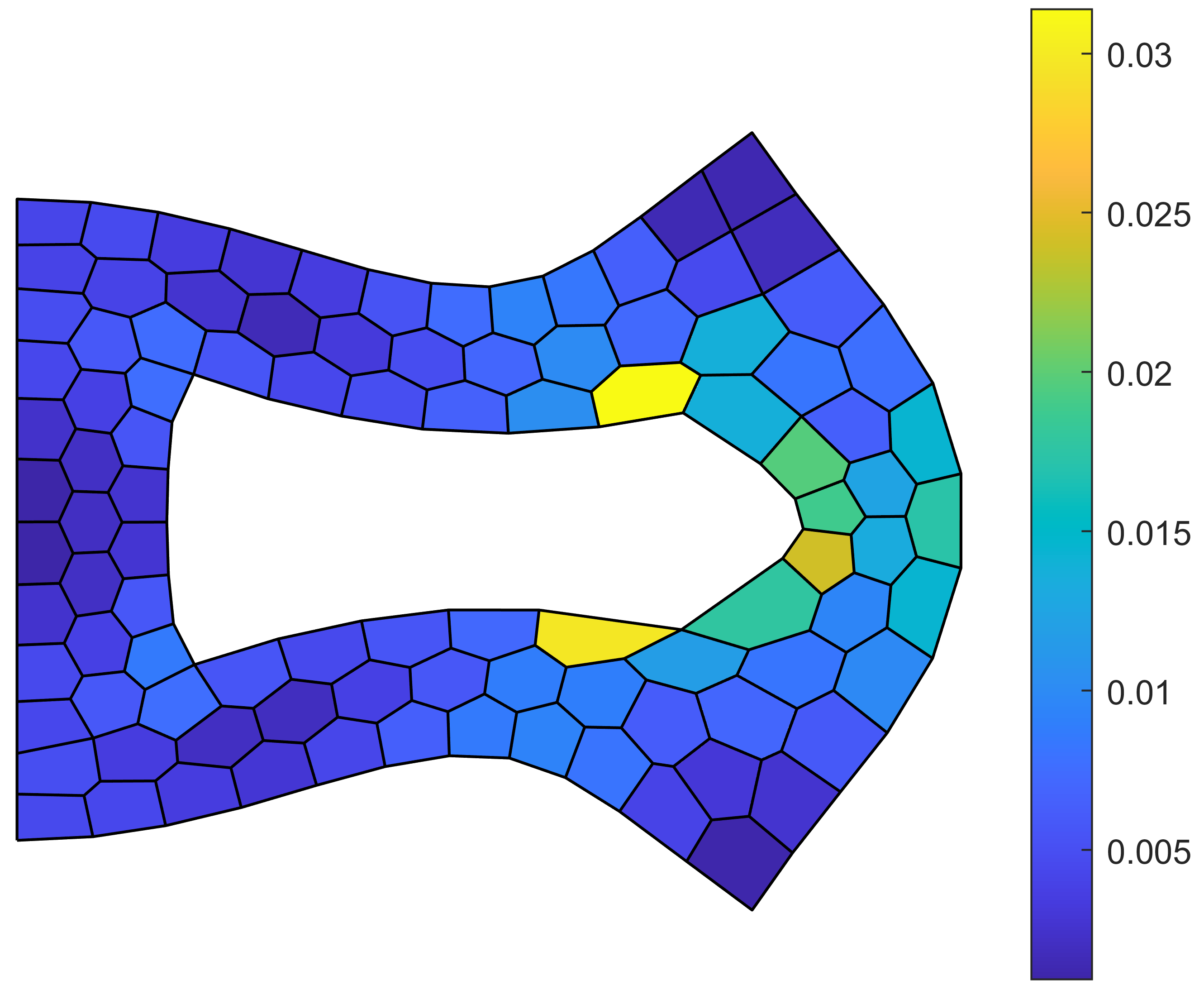}
			\caption{${\Upsilon_{\text{DB}}^{t}}$ - Step 1}
		\end{subfigure}%
		\begin{subfigure}[t]{0.33\textwidth}
			\centering
			\includegraphics[width=0.95\textwidth]{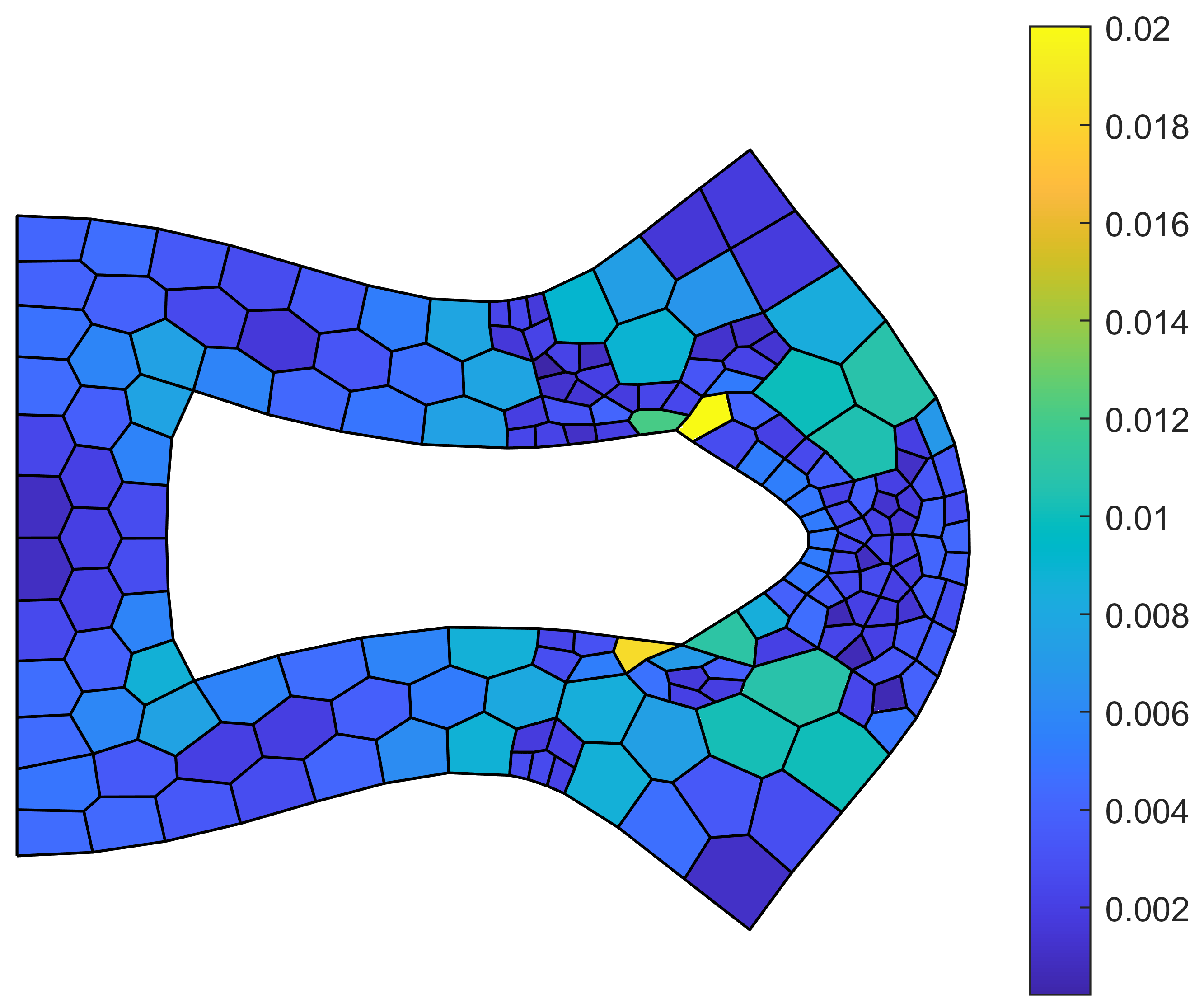}
			\caption{${\Upsilon_{\text{DB}}^{t}}$ - Step 2}
		\end{subfigure}%
		\begin{subfigure}[t]{0.33\textwidth}
			\centering
			\includegraphics[width=0.95\textwidth]{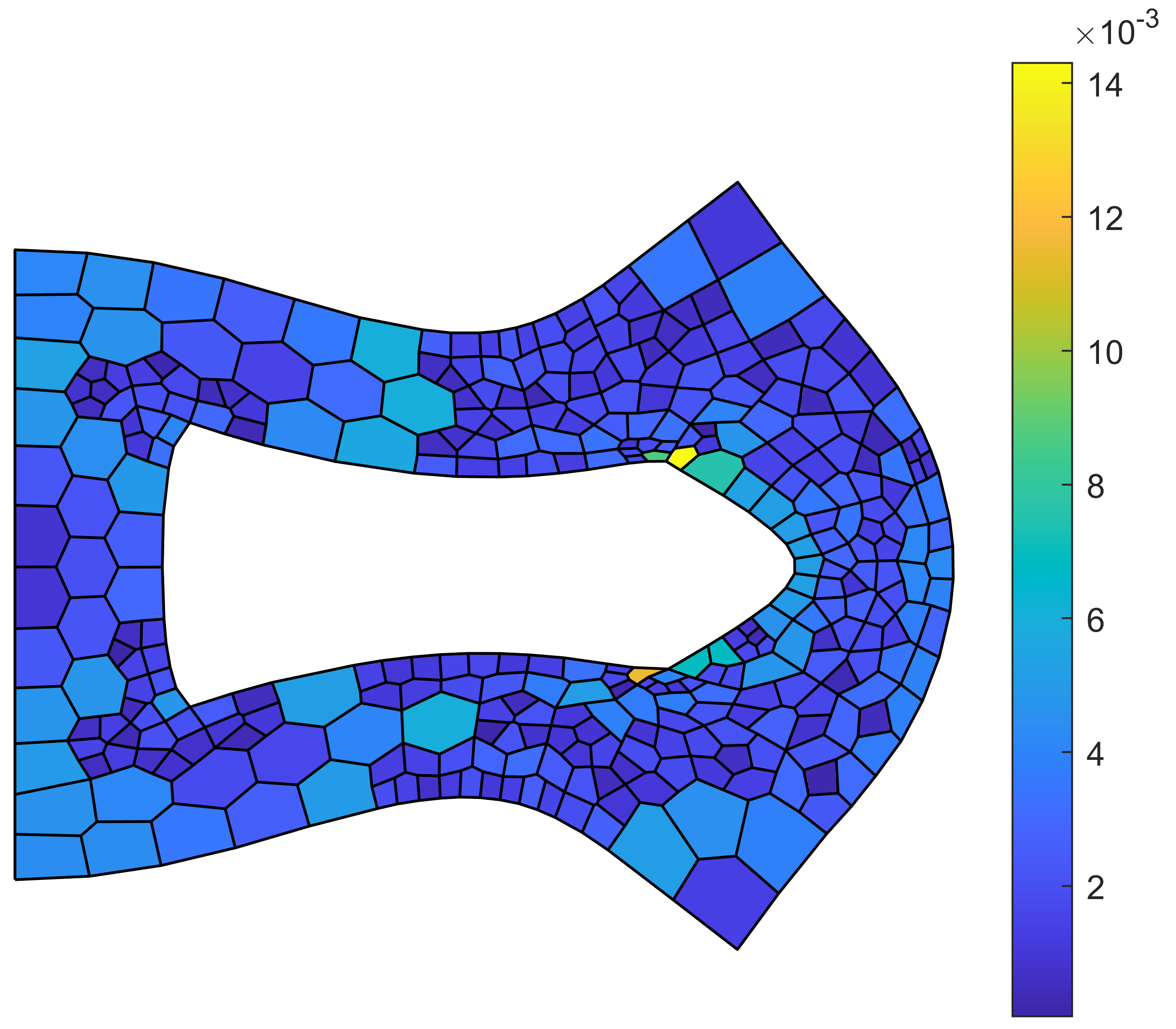}
			\caption{${\Upsilon_{\text{DB}}^{t}}$ - Step 3}
		\end{subfigure}
		\vskip \baselineskip 
		\begin{subfigure}[t]{0.33\textwidth}
			\centering
			\includegraphics[width=0.95\textwidth]{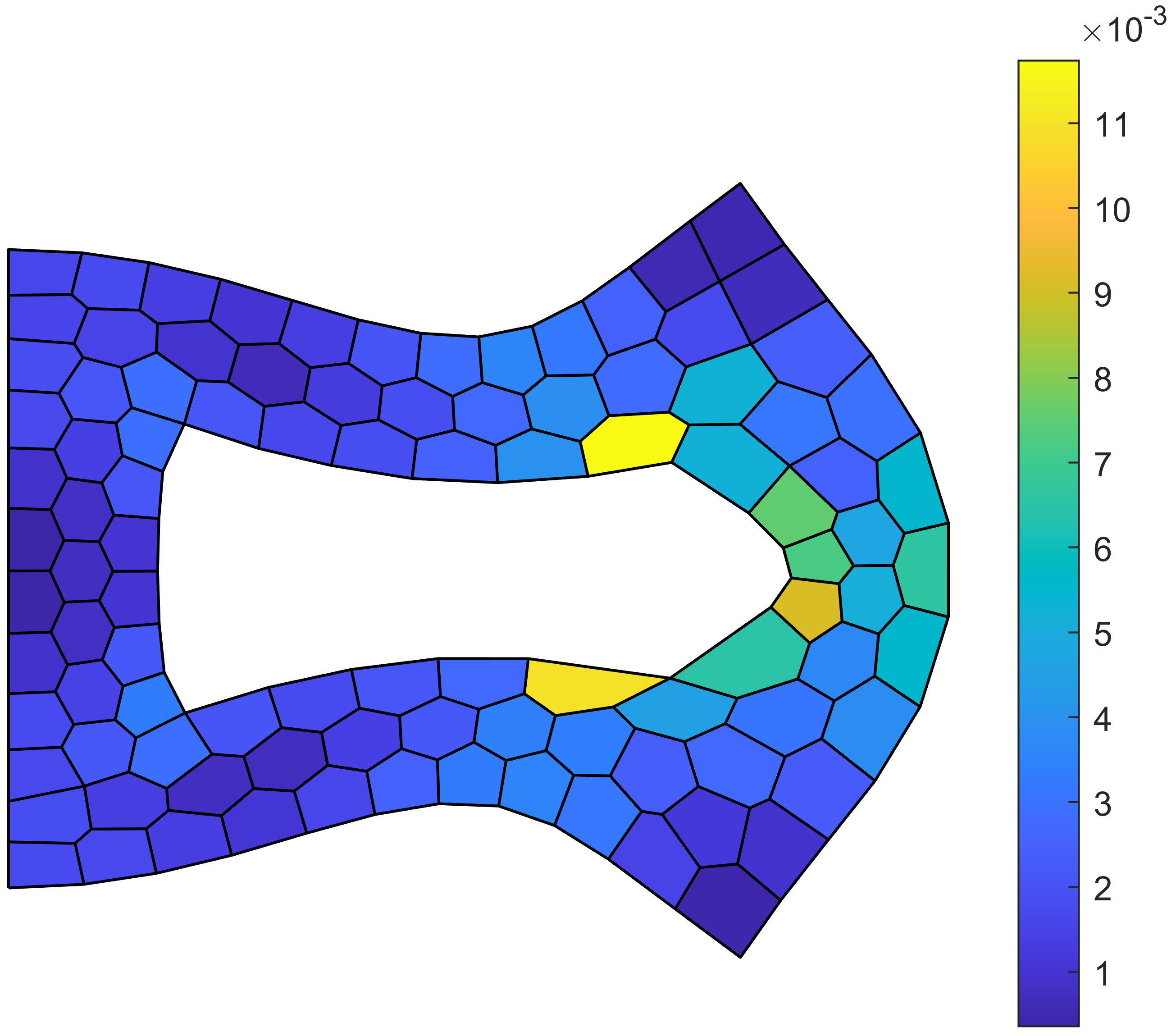}
			\caption{${\| \bK_{\rm s}^{E} \cdot \bd^{E} \|}$ - Step 1}
		\end{subfigure}%
		\begin{subfigure}[t]{0.33\textwidth}
			\centering
			\includegraphics[width=0.95\textwidth]{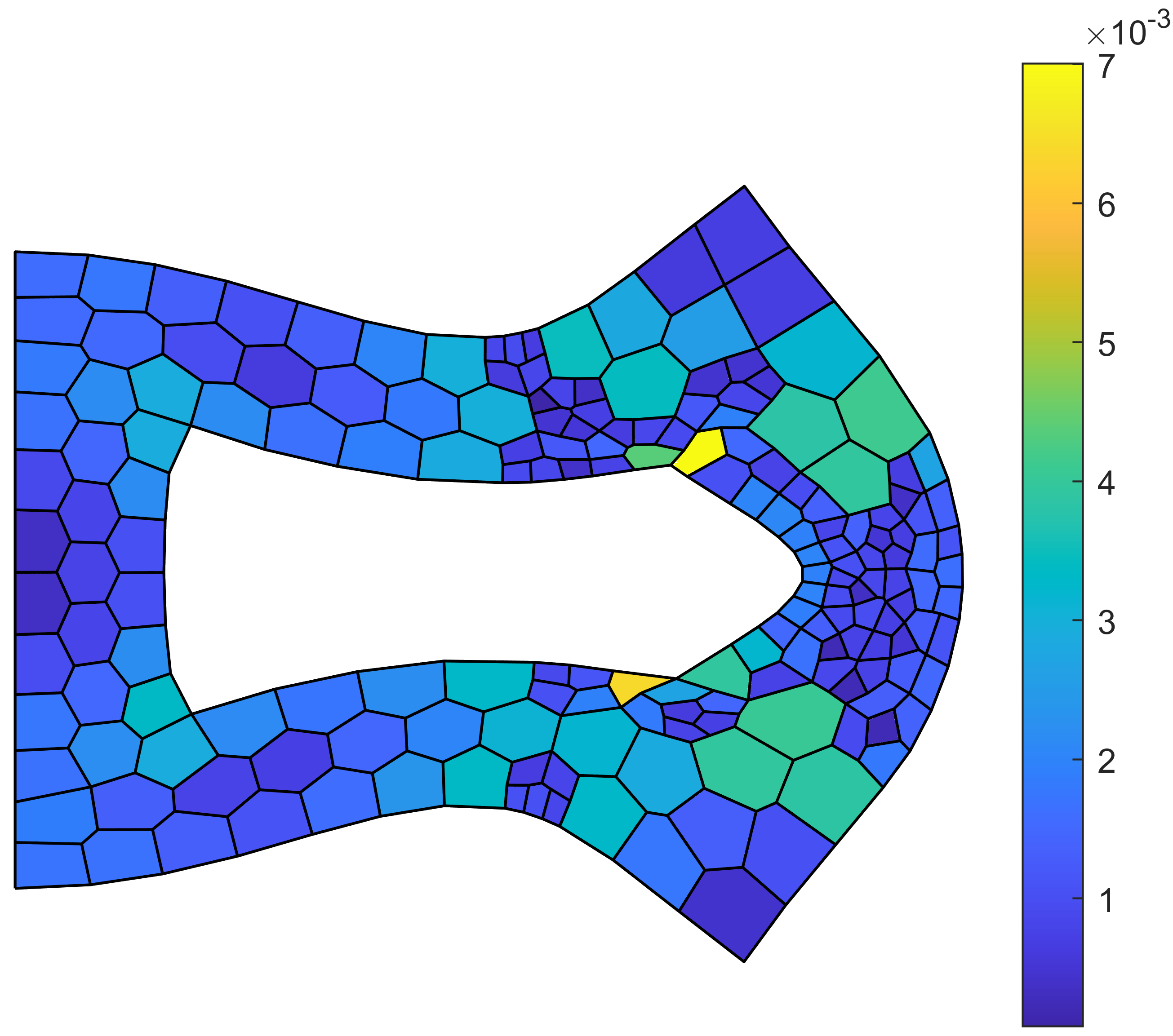}
			\caption{${\| \bK_{\rm s}^{E} \cdot \bd^{E} \|}$ - Step 2}
		\end{subfigure}%
		\begin{subfigure}[t]{0.33\textwidth}
			\centering
			\includegraphics[width=0.95\textwidth]{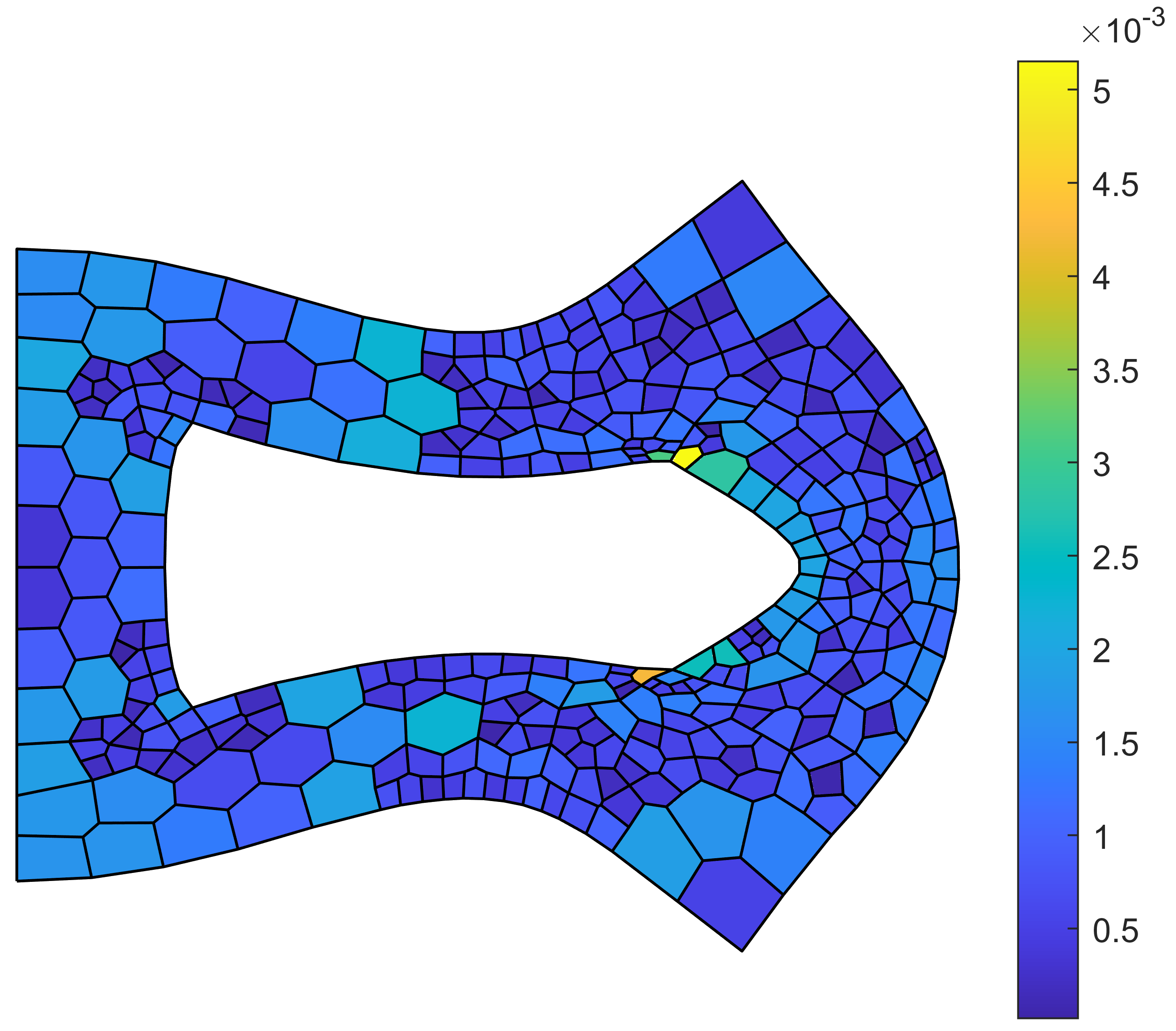}
			\caption{${\| \bK_{\rm s}^{E} \cdot \bd^{E} \|}$ - Step 3}
		\end{subfigure}
		\vskip \baselineskip 
		\caption{Displacement-based refinement indicator ${\Upsilon_{\text{DB}}^{t}}$ and norm of the residual of the stabilization term ${\| \bK_{\rm s}^{E} \cdot \bd^{E} \|}$ for problem~A(1) using the displacement-based refinement procedure with ${T=20\%}$ on Voronoi meshes with a compressible Poisson's ratio.
			\label{fig:PlateWithHoleTractionIndicatorVsStabNormVrn}}
	\end{figure}

	Figure~\ref{fig:PlateWithHoleTractionIndicatorVsStabNormVrn} depicts the evolution of the displacement-based refinement indicator ${\Upsilon_{\text{DB}}^{t}}$ and the norm of the residual of the stabilization term, i.e. ${\| \bK_{\rm s}^{E} \cdot \bd^{E} \|}$, for problem~A(1) using the displacement-based refinement procedure with ${T=20\%}$ on successive unstructured/Voronoi meshes with a compressible Poisson's ratio. 
	It is clear that the distribution of the refinement indicator and the residual term are qualitatively similar over the domain. The only difference between the indicator and the residual term is their respective magnitudes. This behaviour is expected, as mentioned in Section~\ref{subsec:DispBasedIndicator}, and demonstrates that the residual of the stabilization term could be used instead of the displacement-based refinement indicator. For corresponding figures for the case of structured meshes see \cite{vanHuyssteen2022,DanielHuyssteen2023}.
	
	\FloatBarrier
	\FloatBarrier
	
	\subsection{Problem B(4): Plate with notch \\ \cgray Displacement-based and strain jump-based refinement procedures \cblack}
	\label{subsec:PlateWithNotch}
	
\begin{figure}[ht!]
		\centering
		\begin{subfigure}[t]{0.45\textwidth}
			\centering
			\includegraphics[width=0.95\textwidth]{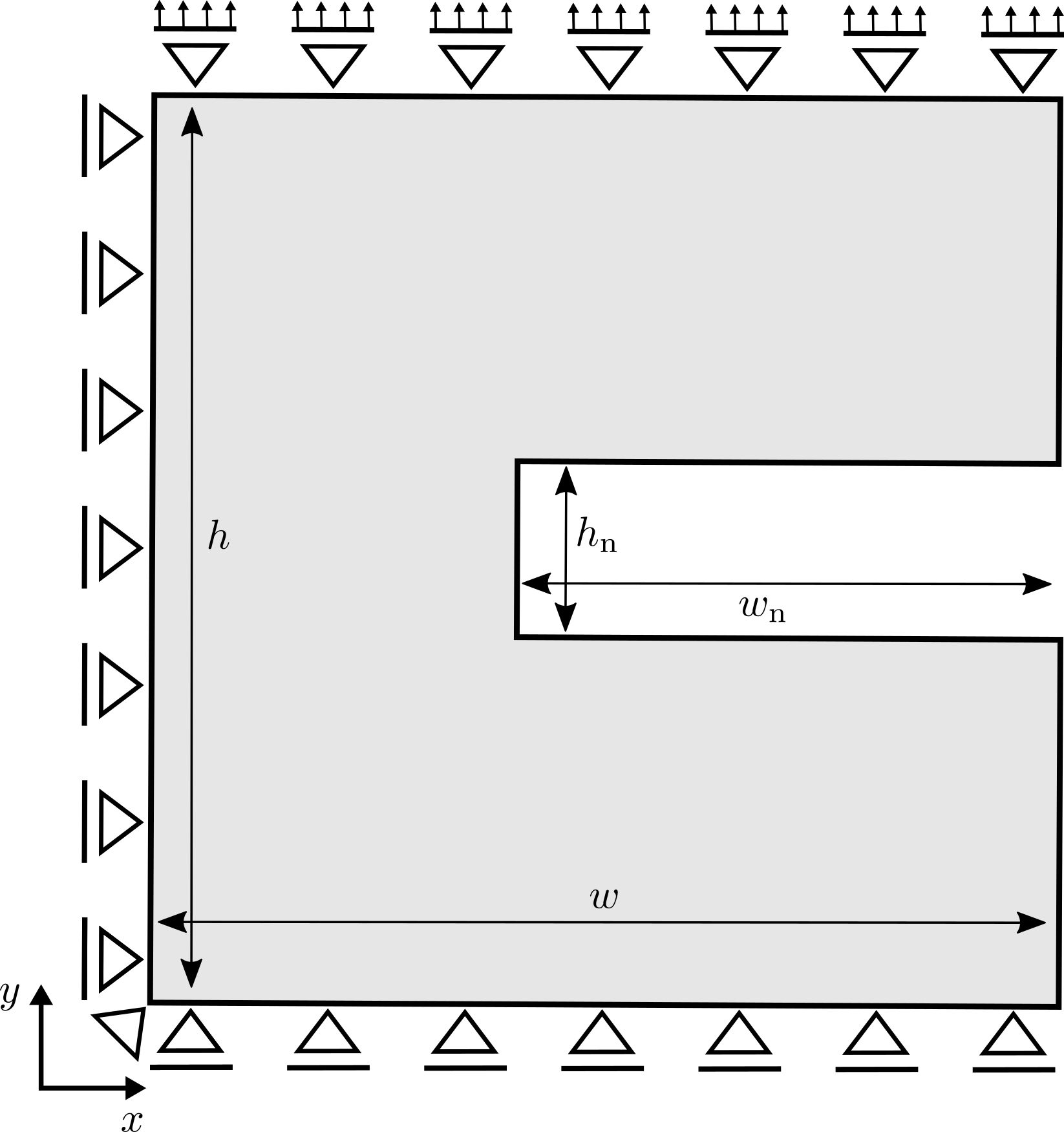}
			\caption{Problem geometry}
		\end{subfigure}%
		\begin{subfigure}[t]{0.55\textwidth}
			\centering
			\includegraphics[width=0.7\textwidth]{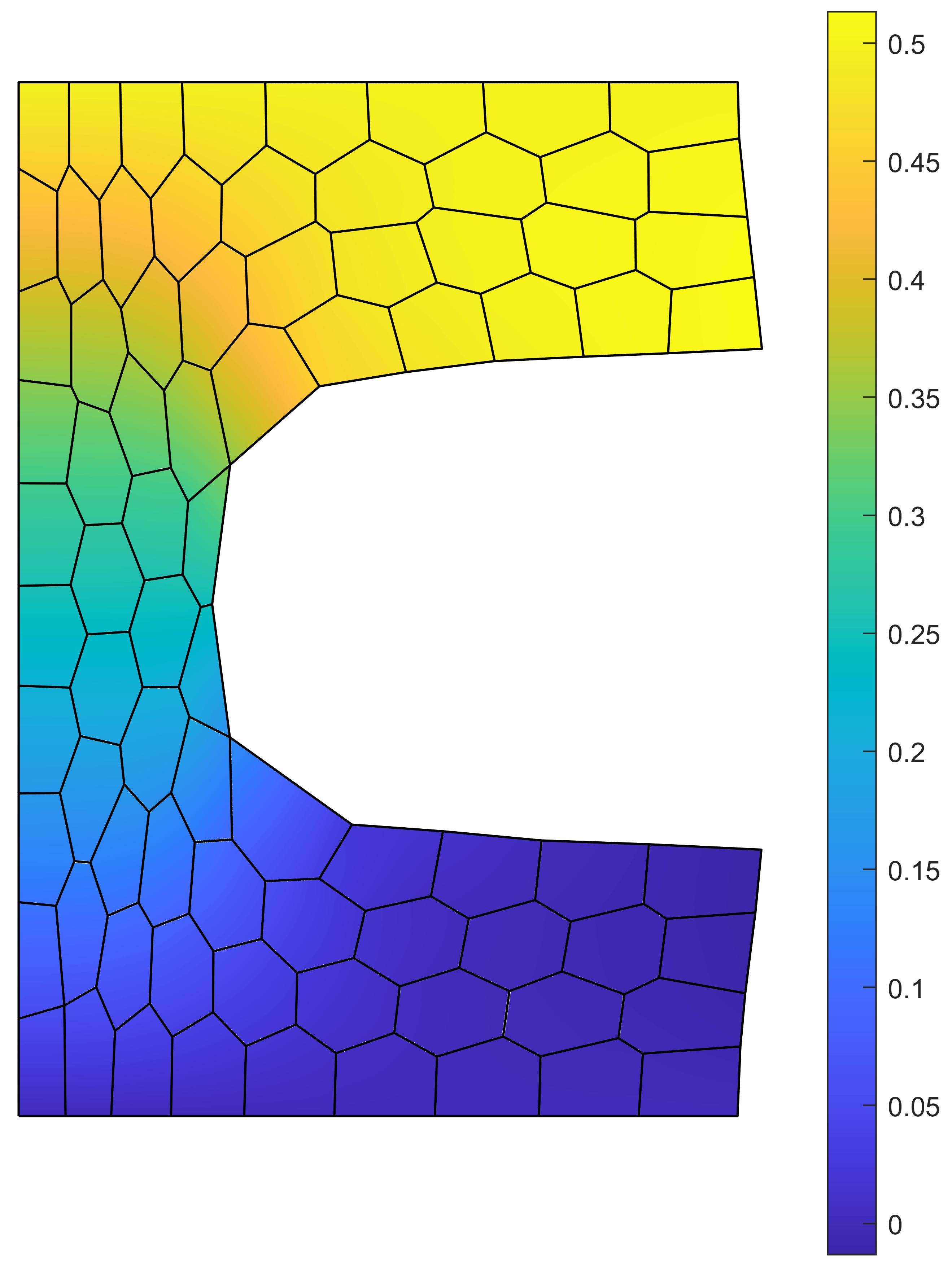}
			\caption{Deformed configuration}
		\end{subfigure}
		\caption{Problem~B(4) (a) geometry, and (b) sample deformed configuration of a Voronoi mesh with ${\nu=0.49995}$. 
			\label{fig:PlateWithNotchGeometry}}
	\end{figure}

	Problem~B(4) comprises a domain of width ${w=1~\rm{m}}$ and height ${h=1~\rm{m}}$ with a notch of width ${w_{\rm n}=0.6~w}$ and height ${w_{\rm n}=0.2~h}$.
	The bottom and left-hand edges of the domain are constrained vertically and horizontally respectively, and the bottom left-hand corner is fully constrained. 
	The top edge is subject to a prescribed displacement in the $y$-direction of ${\bar{u}_{y}=0.5~\rm{m}}$ while the displacement in the $x$-direction is unconstrained (see Figure~\ref{fig:PlateWithNotchGeometry}(a)).
	The results presented for this problem were generated using a combination of the displacement-based and strain jump-based refinement procedures.
	Figure~\ref{fig:PlateWithNotchGeometry}(b) depicts a sample deformed configuration of the plate with a Voronoi mesh and a nearly incompressible Poisson's ratio of ${\nu=0.49995}$. The vertical displacement ${u_{y}}$ is plotted on the colour axis.
	
	\FloatBarrier
	\FloatBarrier
	
\begin{figure}[ht!]
		\centering
		\begin{subfigure}[t]{0.33\textwidth}
			\centering
			\includegraphics[height=0.95\textwidth]{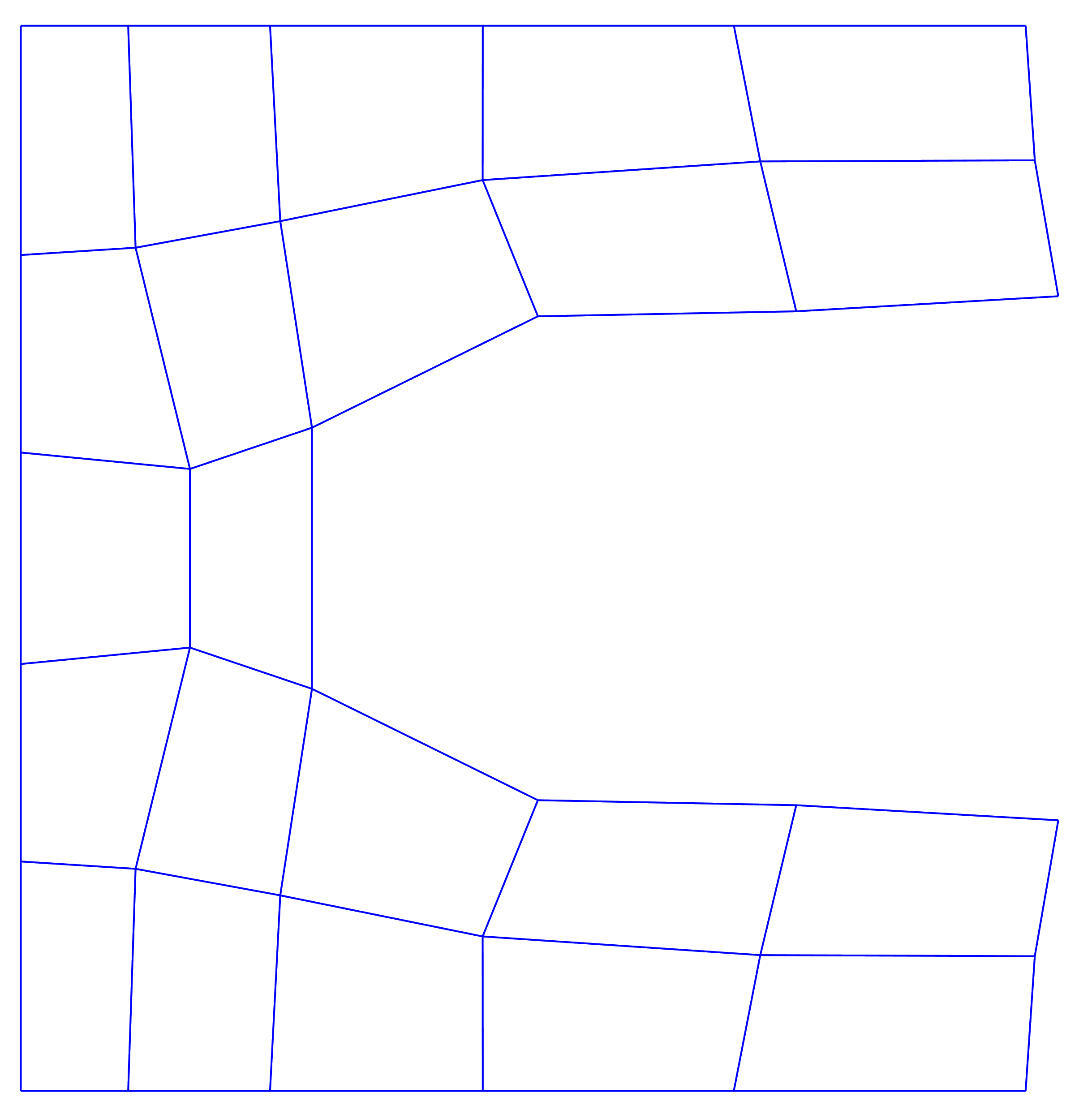}
			\caption{Nearly incompressible - Step 1}
		\end{subfigure}%
		\begin{subfigure}[t]{0.33\textwidth}
			\centering
			\includegraphics[height=0.95\textwidth]{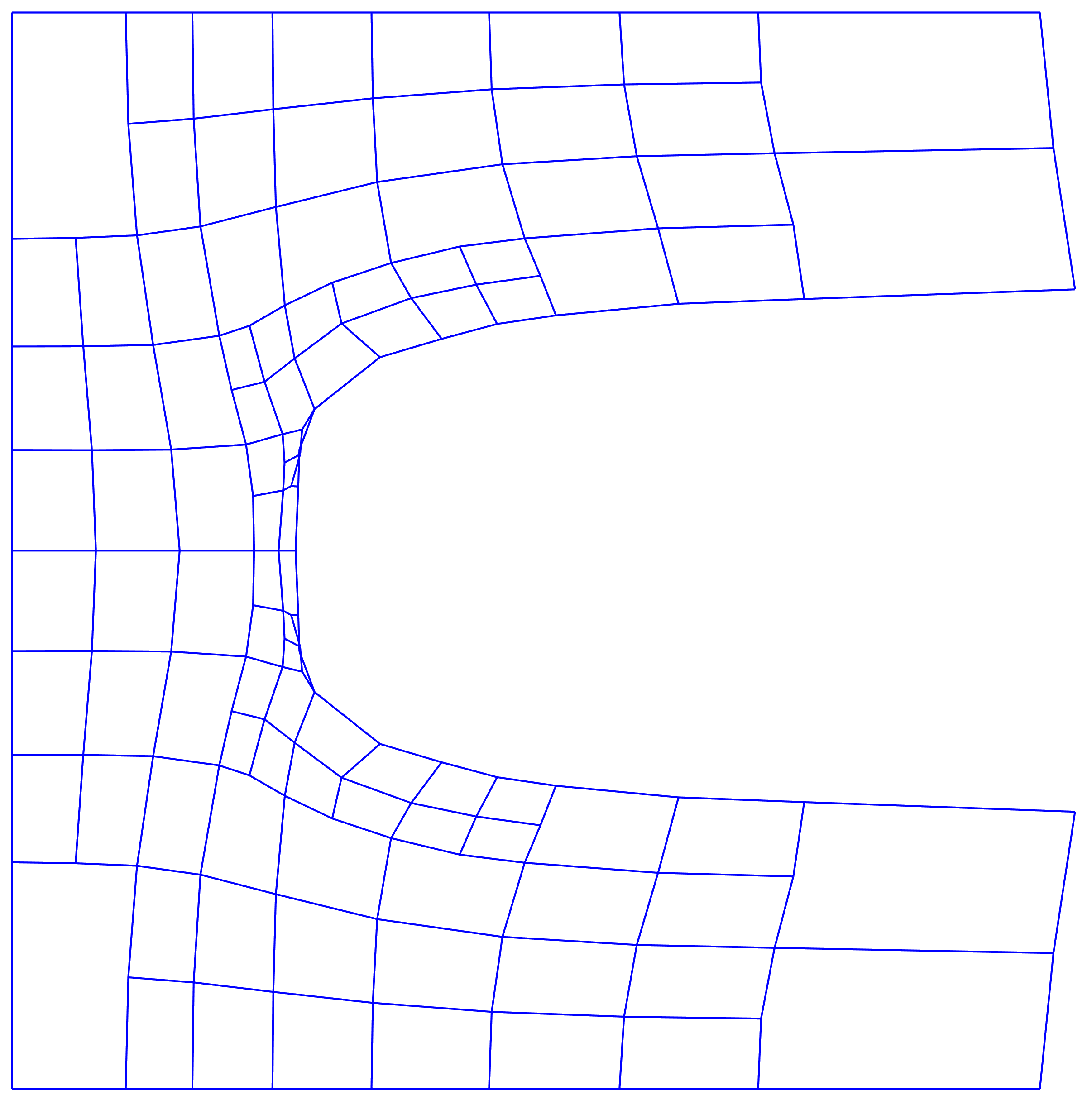}
			\caption{Nearly incompressible - Step 4}
		\end{subfigure}%
		\begin{subfigure}[t]{0.33\textwidth}
			\centering
			\includegraphics[height=0.95\textwidth]{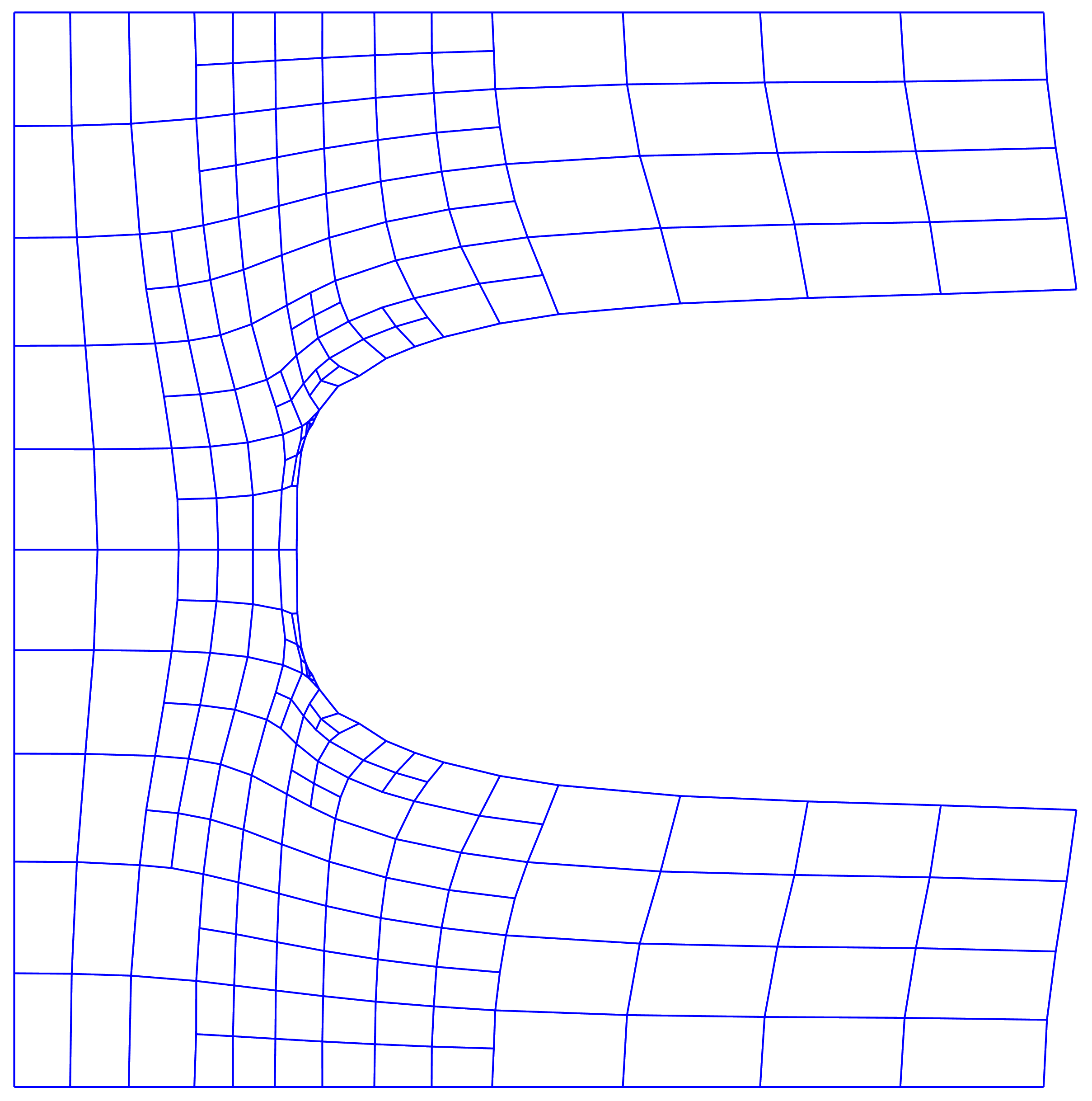}
			\caption{Nearly incompressible - Step 6}
		\end{subfigure}
		\vskip \baselineskip 
		\begin{subfigure}[t]{0.33\textwidth}
			\centering
			\includegraphics[height=0.95\textwidth]{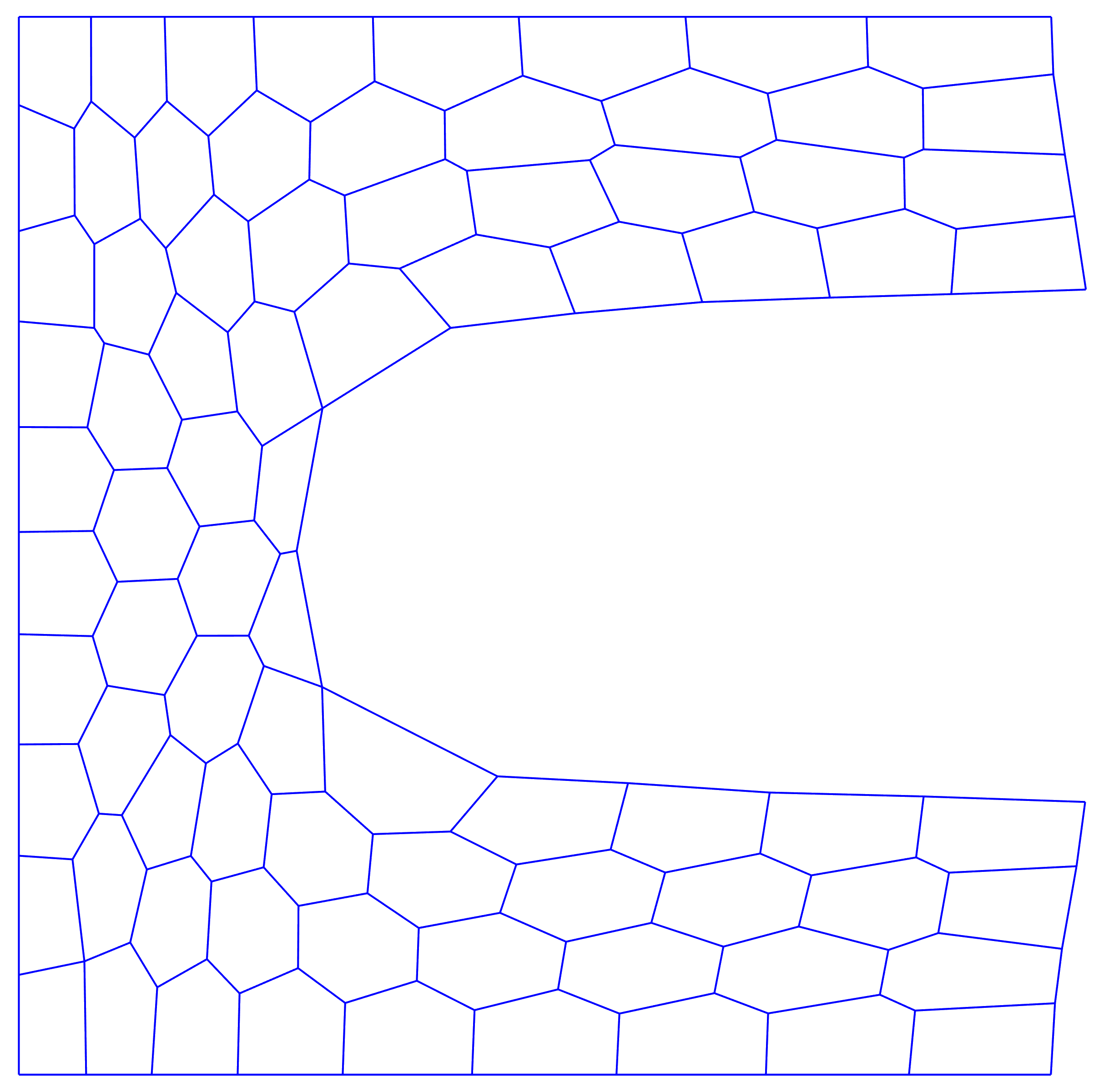}
			\caption{Nearly incompressible - Step 1}
		\end{subfigure}%
		\begin{subfigure}[t]{0.33\textwidth}
			\centering
			\includegraphics[height=0.95\textwidth]{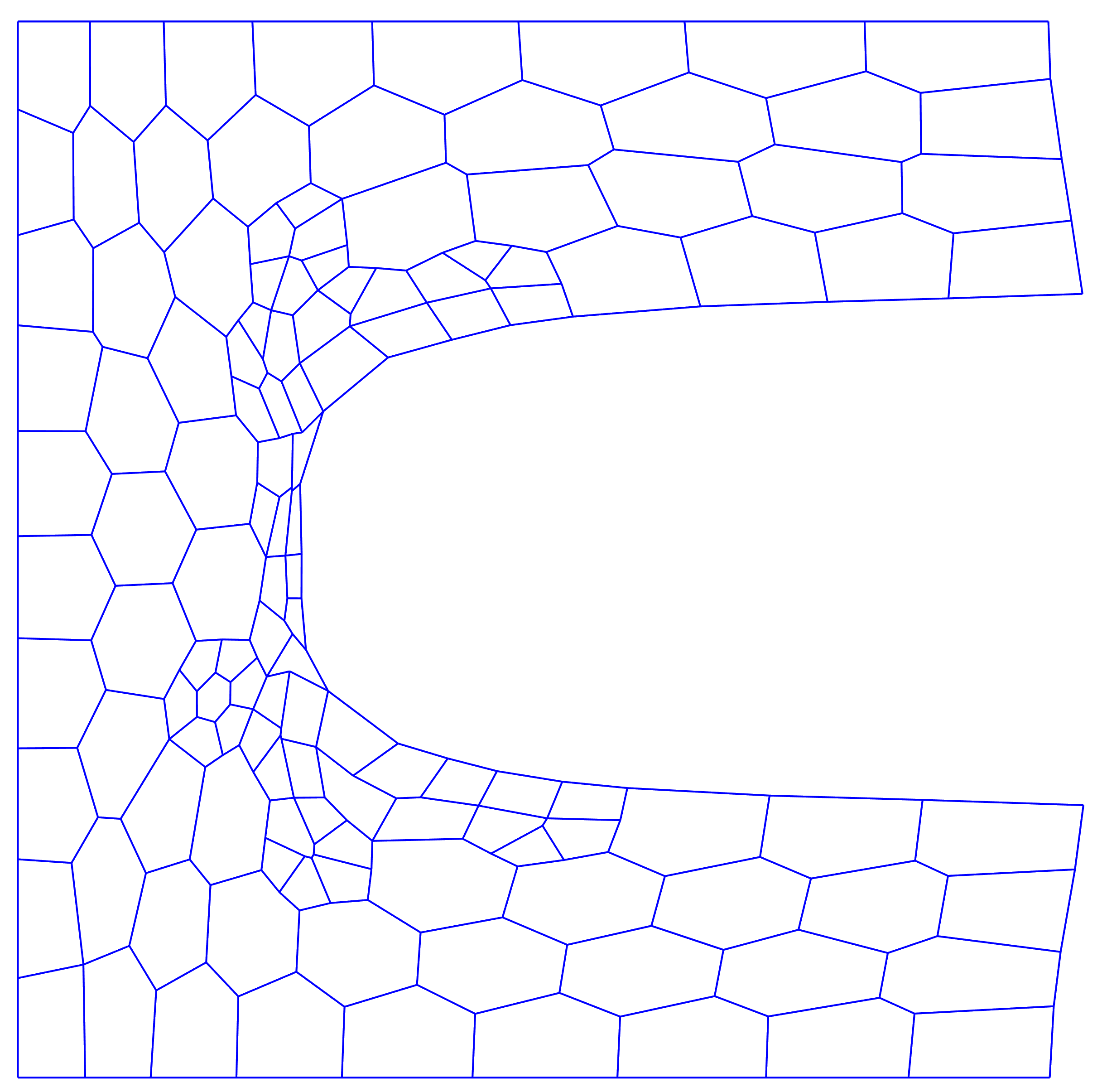}
			\caption{Nearly incompressible - Step 2}
		\end{subfigure}%
		\begin{subfigure}[t]{0.33\textwidth}
			\centering
			\includegraphics[height=0.95\textwidth]{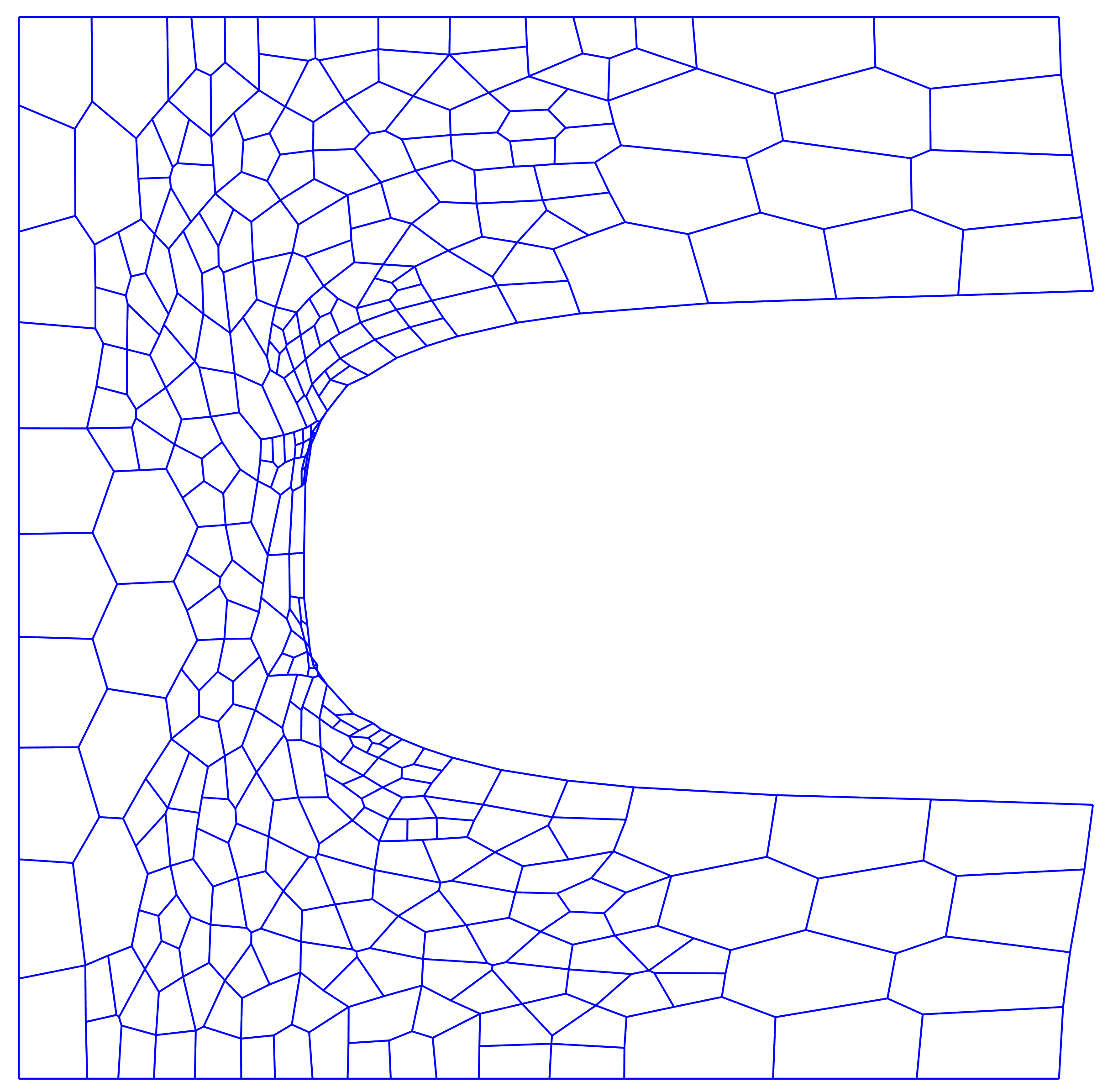}
			\caption{Nearly incompressible - Step 4}
		\end{subfigure}
		\caption{Mesh refinement process for problem~B(4) using a combination of the displacement-based and strain jump-based refinement procedures with ${T=15\%}$ for structured and Voronoi meshes with a nearly incompressible Poisson's ratio.
			\label{fig:PlateWithNotchMeshes}}
	\end{figure}

	Figure~\ref{fig:PlateWithNotchMeshes} depicts the mesh refinement process for problem~B(4) using a combination of the displacement-based and strain jump-based refinement procedures with ${T=15\%}$ for structured and Voronoi meshes with a nearly incompressible Poisson's ratio. 
	Meshes are shown at various refinement steps with step~1 corresponding to the initial mesh. Similar mesh refinement behaviour is observed for both mesh types with increased refinement developing the around the notch as expected. Additionally, the meshes in areas that experience relatively simple deformation, such as the left-hand and far right-hand portions of the domain, remain quite coarse. Thus, the mesh evolution is sensible for this problem.
	
	\FloatBarrier
	\FloatBarrier
	
\begin{figure}[ht!]
		\centering
		\begin{subfigure}[t]{0.5\textwidth}
			\centering
			\includegraphics[width=0.95\textwidth]{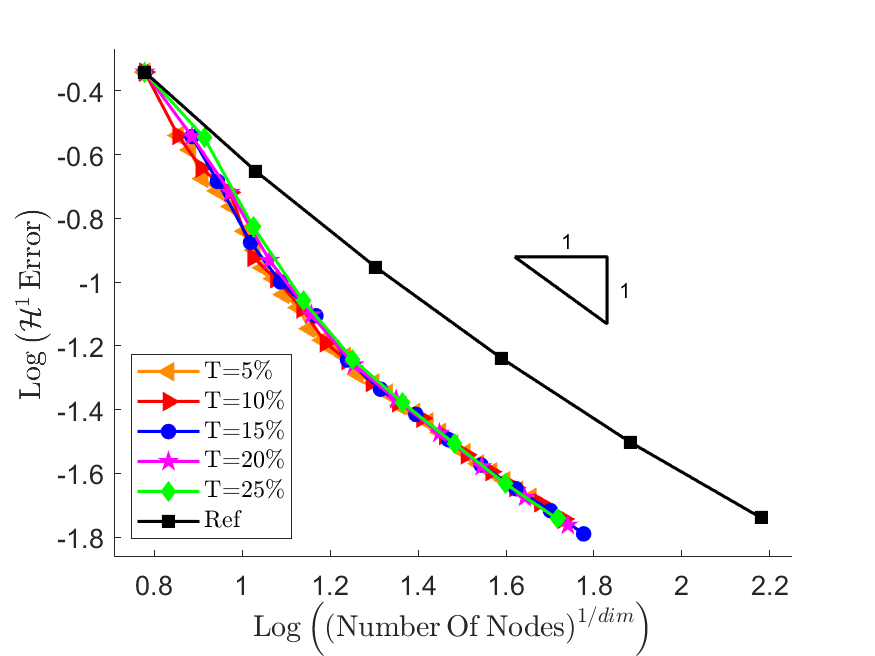}
			\caption{Structured mesh}
		\end{subfigure}%
		\begin{subfigure}[t]{0.5\textwidth}
			\centering
			\includegraphics[width=0.95\textwidth]{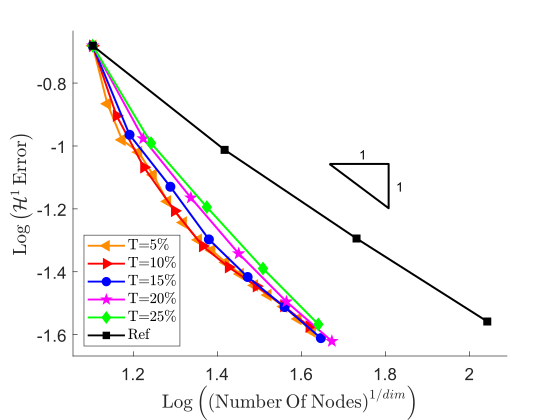}
			\caption{Voronoi mesh}
		\end{subfigure}
		\caption{$\mathcal{H}^{1}$ error vs $n_{\rm v}$ for problem~B(4) using a combination of the displacement-based and strain jump-based refinement procedures with a variety of choices of $T$ for structured and Voronoi meshes with a nearly incompressible Poisson's ratio.
			\label{fig:PlateWithNotchConvergenceNumberOfNodes}}
	\end{figure}

	The convergence behaviour in the ${\mathcal{H}^{1}}$ error norm of the VEM for problem~B(4) using a combination of the displacement-based and strain jump-based refinement procedures is depicted in Figure~\ref{fig:PlateWithNotchConvergenceNumberOfNodes} on a logarithmic scale for a variety of choices of $T$. Here, the ${\mathcal{H}^{1}}$ error is plotted against the number of vertices/nodes in the discretization for structured and unstructured/Voronoi meshes with a nearly incompressible Poisson's ratio. 
	The convergence behaviour is similar to that observed in Figure~\ref{fig:PlateWithHoleTractionConvergenceNumberOfNodes} for problem~A(1). The choice of $T$ has a smaller influence in the case of structured meshes than unstructured/Voronoi meshes where lower choices of $T$ exhibit a slightly faster initial convergence rate. For larger choices of $T$ the convergence rate is more consistent throughout the domain and for all choices of $T$ the adaptive refinement procedure exhibits a superior convergence rate to, and significantly outperforms, the reference refinement procedure for both mesh types.
	
	\FloatBarrier
	\FloatBarrier
	
	The performance of the VEM in terms of its convergence behaviour in the $\mathcal{H}^{1}$ error norm with respect to the number of vertices/nodes in the discretization when using a combination of the displacement-based and strain jump-based refinement procedures for problem~B(4) is summarized in Table~\ref{tab:PerformancePlateWithNotchNumberOfNodes}. Here, the performance, as measured by the PRE, is presented for the cases of structured and Voronoi meshes with compressible and nearly incompressible Poisson's ratios for a variety of choices of $T$. The results are, again, qualitatively similar to those presented in Table~\ref{tab:PerformancePlateWithHoleTractionNumberOfNodes} for problem~A(1). In the case of a compressible Poisson's ratio lower choices of $T$ perform well as a result of the more local/targeted refinement they offer. However, in the case of a nearly incompressible Poisson's ratio, and as discussed previously, the greater distribution of error over the domain requires increased refinement throughout the domain which results in all choices of $T$ exhibiting similar levels of performance. 
	Nevertheless, for both mesh types, and both choices of Poisson's ratio, the adaptive procedure significantly outperforms the reference procedure.
	
	\FloatBarrier
	\begin{table}[ht!]
		\centering 
		\begin{adjustbox}{max width=\textwidth}
			\begin{tabular}{|c|cccc|cccc|}
				\hline
				\multirow{3}{*}{Threshold} & \multicolumn{4}{c|}{Compressible}                                                                & \multicolumn{4}{c|}{Nearly-incompressible}                                                         \\ \cline{2-9} 
				& \multicolumn{2}{c|}{Structured}                           & \multicolumn{2}{c|}{Voronoi}         & \multicolumn{2}{c|}{Structured}                            & \multicolumn{2}{c|}{Voronoi}          \\ \cline{2-9} 
				& \multicolumn{1}{c|}{nNodes}   & \multicolumn{1}{c|}{PRE}  & \multicolumn{1}{c|}{nNodes}   & PRE  & \multicolumn{1}{c|}{nNodes}   & \multicolumn{1}{c|}{PRE}   & \multicolumn{1}{c|}{nNodes}   & PRE   \\ \hline
				T=5\%                      & \multicolumn{1}{c|}{692.67}   & \multicolumn{1}{c|}{3.02} & \multicolumn{1}{c|}{569.97}   & 4.68 & \multicolumn{1}{c|}{2857.62}  & \multicolumn{1}{c|}{12.45} & \multicolumn{1}{c|}{1574.33}  & 12.96 \\ \hline
				T=10\%                     & \multicolumn{1}{c|}{717.61}   & \multicolumn{1}{c|}{3.13} & \multicolumn{1}{c|}{646.88}   & 5.31 & \multicolumn{1}{c|}{2847.76}  & \multicolumn{1}{c|}{12.41} & \multicolumn{1}{c|}{1611.84}  & 13.27 \\ \hline
				T=15\%                     & \multicolumn{1}{c|}{726.99}   & \multicolumn{1}{c|}{3.17} & \multicolumn{1}{c|}{789.82}   & 6.48 & \multicolumn{1}{c|}{2810.36}  & \multicolumn{1}{c|}{12.25} & \multicolumn{1}{c|}{1585.57}  & 13.06 \\ \hline
				T=20\%                     & \multicolumn{1}{c|}{824.22}   & \multicolumn{1}{c|}{3.59} & \multicolumn{1}{c|}{1055.86}  & 8.67 & \multicolumn{1}{c|}{2687.43}  & \multicolumn{1}{c|}{11.71} & \multicolumn{1}{c|}{1726.08}  & 14.21 \\ \hline
				T=25\%                     & \multicolumn{1}{c|}{887.97}   & \multicolumn{1}{c|}{3.87} & \multicolumn{1}{c|}{1216.47}  & 9.99 & \multicolumn{1}{c|}{2669.84}  & \multicolumn{1}{c|}{11.64} & \multicolumn{1}{c|}{1849.84}  & 15.23 \\ \hline
				Ref                        & \multicolumn{1}{c|}{22945.00} & \multicolumn{1}{c|}{}     & \multicolumn{1}{c|}{12182.00} &      & \multicolumn{1}{c|}{22945.00} & \multicolumn{1}{c|}{}      & \multicolumn{1}{c|}{12144.00} &       \\ \hline
			\end{tabular}
		\end{adjustbox}
		\caption{Performance summary of the VEM in terms of its convergence behaviour in the $\mathcal{H}^{1}$ error norm with respect to the number of vertices/nodes in the discretization when using a combination of the displacement-based and strain jump-based refinement procedures for problem~B(4).
			\label{tab:PerformancePlateWithNotchNumberOfNodes}}
	\end{table}
	\FloatBarrier
	
\begin{figure}[ht!]
		\centering
		\begin{subfigure}[t]{0.5\textwidth}
			\centering
			\includegraphics[width=0.95\textwidth]{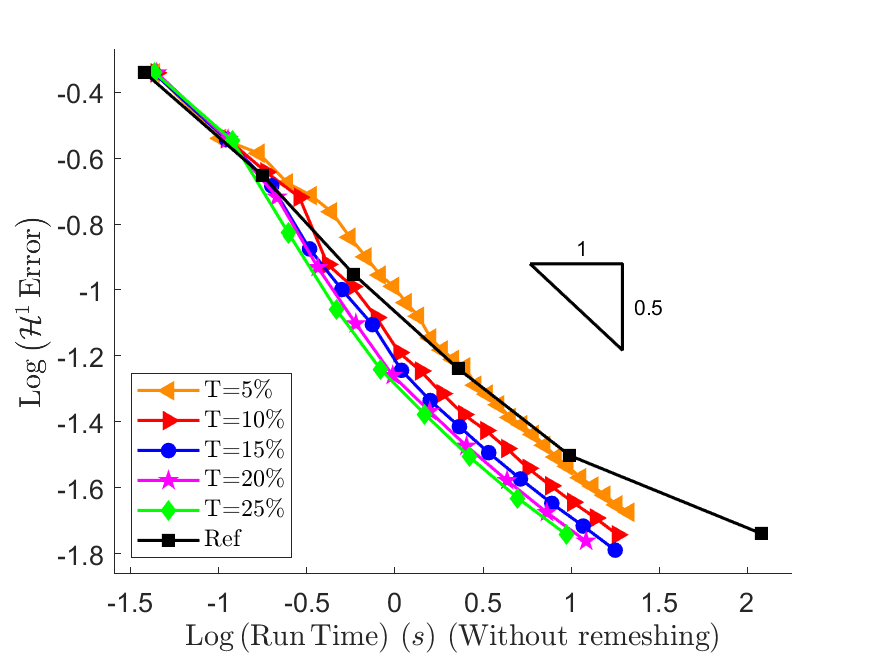}
			\caption{Structured mesh}
		\end{subfigure}%
		\begin{subfigure}[t]{0.5\textwidth}
			\centering
			\includegraphics[width=0.95\textwidth]{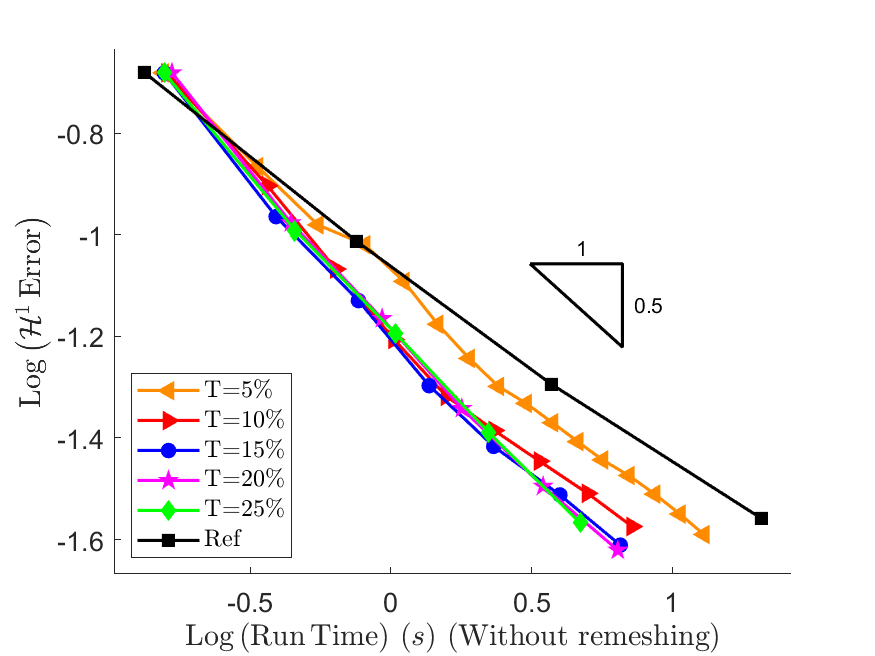}
			\caption{Voronoi mesh}
		\end{subfigure}
		\caption{$\mathcal{H}^{1}$ error vs run time (excluding remeshing time) for problem~B(4) using a combination of the displacement-based and strain jump-based refinement procedures with a variety of choices of $T$ for structured and Voronoi meshes with a nearly incompressible Poisson's ratio.
			\label{fig:PlateWithNotchConvergenceRunTime}}
	\end{figure}

	The convergence behaviour in the ${\mathcal{H}^{1}}$ error norm of the VEM for problem~B(4) using a combination of the displacement-based and strain jump-based refinement procedures is depicted in Figure~\ref{fig:PlateWithNotchConvergenceRunTime} on a logarithmic scale for a variety of choices of $T$. Here, the ${\mathcal{H}^{1}}$ error is plotted against run time (excluding remeshing time) for structured and Voronoi meshes with a nearly incompressible Poisson's ratio. 
	Similar to the behaviour observed in Figure~\ref{fig:PlateWithHoleTractionConvergenceRunTime} for problem~A(1), it is clear that lower choices of $T$ require significantly more remeshing steps, and consequently more run time, than larger values of $T$. Thus, lower choices of $T$ are not particularly efficient in terms of run time.
	However, and particularly in the fine mesh range, for both mesh types and for all choices of $T$ the adaptive procedure exhibits a superior convergence rate to the reference refinement procedure.
	
	\FloatBarrier
	\FloatBarrier
	
	The performance of the VEM in terms of its convergence behaviour in the $\mathcal{H}^{1}$ error norm with respect to run time (excluding remeshing time) when using a combination of the displacement-based and strain jump-based refinement procedures for problem~B(4) is summarized in Table~\ref{tab:PerformancePlateWithNotchRunTime}. Here, the performance, as measured by the PRE, is presented for the cases of structured and Voronoi meshes with compressible and nearly incompressible Poisson's ratios for a variety of choices of $T$. 
	The relative inefficiencies of lower choices of $T$, particularly in the case of near-incompressibility, are again easy to see. This further motivates determination a choice of $T$ that yields a balance of performance in terms of the number of nodes/vertices and run time.
	Nevertheless, independent of the degree of compressibility and mesh type the adaptive procedure, for all choices of $T$, represents a significant improvement in efficiency compared to the reference procedure.
	
	\FloatBarrier
	\begin{table}[ht!]
		\centering 
		\begin{adjustbox}{max width=\textwidth}
			\begin{tabular}{|c|cccc|cccc|}
				\hline
				\multirow{3}{*}{Threshold} & \multicolumn{4}{c|}{Compressible}                                                                 & \multicolumn{4}{c|}{Nearly-incompressible}                                                         \\ \cline{2-9} 
				& \multicolumn{2}{c|}{Structured}                           & \multicolumn{2}{c|}{Voronoi}          & \multicolumn{2}{c|}{Structured}                            & \multicolumn{2}{c|}{Voronoi}          \\ \cline{2-9} 
				& \multicolumn{1}{c|}{Run time} & \multicolumn{1}{c|}{PRE}  & \multicolumn{1}{c|}{Run time} & PRE   & \multicolumn{1}{c|}{Run time} & \multicolumn{1}{c|}{PRE}   & \multicolumn{1}{c|}{Run time} & PRE   \\ \hline
				T=5\%                      & \multicolumn{1}{c|}{6.34}     & \multicolumn{1}{c|}{5.15} & \multicolumn{1}{c|}{2.85}     & 13.74 & \multicolumn{1}{c|}{31.20}    & \multicolumn{1}{c|}{26.09} & \multicolumn{1}{c|}{11.08}    & 53.60 \\ \hline
				T=10\%                     & \multicolumn{1}{c|}{3.68}     & \multicolumn{1}{c|}{2.99} & \multicolumn{1}{c|}{2.03}     & 9.78  & \multicolumn{1}{c|}{17.91}    & \multicolumn{1}{c|}{14.98} & \multicolumn{1}{c|}{6.57}     & 31.81 \\ \hline
				T=15\%                     & \multicolumn{1}{c|}{2.81}     & \multicolumn{1}{c|}{2.28} & \multicolumn{1}{c|}{1.99}     & 9.58  & \multicolumn{1}{c|}{13.31}    & \multicolumn{1}{c|}{11.13} & \multicolumn{1}{c|}{5.01}     & 24.26 \\ \hline
				T=20\%                     & \multicolumn{1}{c|}{2.67}     & \multicolumn{1}{c|}{2.17} & \multicolumn{1}{c|}{2.47}     & 11.89 & \multicolumn{1}{c|}{10.61}    & \multicolumn{1}{c|}{8.87}  & \multicolumn{1}{c|}{4.72}     & 22.83 \\ \hline
				T=25\%                     & \multicolumn{1}{c|}{2.63}     & \multicolumn{1}{c|}{2.13} & \multicolumn{1}{c|}{2.61}     & 12.54 & \multicolumn{1}{c|}{9.20}     & \multicolumn{1}{c|}{7.70}  & \multicolumn{1}{c|}{4.52}     & 21.89 \\ \hline
				Ref                        & \multicolumn{1}{c|}{123.12}   & \multicolumn{1}{c|}{}     & \multicolumn{1}{c|}{20.78}    &       & \multicolumn{1}{c|}{119.58}   & \multicolumn{1}{c|}{}      & \multicolumn{1}{c|}{20.66}    &       \\ \hline
			\end{tabular}
		\end{adjustbox}
		\caption{Performance summary of the VEM in terms of its convergence behaviour in the $\mathcal{H}^{1}$ error norm with respect to run time (excluding remeshing time) when using a combination of the displacement-based and strain jump-based refinement procedures for problem~B(4).
			\label{tab:PerformancePlateWithNotchRunTime}}
	\end{table}
	\FloatBarrier
	
\begin{figure}[ht!]
		\centering
		\begin{subfigure}[t]{0.5\textwidth}
			\centering
			\includegraphics[width=0.95\textwidth]{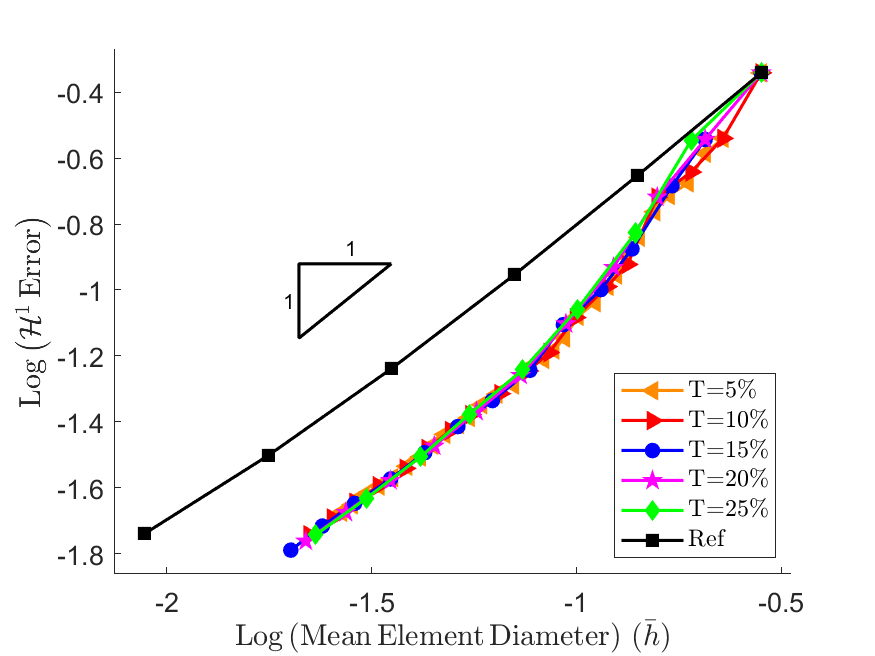}
			\caption{Structured mesh}
		\end{subfigure}%
		\begin{subfigure}[t]{0.5\textwidth}
			\centering
			\includegraphics[width=0.95\textwidth]{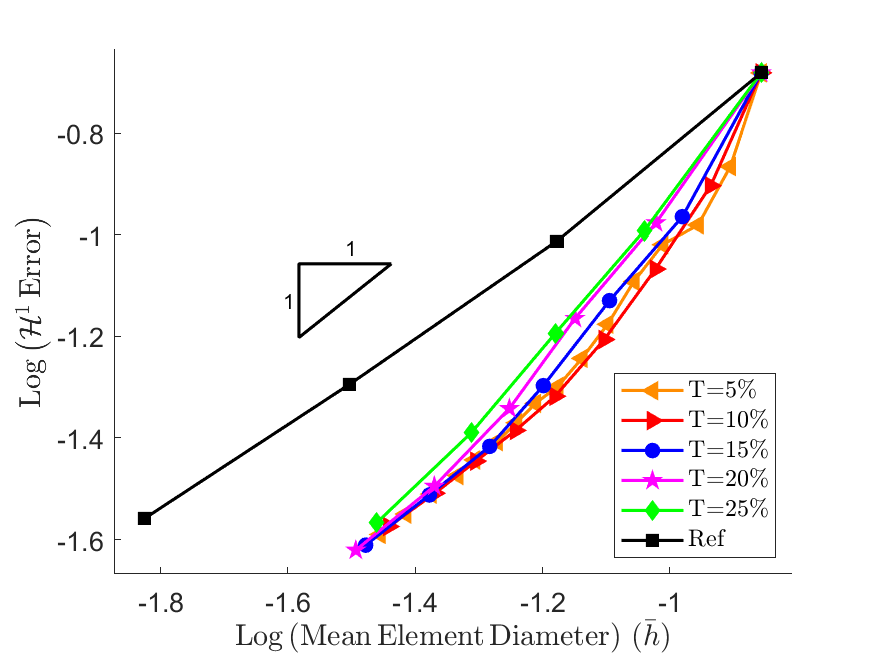}
			\caption{Voronoi mesh}
		\end{subfigure}
		\caption{$\mathcal{H}^{1}$ error vs mesh size for problem~B(4) using a combination of the displacement-based and strain jump-based refinement procedures with a variety of choices of $T$ for structured and Voronoi meshes with a nearly incompressible Poisson's ratio.
			\label{fig:PlateWithNotchConvergenceMeshSize}}
	\end{figure}

	The convergence behaviour in the ${\mathcal{H}^{1}}$ error norm of the VEM for problem~B(4) using a combination of the displacement-based and strain jump-based refinement procedures is depicted in Figure~\ref{fig:PlateWithNotchConvergenceMeshSize} on a logarithmic scale for a variety of choices of $T$. Here, the ${\mathcal{H}^{1}}$ error is plotted against mesh size as measured by the mean element diameter for structured and Voronoi meshes with a nearly incompressible Poisson's ratio. 
	The convergence behaviour is, again, qualitatively similar to that observed in Figure~\ref{fig:PlateWithNotchConvergenceNumberOfNodes} for problem~A(1). Specifically, the choice of $T$ has a smaller influence in the case of structured meshes than unstructured/Voronoi meshes where lower choices of $T$ exhibit a slightly faster initial convergence rate. For larger choices of $T$ the convergence rate is more consistent throughout the domain and for all choices of $T$ the adaptive refinement procedure exhibits a superior convergence rate to, and significantly outperforms, the reference refinement procedure for both mesh types.
	
	\FloatBarrier
	\FloatBarrier
	
	The performance of the VEM in terms of its convergence behaviour in the $\mathcal{H}^{1}$ error norm with respect to mesh size when using a combination of the displacement-based and strain jump-based refinement procedures for problem~B(4) is summarized in Table~\ref{tab:PerformancePlateWithNotchMeshSize}. Here, the performance, as measured by the PRE*, is presented for the cases of structured and Voronoi meshes with compressible and nearly incompressible Poisson's ratios for a variety of choices of $T$. 
	In the case of a compressible Poisson's ratio lower choices of $T$ perform well as a result of the more local/targeted refinement they offer. However, in the case of a nearly incompressible Poisson's ratio the greater distribution of error over the domain means that all choices of $T$ exhibit similar levels of performance. However, for both mesh types and for all choices of $T$ the adaptive procedure significantly outperforms the reference procedure.
	
	\FloatBarrier
	\begin{table}[ht!]
		\centering 
		\begin{adjustbox}{max width=\textwidth}
			\begin{tabular}{|c|cccc|cccc|}
				\hline
				\multirow{3}{*}{Threshold} & \multicolumn{4}{c|}{Compressible}                                                                    & \multicolumn{4}{c|}{Nearly-incompressible}                                                           \\ \cline{2-9} 
				& \multicolumn{2}{c|}{Structured}                             & \multicolumn{2}{c|}{Voronoi}           & \multicolumn{2}{c|}{Structured}                             & \multicolumn{2}{c|}{Voronoi}           \\ \cline{2-9} 
				& \multicolumn{1}{c|}{Mesh size} & \multicolumn{1}{c|}{PRE*}  & \multicolumn{1}{c|}{Mesh size} & PRE*  & \multicolumn{1}{c|}{Mesh size} & \multicolumn{1}{c|}{PRE*}  & \multicolumn{1}{c|}{Mesh size} & PRE*  \\ \hline
				T=5\%                      & \multicolumn{1}{c|}{0.048080}  & \multicolumn{1}{c|}{18.38} & \multicolumn{1}{c|}{0.063072}  & 23.73 & \multicolumn{1}{c|}{0.022569}  & \multicolumn{1}{c|}{39.16} & \multicolumn{1}{c|}{0.037727}  & 39.69 \\ \hline
				T=10\%                     & \multicolumn{1}{c|}{0.047191}  & \multicolumn{1}{c|}{18.73} & \multicolumn{1}{c|}{0.060153}  & 24.88 & \multicolumn{1}{c|}{0.022477}  & \multicolumn{1}{c|}{39.32} & \multicolumn{1}{c|}{0.037698}  & 39.72 \\ \hline
				T=15\%                     & \multicolumn{1}{c|}{0.047098}  & \multicolumn{1}{c|}{18.77} & \multicolumn{1}{c|}{0.055424}  & 27.01 & \multicolumn{1}{c|}{0.022696}  & \multicolumn{1}{c|}{38.94} & \multicolumn{1}{c|}{0.037649}  & 39.77 \\ \hline
				T=20\%                     & \multicolumn{1}{c|}{0.044310}  & \multicolumn{1}{c|}{19.95} & \multicolumn{1}{c|}{0.048939}  & 30.58 & \multicolumn{1}{c|}{0.023186}  & \multicolumn{1}{c|}{38.12} & \multicolumn{1}{c|}{0.036993}  & 40.48 \\ \hline
				T=25\%                     & \multicolumn{1}{c|}{0.042780}  & \multicolumn{1}{c|}{20.66} & \multicolumn{1}{c|}{0.045281}  & 33.05 & \multicolumn{1}{c|}{0.023335}  & \multicolumn{1}{c|}{37.88} & \multicolumn{1}{c|}{0.035225}  & 42.51 \\ \hline
				Ref                        & \multicolumn{1}{c|}{0.008839}  & \multicolumn{1}{c|}{}      & \multicolumn{1}{c|}{0.014967}  &       & \multicolumn{1}{c|}{0.008839}  & \multicolumn{1}{c|}{}      & \multicolumn{1}{c|}{0.014973}  &       \\ \hline
			\end{tabular}
		\end{adjustbox}
		\caption{Performance summary of the VEM in terms of its convergence behaviour in the $\mathcal{H}^{1}$ error norm with respect to mesh size when using a combination of the displacement-based and strain jump-based refinement procedures for problem~B(4).
			\label{tab:PerformancePlateWithNotchMeshSize}}
	\end{table}
	\FloatBarrier
	
\begin{figure}[ht!]
		\centering
		\begin{subfigure}[t]{0.5\textwidth}
			\centering
			\includegraphics[width=0.95\textwidth]{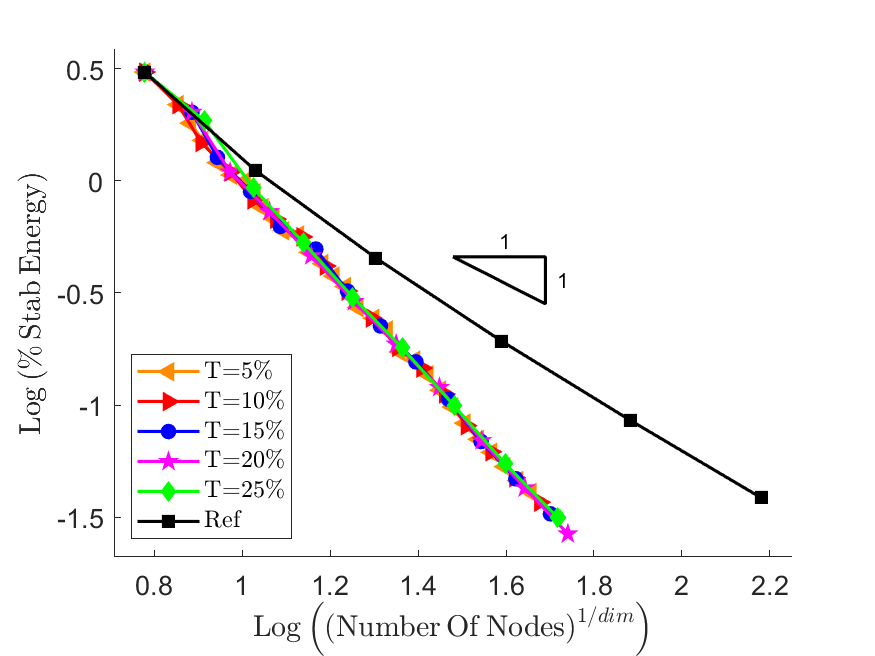}
			\caption{Structured mesh}
		\end{subfigure}%
		\begin{subfigure}[t]{0.5\textwidth}
			\centering
			\includegraphics[width=0.95\textwidth]{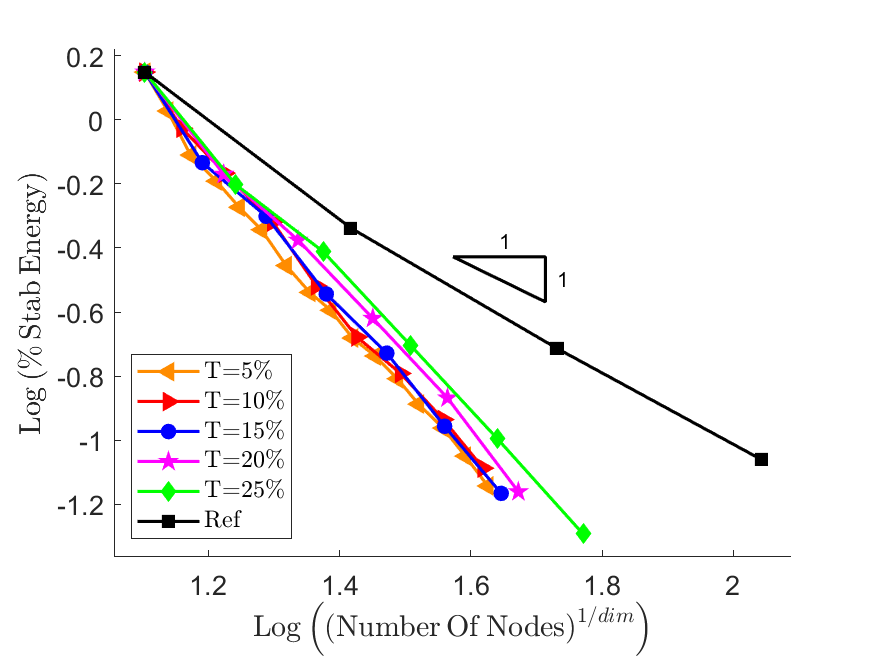}
			\caption{Voronoi mesh}
		\end{subfigure}
		\caption{PSE vs number of nodes/vertices in the discretization for problem~B(4) using a combination of the displacement-based and strain jump-based refinement procedures with a variety of choices of $T$ for structured and Voronoi meshes with a nearly incompressible Poisson's ratio.
			\label{fig:PlateWithNotchConvergencePSE}}
	\end{figure}

	The convergence behaviour in the PSE of the VEM for problem~B(4) using a combination of the displacement-based and strain jump-based refinement procedures is depicted in Figure~\ref{fig:PlateWithNotchConvergencePSE} on a logarithmic scale for a variety of choices of $T$. Here, the PSE is plotted against the number of nodes/vertices in the discretization for structured and Voronoi meshes with a nearly incompressible Poisson's ratio. 
	Similar to the behaviour observed in Figure~\ref{fig:PlateWithHoleTractionConvergencePSE} for problem~A(1), the convergence of the PSE is very similar for all choices of $T$. In the case of structured meshes the adaptive procedure initially exhibits similar convergence behaviour to the reference procedure. However, the performance of the adaptive procedure quickly improves and exhibits superior convergence behaviour to the reference procedure. In the case of unstructured/Voronoi meshes the adaptive procedure exhibits superior convergence behaviour to the reference procedure throughout the domain. 
	
	\FloatBarrier
	\FloatBarrier
	
	The performance of the VEM in terms of its convergence behaviour in the PSE with respect to the number of nodes/vertices in the discretization when using a combination of the displacement-based and strain jump-based refinement procedures for problem~B(4) is summarized in Table~\ref{tab:PerformancePlateWithNotchPSE}. Here, the performance, as measured by the PRE, is presented for the cases of structured and Voronoi meshes with compressible and nearly incompressible Poisson's ratios for a variety of choices of $T$. 
	In the case of unstructured/Voronoi meshes the results show that lower choices of $T$ perform slightly better than larger choices of $T$. However, the influence of the choice of $T$ is relatively small. 
	As observed in Table~\ref{tab:PerformancePlateWithHoleTractionPSE} for problem~A(1), the results do not show any significant dependence on the degree of compressibility. Thus, indicating that the meshes generated by the adaptive procedure are equally well-suited to cases of compressibility and near incompressibility.
	
	\FloatBarrier
	\begin{table}[ht!]
		\centering 
		\begin{adjustbox}{max width=\textwidth}
			\begin{tabular}{|c|cccc|cccc|}
				\hline
				\multirow{3}{*}{Threshold} & \multicolumn{4}{c|}{Compressible}                                                                  & \multicolumn{4}{c|}{Nearly-incompressible}                                                        \\ \cline{2-9} 
				& \multicolumn{2}{c|}{Structured}                            & \multicolumn{2}{c|}{Voronoi}          & \multicolumn{2}{c|}{Structured}                           & \multicolumn{2}{c|}{Voronoi}          \\ \cline{2-9} 
				& \multicolumn{1}{c|}{nNodes}   & \multicolumn{1}{c|}{PRE}   & \multicolumn{1}{c|}{nNodes}   & PRE   & \multicolumn{1}{c|}{nNodes}   & \multicolumn{1}{c|}{PRE}  & \multicolumn{1}{c|}{nNodes}   & PRE   \\ \hline
				T=5\%                      & \multicolumn{1}{c|}{2369.17}  & \multicolumn{1}{c|}{10.33} & \multicolumn{1}{c|}{1608.08}  & 13.20 & \multicolumn{1}{c|}{2168.02}  & \multicolumn{1}{c|}{9.45} & \multicolumn{1}{c|}{1551.00}  & 12.77 \\ \hline
				T=10\%                     & \multicolumn{1}{c|}{2383.95}  & \multicolumn{1}{c|}{10.39} & \multicolumn{1}{c|}{1681.26}  & 13.80 & \multicolumn{1}{c|}{2153.70}  & \multicolumn{1}{c|}{9.39} & \multicolumn{1}{c|}{1642.74}  & 13.53 \\ \hline
				T=15\%                     & \multicolumn{1}{c|}{2347.07}  & \multicolumn{1}{c|}{10.23} & \multicolumn{1}{c|}{1714.82}  & 14.08 & \multicolumn{1}{c|}{2136.65}  & \multicolumn{1}{c|}{9.31} & \multicolumn{1}{c|}{1603.92}  & 13.21 \\ \hline
				T=20\%                     & \multicolumn{1}{c|}{2342.70}  & \multicolumn{1}{c|}{10.21} & \multicolumn{1}{c|}{1933.43}  & 15.87 & \multicolumn{1}{c|}{2125.12}  & \multicolumn{1}{c|}{9.26} & \multicolumn{1}{c|}{1863.62}  & 15.35 \\ \hline
				T=25\%                     & \multicolumn{1}{c|}{2444.67}  & \multicolumn{1}{c|}{10.65} & \multicolumn{1}{c|}{2089.64}  & 17.15 & \multicolumn{1}{c|}{2220.72}  & \multicolumn{1}{c|}{9.68} & \multicolumn{1}{c|}{2175.37}  & 17.91 \\ \hline
				Ref                        & \multicolumn{1}{c|}{22945.00} & \multicolumn{1}{c|}{}      & \multicolumn{1}{c|}{12182.00} &       & \multicolumn{1}{c|}{22945.00} & \multicolumn{1}{c|}{}     & \multicolumn{1}{c|}{12144.00} &       \\ \hline
			\end{tabular}
		\end{adjustbox}
		\caption{Performance summary of the VEM in terms of its convergence behaviour in the PSE with respect to the number of nodes/vertices in the discretization when using a combination of the displacement-based and strain jump-based refinement procedures for problem~B(4).
			\label{tab:PerformancePlateWithNotchPSE}}
	\end{table}
	\FloatBarrier
	
\begin{figure}[ht!]
		\centering
		\begin{subfigure}[t]{0.5\textwidth}
			\centering
			\includegraphics[width=0.8\textwidth]{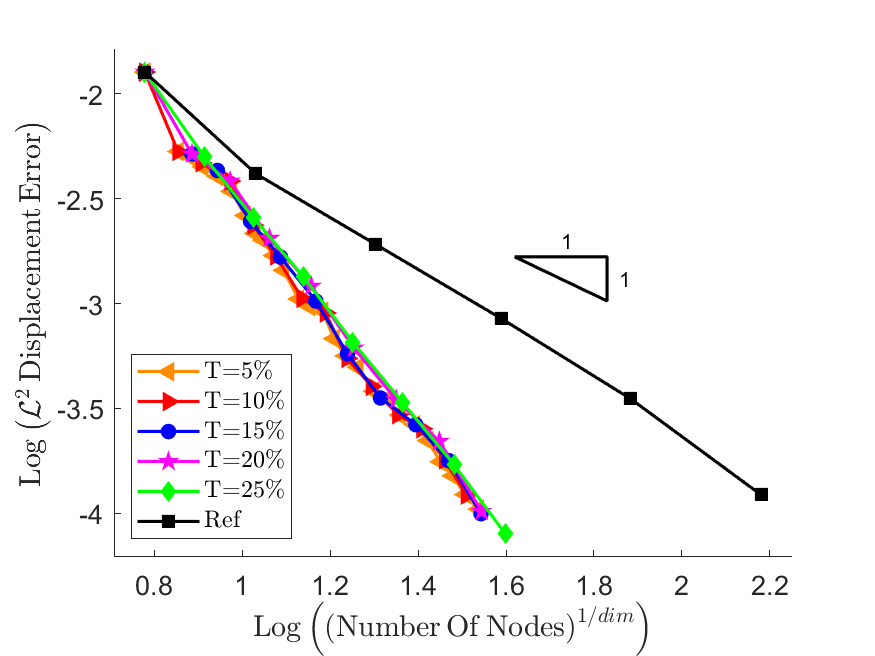}
			\caption{Displacement error - Structured meshes }
		\end{subfigure}%
		\begin{subfigure}[t]{0.5\textwidth}
			\centering
			\includegraphics[width=0.8\textwidth]{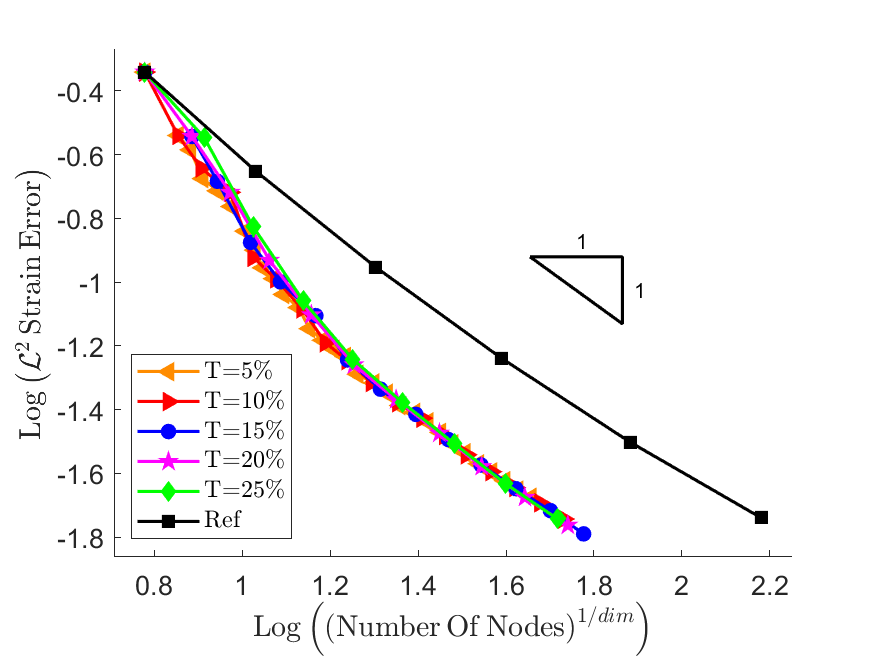}
			\caption{Strain error - Structured meshes }
		\end{subfigure}
		\vskip \baselineskip 
		\vspace*{-3mm}
		\begin{subfigure}[t]{0.5\textwidth}
			\centering
			\includegraphics[width=0.8\textwidth]{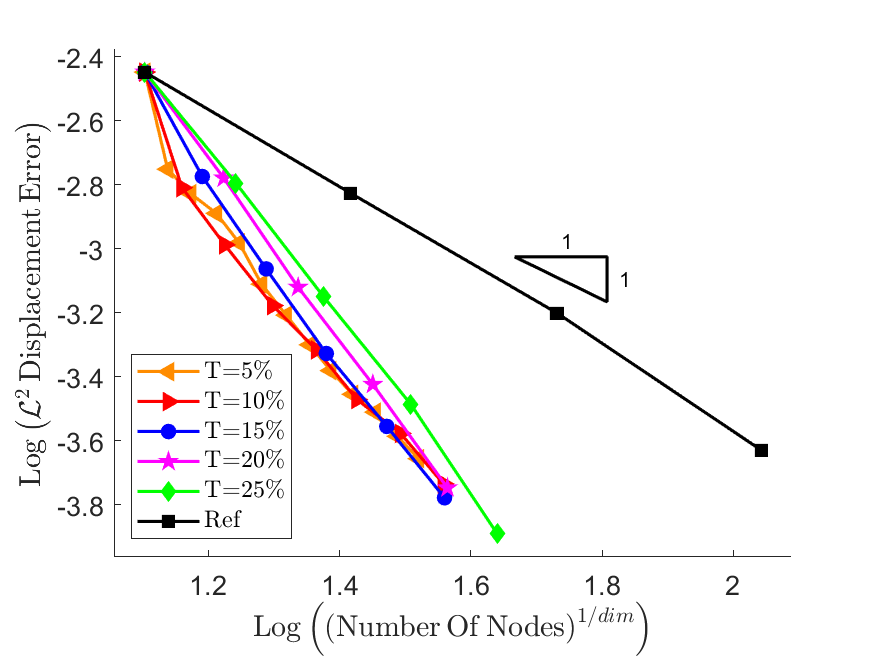}
			\caption{Displacement error - Voronoi meshes }
		\end{subfigure}%
		\begin{subfigure}[t]{0.5\textwidth}
			\centering
			\includegraphics[width=0.8\textwidth]{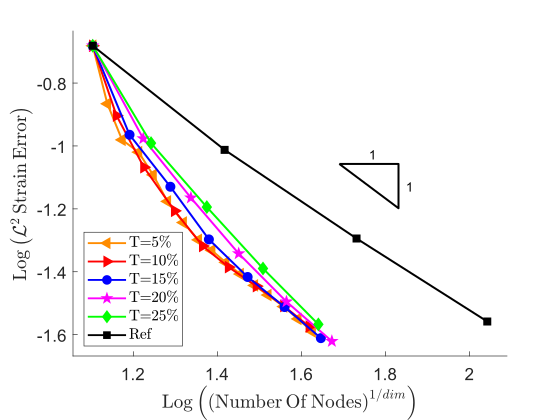}
			\caption{Strain error - Voronoi meshes }
		\end{subfigure}
		\caption{Displacement and strain $\mathcal{L}^{2}$ error components vs $n_{\rm v}$ for problem~B(4) using a combination of the displacement-based and strain jump-based refinement procedures with a variety of choices of $T$ for structured and Voronoi meshes with a nearly incompressible Poisson's ratio.
			\label{fig:PlateWithNotchConvergenceComponents}}
	\end{figure}

	The convergence behaviour in the ${\mathcal{L}^{2}}$ error norm of the displacement and strain field approximations when using a combination of the displacement-based and strain jump-based refinement procedures for problem~B(4) is plotted in Figure~\ref{fig:PlateWithNotchConvergenceComponents} against the number of vertices/nodes in the discretization. Here, plots are presented for the cases of structured and Voronoi meshes with a nearly incompressible Poisson's ratio. 
	As observed in Figure~\ref{fig:PlateWithHoleTractionConvergenceComponents} for problem~A(1), the convergence behaviour in both error components and for both mesh types is similar to the convergence behaviour observed in the ${\mathcal{H}^{1}}$ error norm.
	For structured meshes the choice of $T$ has a comparatively small influence with similar convergence behaviour exhibited by all choices of $T$.
	In the case of unstructured/Voronoi meshes the convergence rate exhibited by lower choices of $T$ is initially slightly faster, while larger values of $T$ exhibit more consistent convergence behaviour throughout the domain.
	Overall, for all choices of $T$ the adaptive procedure exhibits superior convergence behaviour to the reference procedure in both the displacement and strain error components for both mesh types.
	
	\FloatBarrier
	\FloatBarrier
	
	\subsection{Problem C(5): Narrow punch \\ \cgray Displacement-based and ${Z^{2}}$-like refinement procedures \cblack}
	\label{subsec:SquarePlateNarrowPunch}
	
\begin{figure}[ht!]
		\centering
		\begin{subfigure}[t]{0.45\textwidth}
			\centering
			\includegraphics[width=0.95\textwidth]{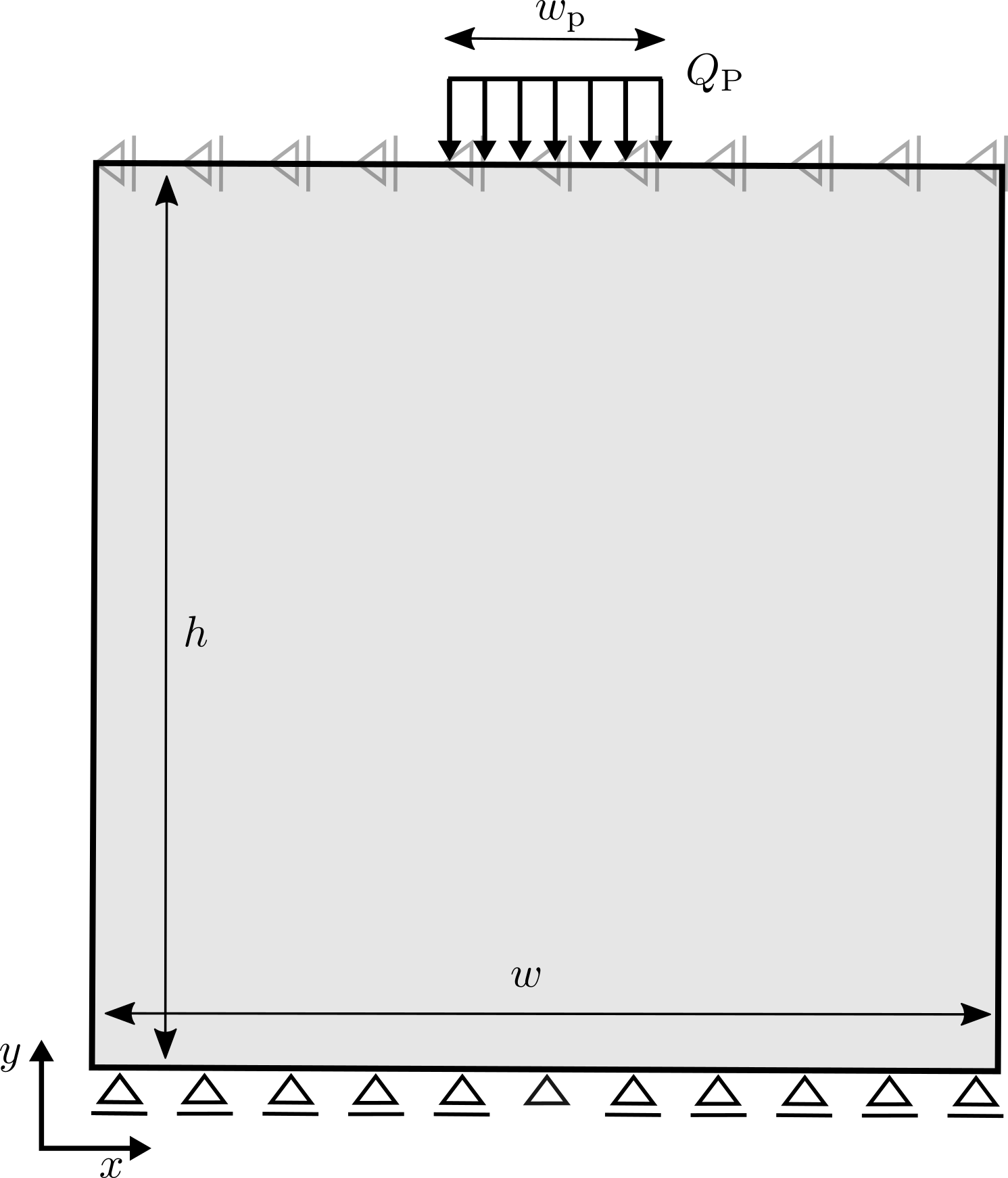}
			\caption{Problem geometry}
		\end{subfigure}%
		\begin{subfigure}[t]{0.55\textwidth}
			\centering
			\includegraphics[width=0.95\textwidth]{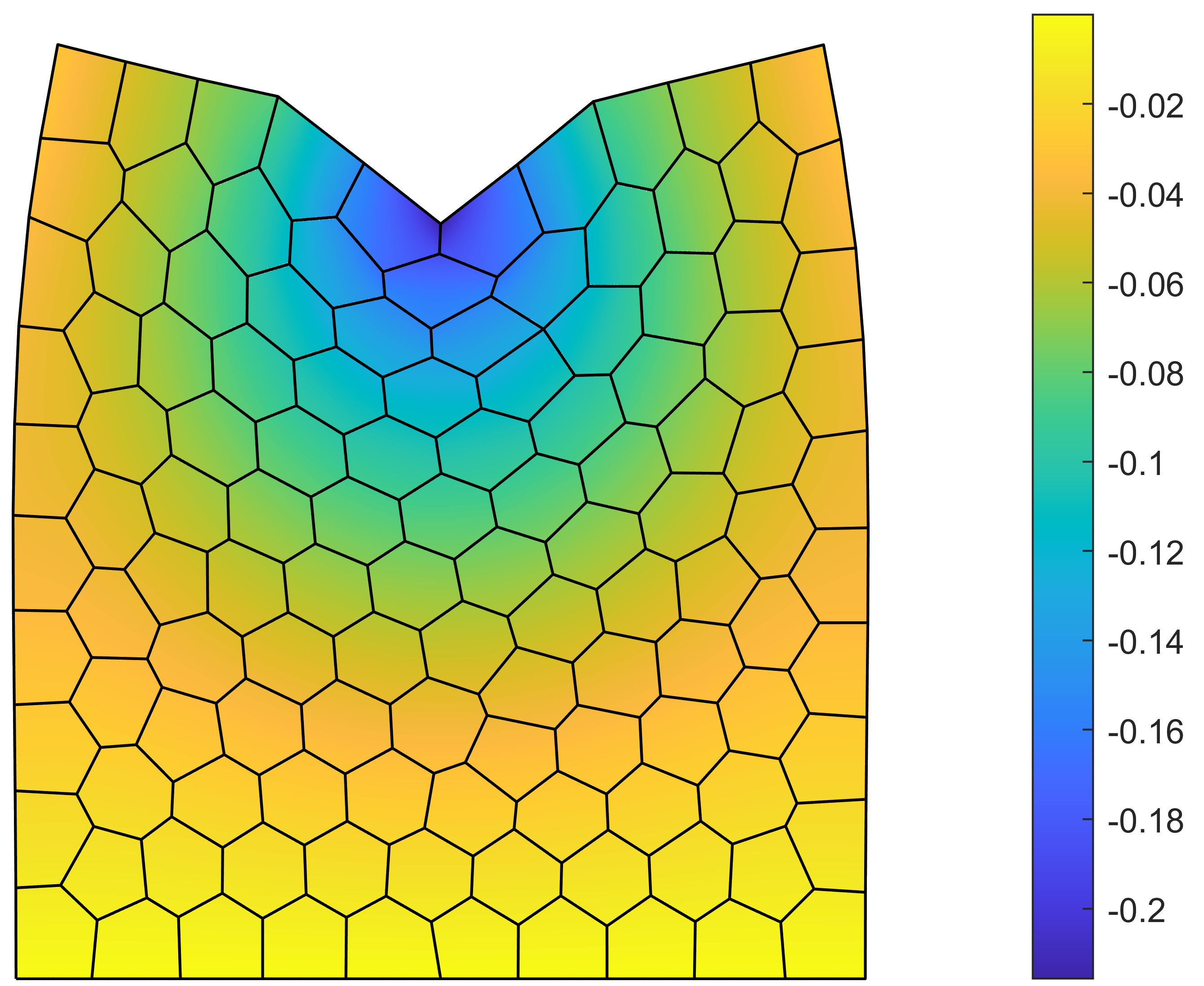}
			\caption{Deformed configuration}
		\end{subfigure}
		\caption{Problem~C(5) (a) geometry, and (b) sample deformed configuration of a Voronoi mesh with ${\nu=0.49995}$. 
			\label{fig:PunchGeometry}}
	\end{figure}

	Problem~C(5) comprises a domain of width ${w=1~\rm{m}}$ and height ${h=1~\rm{m}}$ into which a punch of width ${w_{\rm p}=0.2~w}$ is driven into the middle of the top edge.
	The bottom edge of the domain is constrained vertically and the midpoint of the bottom edge is fully constrained.
	The top edge is constrained horizontally and the punch is modelled as a uniformly distributed load with a magnitude of ${Q_{\rm{P}}=0.675~\frac{\rm N}{\rm m}}$ (see Figure~\ref{fig:PunchGeometry}(a)).
	The results presented for this problem were generated using a combination of the displacement-based and ${Z^{2}}$-like refinement procedures. 
	Figure~\ref{fig:PunchGeometry}(b) depicts a sample deformed configuration of the body with a Voronoi mesh and a nearly incompressible Poisson's ratio of ${\nu=0.49995}$. The vertical displacement ${u_{y}}$ is plotted on the colour axis.
	
	\FloatBarrier
	\FloatBarrier
	
\begin{figure}[ht!]
		\centering
		\begin{subfigure}[t]{0.33\textwidth}
			\centering
			\includegraphics[width=0.95\textwidth]{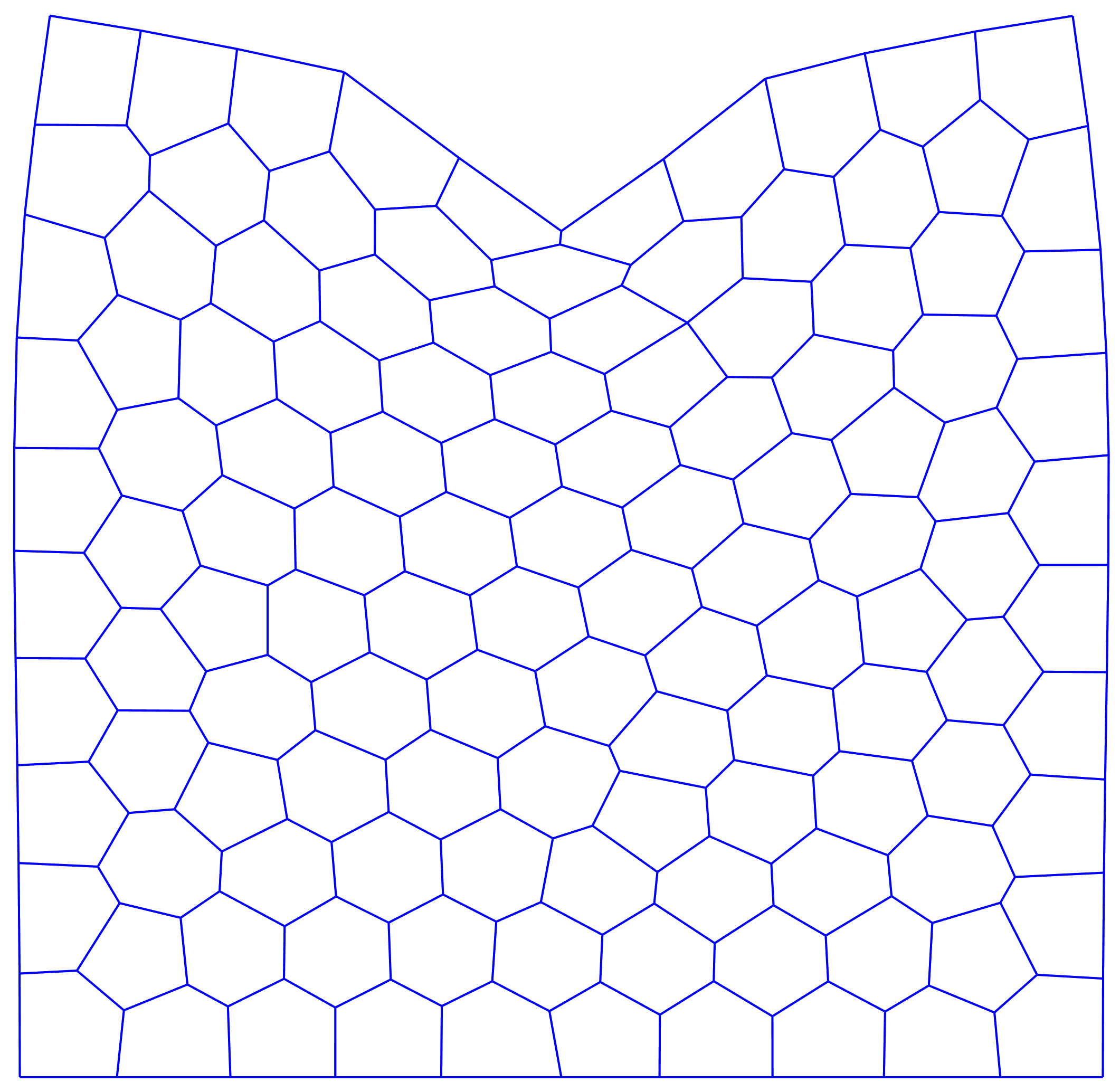}
			\caption{Compressible - Step 1}
		\end{subfigure}%
		\begin{subfigure}[t]{0.33\textwidth}
			\centering
			\includegraphics[width=0.95\textwidth]{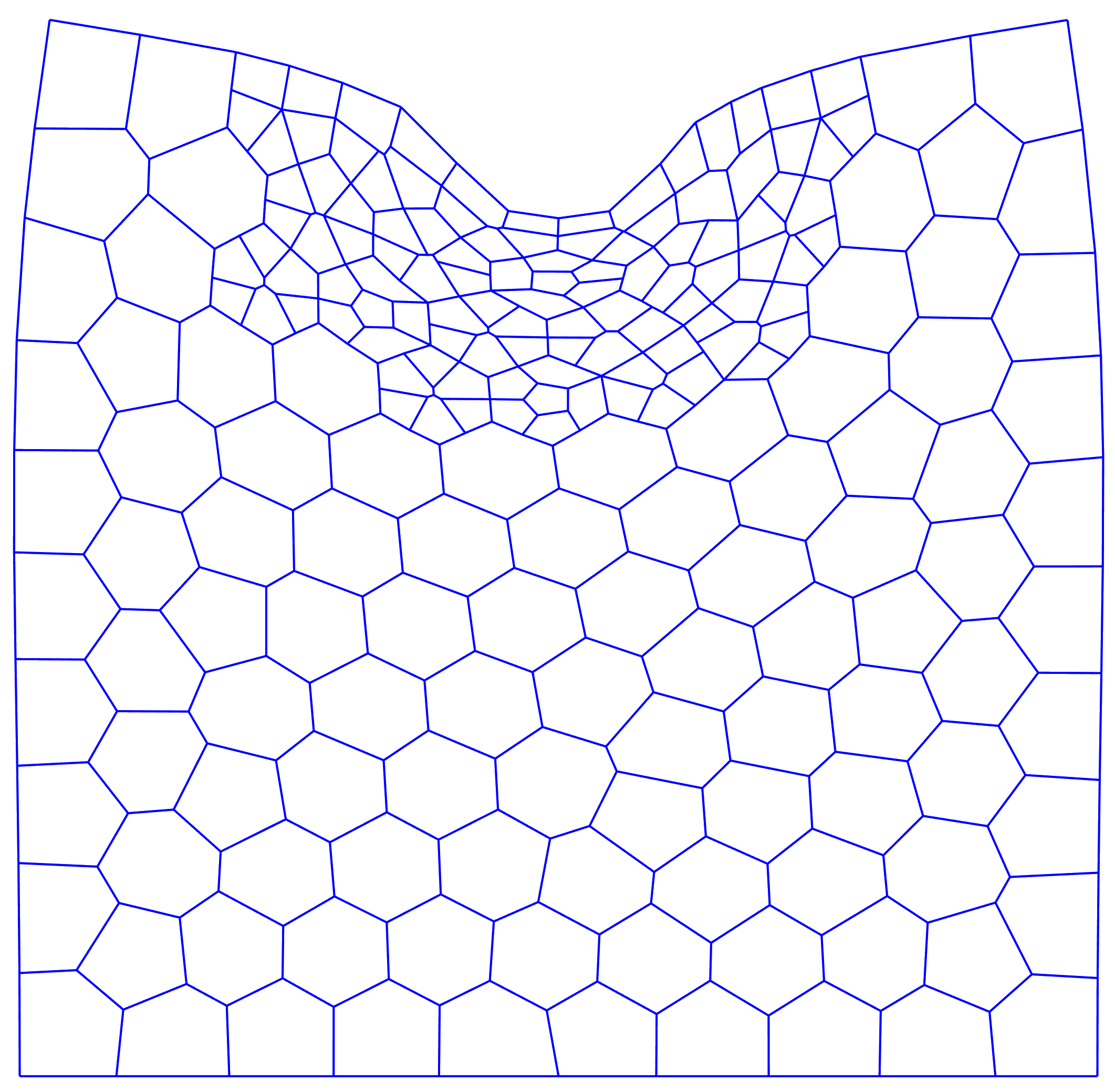}
			\caption{Compressible - Step 2}
		\end{subfigure}%
		\begin{subfigure}[t]{0.33\textwidth}
			\centering
			\includegraphics[width=0.95\textwidth]{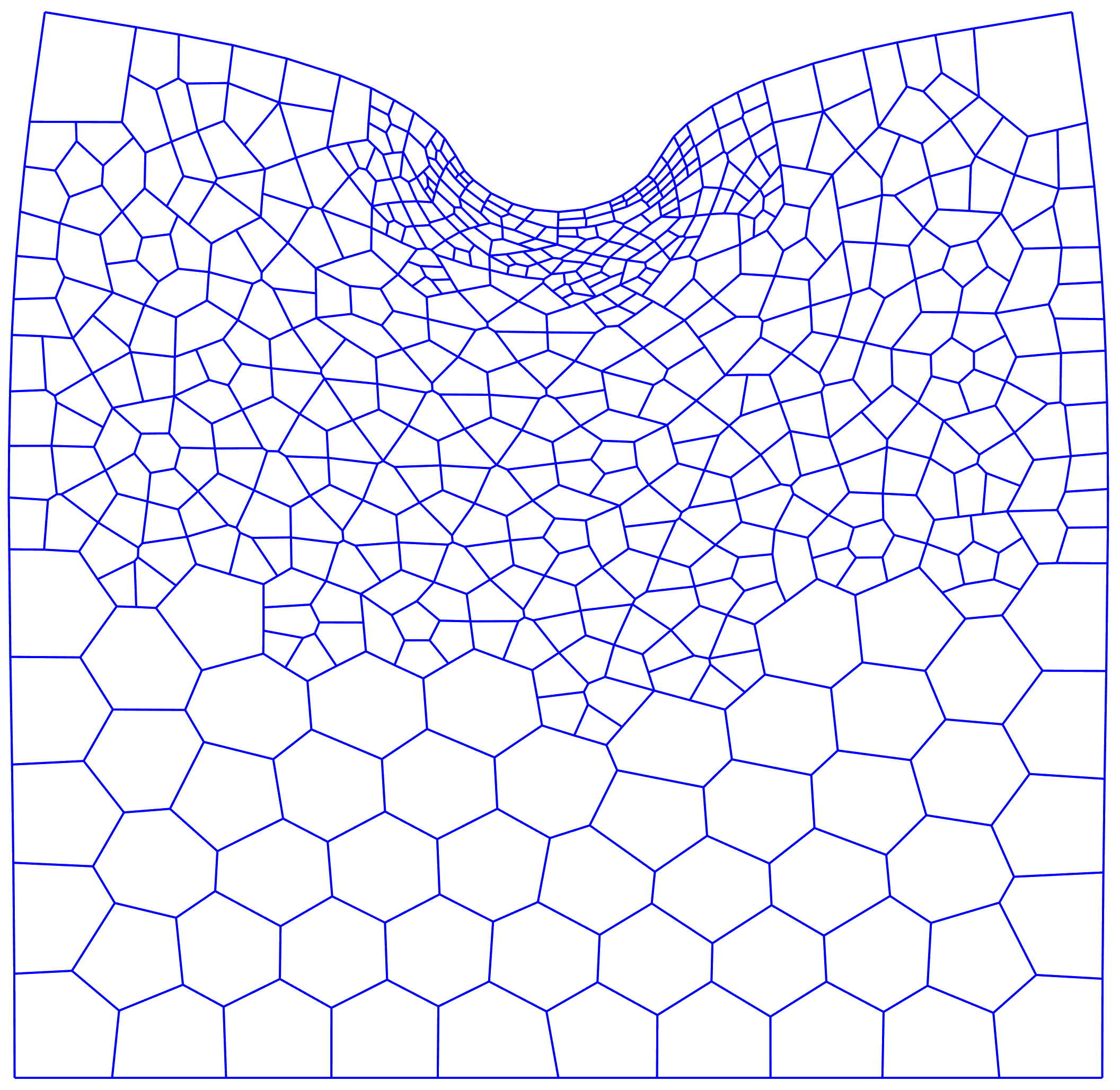}
			\caption{Compressible - Step 4}
		\end{subfigure}
		\vskip \baselineskip 
		\begin{subfigure}[t]{0.33\textwidth}
			\centering
			\includegraphics[width=0.95\textwidth]{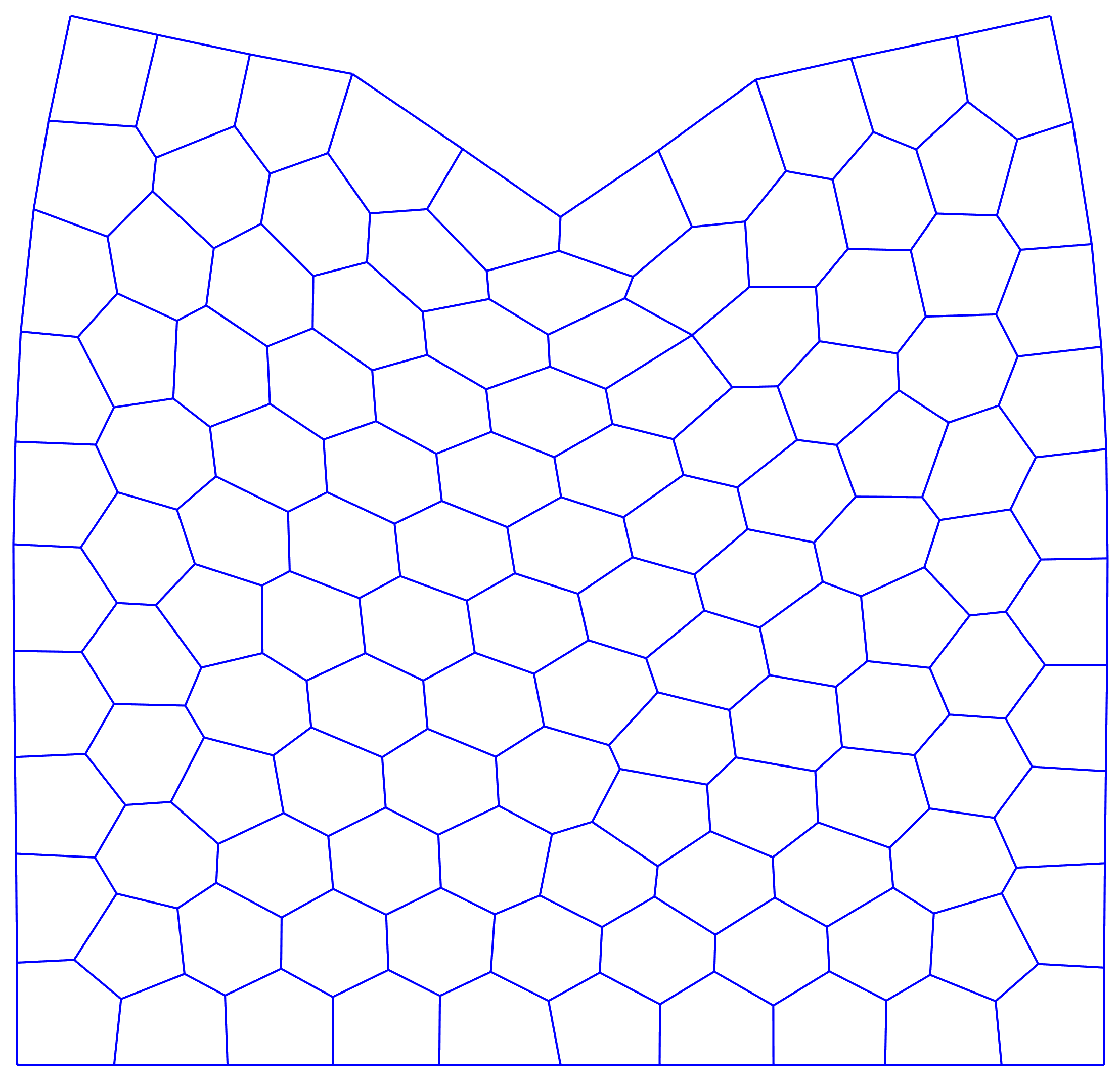}
			\caption{Nearly incompressible - Step 1}
		\end{subfigure}%
		\begin{subfigure}[t]{0.33\textwidth}
			\centering
			\includegraphics[width=0.95\textwidth]{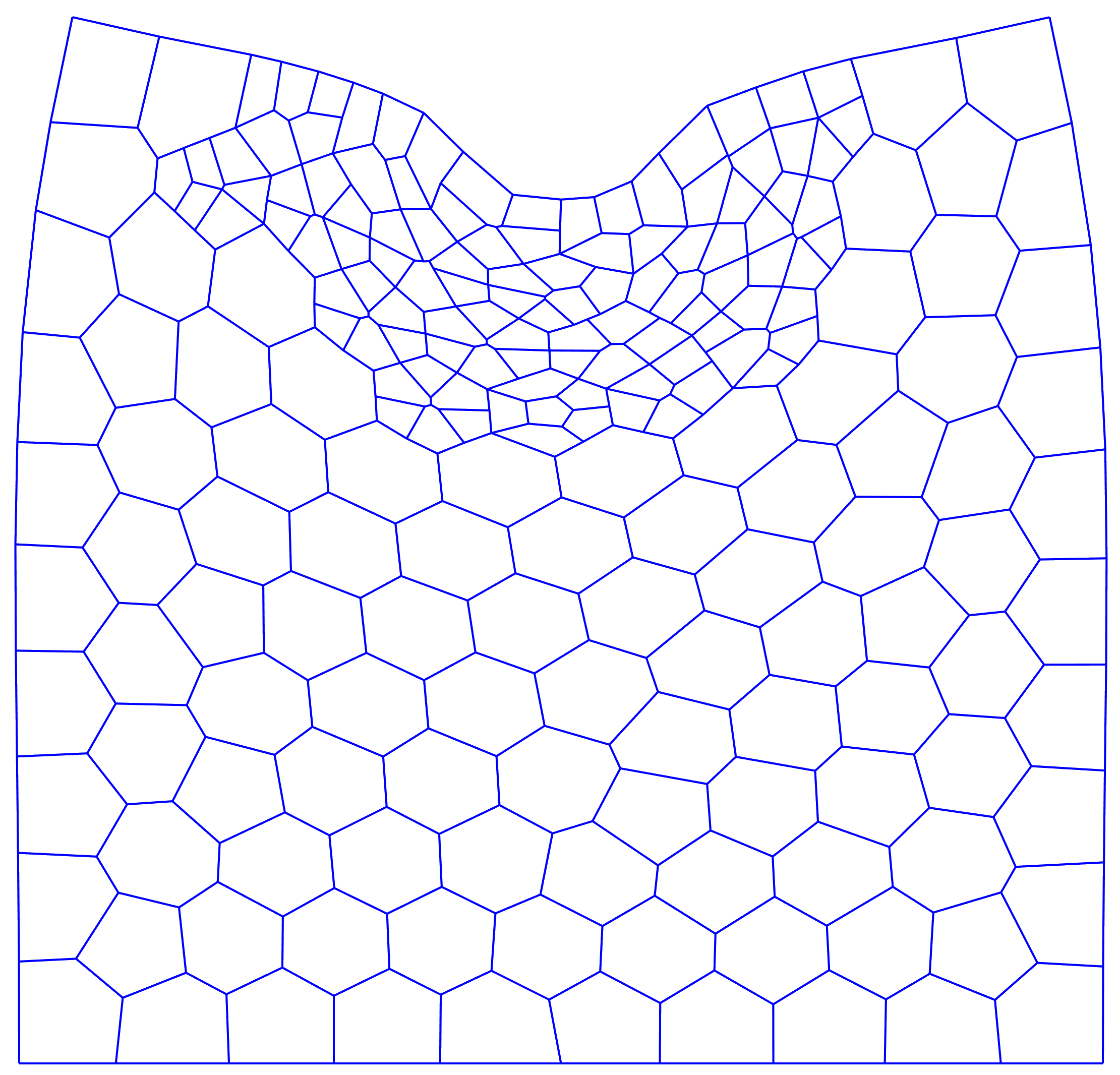}
			\caption{Nearly incompressible - Step 2}
		\end{subfigure}%
		\begin{subfigure}[t]{0.33\textwidth}
			\centering
			\includegraphics[width=0.95\textwidth]{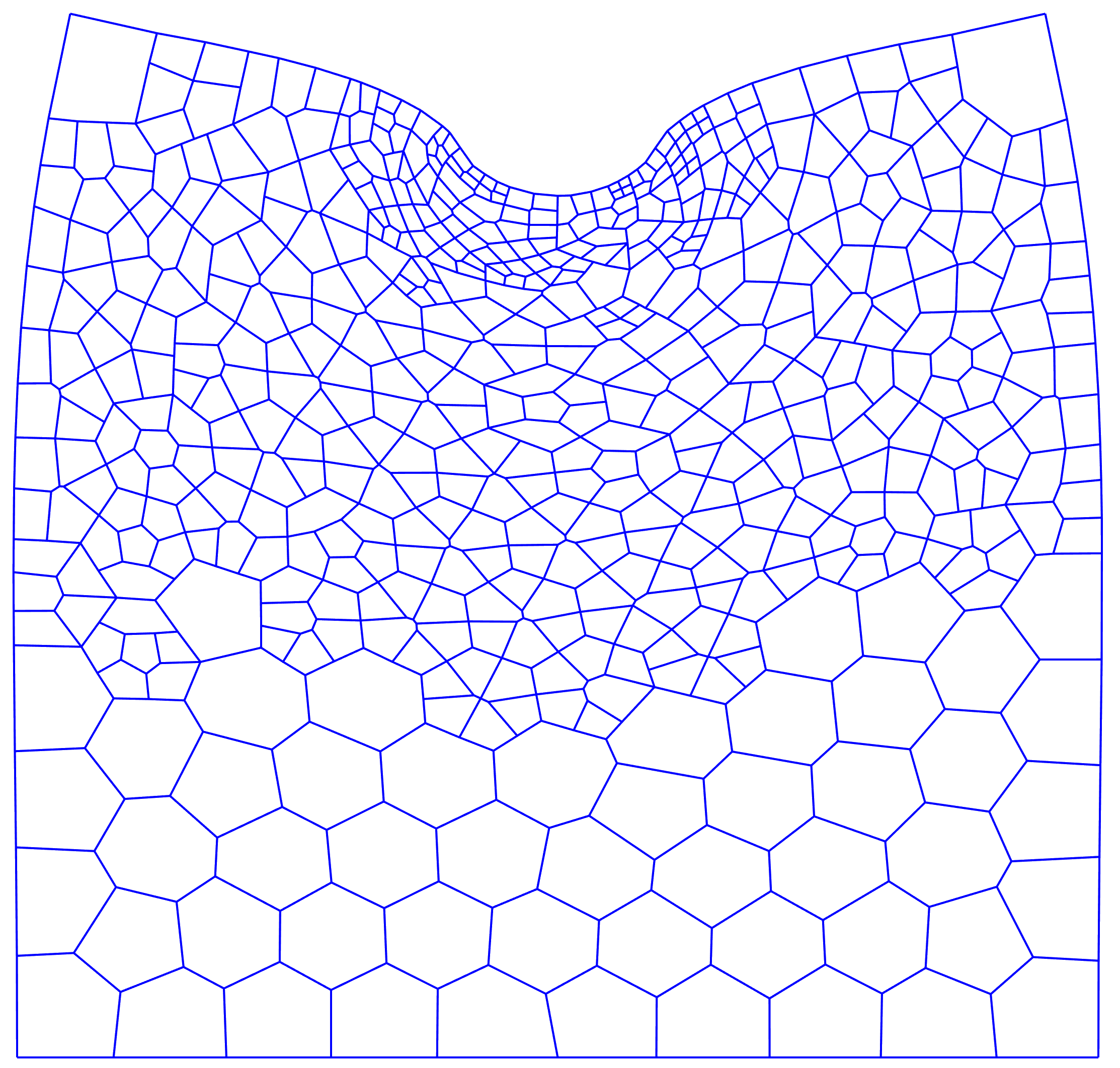}
			\caption{Nearly incompressible - Step 4}
		\end{subfigure}
		\caption{Mesh refinement process for problem~C(5) using a combination of the displacement-based and ${Z^{2}}$-like refinement procedures with ${T=15\%}$ for Voronoi meshes with compressible and nearly incompressible Poisson's ratios.
			\label{fig:SquarePlateNarrowPunchMeshesVrn}}
	\end{figure}

	Figure~\ref{fig:SquarePlateNarrowPunchMeshesVrn} depicts the mesh refinement process for problem~C(5) using a combination of the displacement-based and ${Z^{2}}$-like refinement procedures with ${T=15\%}$ for Voronoi meshes with compressible and nearly incompressible Poisson's ratios. 
	Meshes are shown at various refinement steps with step~1 corresponding to the initial mesh. In this problem most of the deformation occurs around the location of the punch. The rest of the body experiences comparatively little deformation and the magnitude of the deformation decreases towards the bottom of the domain.
	This behaviour is reflected in the mesh refinement process with the refinement strongly focused in the area around the punch with decreased refinement moving away from the punch. Thus, the mesh evolution is sensible for this problem.
	
	\FloatBarrier
	\FloatBarrier
	
\begin{figure}[ht!]
		\centering
		\begin{subfigure}[t]{0.5\textwidth}
			\centering
			\includegraphics[width=0.95\textwidth]{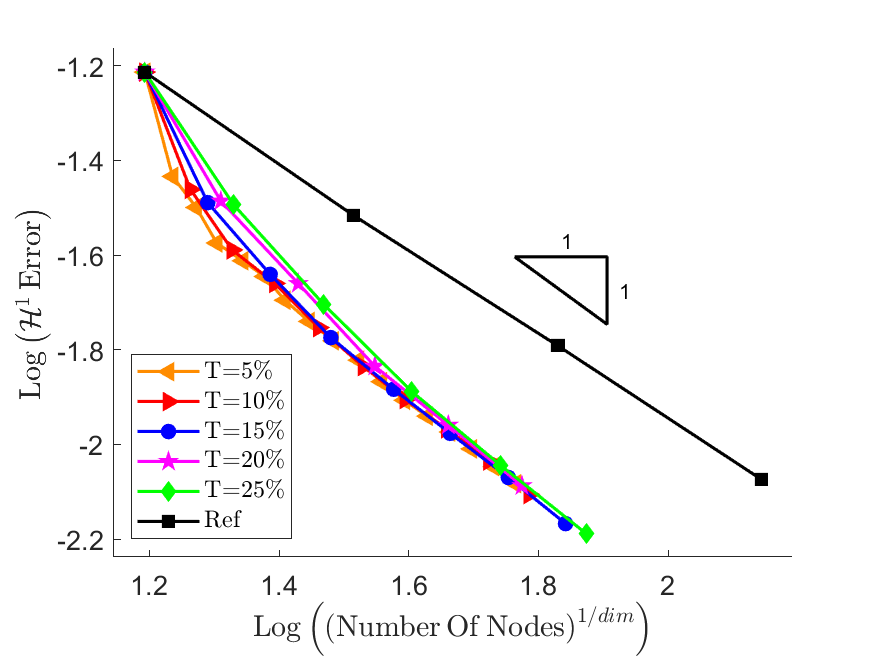}
			\caption{Compressible}
		\end{subfigure}%
		\begin{subfigure}[t]{0.5\textwidth}
			\centering
			\includegraphics[width=0.95\textwidth]{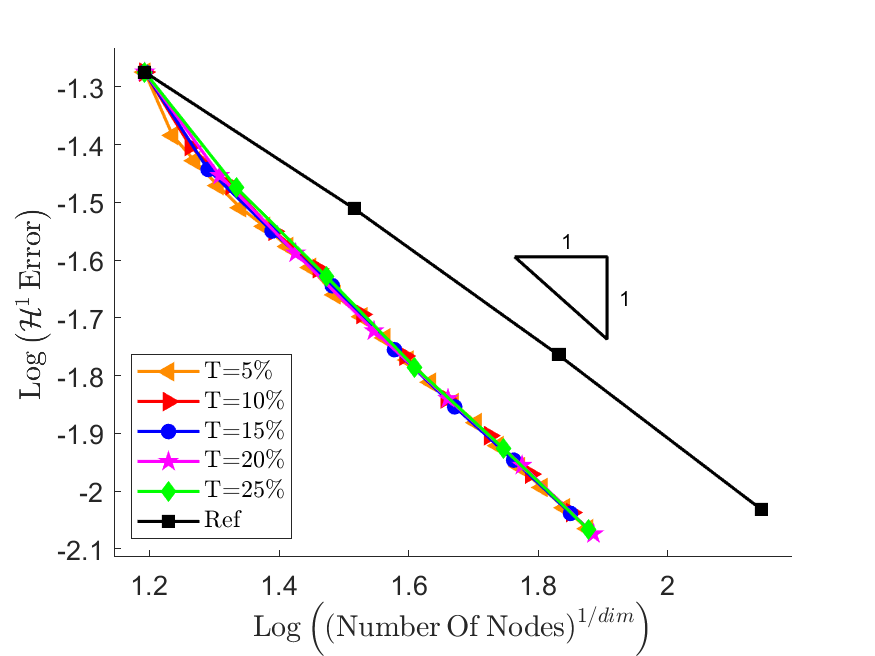}
			\caption{Nearly incompressible}
		\end{subfigure}
		\caption{$\mathcal{H}^{1}$ error vs $n_{\rm v}$ for problem~C(5) using a combination of the displacement-based and ${Z^{2}}$-like refinement procedures with a variety of choices of $T$ for Voronoi meshes with compressible and nearly incompressible Poisson's ratios.
			\label{fig:SquarePlateNarrowPunchConvergenceNumberOfNodes}}
	\end{figure}

	The convergence behaviour in the ${\mathcal{H}^{1}}$ error norm of the VEM for problem~C(5) using a combination of the displacement-based and ${Z^{2}}$-like refinement procedures is depicted in Figure~\ref{fig:SquarePlateNarrowPunchConvergenceNumberOfNodes} on a logarithmic scale for a variety of choices of $T$. Here, the ${\mathcal{H}^{1}}$ error is plotted against the number of vertices/nodes in the discretization for Voronoi meshes with compressible and nearly incompressible Poisson's ratios. The adaptive procedure exhibits the expected behaviour with lower choices of $T$ initially exhibiting slightly faster convergence rates. However, in the fine mesh range all choices of $T$ exhibit similar convergence behaviour and significantly outperform the reference procedure. While similar behaviour is observed for the cases of compressibility and near-incompressibility, the influence of $T$ is expectedly smaller in the nearly incompressible case. 
	
	\FloatBarrier
	\FloatBarrier
	
	The performance of the VEM in terms of its convergence behaviour in the $\mathcal{H}^{1}$ error norm with respect to the number of vertices/nodes in the discretization when using a combination of the displacement-based and ${Z^{2}}$-like refinement procedures for problem~C(5) is summarized in Table~\ref{tab:PerformanceSquarePlateNarrowPunchNumberOfNodes}. Here, the performance, as measured by the PRE, is presented for the cases of structured and Voronoi meshes with compressible and nearly incompressible Poisson's ratios for a variety of choices of $T$. 
	For this example for both mesh types and for both choices of Poisson's ratio the choice of $T$ has a small effect on the PRE. However, in general lower choices of $T$ perform slightly better than larger choices.
	
	\FloatBarrier
	\begin{table}[ht!]
		\centering 
		\begin{adjustbox}{max width=\textwidth}
			\begin{tabular}{|c|cccc|cccc|}
				\hline
				\multirow{3}{*}{Threshold} & \multicolumn{4}{c|}{Compressible}                                                                 & \multicolumn{4}{c|}{Nearly-incompressible}                                                        \\ \cline{2-9} 
				& \multicolumn{2}{c|}{Structured}                           & \multicolumn{2}{c|}{Voronoi}          & \multicolumn{2}{c|}{Structured}                           & \multicolumn{2}{c|}{Voronoi}          \\ \cline{2-9} 
				& \multicolumn{1}{c|}{nNodes}  & \multicolumn{1}{c|}{PRE}   & \multicolumn{1}{c|}{nNodes}   & PRE   & \multicolumn{1}{c|}{nNodes}  & \multicolumn{1}{c|}{PRE}   & \multicolumn{1}{c|}{nNodes}   & PRE   \\ \hline
				T=5\%                      & \multicolumn{1}{c|}{1514.66} & \multicolumn{1}{c|}{23.09} & \multicolumn{1}{c|}{3245.14}  & 16.76 & \multicolumn{1}{c|}{2048.90} & \multicolumn{1}{c|}{31.23} & \multicolumn{1}{c|}{4854.09}  & 24.99 \\ \hline
				T=10\%                     & \multicolumn{1}{c|}{1522.64} & \multicolumn{1}{c|}{23.21} & \multicolumn{1}{c|}{3271.67}  & 16.90 & \multicolumn{1}{c|}{2059.60} & \multicolumn{1}{c|}{31.39} & \multicolumn{1}{c|}{4933.02}  & 25.39 \\ \hline
				T=15\%                     & \multicolumn{1}{c|}{1528.55} & \multicolumn{1}{c|}{23.30} & \multicolumn{1}{c|}{3264.21}  & 16.86 & \multicolumn{1}{c|}{2079.66} & \multicolumn{1}{c|}{31.70} & \multicolumn{1}{c|}{4858.00}  & 25.01 \\ \hline
				T=20\%                     & \multicolumn{1}{c|}{1545.05} & \multicolumn{1}{c|}{23.55} & \multicolumn{1}{c|}{3346.74}  & 17.29 & \multicolumn{1}{c|}{2087.51} & \multicolumn{1}{c|}{31.82} & \multicolumn{1}{c|}{4916.83}  & 25.31 \\ \hline
				T=25\%                     & \multicolumn{1}{c|}{1552.66} & \multicolumn{1}{c|}{23.67} & \multicolumn{1}{c|}{3443.74}  & 17.79 & \multicolumn{1}{c|}{2078.83} & \multicolumn{1}{c|}{31.68} & \multicolumn{1}{c|}{4897.56}  & 25.21 \\ \hline
				Ref                        & \multicolumn{1}{c|}{6561.00} & \multicolumn{1}{c|}{}      & \multicolumn{1}{c|}{19362.00} &       & \multicolumn{1}{c|}{6561.00} & \multicolumn{1}{c|}{}      & \multicolumn{1}{c|}{19427.00} &       \\ \hline
			\end{tabular}
		\end{adjustbox}
		\caption{Performance summary of the VEM in terms of its convergence behaviour in the $\mathcal{H}^{1}$ error norm with respect to the number of vertices/nodes in the discretization when using a combination of the displacement-based and ${Z^{2}}$-like refinement procedures for problem~C(5).
			\label{tab:PerformanceSquarePlateNarrowPunchNumberOfNodes}}
	\end{table}
	\FloatBarrier
	
\begin{figure}[ht!]
		\centering
		\begin{subfigure}[t]{0.5\textwidth}
			\centering
			\includegraphics[width=0.95\textwidth]{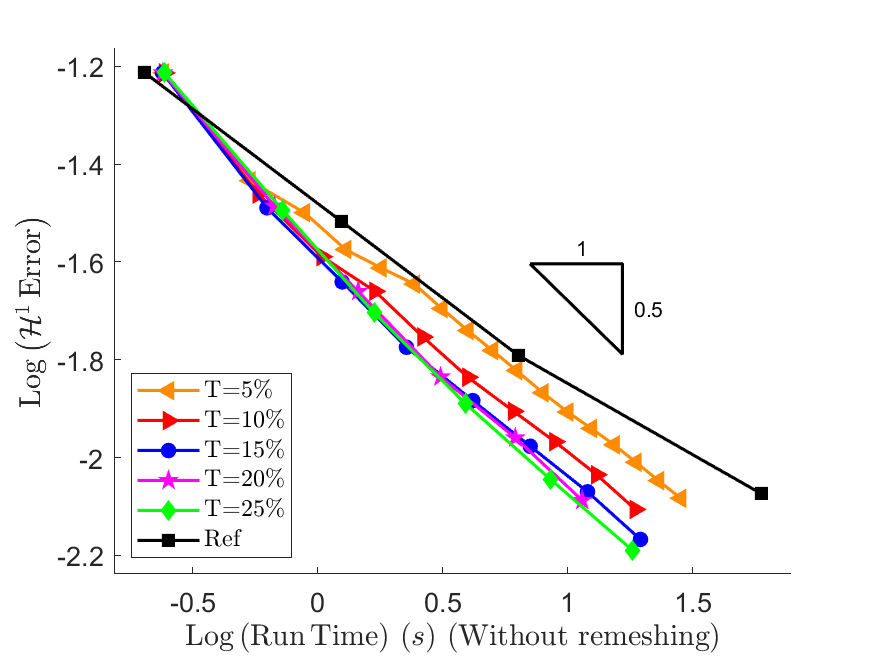}
			\caption{Voronoi mesh - Compressible}
		\end{subfigure}%
		\begin{subfigure}[t]{0.5\textwidth}
			\centering
			\includegraphics[width=0.95\textwidth]{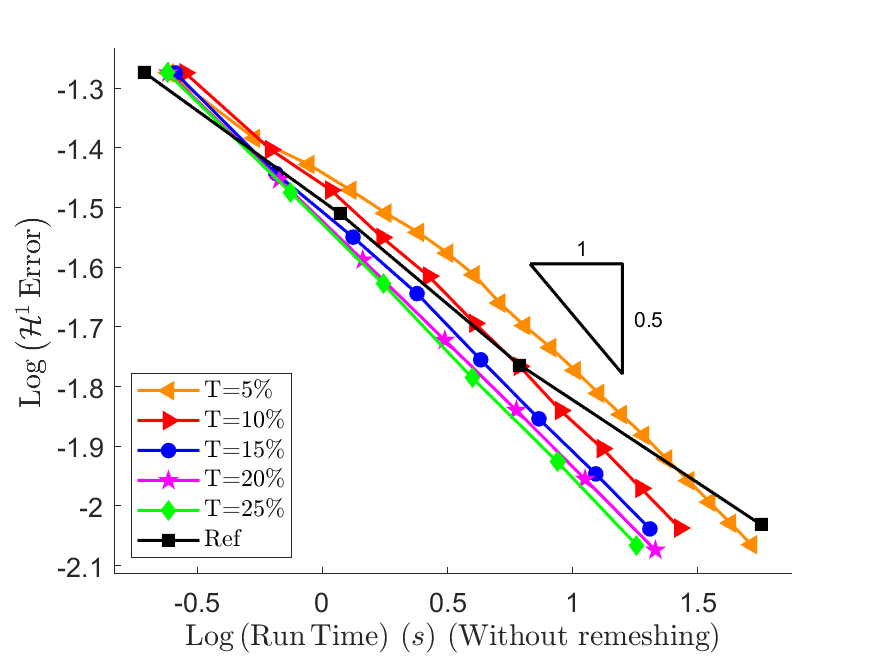}
			\caption{Voronoi mesh - Nearly incompressible}
		\end{subfigure}
		\caption{$\mathcal{H}^{1}$ error vs run time (excluding remeshing time) for problem~C(5) using a combination of the displacement-based and ${Z^{2}}$-like refinement procedures with a variety of choices of $T$ for Voronoi meshes with compressible and nearly incompressible Poisson's ratios.
			\label{fig:SquarePlateNarrowPunchConvergenceRunTime}}
	\end{figure}

	The convergence behaviour in the ${\mathcal{H}^{1}}$ error norm of the VEM for problem~C(5) using a combination of the displacement-based and ${Z^{2}}$-like refinement procedures is depicted in Figure~\ref{fig:SquarePlateNarrowPunchConvergenceRunTime} on a logarithmic scale for a variety of choices of $T$. Here, the ${\mathcal{H}^{1}}$ error is plotted against run time (excluding remeshing time) for Voronoi meshes with compressible and nearly incompressible Poisson's ratios. 
	As observed in previous examples, lower choices of $T$ require significantly more remeshing steps, and consequently more run time, than larger values of $T$ to reach the desired level of accuracy. Thus, lower choices of $T$ are not particularly efficient in terms of run time.
	However, for ${T>5\%}$, for both mesh types and for all choices of $T$ the adaptive procedure exhibits a superior convergence rate to, and significantly outperforms, the reference refinement procedure.
	
	\FloatBarrier
	\FloatBarrier
	
	The performance of the VEM in terms of its convergence behaviour in the $\mathcal{H}^{1}$ error norm with respect to run time (excluding remeshing time) when using a combination of the displacement-based and ${Z^{2}}$-like refinement procedures for problem~C(5) is summarized in Table~\ref{tab:PerformanceSquarePlateNarrowPunchRunTime}. Here, the performance, as measured by the PRE, is presented for the cases of structured and Voronoi meshes with compressible and nearly incompressible Poisson's ratios for a variety of choices of $T$. 
	The relative inefficiencies of lower choices of $T$ are again easy to see. 
	Nevertheless, independent of the degree of compressibility and mesh type the adaptive procedure, for larger values of $T$, represents a significant improvement in efficiency compared to the reference procedure.
	
	\FloatBarrier
	\begin{table}[ht!]
		\centering 
		\begin{adjustbox}{max width=\textwidth}
			\begin{tabular}{|c|cccc|cccc|}
				\hline
				\multirow{3}{*}{Threshold} & \multicolumn{4}{c|}{Compressible}                                                                   & \multicolumn{4}{c|}{Nearly-incompressible}                                                          \\ \cline{2-9} 
				& \multicolumn{2}{c|}{Structured}                             & \multicolumn{2}{c|}{Voronoi}          & \multicolumn{2}{c|}{Structured}                             & \multicolumn{2}{c|}{Voronoi}          \\ \cline{2-9} 
				& \multicolumn{1}{c|}{Run time} & \multicolumn{1}{c|}{PRE}    & \multicolumn{1}{c|}{Run time} & PRE   & \multicolumn{1}{c|}{Run time} & \multicolumn{1}{c|}{PRE}    & \multicolumn{1}{c|}{Run time} & PRE   \\ \hline
				T=5\%                      & \multicolumn{1}{c|}{14.69}    & \multicolumn{1}{c|}{131.43} & \multicolumn{1}{c|}{26.80}    & 45.14 & \multicolumn{1}{c|}{21.59}    & \multicolumn{1}{c|}{191.87} & \multicolumn{1}{c|}{43.32}    & 76.45 \\ \hline
				T=10\%                     & \multicolumn{1}{c|}{8.68}     & \multicolumn{1}{c|}{77.70}  & \multicolumn{1}{c|}{15.84}    & 26.68 & \multicolumn{1}{c|}{12.56}    & \multicolumn{1}{c|}{111.61} & \multicolumn{1}{c|}{25.99}    & 45.86 \\ \hline
				T=15\%                     & \multicolumn{1}{c|}{6.56}     & \multicolumn{1}{c|}{58.69}  & \multicolumn{1}{c|}{12.23}    & 20.59 & \multicolumn{1}{c|}{9.65}     & \multicolumn{1}{c|}{85.74}  & \multicolumn{1}{c|}{19.59}    & 34.58 \\ \hline
				T=20\%                     & \multicolumn{1}{c|}{5.62}     & \multicolumn{1}{c|}{50.33}  & \multicolumn{1}{c|}{10.69}    & 18.01 & \multicolumn{1}{c|}{8.31}     & \multicolumn{1}{c|}{73.87}  & \multicolumn{1}{c|}{16.98}    & 29.97 \\ \hline
				T=25\%                     & \multicolumn{1}{c|}{4.87}     & \multicolumn{1}{c|}{43.61}  & \multicolumn{1}{c|}{9.93}     & 16.72 & \multicolumn{1}{c|}{7.30}     & \multicolumn{1}{c|}{64.91}  & \multicolumn{1}{c|}{15.06}    & 26.57 \\ \hline
				Ref                        & \multicolumn{1}{c|}{11.17}    & \multicolumn{1}{c|}{}       & \multicolumn{1}{c|}{59.36}    &       & \multicolumn{1}{c|}{11.25}    & \multicolumn{1}{c|}{}       & \multicolumn{1}{c|}{56.66}    &       \\ \hline
			\end{tabular}
		\end{adjustbox}
		\caption{Performance summary of the VEM in terms of its convergence behaviour in the $\mathcal{H}^{1}$ error norm with respect to run time (excluding remeshing time) when using a combination of the displacement-based and ${Z^{2}}$-like refinement procedures for problem~C(5).
			\label{tab:PerformanceSquarePlateNarrowPunchRunTime}}
	\end{table}
	\FloatBarrier
	
\begin{figure}[ht!]
		\centering
		\begin{subfigure}[t]{0.5\textwidth}
			\centering
			\includegraphics[width=0.95\textwidth]{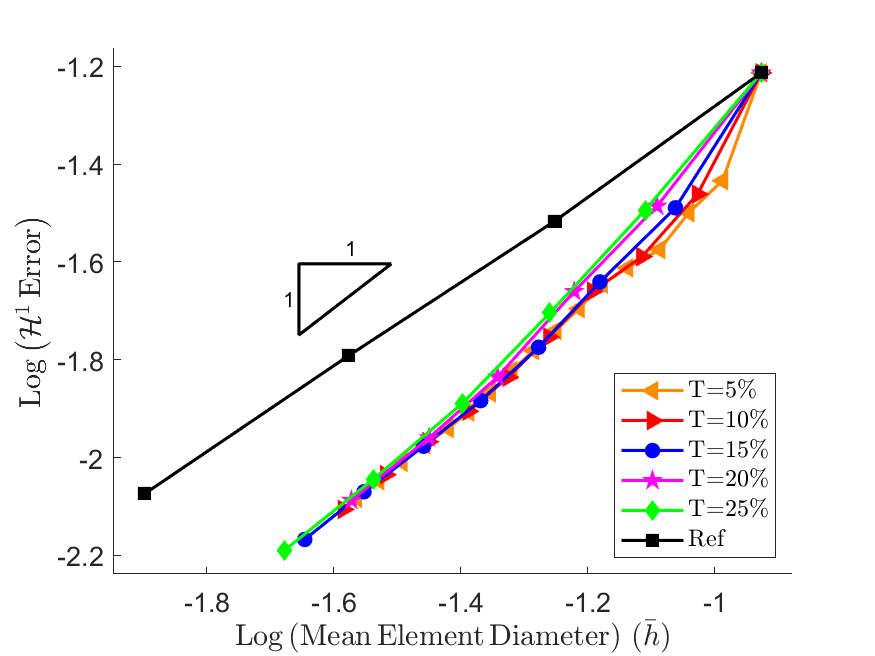}
			\caption{Voronoi mesh - Compressible}
		\end{subfigure}%
		\begin{subfigure}[t]{0.5\textwidth}
			\centering
			\includegraphics[width=0.95\textwidth]{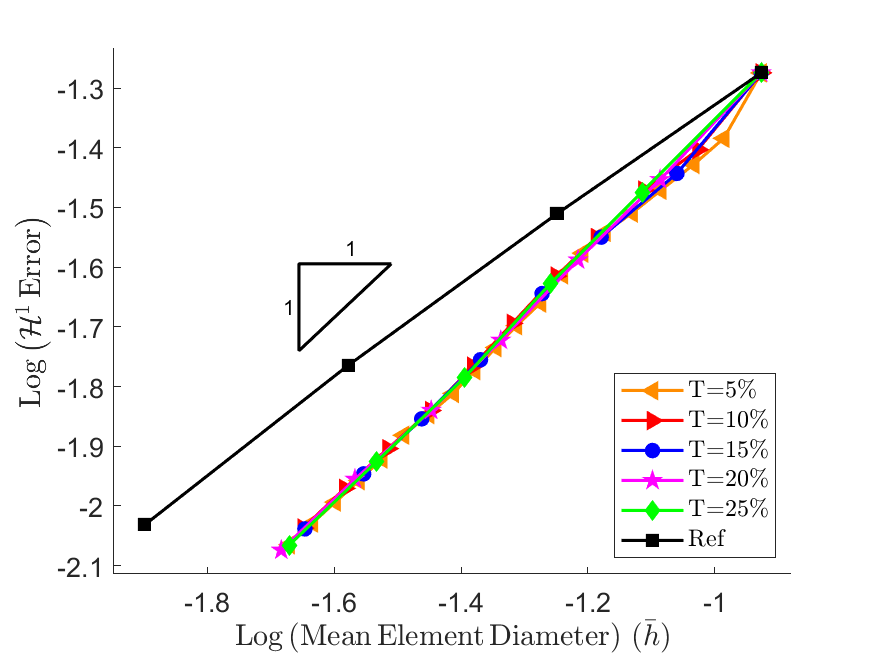}
			\caption{Voronoi mesh - Nearly incompressible}
		\end{subfigure}
		\caption{$\mathcal{H}^{1}$ error vs mesh size for problem~C(5) using a combination of the displacement-based and ${Z^{2}}$-like refinement procedures with a variety of choices of $T$.
			\label{fig:SquarePlateNarrowPunchConvergenceMeshSize}}
	\end{figure}

	The convergence behaviour in the ${\mathcal{H}^{1}}$ error norm of the VEM for problem~C(5) using a combination of the displacement-based and ${Z^{2}}$-like refinement procedures is depicted in Figure~\ref{fig:SquarePlateNarrowPunchConvergenceMeshSize} on a logarithmic scale for a variety of choices of $T$. Here, the ${\mathcal{H}^{1}}$ error is plotted against mesh size as measured by the mean element diameter for Voronoi meshes with compressible and nearly incompressible Poisson's ratios. 
	The convergence behaviour is very similar to that observed in Figure~\ref{fig:SquarePlateNarrowPunchConvergenceNumberOfNodes} with respect to the number of nodes. Specifically, lower choices of $T$ initially exhibit slightly faster convergence rates. However, in the fine mesh range all choices of $T$ exhibit similar convergence behaviour and significantly outperform the reference procedure.
	Additionally, as observed in Figure~\ref{fig:SquarePlateNarrowPunchConvergenceNumberOfNodes}, while similar behaviour is exhibited in the cases of compressibility and near-incompressibility, the influence of $T$ is smaller in the nearly incompressible case. 
	
	\FloatBarrier
	\FloatBarrier
	
	The performance of the VEM in terms of its convergence behaviour in the $\mathcal{H}^{1}$ error norm with respect to mesh size when using a combination of the displacement-based and ${Z^{2}}$-like refinement procedures for problem~C(5) is summarized in Table~\ref{tab:PerformanceSquarePlateNarrowPunchMeshSize}. Here, the performance, as measured by the PRE*, is presented for the cases of structured and Voronoi meshes with compressible and nearly incompressible Poisson's ratios for a variety of choices of $T$. 
	It is clear, again, that for this problem the choice of $T$ has a relatively small influence on the performance of the adaptive procedure. For all choices of $T$, for both mesh types, and for the cases of compressibility and near-incompressibility, the adaptive procedure significantly outperforms the reference procedure.
	
	\FloatBarrier
	\begin{table}[ht!]
		\centering 
		\begin{adjustbox}{max width=\textwidth}
			\begin{tabular}{|c|cccc|cccc|}
				\hline
				\multirow{3}{*}{Threshold} & \multicolumn{4}{c|}{Compressible}                                                                    & \multicolumn{4}{c|}{Nearly-incompressible}                                                           \\ \cline{2-9} 
				& \multicolumn{2}{c|}{Structured}                             & \multicolumn{2}{c|}{Voronoi}           & \multicolumn{2}{c|}{Structured}                             & \multicolumn{2}{c|}{Voronoi}           \\ \cline{2-9} 
				& \multicolumn{1}{c|}{Mesh size} & \multicolumn{1}{c|}{PRE*}  & \multicolumn{1}{c|}{Mesh size} & PRE*  & \multicolumn{1}{c|}{Mesh size} & \multicolumn{1}{c|}{PRE*}  & \multicolumn{1}{c|}{Mesh size} & PRE*  \\ \hline
				T=5\%                      & \multicolumn{1}{c|}{0.032783}  & \multicolumn{1}{c|}{53.92} & \multicolumn{1}{c|}{0.027832}  & 45.46 & \multicolumn{1}{c|}{0.028504}  & \multicolumn{1}{c|}{62.02} & \multicolumn{1}{c|}{0.023070}  & 54.75 \\ \hline
				T=10\%                     & \multicolumn{1}{c|}{0.032648}  & \multicolumn{1}{c|}{54.15} & \multicolumn{1}{c|}{0.027885}  & 45.37 & \multicolumn{1}{c|}{0.028440}  & \multicolumn{1}{c|}{62.16} & \multicolumn{1}{c|}{0.022735}  & 55.56 \\ \hline
				T=15\%                     & \multicolumn{1}{c|}{0.032760}  & \multicolumn{1}{c|}{53.96} & \multicolumn{1}{c|}{0.027785}  & 45.54 & \multicolumn{1}{c|}{0.028327}  & \multicolumn{1}{c|}{62.40} & \multicolumn{1}{c|}{0.022953}  & 55.03 \\ \hline
				T=20\%                     & \multicolumn{1}{c|}{0.032413}  & \multicolumn{1}{c|}{54.54} & \multicolumn{1}{c|}{0.027588}  & 45.86 & \multicolumn{1}{c|}{0.028230}  & \multicolumn{1}{c|}{62.62} & \multicolumn{1}{c|}{0.022810}  & 55.38 \\ \hline
				T=25\%                     & \multicolumn{1}{c|}{0.032281}  & \multicolumn{1}{c|}{54.76} & \multicolumn{1}{c|}{0.027143}  & 46.61 & \multicolumn{1}{c|}{0.028304}  & \multicolumn{1}{c|}{62.46} & \multicolumn{1}{c|}{0.023053}  & 54.79 \\ \hline
				Ref                        & \multicolumn{1}{c|}{0.017678}  & \multicolumn{1}{c|}{}      & \multicolumn{1}{c|}{0.012652}  &       & \multicolumn{1}{c|}{0.017678}  & \multicolumn{1}{c|}{}      & \multicolumn{1}{c|}{0.012632}  &       \\ \hline
			\end{tabular}
		\end{adjustbox}
		\caption{Performance summary of the VEM in terms of its convergence behaviour in the $\mathcal{H}^{1}$ error norm with respect to mesh size when using a combination of the displacement-based and ${Z^{2}}$-like refinement procedures for problem~C(5).
			\label{tab:PerformanceSquarePlateNarrowPunchMeshSize}}
	\end{table}
	\FloatBarrier
	
\begin{figure}[ht!]
		\centering
		\begin{subfigure}[t]{0.5\textwidth}
			\centering
			\includegraphics[width=0.95\textwidth]{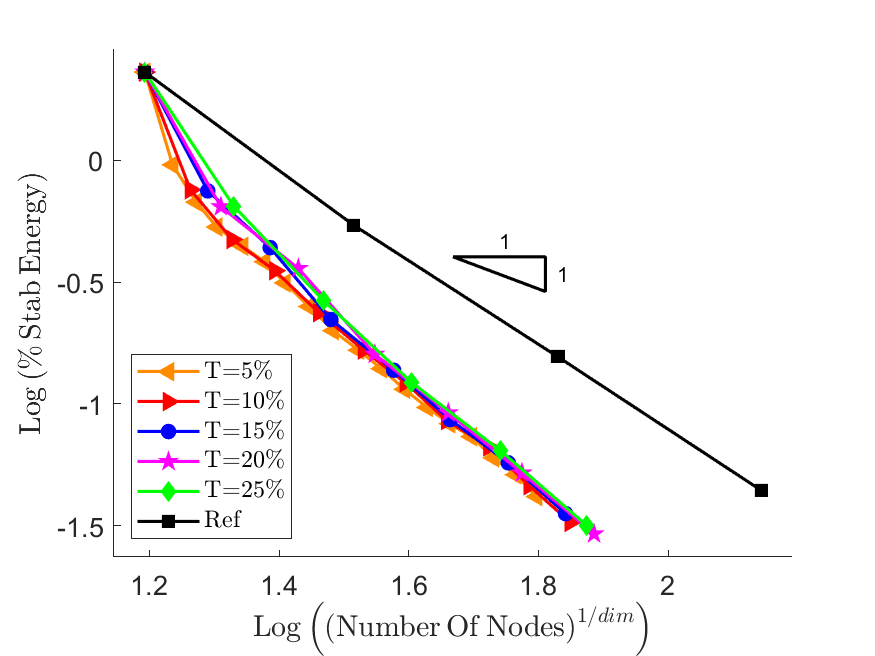}
			\caption{Voronoi mesh - Compressible}
		\end{subfigure}%
		\begin{subfigure}[t]{0.5\textwidth}
			\centering
			\includegraphics[width=0.95\textwidth]{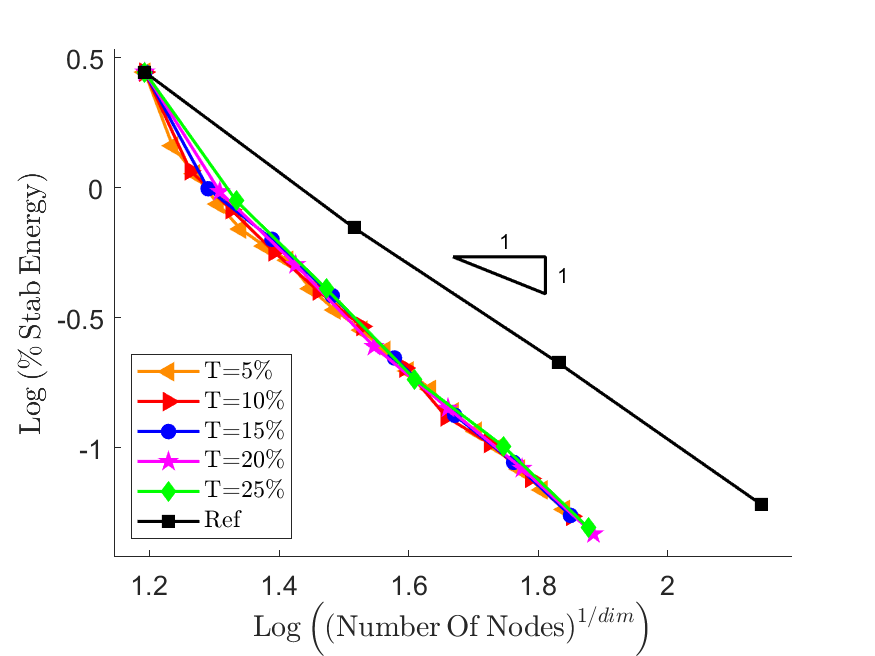}
			\caption{Voronoi mesh - Nearly incompressible}
		\end{subfigure}
		\caption{PSE vs number of nodes/vertices in the discretization for problem~C(5) using a combination of the displacement-based and ${Z^{2}}$-like refinement procedures with a variety of choices of $T$ for Voronoi meshes with compressible and nearly incompressible Poisson's ratios.
			\label{fig:SquarePlateNarrowPunchConvergencePSE}}
	\end{figure}

	The convergence behaviour in the PSE of the VEM for problem~C(5) using a combination of the displacement-based and ${Z^{2}}$-like refinement procedures is depicted in Figure~\ref{fig:SquarePlateNarrowPunchConvergencePSE} on a logarithmic scale for a variety of choices of $T$. Here, the PSE is plotted against the number of nodes/vertices in the discretization for Voronoi meshes with compressible and nearly incompressible Poisson's ratios.
	Similar to the previous results for this problem, lower choices of $T$ initially exhibit slightly faster convergence rates. However, in the fine mesh range all choices of $T$ exhibit similar convergence behaviour and significantly outperform the reference procedure. 
	
	\FloatBarrier
	\FloatBarrier
	
	The performance of the VEM in terms of its convergence behaviour in the PSE with respect to the number of nodes/vertices in the discretization when using a combination of the displacement-based and ${Z^{2}}$-like refinement procedures for problem~C(5) is summarized in Table~\ref{tab:PerformanceSquarePlateNarrowPunchPSE}. Here, the performance, as measured by the PRE, is presented for the cases of structured and Voronoi meshes with compressible and nearly incompressible Poisson's ratios for a variety of choices of $T$. As observed in the previous examples, for both mesh types and for the cases of both compressibility and near-incompressibility, the choice of $T$ has very little influence on the PRE. Additionally, similar performance is observed in cases of both compressibility and near-incompressibility indicating that the degree of compressibility doesn't strongly effect the convergence behaviour of the PSE.
	
	\FloatBarrier
	\begin{table}[ht!]
		\centering 
		\begin{adjustbox}{max width=\textwidth}
			\begin{tabular}{|c|cccc|cccc|}
				\hline
				\multirow{3}{*}{Threshold} & \multicolumn{4}{c|}{Compressible}                                                                 & \multicolumn{4}{c|}{Nearly-incompressible}                                                        \\ \cline{2-9} 
				& \multicolumn{2}{c|}{Structured}                           & \multicolumn{2}{c|}{Voronoi}          & \multicolumn{2}{c|}{Structured}                           & \multicolumn{2}{c|}{Voronoi}          \\ \cline{2-9} 
				& \multicolumn{1}{c|}{nNodes}  & \multicolumn{1}{c|}{PRE}   & \multicolumn{1}{c|}{nNodes}   & PRE   & \multicolumn{1}{c|}{nNodes}  & \multicolumn{1}{c|}{PRE}   & \multicolumn{1}{c|}{nNodes}   & PRE   \\ \hline
				T=5\%                      & \multicolumn{1}{c|}{1523.57} & \multicolumn{1}{c|}{23.22} & \multicolumn{1}{c|}{3764.34}  & 19.44 & \multicolumn{1}{c|}{2037.64} & \multicolumn{1}{c|}{31.06} & \multicolumn{1}{c|}{4604.98}  & 23.70 \\ \hline
				T=10\%                     & \multicolumn{1}{c|}{1518.12} & \multicolumn{1}{c|}{23.14} & \multicolumn{1}{c|}{3868.79}  & 19.98 & \multicolumn{1}{c|}{2059.23} & \multicolumn{1}{c|}{31.39} & \multicolumn{1}{c|}{4613.85}  & 23.75 \\ \hline
				T=15\%                     & \multicolumn{1}{c|}{1517.86} & \multicolumn{1}{c|}{23.13} & \multicolumn{1}{c|}{4004.96}  & 20.68 & \multicolumn{1}{c|}{2080.66} & \multicolumn{1}{c|}{31.71} & \multicolumn{1}{c|}{4606.02}  & 23.71 \\ \hline
				T=20\%                     & \multicolumn{1}{c|}{1538.93} & \multicolumn{1}{c|}{23.46} & \multicolumn{1}{c|}{4096.60}  & 21.16 & \multicolumn{1}{c|}{2082.35} & \multicolumn{1}{c|}{31.74} & \multicolumn{1}{c|}{4694.64}  & 24.17 \\ \hline
				T=25\%                     & \multicolumn{1}{c|}{1553.29} & \multicolumn{1}{c|}{23.67} & \multicolumn{1}{c|}{4205.20}  & 21.72 & \multicolumn{1}{c|}{2074.74} & \multicolumn{1}{c|}{31.62} & \multicolumn{1}{c|}{4785.93}  & 24.64 \\ \hline
				Ref                        & \multicolumn{1}{c|}{6561.00} & \multicolumn{1}{c|}{}      & \multicolumn{1}{c|}{19362.00} &       & \multicolumn{1}{c|}{6561.00} & \multicolumn{1}{c|}{}      & \multicolumn{1}{c|}{19427.00} &       \\ \hline
			\end{tabular}
		\end{adjustbox}
		\caption{Performance summary of the VEM in terms of its convergence behaviour in the PSE with respect to the number of nodes/vertices in the discretization when using a combination of the displacement-based and ${Z^{2}}$-like refinement procedures for problem~C(5).
			\label{tab:PerformanceSquarePlateNarrowPunchPSE}}
	\end{table}
	\FloatBarrier
	
\begin{figure}[ht!]
		\centering
		\begin{subfigure}[t]{0.5\textwidth}
			\centering
			\includegraphics[width=0.8\textwidth]{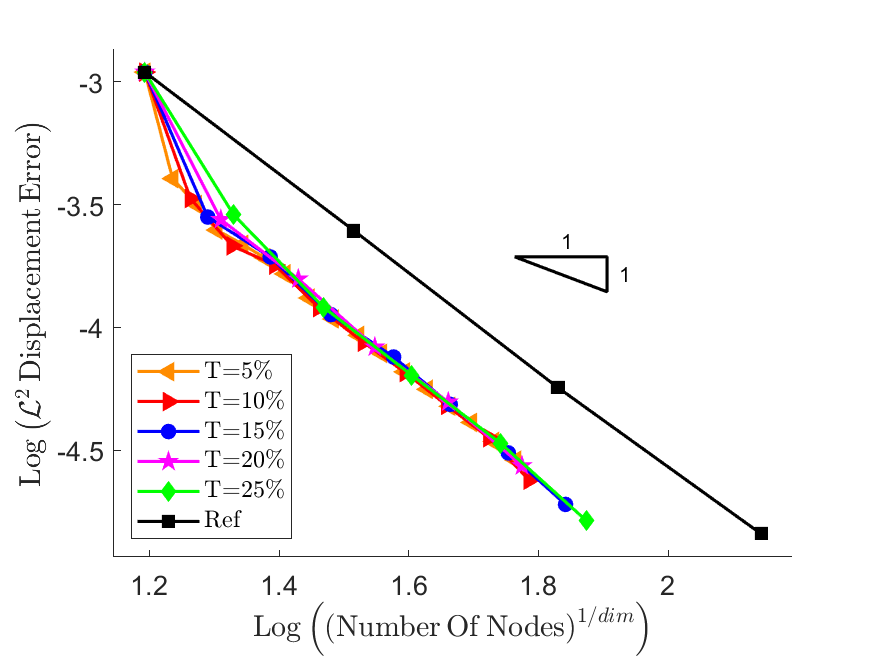}
			\caption{Displacement error - Compressible}
		\end{subfigure}%
		\begin{subfigure}[t]{0.5\textwidth}
			\centering
			\includegraphics[width=0.8\textwidth]{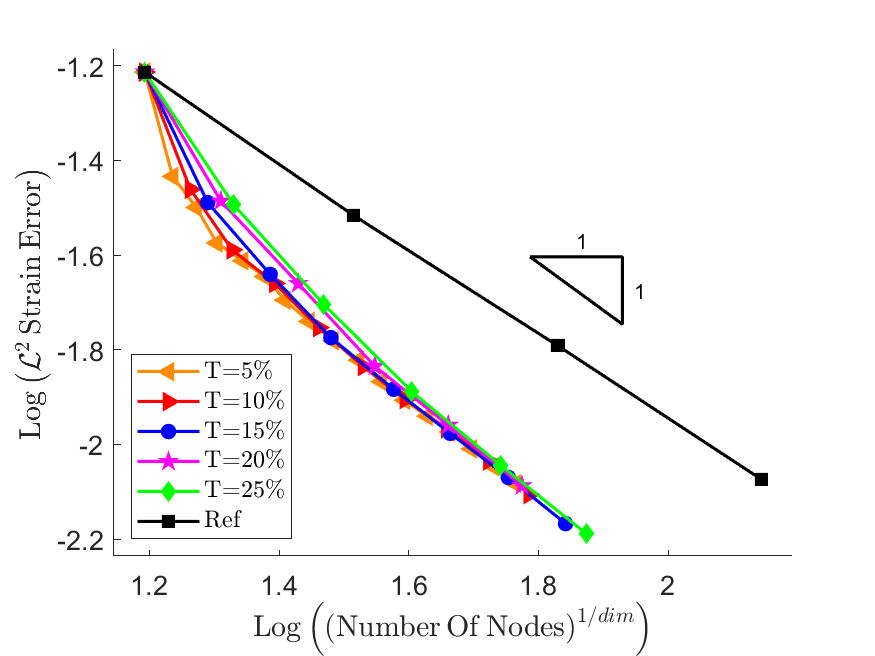}
			\caption{Strain error - Compressible}
		\end{subfigure}
		\vskip \baselineskip 
		\vspace*{-3mm}
		\begin{subfigure}[t]{0.5\textwidth}
			\centering
			\includegraphics[width=0.8\textwidth]{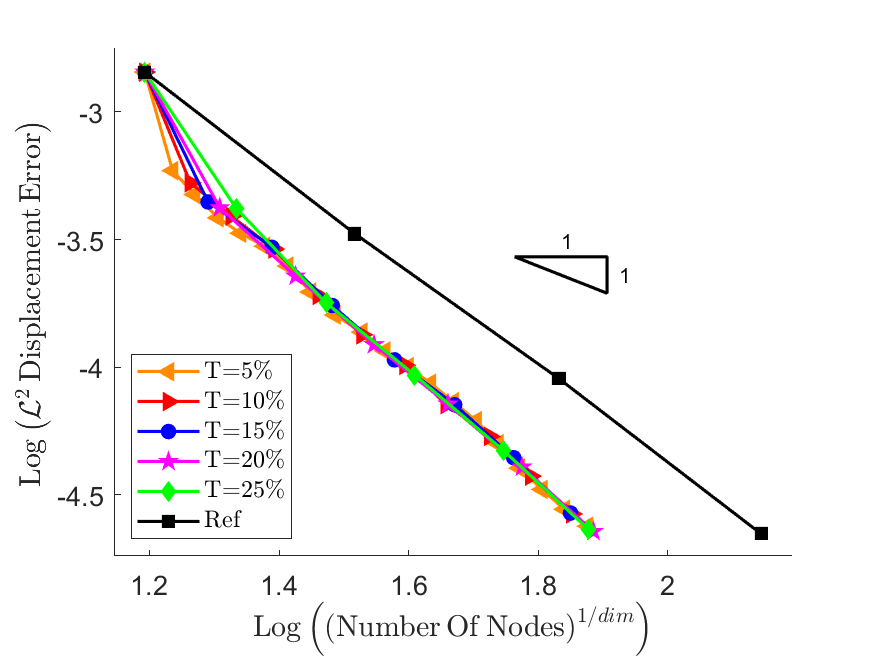}
			\caption{Displacement error - Nearly incompressible}
		\end{subfigure}%
		\begin{subfigure}[t]{0.5\textwidth}
			\centering
			\includegraphics[width=0.8\textwidth]{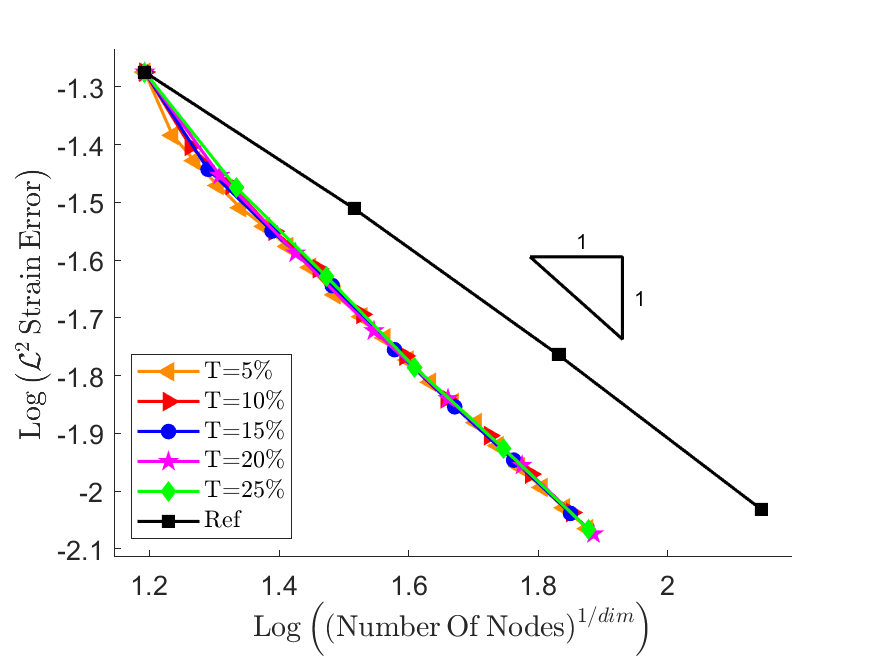}
			\caption{Strain error - Nearly incompressible}
		\end{subfigure}
		\caption{Displacement and strain $\mathcal{L}^{2}$ error components vs $n_{\rm v}$ for problem~C(5) using a combination of the displacement-based and ${Z^{2}}$-like refinement procedures with a variety of choices of $T$ for Voronoi meshes with compressible and nearly incompressible Poisson's ratios.
			\label{fig:SquarePlateNarrowPunchConvergenceComponents}}
	\end{figure}

	The convergence behaviour in the ${\mathcal{L}^{2}}$ error norm of the displacement and strain field approximations when using a combination of the displacement-based and ${Z^{2}}$-like refinement procedures for problem~C(5) is plotted in Figure~\ref{fig:SquarePlateNarrowPunchConvergenceComponents} against the number of vertices/nodes in the discretization for Voronoi meshes with compressible and nearly incompressible Poisson's ratios. 
	Similar to the convergence behaviour observed in the ${\mathcal{H}^{1}}$ error norm, lower choices of $T$ initially exhibit slightly faster convergence rates. However, in the fine mesh range for both mesh types and in the cases of both compressibility and near-incompressibility all choices of $T$ exhibit similar convergence behaviour and significantly outperform the reference procedure in both error components.
	
	\FloatBarrier
	\FloatBarrier
	
	\subsection{Threshold comparison}
	\label{subsec:ThresholdComp}
	
	To determine the best refinement threshold percentage for a particular choice, or particular combination, of refinement indicator(s) the combined performance, as measured by the PRE, in terms of the number of nodes/vertices and run time (without remeshing time) is considered. Furthermore, the performance over the full range of mesh types, Poisson's ratios, and sample problems is considered. 
	
	It is well-known that computational run time is strongly influenced by the specific implementation of the necessary algorithms Thus, to reduce the influence of implementation on the determination of the best choice of $T$ the time taken to perform element refinement is excluded. However, the time taken to compute the refinement indicators and mark elements for refinement is included. This approach results in a choice of $T$ that is optimal when the implementation of the remeshing is instantaneous. 
	
	The threshold comparison for problem~A(1) using the displacement-based refinement procedure is presented in Table~\ref{tab:ThreshComp_ProblemOne}. Here, the PRE in terms of the convergence behaviour in the $\mathcal{H}^{1}$ error norm with respect to the number of vertices/nodes in the discretization and with respect to run time (without remeshing) is summarized. Additionally, a combined average PRE of both metrics (number of nodes/vertices and run time) over the range of mesh types and Poisson's ratios considered is computed. This average PRE is presented in the right-hand portion of the table in both unsorted and sorted formats. In the sorted format, the various refinement thresholds considered are sorted by their average PRE in ascending order. For clarity, a lower average PRE corresponds to better/more efficient overall performance. For problem~A(1) better performance is exhibited by larger values of $T$. The choice of ${T=25\%}$ exhibits the best performance, with ${T=20\%}$ exhibiting similarly good performance. The corresponding tables for the full range of example problems can be found in the supplementary material.
	
	\FloatBarrier
	\begin{table}[ht!]
		\centering 
		\begin{adjustbox}{max width=\textwidth}
			\begin{tabular}{|c|cccc|cccc|c|cc|}
				\hline
				\multirow{3}{*}{Threshold} & \multicolumn{4}{c|}{Compressible}                                                                    & \multicolumn{4}{c|}{Nearly-incompressible}                                                           & \multirow{2}{*}{Unsorted} & \multicolumn{2}{c|}{\multirow{2}{*}{Sorted}} \\ \cline{2-9}
				& \multicolumn{2}{c|}{Structured}                             & \multicolumn{2}{c|}{Voronoi}           & \multicolumn{2}{c|}{Structured}                             & \multicolumn{2}{c|}{Voronoi}           &                           & \multicolumn{2}{c|}{}                        \\ \cline{2-12} 
				& \multicolumn{1}{c|}{nNodes} & \multicolumn{1}{c|}{Run time} & \multicolumn{1}{c|}{nNodes} & Run time & \multicolumn{1}{c|}{nNodes} & \multicolumn{1}{c|}{Run time} & \multicolumn{1}{c|}{nNodes} & Run time & Avg PRE                   & \multicolumn{1}{c|}{Threshold}   & Avg. PRE  \\ \hline
				T=5\%                      & \multicolumn{1}{c|}{8.92}   & \multicolumn{1}{c|}{34.68}    & \multicolumn{1}{c|}{11.56}  & 35.83    & \multicolumn{1}{c|}{32.43}  & \multicolumn{1}{c|}{154.21}   & \multicolumn{1}{c|}{26.67}  & 105.42   & 51.22                     & \multicolumn{1}{c|}{T=25\%}      & 25.88     \\ \hline
				T=10\%                     & \multicolumn{1}{c|}{9.73}   & \multicolumn{1}{c|}{22.41}    & \multicolumn{1}{c|}{11.83}  & 23.01    & \multicolumn{1}{c|}{32.84}  & \multicolumn{1}{c|}{89.21}    & \multicolumn{1}{c|}{26.42}  & 62.58    & 34.75                     & \multicolumn{1}{c|}{T=20\%}      & 26.32     \\ \hline
				T=15\%                     & \multicolumn{1}{c|}{9.28}   & \multicolumn{1}{c|}{16.01}    & \multicolumn{1}{c|}{12.72}  & 19.75    & \multicolumn{1}{c|}{32.63}  & \multicolumn{1}{c|}{67.65}    & \multicolumn{1}{c|}{26.46}  & 46.69    & 28.90                     & \multicolumn{1}{c|}{T=15\%}      & 28.90     \\ \hline
				T=20\%                     & \multicolumn{1}{c|}{9.02}   & \multicolumn{1}{c|}{12.99}    & \multicolumn{1}{c|}{14.74}  & 19.90    & \multicolumn{1}{c|}{32.21}  & \multicolumn{1}{c|}{55.08}    & \multicolumn{1}{c|}{26.56}  & 40.04    & 26.32                     & \multicolumn{1}{c|}{T=10\%}      & 34.75     \\ \hline
				T=25\%                     & \multicolumn{1}{c|}{10.56}  & \multicolumn{1}{c|}{13.68}    & \multicolumn{1}{c|}{17.44}  & 21.82    & \multicolumn{1}{c|}{31.77}  & \multicolumn{1}{c|}{48.32}    & \multicolumn{1}{c|}{26.66}  & 36.80    & 25.88                     & \multicolumn{1}{c|}{T=5\%}       & 51.22     \\ \hline
			\end{tabular}
		\end{adjustbox}
		\caption{Threshold comparison for problem~A(1) using the displacement-based refinement procedure.
			\label{tab:ThreshComp_ProblemOne}}
	\end{table}
	\FloatBarrier
	
	\subsubsection{Overall}
	\label{subsubsec:ThreshComp_Overall}
	
	The overall performance of the displacement-based refinement procedure (measured by the combined average PRE as computed in Table~\ref{tab:ThreshComp_ProblemOne}) for the full range of sample problems is presented in Table~\ref{tab:ThreshComp_Overall}. Here, an overall average PRE over the range of problems is computed and presented in the right-hand portion of the table in both unsorted and sorted formats. For the displacement-based refinement indicator, the overall best choice of the refinement threshold percentage is ${T=20\%}$.
	
	\FloatBarrier
	\begin{table}[ht!]
		\centering 
		\begin{adjustbox}{max width=\textwidth}
			\begin{tabular}{|c|c|c|c|c|c|c|c|cc|}
				\hline
				\multirow{2}{*}{Threshold} & \multirow{2}{*}{Prob. 1} & \multirow{2}{*}{Prob. 2} & \multirow{2}{*}{Prob. 3} & \multirow{2}{*}{Prob. 4} & \multirow{2}{*}{Prob. 5} & \multirow{2}{*}{Prob. 6} & Unsorted & \multicolumn{2}{c|}{Sorted}               \\ \cline{8-10} 
				&                          &                          &                          &                          &                          &                          & Avg PRE  & \multicolumn{1}{c|}{Threshold} & Avg. PRE \\ \hline
				T=5\%                      & 51.22                    & 32.17                    & 10.61                    & 13.53                    & 56.01                    & 9.25                     & 28.80    & \multicolumn{1}{c|}{T=20\%}    & 17.52    \\ \hline
				T=10\%                     & 34.75                    & 22.52                    & 7.54                     & 9.75                     & 38.59                    & 9.91                     & 20.51    & \multicolumn{1}{c|}{T=25\%}    & 17.90    \\ \hline
				T=15\%                     & 28.90                    & 18.97                    & 7.83                     & 9.04                     & 32.72                    & 11.79                    & 18.21    & \multicolumn{1}{c|}{T=15\%}    & 18.21    \\ \hline
				T=20\%                     & 26.32                    & 17.68                    & 8.22                     & 9.00                     & 30.16                    & 13.76                    & 17.52    & \multicolumn{1}{c|}{T=10\%}    & 20.51    \\ \hline
				T=25\%                     & 25.88                    & 17.97                    & 9.83                     & 9.58                     & 28.25                    & 15.88                    & 17.90    & \multicolumn{1}{c|}{T=5\%}     & 28.80    \\ \hline
			\end{tabular}
		\end{adjustbox}
		\caption{Overall threshold comparison for the displacement-based refinement procedure.
			\label{tab:ThreshComp_Overall}}
	\end{table}
	\FloatBarrier
	
	The procedure demonstrated above to determine the best refinement threshold percentage for the displacement-based refinement procedure was repeated for each of the refinement indicators, and combinations of refinement indicators, considered. The results are summarized in Table~\ref{tab:ThreshComp_Summary} and are qualitatively similar to those presented in \cite{vanHuyssteen2022} where it was also found that the best threshold in the case of combined indicators is slightly lower than that of individual indicators. Additionally, it is noted that the best choices of the thresholds presented in \cite{vanHuyssteen2022} are slightly lower than those presented here because the run time of the remeshing process was not excluded from the calculations.
	
	\FloatBarrier
	\begin{table}[ht!]
		\centering 
		\begin{adjustbox}{max width=\textwidth}
			\begin{tabular}{|c|c|}
				\hline
				Ref.   Indicator   & Best threshold \\ \hline
				Disp               & $T=20\%$         \\ \hline
				SJ                 & $T=25\%$         \\ \hline
				Z2                 & $T=25\%$         \\ \hline
				Disp \& SJ         & $T=15\%$         \\ \hline
				Disp \& $Z2$       & $T=15\%$         \\ \hline
			\end{tabular}
		\end{adjustbox}
		\caption{Threshold comparison summary.
			\label{tab:ThreshComp_Summary}}
	\end{table}
	\FloatBarrier
	
	\subsection{Method comparison}
	\label{subsec:MethodComp}
	
	In this section the performance of the various mesh refinement indicators, each with its best choice of refinement threshold percentage, are comparatively assessed to determine the best overall refinement procedure. For brevity, results are only presented for problem~A(1). The results for the full range of problems can be found in the supplementary material.
	
\begin{figure}[ht!]
		\centering
		\begin{subfigure}[t]{0.5\textwidth}
			\centering
			\includegraphics[width=0.95\textwidth]{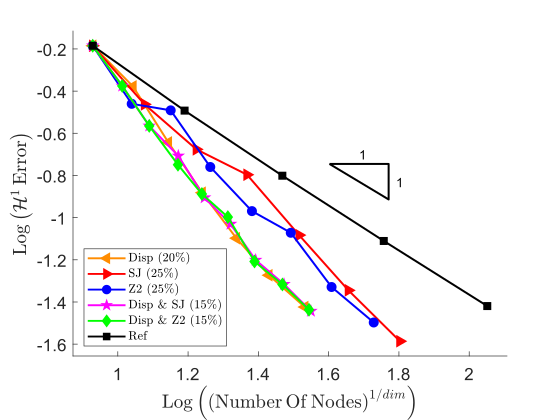}
			\caption{Structured mesh}
		\end{subfigure}%
		\begin{subfigure}[t]{0.5\textwidth}
			\centering
			\includegraphics[width=0.95\textwidth]{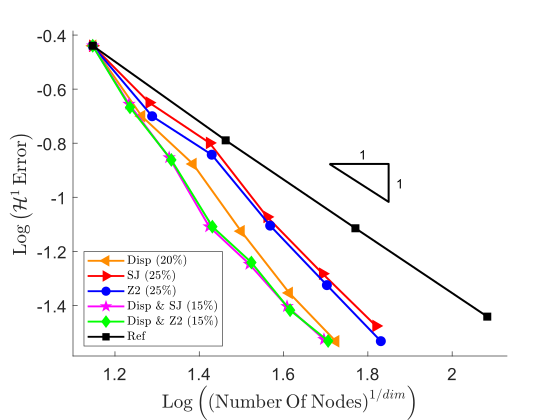}
			\caption{Voronoi mesh}
		\end{subfigure}
		\caption{$\mathcal{H}^{1}$ error vs $n_{\rm v}$ for problem~A(1) using the various refinement procedures considered for structured and unstructured/Voronoi meshes with a compressible Poisson's ratio.
			\label{fig:MethComp_PlateWithHoleTractionConvergenceNumberOfNodes}}
	\end{figure}

	The convergence behaviour in the ${\mathcal{H}^{1}}$ error norm of the VEM for problem~A(1) for the various refinement procedures considered is depicted in Figure~\ref{fig:MethComp_PlateWithHoleTractionConvergenceNumberOfNodes} on a logarithmic scale. Here, the ${\mathcal{H}^{1}}$ error is plotted against the number of vertices/nodes in the discretization for structured and unstructured/Voronoi meshes with a compressible Poison's ratio. For both mesh types it is clear that the displacement-based, combined displacement-based and strain jump-based, and combined displacement-based and ${Z^{2}}$-like procedures exhibit similarly good convergence behaviour and converge much faster than the strain jump-based and ${Z^{2}}$-like procedures. In the case of unstructured/Voronoi meshes the combined displacement-based and strain jump-based, and combined displacement-based and ${Z^{2}}$-like procedures exhibit slightly superior convergence behaviour to the displacement-based procedure.
	However, for both mesh types, the proposed refinement procedures all exhibit significantly improved convergence rates compared to the reference procedure.
	
	\FloatBarrier
	\FloatBarrier
	
	The performance of the VEM in terms of its convergence behaviour in the $\mathcal{H}^{1}$ error norm with respect to the number of vertices/nodes in the discretization for the various refinement procedures considered for problem~A(1) is summarized in Table~\ref{tab:MethComp_PerformancePlateWithHoleTractionNumberOfNodes}. Here, the performance, as measured by the PRE, is presented for the cases of structured and Voronoi meshes with compressible and nearly incompressible Poisson's ratios. Additionally, a combined average PRE over the range of mesh types and Poisson's ratios considered is computed and presented in the right-hand portion of the table. As observed in Figure~\ref{fig:MethComp_PlateWithHoleTractionConvergenceNumberOfNodes} good overall performance is exhibited by the displacement-based, combined displacement-based and strain jump-based, and combined displacement-based and ${Z^{2}}$-like procedures. Furthermore, the superior performance of the combined displacement-based and strain jump-based, and combined displacement-based and ${Z^{2}}$-like procedures in the case of unstructured/Voronoi meshes means that these two procedures exhibit the best overall performance. However, all proposed procedures represent significant improvements in performance when compared to the reference procedure.
	
	\FloatBarrier
	\begin{table}[ht!]
		\centering 
		\begin{adjustbox}{max width=\textwidth}
			\begin{tabular}{|c|cccc|cccc|ccc|}
				\hline
				\multirow{3}{*}{Method} & \multicolumn{4}{c|}{Compressible}                                                                  & \multicolumn{4}{c|}{Nearly-incompressible}                                                         & \multicolumn{1}{c|}{\multirow{2}{*}{Unsorted}} & \multicolumn{2}{c|}{\multirow{2}{*}{Sorted}}  \\ \cline{2-9}
				& \multicolumn{2}{c|}{Structured}                            & \multicolumn{2}{c|}{Voronoi}          & \multicolumn{2}{c|}{Structured}                            & \multicolumn{2}{c|}{Voronoi}          & \multicolumn{1}{c|}{}                          & \multicolumn{2}{c|}{}                         \\ \cline{2-12} 
				& \multicolumn{1}{c|}{nNodes}   & \multicolumn{1}{c|}{PRE}   & \multicolumn{1}{c|}{nNodes}   & PRE   & \multicolumn{1}{c|}{nNodes}   & \multicolumn{1}{c|}{PRE}   & \multicolumn{1}{c|}{nNodes}   & PRE   & \multicolumn{1}{c|}{Avg PRE}                   & \multicolumn{1}{c|}{Method}        & Avg. PRE \\ \hline
				Disp-20\%               & \multicolumn{1}{c|}{1142.64}  & \multicolumn{1}{c|}{9.02}  & \multicolumn{1}{c|}{2162.95}  & 14.74 & \multicolumn{1}{c|}{4081.88}  & \multicolumn{1}{c|}{32.21} & \multicolumn{1}{c|}{3862.79}  & 26.56 & \multicolumn{1}{c|}{20.63}                     & \multicolumn{1}{c|}{Disp\&SJ-15\%} & 20.20    \\ \hline
				SJ-25\%                 & \multicolumn{1}{c|}{2530.22}  & \multicolumn{1}{c|}{19.97} & \multicolumn{1}{c|}{3919.55}  & 26.71 & \multicolumn{1}{c|}{4517.62}  & \multicolumn{1}{c|}{35.65} & \multicolumn{1}{c|}{4139.06}  & 28.46 & \multicolumn{1}{c|}{27.70}                     & \multicolumn{1}{c|}{Disp\&Z2-15\%} & 20.20    \\ \hline
				Z2-25\%                 & \multicolumn{1}{c|}{2216.05}  & \multicolumn{1}{c|}{17.49} & \multicolumn{1}{c|}{3543.40}  & 24.15 & \multicolumn{1}{c|}{4713.14}  & \multicolumn{1}{c|}{37.19} & \multicolumn{1}{c|}{4768.48}  & 32.79 & \multicolumn{1}{c|}{27.90}                     & \multicolumn{1}{c|}{Disp-20\%}     & 20.63    \\ \hline
				Disp\&SJ-15\%           & \multicolumn{1}{c|}{1175.82}  & \multicolumn{1}{c|}{9.28}  & \multicolumn{1}{c|}{1865.98}  & 12.72 & \multicolumn{1}{c|}{4164.37}  & \multicolumn{1}{c|}{32.86} & \multicolumn{1}{c|}{3772.34}  & 25.94 & \multicolumn{1}{c|}{20.20}                     & \multicolumn{1}{c|}{SJ-25\%}       & 27.70    \\ \hline
				Disp\&Z2-15\%           & \multicolumn{1}{c|}{1164.27}  & \multicolumn{1}{c|}{9.19}  & \multicolumn{1}{c|}{1855.68}  & 12.65 & \multicolumn{1}{c|}{4153.59}  & \multicolumn{1}{c|}{32.78} & \multicolumn{1}{c|}{3808.53}  & 26.19 & \multicolumn{1}{c|}{20.20}                     & \multicolumn{1}{c|}{Z2-25\%}       & 27.90    \\ \hline
				Ref                     & \multicolumn{1}{c|}{12672.00} & \multicolumn{1}{c|}{}      & \multicolumn{1}{c|}{14675.00} &       & \multicolumn{1}{c|}{12672.00} & \multicolumn{1}{c|}{}      & \multicolumn{1}{c|}{14542.00} &       & \multicolumn{3}{c|}{}                                                                          \\ \hline
			\end{tabular}
		\end{adjustbox}
		\caption{Performance summary of the VEM in terms of its convergence behaviour in the $\mathcal{H}^{1}$ error norm with respect to the number of vertices/nodes in the discretization for the various refinement procedures considered for problem~A(1).
			\label{tab:MethComp_PerformancePlateWithHoleTractionNumberOfNodes}}
	\end{table}
	\FloatBarrier
	
\begin{figure}[ht!]
		\centering
		\begin{subfigure}[t]{0.5\textwidth}
			\centering
			\includegraphics[width=0.95\textwidth]{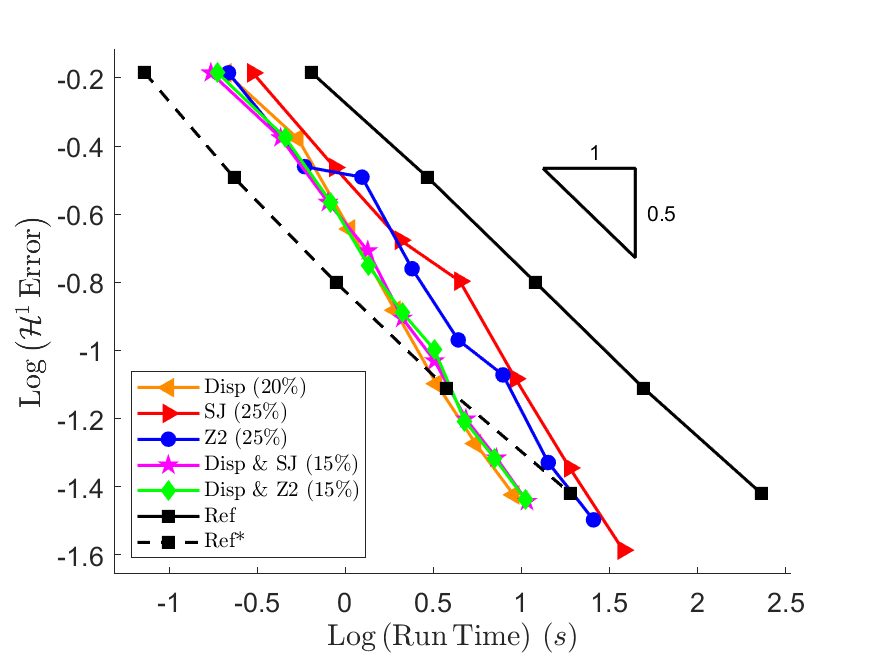}
			\caption{Structured mesh}
		\end{subfigure}%
		\begin{subfigure}[t]{0.5\textwidth}
			\centering
			\includegraphics[width=0.95\textwidth]{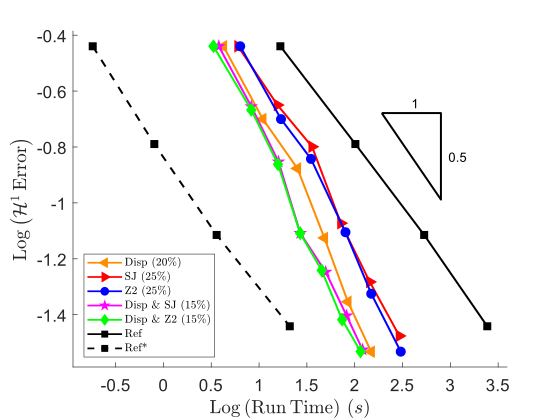}
			\caption{Voronoi mesh}
		\end{subfigure}
		\caption{$\mathcal{H}^{1}$ error vs run time (including remeshing time) for problem~A(1) for the various refinement procedures considered for structured and unstructured/Voronoi meshes with a compressible Poisson's ratio.
			\label{fig:MethComp_PlateWithHoleTractionConvergenceRunTime}}
	\end{figure}

	The convergence behaviour in the ${\mathcal{H}^{1}}$ error norm of the VEM for problem~A(1) for the various refinement procedures considered is depicted in Figure~\ref{fig:MethComp_PlateWithHoleTractionConvergenceRunTime} on a logarithmic scale. Here, the ${\mathcal{H}^{1}}$ error is plotted against run time (including remeshing time) for structured and unstructured/Voronoi meshes with a compressible Poisson's ratio. 
	Since the best choice of $T$ has been determined independently of the implementation of the mesh refinement algorithms, it is chosen now to include the remeshing time in the comparative evaluation of the various refinement procedures. This allows for an objective comparison of the various procedures against the reference procedure with the current implementation. 
	Hereinafter, all run time related computations will include remeshing time.
	The run time of the reference procedure is presented in two ways. Firstly, and similar to that described in Section~\ref{subsec:PlateWithHoleTraction}, a cumulative run time including the remeshing process is presented and is indicated by a solid black curve with square markers.
	Secondly, a non-cumulative run time which does not include any remeshing processes is presented and is indicated by a dashed black curve with square markers. Here, each marker represents the run time required to solve the problem for a particular level of mesh refinement.
	Similar to the convergence behaviour in Figure~\ref{fig:MethComp_PlateWithHoleTractionConvergenceNumberOfNodes} good overall performance is exhibited by the displacement-based, combined displacement-based and strain jump-based, and combined displacement-based and ${Z^{2}}$-like procedures. In the case of unstructured/Voronoi meshes the best performance is again exhibited by the combined displacement-based and strain jump-based, and combined displacement-based and ${Z^{2}}$-like procedures. Furthermore, all proposed procedures represent significant improvements in performance compared to the cumulative run time of the reference procedure.
	
	A noteworthy feature is the difference in convergence rates of the proposed adaptive refinement procedures and that of the reference procedure. As expected, the run time of the first step for all adaptive procedures is somewhere between the run time of the first step for the non-cumulative and cumulative reference procedures. This is expected because the cumulative procedure includes the time taken to refine every element in the domain, while the adaptive procedure computes the refinement indicators and then refines an identified subset of the elements. Furthermore, the non-cumulative procedure involves no remeshing and is thus the fastest procedure for the first refinement step. Since the adaptive procedures exhibit superior convergence rates to the reference procedure there will be some level of mesh refinement after which the cumulative run time of the adaptive procedures is less than that of the non-cumulative reference procedure. Therefore, after this particular level of refinement is reached the adaptive procedures, even starting from a coarse mesh, will generate solutions of the same accuracy as the reference procedure in less time than it takes the non-cumulative procedure to solve the problem using a `pre-generated´ mesh. For this example problem, in the case of structured meshes, the displacement-based, combined displacement-based and strain jump-based, and combined displacement-based and ${Z^{2}}$-like procedures outperform the non-cumulative reference procedure after the fourth reference mesh considered (indicated by the fourth square marker). In the range of mesh refinement levels considered for the unstructured/Voronoi meshes the level of refinement at which the adaptive procedures outperform the non-cumulative reference procedure is not reached. However, it is easy to see that for sufficiently fine (reference) meshes the adaptive procedures will outperform the non-cumulative procedure.
	
	\FloatBarrier
	\FloatBarrier
	
	The performance of the VEM in terms of its convergence behaviour in the $\mathcal{H}^{1}$ error norm with respect to run time (including remeshing time) for the various refinement procedures considered for problem~A(1) is summarized in Table~\ref{tab:MethComp_PerformancePlateWithHoleTractionRunTime}. Here, the performance, as measured by the PRE, is presented for the cases of structured and Voronoi meshes with compressible and nearly incompressible Poisson's ratios. Additionally, a combined average PRE over the range of mesh types and Poisson's ratios considered is computed and presented in the right-hand portion of the table.
	As expected, good overall performance is exhibited by the displacement-based, combined displacement-based and strain jump-based, and combined displacement-based and ${Z^{2}}$-like procedures. For this example the best overall performance is exhibited by the displacement-based procedure. However, all proposed procedures significantly outperform the (cumulative) reference procedure.
	
	\FloatBarrier
	\begin{table}[ht!]
		\centering 
		\begin{adjustbox}{max width=\textwidth}
			\begin{tabular}{|c|cccc|cccc|ccc|}
				\hline
				\multirow{3}{*}{Method} & \multicolumn{4}{c|}{Compressible}                                                                  & \multicolumn{4}{c|}{Nearly-incompressible}                                                         & \multicolumn{1}{c|}{\multirow{2}{*}{Unsorted}} & \multicolumn{2}{c|}{\multirow{2}{*}{Sorted}}  \\ \cline{2-9}
				& \multicolumn{2}{c|}{Structured}                            & \multicolumn{2}{c|}{Voronoi}          & \multicolumn{2}{c|}{Structured}                            & \multicolumn{2}{c|}{Voronoi}          & \multicolumn{1}{c|}{}                          & \multicolumn{2}{c|}{}                         \\ \cline{2-12} 
				& \multicolumn{1}{c|}{Run time} & \multicolumn{1}{c|}{PRE}   & \multicolumn{1}{c|}{Run time} & PRE   & \multicolumn{1}{c|}{Run time} & \multicolumn{1}{c|}{PRE}   & \multicolumn{1}{c|}{Run time} & PRE   & \multicolumn{1}{c|}{Avg PRE}                   & \multicolumn{1}{c|}{Method}        & Avg. PRE \\ \hline
				Disp-20\%               & \multicolumn{1}{c|}{8.91}     & \multicolumn{1}{c|}{3.87}  & \multicolumn{1}{c|}{111.77}   & 4.59  & \multicolumn{1}{c|}{35.81}    & \multicolumn{1}{c|}{15.26} & \multicolumn{1}{c|}{220.67}   & 9.36  & \multicolumn{1}{c|}{8.27}                      & \multicolumn{1}{c|}{Disp-20\%}     & 8.27     \\ \hline
				SJ-25\%                 & \multicolumn{1}{c|}{23.32}    & \multicolumn{1}{c|}{10.14} & \multicolumn{1}{c|}{258.25}   & 10.61 & \multicolumn{1}{c|}{43.39}    & \multicolumn{1}{c|}{18.49} & \multicolumn{1}{c|}{287.77}   & 12.20 & \multicolumn{1}{c|}{12.86}                     & \multicolumn{1}{c|}{Disp\&Z2-15\%} & 8.55     \\ \hline
				Z2-25\%                 & \multicolumn{1}{c|}{19.51}    & \multicolumn{1}{c|}{8.48}  & \multicolumn{1}{c|}{222.07}   & 9.13  & \multicolumn{1}{c|}{44.81}    & \multicolumn{1}{c|}{19.09} & \multicolumn{1}{c|}{322.21}   & 13.66 & \multicolumn{1}{c|}{12.59}                     & \multicolumn{1}{c|}{Disp\&SJ-15\%} & 8.66     \\ \hline
				Disp\&SJ-15\%           & \multicolumn{1}{c|}{9.98}     & \multicolumn{1}{c|}{4.34}  & \multicolumn{1}{c|}{92.72}    & 3.81  & \multicolumn{1}{c|}{43.52}    & \multicolumn{1}{c|}{18.55} & \multicolumn{1}{c|}{187.61}   & 7.96  & \multicolumn{1}{c|}{8.66}                      & \multicolumn{1}{c|}{Z2-25\%}       & 12.59    \\ \hline
				Disp\&Z2-15\%           & \multicolumn{1}{c|}{9.89}     & \multicolumn{1}{c|}{4.30}  & \multicolumn{1}{c|}{81.05}    & 3.33  & \multicolumn{1}{c|}{43.23}    & \multicolumn{1}{c|}{18.42} & \multicolumn{1}{c|}{191.95}   & 8.14  & \multicolumn{1}{c|}{8.55}                      & \multicolumn{1}{c|}{SJ-25\%}       & 12.86    \\ \hline
				Ref                     & \multicolumn{1}{c|}{229.93}   & \multicolumn{1}{c|}{}      & \multicolumn{1}{c|}{2433.15}  &       & \multicolumn{1}{c|}{234.67}   & \multicolumn{1}{c|}{}      & \multicolumn{1}{c|}{2358.32}  &       & \multicolumn{3}{c|}{}                                                                          \\ \hline
			\end{tabular}
		\end{adjustbox}
		\caption{Performance summary of the VEM in terms of its convergence behaviour in the $\mathcal{H}^{1}$ error norm with respect to run time (including remeshing time) for the various refinement procedures considered for problem~A(1).
			\label{tab:MethComp_PerformancePlateWithHoleTractionRunTime}}
	\end{table}
	\FloatBarrier
	
\begin{figure}[ht!]
		\centering
		\begin{subfigure}[t]{0.5\textwidth}
			\centering
			\includegraphics[width=0.95\textwidth]{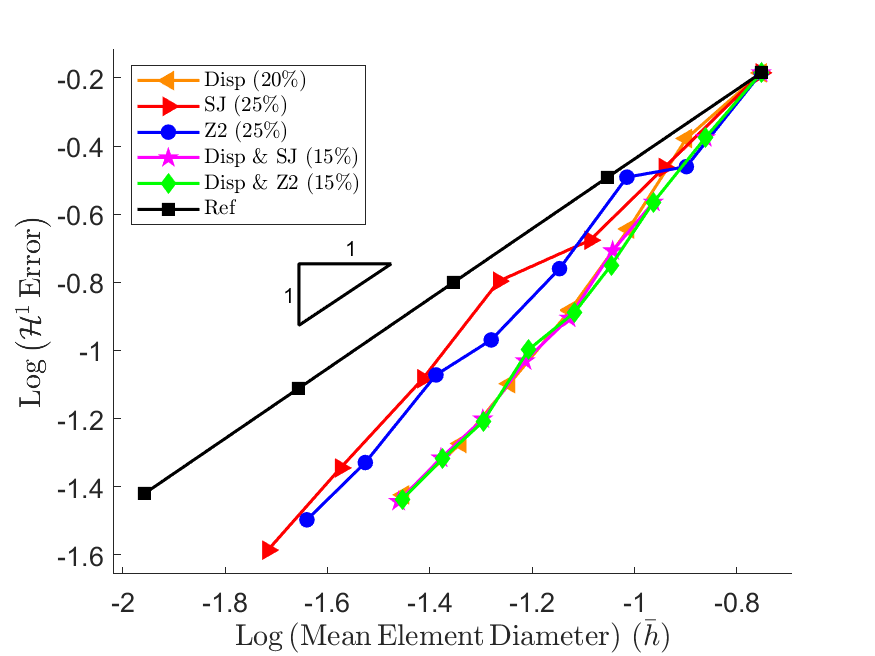}
			\caption{Structured mesh}
		\end{subfigure}%
		\begin{subfigure}[t]{0.5\textwidth}
			\centering
			\includegraphics[width=0.95\textwidth]{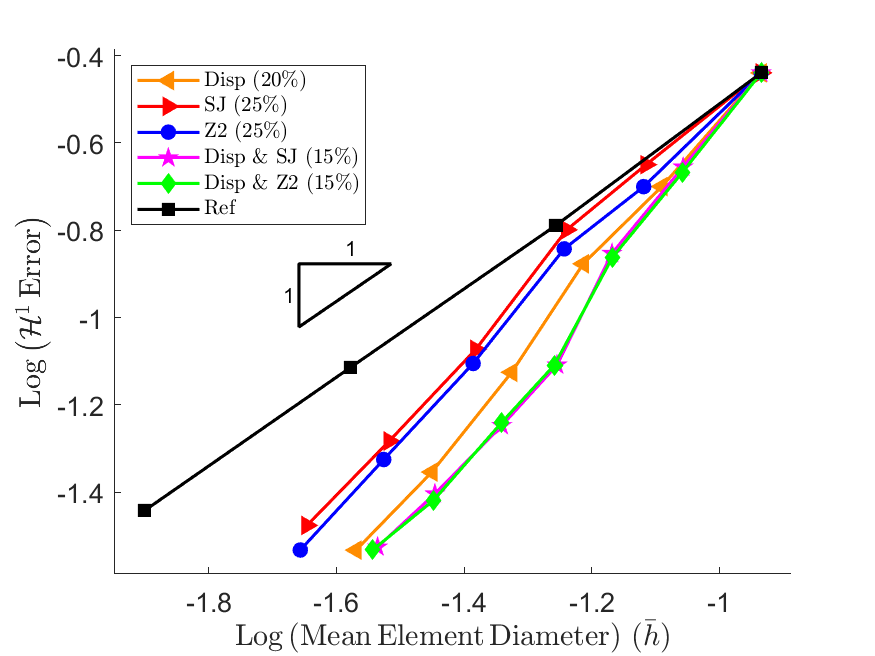}
			\caption{Voronoi mesh}
		\end{subfigure}
		\caption{$\mathcal{H}^{1}$ error vs mesh size for problem~A(1) for the various refinement procedures considered for structured and unstructured/Voronoi meshes with a compressible Poisson's ratio.
			\label{fig:MethComp_PlateWithHoleTractionConvergenceMeshSize}}
	\end{figure}

	The convergence behaviour in the ${\mathcal{H}^{1}}$ error norm of the VEM for problem~A(1) for the various refinement procedures considered is depicted in Figure~\ref{fig:MethComp_PlateWithHoleTractionConvergenceMeshSize} on a logarithmic scale. Here, the ${\mathcal{H}^{1}}$ error is plotted against mesh size as measured by the mean element diameter for structured and unstructured/Voronoi meshes with a compressible Poisson's ratio. 
	The convergence behaviour is qualitatively similar to that observed in Figure~\ref{fig:MethComp_PlateWithHoleTractionConvergenceNumberOfNodes} with respect to the number of nodes/vertices.
	For both mesh types the displacement-based, combined displacement-based and strain jump-based, and combined displacement-based and ${Z^{2}}$-like procedures exhibit similarly good convergence behaviour and converge much faster than the strain jump-based and ${Z^{2}}$-like procedures. In the case of unstructured/Voronoi meshes the combined displacement-based and strain jump-based, and combined displacement-based and ${Z^{2}}$-like procedures exhibit slightly superior convergence behaviour to the displacement-based procedure.
	However, for both mesh types, the proposed refinement procedures all exhibit superior convergence behaviour to, and significantly outperform, the reference procedure.
	
	\FloatBarrier
	\FloatBarrier
	
	The performance of the VEM in terms of its convergence behaviour in the $\mathcal{H}^{1}$ error norm with respect to mesh size for the various refinement procedures considered for problem~A(1) is summarized in Table~\ref{tab:MethComp_PerformancePlateWithHoleTractionMeshSize}. Here, the performance, as measured by the PRE*, is presented for the cases of structured and Voronoi meshes with compressible and nearly incompressible Poisson's ratios. Additionally, a combined average PRE* over the range of mesh types and Poisson's ratios considered is computed and presented in the right-hand portion of the table.
	As observed in Figure~\ref{fig:MethComp_PlateWithHoleTractionConvergenceMeshSize} the displacement-based, combined displacement-based and strain jump-based, and combined displacement-based and ${Z^{2}}$-like procedures exhibit similarly good convergence performance. Overall, however, the best performance is exhibited by the combined displacement-based and ${Z^{2}}$-like procedures.
	
	\FloatBarrier
	\begin{table}[ht!]
		\centering 
		\begin{adjustbox}{max width=\textwidth}
			\begin{tabular}{|c|cccc|cccc|ccc|}
				\hline
				\multirow{3}{*}{Method} & \multicolumn{4}{c|}{Compressible}                                                                    & \multicolumn{4}{c|}{Nearly-incompressible}                                                           & \multicolumn{1}{c|}{\multirow{2}{*}{Unsorted}} & \multicolumn{2}{c|}{\multirow{2}{*}{Sorted}}  \\ \cline{2-9}
				& \multicolumn{2}{c|}{Structured}                             & \multicolumn{2}{c|}{Voronoi}           & \multicolumn{2}{c|}{Structured}                             & \multicolumn{2}{c|}{Voronoi}           & \multicolumn{1}{c|}{}                          & \multicolumn{2}{c|}{}                         \\ \cline{2-12} 
				& \multicolumn{1}{c|}{Mesh size} & \multicolumn{1}{c|}{PRE*}  & \multicolumn{1}{c|}{Mesh size} & PRE*  & \multicolumn{1}{c|}{Mesh size} & \multicolumn{1}{c|}{PRE*}  & \multicolumn{1}{c|}{Mesh size} & PRE*  & \multicolumn{1}{c|}{Avg. PRE*}                   & \multicolumn{1}{c|}{Method}        & Avg. PRE* \\ \hline
				Disp-20\%               & \multicolumn{1}{c|}{0.035655}  & \multicolumn{1}{c|}{30.99} & \multicolumn{1}{c|}{0.030985}  & 40.62 & \multicolumn{1}{c|}{0.018041}  & \multicolumn{1}{c|}{61.24} & \multicolumn{1}{c|}{0.022757}  & 55.34 & \multicolumn{1}{c|}{47.05}                     & \multicolumn{1}{c|}{Disp\&Z2-15\%} & 45.29    \\ \hline
				SJ-25\%                 & \multicolumn{1}{c|}{0.024089}  & \multicolumn{1}{c|}{45.87} & \multicolumn{1}{c|}{0.023812}  & 52.85 & \multicolumn{1}{c|}{0.017906}  & \multicolumn{1}{c|}{61.70} & \multicolumn{1}{c|}{0.023054}  & 54.62 & \multicolumn{1}{c|}{53.76}                     & \multicolumn{1}{c|}{Disp\&SJ-15\%} & 45.39    \\ \hline
				Z2-25\%                 & \multicolumn{1}{c|}{0.025869}  & \multicolumn{1}{c|}{42.71} & \multicolumn{1}{c|}{0.025154}  & 50.03 & \multicolumn{1}{c|}{0.017472}  & \multicolumn{1}{c|}{63.24} & \multicolumn{1}{c|}{0.021670}  & 58.11 & \multicolumn{1}{c|}{53.52}                     & \multicolumn{1}{c|}{Disp-20\%}     & 47.05    \\ \hline
				Disp\&SJ-15\%           & \multicolumn{1}{c|}{0.035854}  & \multicolumn{1}{c|}{30.82} & \multicolumn{1}{c|}{0.033592}  & 37.46 & \multicolumn{1}{c|}{0.018339}  & \multicolumn{1}{c|}{60.25} & \multicolumn{1}{c|}{0.023738}  & 53.05 & \multicolumn{1}{c|}{45.39}                     & \multicolumn{1}{c|}{Z2-25\%}       & 53.52    \\ \hline
				Disp\&Z2-15\%           & \multicolumn{1}{c|}{0.036150}  & \multicolumn{1}{c|}{30.56} & \multicolumn{1}{c|}{0.034096}  & 36.91 & \multicolumn{1}{c|}{0.018319}  & \multicolumn{1}{c|}{60.31} & \multicolumn{1}{c|}{0.023591}  & 53.38 & \multicolumn{1}{c|}{45.29}                     & \multicolumn{1}{c|}{SJ-25\%}       & 53.76    \\ \hline
				Ref                     & \multicolumn{1}{c|}{0.011049}  & \multicolumn{1}{c|}{}      & \multicolumn{1}{c|}{0.012585}  &       & \multicolumn{1}{c|}{0.011049}  & \multicolumn{1}{c|}{}      & \multicolumn{1}{c|}{0.012593}  &       & \multicolumn{3}{c|}{}                                                                          \\ \hline
			\end{tabular}
		\end{adjustbox}
		\caption{Performance summary of the VEM in terms of its convergence behaviour in the $\mathcal{H}^{1}$ error norm with respect to mesh size for the various refinement procedures considered for problem~A(1).
			\label{tab:MethComp_PerformancePlateWithHoleTractionMeshSize}}
	\end{table}
	\FloatBarrier
	
\begin{figure}[ht!]
		\centering
		\begin{subfigure}[t]{0.5\textwidth}
			\centering
			\includegraphics[width=0.95\textwidth]{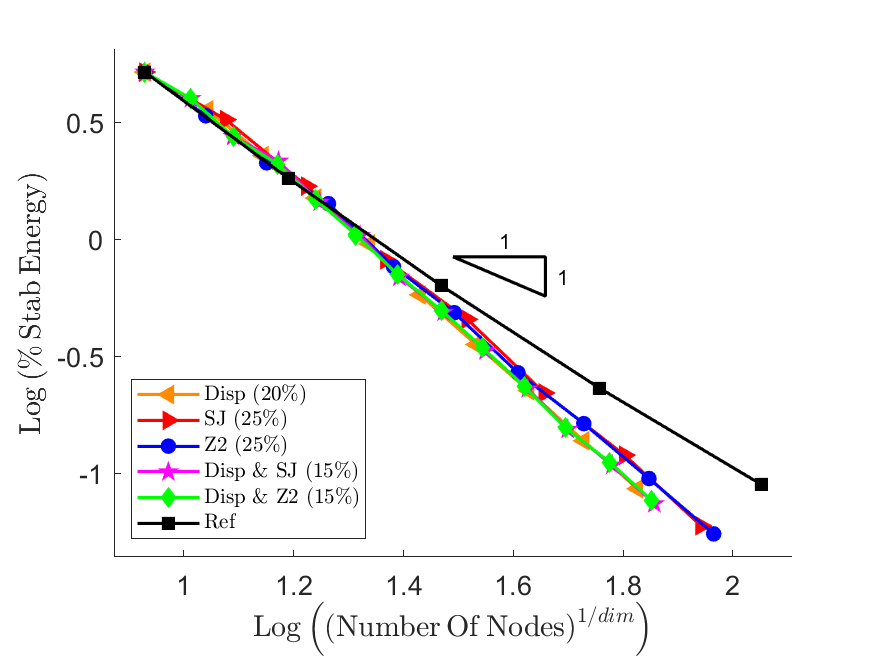}
			\caption{Structured mesh}
		\end{subfigure}%
		\begin{subfigure}[t]{0.5\textwidth}
			\centering
			\includegraphics[width=0.95\textwidth]{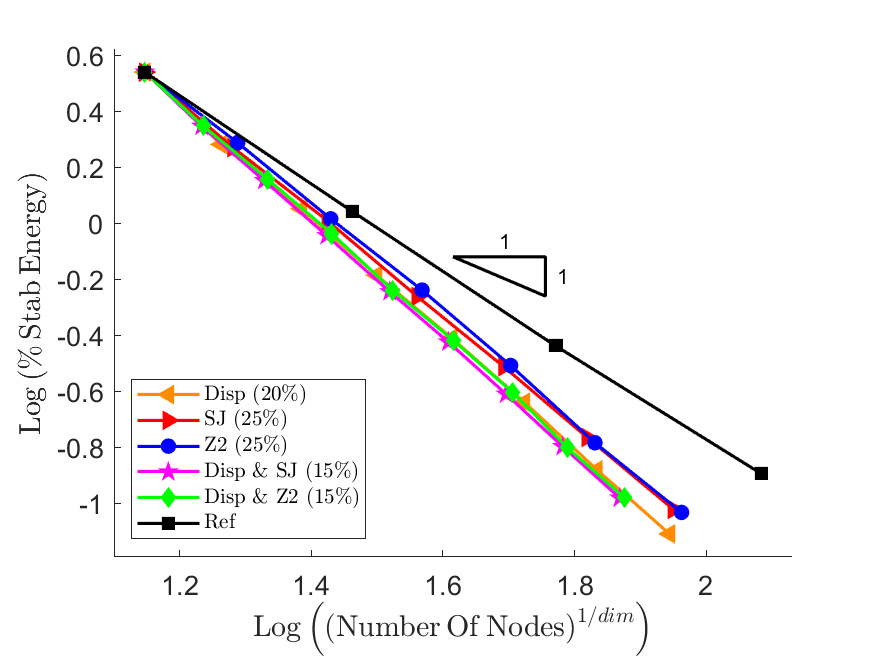}
			\caption{Voronoi mesh}
		\end{subfigure}
		\caption{PSE vs number of nodes/vertices in the discretization for problem~A(1) for the various refinement procedures considered for structured and unstructured meshes with a compressible Poisson's ratio.
			\label{fig:MethComp_PlateWithHoleTractionConvergencePSE}}
	\end{figure}

	The convergence behaviour in the PSE of the VEM for problem~A(1) for the various refinement procedures considered is depicted in Figure~\ref{fig:MethComp_PlateWithHoleTractionConvergencePSE} on a logarithmic scale. Here, the PSE is plotted against the number of nodes/vertices in the discretization for structured and unstructured/Voronoi meshes with a compressible Poisson's ratio. 
	The convergence behaviour is similar for all adaptive procedures. For structured meshes in the coarse mesh range the adaptive procedures exhibit similar convergence behaviour to the reference procedure. However, as the level of mesh refinement increases the adaptive procedures exhibit superior convergence behaviour and outperform the reference procedure. For unstructured/Voronoi meshes the adaptive procedures exhibit superior convergence behaviour and outperform the reference procedure throughout the domain. Additionally, the displacement-based, combined displacement-based and strain jump-based, and combined displacement-based and ${Z^{2}}$-like procedures exhibit slightly superior convergence behaviour to the strain jump-based and ${Z^{2}}$-like procedures.

	\FloatBarrier
	\FloatBarrier
	
	The performance of the VEM in terms of its convergence behaviour in the PSE with respect to the number of nodes/vertices in the discretization for the various refinement procedures considered for problem~A(1) is summarized in Table~\ref{tab:MethComp_PerformancePlateWithHoleTractionPSE}. Here, the performance, as measured by the PRE, is presented for the cases of structured and Voronoi meshes with compressible and nearly incompressible Poisson's ratios. Additionally, a combined average PRE over the range of mesh types and Poisson's ratios considered is computed and presented in the right-hand portion of the table.
	The good performance of the combined displacement-based and strain jump-based, and combined displacement-based and ${Z^{2}}$-like procedures is again evident, with these procedures slightly outperforming the displacement-based procedure.
	
	\FloatBarrier
	\begin{table}[ht!]
		\centering 
		\begin{adjustbox}{max width=\textwidth}
			\begin{tabular}{|c|cccc|cccc|ccc|}
				\hline
				\multirow{3}{*}{Method} & \multicolumn{4}{c|}{Compressible}                                                                  & \multicolumn{4}{c|}{Nearly-incompressible}                                                         & \multicolumn{1}{c|}{\multirow{2}{*}{Unsorted}} & \multicolumn{2}{c|}{\multirow{2}{*}{Sorted}}  \\ \cline{2-9}
				& \multicolumn{2}{c|}{Structured}                            & \multicolumn{2}{c|}{Voronoi}          & \multicolumn{2}{c|}{Structured}                            & \multicolumn{2}{c|}{Voronoi}          & \multicolumn{1}{c|}{}                          & \multicolumn{2}{c|}{}                         \\ \cline{2-12} 
				& \multicolumn{1}{c|}{nNodes}   & \multicolumn{1}{c|}{PRE}   & \multicolumn{1}{c|}{nNodes}   & PRE   & \multicolumn{1}{c|}{nNodes}   & \multicolumn{1}{c|}{PRE}   & \multicolumn{1}{c|}{nNodes}   & PRE   & \multicolumn{1}{c|}{Avg PRE}                   & \multicolumn{1}{c|}{Method}        & Avg. PRE \\ \hline
				Disp-20\%               & \multicolumn{1}{c|}{4303.95}  & \multicolumn{1}{c|}{33.96} & \multicolumn{1}{c|}{4795.70}  & 32.68 & \multicolumn{1}{c|}{3799.42}  & \multicolumn{1}{c|}{29.98} & \multicolumn{1}{c|}{4597.87}  & 31.62 & \multicolumn{1}{c|}{32.06}                     & \multicolumn{1}{c|}{Disp\&SJ-15\%} & 31.22    \\ \hline
				SJ-25\%                 & \multicolumn{1}{c|}{5257.28}  & \multicolumn{1}{c|}{41.49} & \multicolumn{1}{c|}{5885.60}  & 40.11 & \multicolumn{1}{c|}{4606.90}  & \multicolumn{1}{c|}{36.35} & \multicolumn{1}{c|}{5189.88}  & 35.69 & \multicolumn{1}{c|}{38.41}                     & \multicolumn{1}{c|}{Disp\&Z2-15\%} & 31.49    \\ \hline
				Z2-25\%                 & \multicolumn{1}{c|}{5238.11}  & \multicolumn{1}{c|}{41.34} & \multicolumn{1}{c|}{5999.77}  & 40.88 & \multicolumn{1}{c|}{4712.94}  & \multicolumn{1}{c|}{37.19} & \multicolumn{1}{c|}{5727.16}  & 39.38 & \multicolumn{1}{c|}{39.70}                     & \multicolumn{1}{c|}{Disp-20\%}     & 32.06    \\ \hline
				Disp\&SJ-15\%           & \multicolumn{1}{c|}{4337.02}  & \multicolumn{1}{c|}{34.23} & \multicolumn{1}{c|}{4516.62}  & 30.78 & \multicolumn{1}{c|}{3855.01}  & \multicolumn{1}{c|}{30.42} & \multicolumn{1}{c|}{4282.17}  & 29.45 & \multicolumn{1}{c|}{31.22}                     & \multicolumn{1}{c|}{SJ-25\%}       & 38.41    \\ \hline
				Disp\&Z2-15\%           & \multicolumn{1}{c|}{4365.86}  & \multicolumn{1}{c|}{34.45} & \multicolumn{1}{c|}{4654.08}  & 31.71 & \multicolumn{1}{c|}{3846.33}  & \multicolumn{1}{c|}{30.35} & \multicolumn{1}{c|}{4281.16}  & 29.44 & \multicolumn{1}{c|}{31.49}                     & \multicolumn{1}{c|}{Z2-25\%}       & 39.70    \\ \hline
				Ref                     & \multicolumn{1}{c|}{12672.00} & \multicolumn{1}{c|}{}      & \multicolumn{1}{c|}{14675.00} &       & \multicolumn{1}{c|}{12672.00} & \multicolumn{1}{c|}{}      & \multicolumn{1}{c|}{14542.00} &       & \multicolumn{3}{c|}{}                                                                          \\ \hline
			\end{tabular}
		\end{adjustbox}
		\caption{Performance summary of the VEM in terms of its convergence behaviour in the PSE with respect to the number of nodes/vertices in the discretization for the various refinement procedures considered for problem~A(1).
			\label{tab:MethComp_PerformancePlateWithHoleTractionPSE}}
	\end{table}
	\FloatBarrier
	
\begin{figure}[ht!]
		\centering
		\begin{subfigure}[t]{0.5\textwidth}
			\centering
			\includegraphics[width=0.8\textwidth]{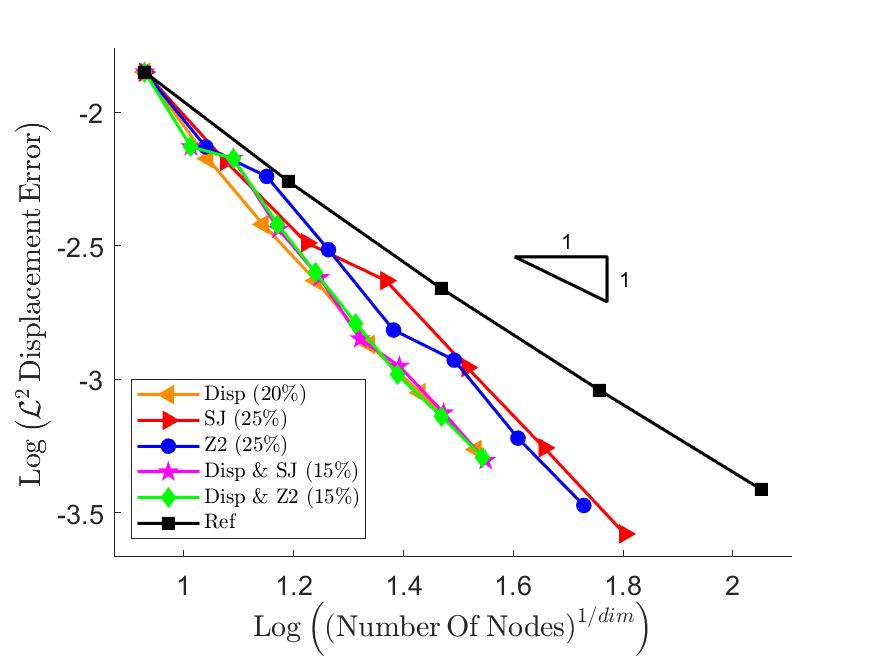}
			\caption{Displacement error - Structured mesh}
		\end{subfigure}%
		\begin{subfigure}[t]{0.5\textwidth}
			\centering
			\includegraphics[width=0.8\textwidth]{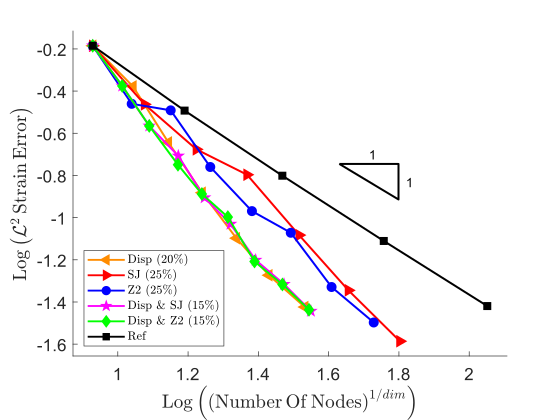}
			\caption{Strain error - Structured mesh}
		\end{subfigure}
		\vskip \baselineskip 
		\vspace*{-3mm}
		\begin{subfigure}[t]{0.5\textwidth}
			\centering
			\includegraphics[width=0.8\textwidth]{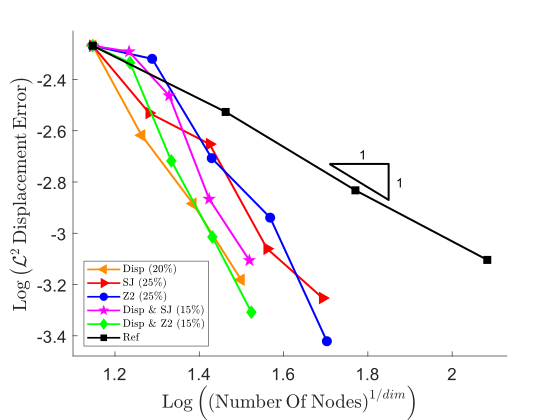}
			\caption{Displacement error - Voronoi mesh}
		\end{subfigure}%
		\begin{subfigure}[t]{0.5\textwidth}
			\centering
			\includegraphics[width=0.8\textwidth]{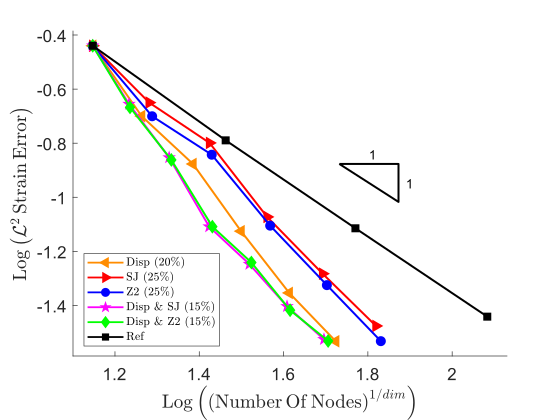}
			\caption{Strain error - Voronoi mesh}
		\end{subfigure}
		\caption{Displacement and strain $\mathcal{L}^{2}$ error components vs $n_{\rm v}$ for problem~A(1) for the various refinement procedures considered for structured and unstructured/Voronoi meshes with a compressible Poisson's ratio.
			\label{fig:MethComp_PlateWithHoleTractionConvergenceComponents}}
	\end{figure}

	The convergence behaviour in the ${\mathcal{L}^{2}}$ error norm of the displacement and strain field approximations for the various refinement procedures considered is plotted in Figure~\ref{fig:MethComp_PlateWithHoleTractionConvergenceComponents} against the number of vertices/nodes in the discretization for problem~A(1) for structured and unstructured/Voronoi meshes with a compressible Poisson's ratio.
	In general, the expected behaviour is observed with good convergence behaviour exhibited by all adaptive refinement procedures, while the best overall behaviour is exhibited by the combined displacement-based and strain jump-based, and combined displacement-based and ${Z^{2}}$-like procedures. 
	In the case of the displacement error for Voronoi meshes the convergence behaviour is slightly more erratic than observed in the other cases. This is most likely due to a combination of the random nature of the Voronoi mesh refinement coupled with the relatively low magnitude of the displacement error. These two features would cause small fluctuations in the error resulting from the randomized remeshing procedure to be exaggerated by the logarithmic scale of the plot.
	Interestingly, there is no clear correlation between the type of refinement procedure used and relative performance in a specific error contribution. For example, the displacement-based procedure does not exhibit clearly superior convergence behaviour in the displacement error, nor clearly worse convergence behaviour in the strain error. 
	Rather, all proposed adaptive refinement procedures improve the approximation of the displacement and strain fields very effectively. Overall, for both mesh types and both error components the best, and most consistently good, performance is exhibited by the combined displacement-based and ${Z^{2}}$-like procedures.
	
	\FloatBarrier
	\FloatBarrier
	
	\subsubsection{Overall}
	\label{subsubsec:MethodComp_Overall}
	
	To determine the best overall refinement procedure the performance over the full range of problems with respect to the number of nodes/vertices, run time, mesh size, and PSE is considered for both mesh types and for compressible and nearly incompressible Poisson's ratios.
	
	The performance of the VEM in terms of its convergence behaviour in the $\mathcal{H}^{1}$ error norm with respect to the number of vertices/nodes in the discretization for the various refinement procedures considered and over the range of sample problems considered is summarized in Table~\ref{tab:MethComp_PerformanceNumberOfNodes}. Here, the combined average performance, as described and computed in Table~\ref{tab:MethComp_PerformancePlateWithHoleTractionNumberOfNodes}, is presented for each sample problem. Additionally, an overall average PRE over the range of sample problems considered is computed and presented in the right-hand portion of the table. 
	All of the proposed adaptive refinement procedures exhibit good performance over the full range of example problems and represent significant improvements in computational efficiency compared to the reference procedure. 
	The worst performing of the proposed procedures, i.e. the ${Z^{2}}$-like procedure with ${T=25\%}$, requires on average ${22.29\%}$ of the number of modes/vertices required by the reference procedure to achieve the same level of accuracy. This improvement in performance is most important in the context of the global stiffness matrix whose size is the square of the number of degrees of freedom. Thus, the total number of matrix entries decreases significantly and the size of the global stiffness matrix required when using the ${Z^{2}}$-like procedure is ${\approx 5\%}$ of the size required by the reference procedure.
	The best performing of the proposed procedures is the combination of the displacement-based and ${Z^{2}}$-like procedures with ${T=15\%}$. This procedure requires on average ${16.79\%}$ of the number of nodes/vertices required by the reference procedure, which corresponds to a global stiffness matrix that is  ${<3\%}$ of the size required by the reference procedure.
	
	\FloatBarrier
	\begin{table}[ht!]
		\centering 
		\begin{adjustbox}{max width=\textwidth}
			\begin{tabular}{|c|c|c|c|c|c|c|c|cc|}
				\hline
				\multirow{2}{*}{Method} & \multirow{2}{*}{Prob. 1} & \multirow{2}{*}{Prob. 2} & \multirow{2}{*}{Prob. 3} & \multirow{2}{*}{Prob. 4} & \multirow{2}{*}{Prob. 5} & \multirow{2}{*}{Prob. 6} & Unsorted & \multicolumn{2}{c|}{Sorted}                   \\ \cline{8-10} 
				&                          &                          &                          &                          &                          &                          & Avg PRE  & \multicolumn{1}{c|}{Method}        & Avg. PRE \\ \hline
				Disp-20\%               & 20.63                    & 15.62                    & 7.03                     & 8.98                     & 24.39                    & 26.56                    & 17.20    & \multicolumn{1}{c|}{Disp\&Z2-15\%} & 16.79    \\ \hline
				SJ-25\%                 & 27.70                    & 24.44                    & 9.29                     & 13.28                    & 28.33                    & 28.46                    & 21.91    & \multicolumn{1}{c|}{Disp\&SJ-15\%} & 16.85    \\ \hline
				Z2-25\%                 & 27.90                    & 22.49                    & 8.74                     & 13.14                    & 28.67                    & 32.79                    & 22.29    & \multicolumn{1}{c|}{Disp-20\%}     & 17.20    \\ \hline
				Disp\&SJ-15\%           & 20.20                    & 16.06                    & 5.97                     & 8.74                     & 24.16                    & 25.94                    & 16.85    & \multicolumn{1}{c|}{SJ-25\%}       & 21.91    \\ \hline
				Disp\&Z2-15\%           & 20.20                    & 15.76                    & 5.77                     & 8.62                     & 24.22                    & 26.19                    & 16.79    & \multicolumn{1}{c|}{Z2-25\%}       & 22.29    \\ \hline
			\end{tabular}
		\end{adjustbox}
		\caption{Performance summary of the VEM in terms of its convergence behaviour in the $\mathcal{H}^{1}$ error norm with respect to the number of vertices/nodes in the discretization for the various refinement procedures and range of sample problems considered.
			\label{tab:MethComp_PerformanceNumberOfNodes}}
	\end{table}
	\FloatBarrier
	
	The performance of the VEM in terms of its convergence behaviour in the $\mathcal{H}^{1}$ error norm with respect to run time (including remeshing time) for the various refinement procedures considered and over the range of sample problems considered is summarized in Table~\ref{tab:MethComp_PerformanceRunTime}. Here, the combined average PRE is presented for each sample problem along with an overall average PRE. As observed in Table~\ref{tab:MethComp_PerformanceNumberOfNodes} the best performance is achieved when using a combination of refinement procedures. The combined displacement-based and strain jump-based procedures and the combined displacement-based and ${Z^{2}}$-like procedures exhibit very similar levels of performance. However, the best performance is exhibited by the combined displacement-based and strain jump-based procedures with ${T=15\%}$. This procedure requires on average ${6.44\%}$ of the run time required by the (cumulative) reference procedure.
	
	
	\FloatBarrier
	\begin{table}[ht!]
		\centering 
		\begin{adjustbox}{max width=\textwidth}
			\begin{tabular}{|c|c|c|c|c|c|c|c|cc|}
				\hline
				\multirow{2}{*}{Method} & \multirow{2}{*}{Prob. 1} & \multirow{2}{*}{Prob. 2} & \multirow{2}{*}{Prob. 3} & \multirow{2}{*}{Prob. 4} & \multirow{2}{*}{Prob. 5} & \multirow{2}{*}{Prob. 6} & Unsorted & \multicolumn{2}{c|}{Sorted}                   \\ \cline{8-10} 
				&                          &                          &                          &                          &                          &                          & Avg PRE  & \multicolumn{1}{c|}{Method}        & Avg. PRE \\ \hline
				Disp-20\%               & 8.27                     & 5.89                     & 2.50                     & 3.03                     & 10.93                    & 9.36                     & 6.66     & \multicolumn{1}{c|}{Disp\&SJ-15\%} & 6.44     \\ \hline
				SJ-25\%                 & 12.86                    & 10.69                    & 3.65                     & 4.83                     & 14.23                    & 12.20                    & 9.74     & \multicolumn{1}{c|}{Disp\&Z2-15\%} & 6.46     \\ \hline
				Z2-25\%                 & 12.59                    & 9.56                     & 3.30                     & 4.99                     & 14.20                    & 13.66                    & 9.72     & \multicolumn{1}{c|}{Disp-20\%}     & 6.66     \\ \hline
				Disp\&SJ-15\%           & 8.66                     & 6.15                     & 2.01                     & 2.91                     & 10.99                    & 7.96                     & 6.44     & \multicolumn{1}{c|}{Z2-25\%}       & 9.72     \\ \hline
				Disp\&Z2-15\%           & 8.55                     & 6.11                     & 1.95                     & 2.93                     & 11.05                    & 8.14                     & 6.46     & \multicolumn{1}{c|}{SJ-25\%}       & 9.74     \\ \hline
			\end{tabular}
		\end{adjustbox}
		\caption{Performance summary of the VEM in terms of its convergence behaviour in the $\mathcal{H}^{1}$ error norm with respect to run time (including remeshing time) for the various refinement procedures and range of sample problems considered.
			\label{tab:MethComp_PerformanceRunTime}}
	\end{table}
	\FloatBarrier
	
	The performance of the VEM in terms of its convergence behaviour in the $\mathcal{H}^{1}$ error norm with respect to mesh size for the various refinement procedures considered and over the range of sample problems considered is summarized in Table~\ref{tab:MethComp_PerformanceMeshSize}. Here, the combined average PRE* is presented for each sample problem along with an overall average PRE*.
	The qualitative behaviour of the overall average PRE* is very similar to that observed in Table~\ref{tab:MethComp_PerformanceNumberOfNodes} with respect to the number of nodes. 
	All of the proposed adaptive refinement procedures exhibit good performance over the full range of example problems and represent significant improvements in computational efficiency compared to the reference procedure. Particularly good performance is achieved when using a combination of refinement procedures. The combined displacement-based and strain jump-based procedures and the combined displacement-based and ${Z^{2}}$-like procedures exhibit very similar levels of performance. However, the best performance is exhibited by the combined displacement-based and ${Z^{2}}$-like procedures with ${T=15\%}$. This procedure has an overall average PRE* of ${41.88\%}$ which corresponds to requiring elements that are on average ${\approx 2.4}$ times larger than those used by the reference procedure to achieve the same level of accuracy.
	
	\FloatBarrier
	\begin{table}[ht!]
		\centering 
		\begin{adjustbox}{max width=\textwidth}
			\begin{tabular}{|c|c|c|c|c|c|c|c|cc|}
				\hline
				\multirow{2}{*}{Method} & \multirow{2}{*}{Prob. 1} & \multirow{2}{*}{Prob. 2} & \multirow{2}{*}{Prob. 3} & \multirow{2}{*}{Prob. 4} & \multirow{2}{*}{Prob. 5} & \multirow{2}{*}{Prob. 6} & Unsorted & \multicolumn{2}{c|}{Sorted}                   \\ \cline{8-10} 
				&                          &                          &                          &                          &                          &                          & Avg PRE*  & \multicolumn{1}{c|}{Method}        & Avg. PRE* \\ \hline
				Disp-20\%               & 47.05                    & 41.37                    & 32.52                    & 33.19                    & 58.19                    & 55.34                    & 44.61    & \multicolumn{1}{c|}{Disp\&Z2-15\%} & 41.88    \\ \hline
				SJ-25\%                 & 53.76                    & 50.33                    & 34.14                    & 38.09                    & 56.95                    & 54.62                    & 47.98    & \multicolumn{1}{c|}{Disp\&SJ-15\%} & 42.05    \\ \hline
				Z2-25\%                 & 53.52                    & 48.05                    & 33.17                    & 37.65                    & 56.90                    & 58.11                    & 47.90    & \multicolumn{1}{c|}{Disp-20\%}     & 44.61    \\ \hline
				Disp\&SJ-15\%           & 45.39                    & 40.65                    & 27.80                    & 31.12                    & 54.26                    & 53.05                    & 42.05    & \multicolumn{1}{c|}{Z2-25\%}       & 47.90    \\ \hline
				Disp\&Z2-15\%           & 45.29                    & 40.51                    & 27.08                    & 30.81                    & 54.23                    & 53.38                    & 41.88    & \multicolumn{1}{c|}{SJ-25\%}       & 47.98    \\ \hline
			\end{tabular}
		\end{adjustbox}
		\caption{Performance summary of the VEM in terms of its convergence behaviour in the $\mathcal{H}^{1}$ error norm with respect to mesh size for the various refinement procedures and range of sample problems considered.
			\label{tab:MethComp_PerformanceMeshSize}}
	\end{table}
	\FloatBarrier
	
	The performance of the VEM in terms of its convergence behaviour in the PSE with respect to the number of nodes/vertices in the discretization for the various refinement procedures considered and over the range of sample problems considered is summarized in Table~\ref{tab:MethComp_PerformancePSE}. Here, the combined average PRE is presented for each sample problem along with an overall average PRE. Interestingly, the best performance is again exhibited by the cases of combined refinement procedures even though the formulation of the displacement-based indicator is similar to that of the stabilization term. This further demonstrates the efficacy of using refinement procedures in combination. 
	The best performance is exhibited by the combined displacement-based and strain jump-based procedures with ${T=15\%}$. This procedure requires on average ${21.62\%}$ of the number of nodes/vertices as the reference procedure to achieve the same PSE.
	
	\FloatBarrier
	\begin{table}[ht!]
		\centering 
		\begin{adjustbox}{max width=\textwidth}
			\begin{tabular}{|c|c|c|c|c|c|c|c|cc|}
				\hline
				\multirow{2}{*}{Method} & \multirow{2}{*}{Prob. 1} & \multirow{2}{*}{Prob. 2} & \multirow{2}{*}{Prob. 3} & \multirow{2}{*}{Prob. 4} & \multirow{2}{*}{Prob. 5} & \multirow{2}{*}{Prob. 6} & Unsorted & \multicolumn{2}{c|}{Sorted}                   \\ \cline{8-10} 
				&                          &                          &                          &                          &                          &                          & Avg PRE  & \multicolumn{1}{c|}{Method}        & Avg. PRE \\ \hline
				Disp-20\%               & 32.06                    & 24.41                    & 8.53                     & 12.24                    & 25.03                    & 31.62                    & 22.32    & \multicolumn{1}{c|}{Disp\&SJ-15\%} & 21.62    \\ \hline
				SJ-25\%                 & 38.41                    & 32.20                    & 11.72                    & 17.03                    & 29.39                    & 35.69                    & 27.41    & \multicolumn{1}{c|}{Disp\&Z2-15\%} & 21.65    \\ \hline
				Z2-25\%                 & 39.70                    & 31.32                    & 11.28                    & 17.13                    & 30.35                    & 39.38                    & 28.19    & \multicolumn{1}{c|}{Disp-20\%}     & 22.32    \\ \hline
				Disp\&SJ-15\%           & 31.22                    & 25.01                    & 7.73                     & 11.71                    & 24.60                    & 29.45                    & 21.62    & \multicolumn{1}{c|}{SJ-25\%}       & 27.41    \\ \hline
				Disp\&Z2-15\%           & 31.49                    & 24.40                    & 8.02                     & 11.74                    & 24.81                    & 29.44                    & 21.65    & \multicolumn{1}{c|}{Z2-25\%}       & 28.19    \\ \hline
			\end{tabular}
		\end{adjustbox}
		\caption{Performance summary of the VEM in terms of its convergence behaviour in the PSE with respect to the number of nodes/vertices in the discretization for the various refinement procedures and range of sample problems considered.
			\label{tab:MethComp_PerformancePSE}}
	\end{table}
	\FloatBarrier
	
	The performance of the VEM over all metrics considered, for the various refinement procedures considered, and over the range of sample problems considered, is summarized in Table~\ref{tab:MethComp_Overall}. Here, the overall average PRE, or PRE* in the case of mesh size, is presented for each metric along with a total average PRE over the range of performance metrics. For clarity the best overall average PRE in each metric is indicated in bold. The superior performance achieved when using a combination of refinement procedures is again evident with the best performance in each metric exhibited by either the combined displacement-based and strain jump-based procedures or the combined displacement-based and ${Z^{2}}$-like procedures.
	Over the range of performance metrics it is found, consistently with \cite{vanHuyssteen2022}, that the best performance is achieved when using a combination of the displacement-based and ${Z^{2}}$-like procedures.
	
	\FloatBarrier
	\begin{table}[ht!]
		\centering 
		\begin{adjustbox}{max width=\textwidth}
			\begin{tabular}{|c|c|c|c|c|c|cc|}
				\hline
				\multirow{2}{*}{Method} & \multirow{2}{*}{nNodes} & \multirow{2}{*}{Run time} & \multirow{2}{*}{Mesh size} & \multirow{2}{*}{PSE} & Unsorted & \multicolumn{2}{c|}{Sorted}                   \\ \cline{6-8} 
				&                         &                           &                            &                      & Avg PRE  & \multicolumn{1}{c|}{Method}        & Avg. PRE \\ \hline
				Disp-20\%               & 17.20                   & 6.66                      & 44.61                      & 22.32                & 22.70    & \multicolumn{1}{c|}{Disp\&Z2-15\%} & 21.70    \\ \hline
				SJ-25\%                 & 21.91                   & 9.74                      & 47.98                      & 27.41                & 26.76    & \multicolumn{1}{c|}{Disp\&SJ-15\%} & 21.74    \\ \hline
				Z2-25\%                 & 22.29                   & 9.72                      & 47.90                      & 28.19                & 27.03    & \multicolumn{1}{c|}{Disp-20\%}     & 22.70    \\ \hline
				Disp\&SJ-15\%           & 16.85                   & \textbf{6.44}                      & 42.05                      & \textbf{21.62}                & 21.74    & \multicolumn{1}{c|}{SJ-25\%}       & 26.76    \\ \hline
				Disp\&Z2-15\%           & \textbf{16.79}                   & 6.46                      & \textbf{41.88}                      & 21.65                & \textbf{21.70}    & \multicolumn{1}{c|}{Z2-25\%}       & 27.03    \\ \hline
			\end{tabular}
		\end{adjustbox}
		\caption{Performance summary of the VEM over all metrics considered and for the various refinement procedures considered.
			\label{tab:MethComp_Overall}}
	\end{table}
	\FloatBarrier
	
	\section{Discussion and conclusion} 
	\label{sec:Conclusion}
	%
	In this work a range of mesh refinement procedures for meshes of arbitrarily shaped virtual elements have been formulated and implemented for isotropic linear elastic problems in two-dimensions. Specifically, the three mesh refinement indicators originally proposed in \cite{vanHuyssteen2022} for structured meshes have been presented for the case of arbitrary polygonal element geometries. Additionally, a simple procedure for selecting elements qualifying for refinement has been presented.
	
	The proposed procedures comprise displacement-based, strain jump-based, and $Z^{2}$-like refinement indicators. The indicators were motivated respectively by seeking to; improve the accuracy of the displacement field approximation, smooth the approximation of the strain field over the problem domain, and formulate a variation of the well-known $Z^{2}$ error estimator for finite elements that is suitable for implementation in a virtual element setting. 
	
	The various mesh refinement procedures were studied both individually and in combination for problems involving structured and unstructured/Voronoi meshes and with Poisson's ratios corresponding to compressible and nearly incompressible materials.
	The efficacy and efficiency of the VEM using the proposed procedures was studied numerically through a wide range of example problems. The efficacy was investigated in terms of the ${\mathcal{H}^{1}}$ error norm and the efficiency was studied in terms of the number of degrees of freedom, run time, and mesh size, compared to a reference procedure. Additionally, the convergence of the contribution of the stabilization term, and the convergence behaviour of the displacement and strain error components in the ${\mathcal{L}^{2}}$ error norm was investigated.
	
	Through the numerical investigation it was found that all of the proposed adaptive refinement procedures significantly outperformed the reference procedure. That is, the adaptive procedures generated solutions of equivalent accuracy to the reference procedure while using significantly fewer nodes and less run time. Thus, representing a significant improvement in computational efficiency. The best overall performance was achieved using a combination of the displacement-based and ${Z^{2}}$-like refinement procedures, when refining ${15\%}$ of the elements in a mesh at each refinement step. It was found, on average, that this combined procedure was able to generate solutions of equivalent accuracy to the reference procedure while using ${16.79\%}$ of the number of nodes and ${6.46\%}$ of the run time used by the reference procedure.
	The good performance exhibited by this adaptive approach over a wide range of challenging example problems, for a variety of mesh types, and for cases of compressible and nearly incompressible materials, demonstrates its versatility, efficiency and suitability for application to elastic analyses using the virtual element method.

	Future work of interest would be the extension to non-linear problems, higher-order formulations, and problems in three-dimensions. Additionally, the investigation of approaches for combined adaptive mesh coarsening and refinement would be of significant interest.
	
	\section*{Supporting information}
	This work is supported by supplementary material comprising additional example problems\cite{DanielHuyssteen2023}.
		
	\section*{Declaration of competing interest} The authors declare that they have no known competing financial interests or personal relationships that could have appeared to influence the work reported in this paper.
	
	\section*{Acknowledgements} 
	This work was carried out with support from the German Science Foundation (DFG) and the National Scientific and Technical Research Council of Argentina (CONICET) through project number DFG 544/68-1 (431843479).
	The authors acknowledge with thanks this support.
	
	\bibliographystyle{elsarticle-num}
	\bibliography{VEM_References}

\end{document}